\documentclass{article}

\usepackage{silence}
\usepackage {amsmath, amsfonts, amssymb, mathrsfs, amsthm,stmaryrd, sectsty,hyphenat,url,accents,graphicx,comment, tocloft}

\usepackage{hyperref}  
\usepackage{xcolor}
\usepackage[all]{xy}
\usepackage[toc,page]{appendix}

\usepackage[shortlabels]{enumitem} 

\setlist[enumerate,1]{label=(\arabic*)} 
\setlist[enumerate,2]{label=(\alph*)}
\setlist[enumerate,3]{label=(\roman*)}

\usepackage[nonewpage]{imakeidx}
\makeindex[intoc,name=notation,title=Selected notation]
\makeindex[intoc,name=terminology,title=Selected terminology]

\def\mf#1{\mathfrak{#1}}
\def\mc#1{\mathcal{#1}}
\def\mb#1{\mathbb{#1}}
\def\tx#1{\textrm{#1}}

\def\R{\mathbb{R}}

\def\C{\mathbb{C}}
\def\Q{\mathbb{Q}}

\def\Z{\mathbb{Z}}

\def\F{\mathbb{F}}
\def\lmod{\backslash}

\def\ol#1{\overline{#1}}
\def\ul#1{\underline{#1}}

\def\hat{\widehat}

\def\sm{\smallsetminus}

\def\<{\langle}
\def\>{\rangle}
\def\Wedge{\bigwedge}

\def\cL{\mathcal{L}}
\def\cO{\mathcal{O}}
\def\sB{\mathcal{B}} \def\ra{\rightarrow}

\def\simto{\buildrel\hbox{\tiny{$\sim$}}\over\longrightarrow}

\newenvironment{mytitle}
{\begin{center}\large\sc}
{\end{center}}

\newcommand{\cB}{\mathcal{B}}

\newcommand{\bF}{\mathbb{F}}
\newcommand{\bC}{\mathbb{C}}
\newcommand{\bR}{\mathbb{R}}
\newcommand{\bZ}{\mathbb{Z}}

\newcommand{\fp}{\mathfrak{p}}

\newcommand{\eps}{\epsilon}
\DeclareMathOperator{\Sp}{Sp}	
\DeclareMathOperator{\GL}{GL}

\DeclareMathOperator{\Inn}{Inn}
\DeclareMathOperator{\cind}{c-ind}
 
\newtheorem{thm}{Theorem}[subsection]
\newtheorem{lem}[thm]{Lemma}

\newtheorem{pro}[thm]{Proposition}
\newtheorem{cor}[thm]{Corollary}
\newtheorem{fct}[thm]{Fact}

\theoremstyle{definition}
\newtheorem{rem}[thm]{Remark}
\newtheorem{dfn}[thm]{Definition}

\newtheorem{exa}[thm]{Example}

\sectionfont{\center\sc\normalsize}
\subsectionfont{\bf\normalsize}

\numberwithin{equation}{subsection}
\renewcommand{\-}{\hyp{}}

\newcommand{\diagblock}[2]{\vcenter{\hbox{\parbox[c][#1][c]{#1}{\centering $#2$}}}}

\DeclareMathOperator{\SO}{SO} 
\DeclareMathOperator{\Mp}{Mp}
\DeclareMathOperator{\SU}{SU} 
\DeclareMathOperator{\U}{U}

\DeclareMathOperator{\WM}{WM}
\newcommand{\WMuniv}{\mathrm{WM}}
\DeclareMathOperator{\BW}{BW}
\DeclareMathOperator{\HW}{HW}
\DeclareMathOperator{\QW}{QW}

\DeclareMathOperator{\Aut}{Aut}

\newcommand{\rf}{\mathfrak{f}}
\newcommand{\oV}{\ol{V}}
\newcommand{\oL}{\overline{L}}
\newcommand{\oVS}{\overline{V}}
\newcommand{\oVO}{\overline{V_a}}

\newcommand{\rN}{\mathrm{N}}

\AtEndDocument{
	{\footnotesize
		\vspace{.5cm}

\textsc{Universität Bonn, Mathematisches Institut,}
Endenicher Allee 60, 53115 Bonn, Germany.

\textit{E-mail addresses}: \texttt{fintzen@math.uni-bonn.de} and  \texttt{kaletha@math.uni-bonn.de}

}}

\begin{document}

\begin{mytitle}
    Supercuspidal representations and types via double covers and the Fourier--Weil representation
\end{mytitle}

\begin{center}
Jessica Fintzen and Tasho Kaletha
\end{center}

\begin{abstract}
	We present a construction of supercuspidal representations and types for tame $p$-adic groups whose input is a generic essentially depth-zero representation of a double covers of twisted Levi subgroups. This clarifies and simplifies existing constructions by automatically producing and explaining the twisting quadratic character of Fintzen--Kaletha--Spice (\cite{FKS23}), the roots of unity in the Harish-Chandra character formula, and making the connection to the Galois side of the local Langlands correspondence immediate. The main innovation is the use of a Fourier model of the Weil representation defined on an extension of the symplectic group by the Witt group. \end{abstract}

{
	\renewcommand{\thefootnote}{}  
	\footnotetext{The first-named author was partially supported by 
		the European Research Council (ERC) under the European Union's Horizon 2020 research
		and innovation programme (grant agreement n° 950326).}          
}

\tableofcontents

\section{Introduction} \label{sec:intro}

The construction of supercuspidal representations of tame $p$-adic groups by J.-K.\ Yu (\cite{Yu01}) and its extension (\cite{KimYu}) to types has been widely used in representation theory and number theory. However, in order to use this construction for many natural problems, such as an explicit construction of a local Langlands correspondence, a reduction to depth zero, or to calculate character formulas, the construction needed to be twisted by a subtle quadratic character (\cite{FKS23}) and certain auxiliary choices had to be made, the conceptual origin of which was unclear. In the present paper, we provide a different construction of the same supercuspidal representations and types (which provides all supercuspidal representations and types for all Bernstein blocks if $p$ is not too small) that starts from a slightly different input, more precisely that takes a genuine essentially depth-zero representation of a double cover of a twisted Levi subgroup as an input. 
This new construction fits very easily into the Langlands program and does not require an auxiliary twist for a reduction-to-depth-zero statement. In fact, the quadratic twist of \cite{FKS23} arises naturally in the comparison of this construction with the traditional construction. The new approach has been made possible by the use and careful study of the Fourier--Weil representation of an extension of the (finite and $p$-adic) symplectic group by the Witt group, which has been classically studied in \cite{LionVergne}, \cite{RR93}, and more recently \cite{Thomas2013}.

\subsection{A brief summary of the result}\label{sub:introoverview}
We start with a brief overview of our main results before providing more motivation and discussion and a more detailed overview of the results in a special case. 

Let $F$ be a non-archimedean local field with residue field $\rf$ of odd characteristic $p$.
Let $G$ be a connected reductive group over $F$ that splits over a tamely ramified field extension of $F$.

Yu's construction of supercuspidal representations of $G(F)$ in \cite{Yu01} starts with a supercuspidal $G$-datum, which is a tuple
\[((G_i)_{1\leq i\leq n+1},x,(r_i)_{1\leq i\leq n},\rho^0,(\phi_i)_{1\leq i\leq n}) \]
consisting of a sequence of twisted Levi subgroups $G=G_1 \supseteq G_2 \supsetneq G_3 \supsetneq \hdots \supsetneq G_{n+1}=: G^0$, a point $x$ in the Bruhat--Tits building of $G^0$, characters $\phi_i$ of $G_{i+1}(F)$ of depth $r_i$, and a depth-zero representation $\rho^0$ of the stabilizer $G^0_x$ of $x$ in $G^0(F)$, which satisfy some conditions recalled in Definition \ref{dfn:supercuspidalGdatum}.
From such a tuple, Yu constructs an open and compact modulo center subgroup $K \subset G(F)$ and a representation $\tilde \rho$ of $K$ such that the compact induction $\cind_K^{G(F)} \tilde \rho$ is a supercuspidal irreducible representation of $G(F)$ (\cite{Yu01, Fin21Comp}).
The subgroup $K$ is of the form $G^0_xK^\dagger$ for a pro-$p$ group $K^\dagger \subseteq G(F)_{x,0+}$. The representation $\tilde \rho$ is the tensor product of the extension of $\rho^0$ to $G^0_xK^\dagger$ which is trivial on $K^\dagger$ and a representation $\kappa^{\mathrm{nt}}$. To describe the representation $\kappa^{\mathrm{nt}}$, we note that there is a slightly smaller pro-$p$ subgroup $K^+ \subseteq K^\dagger$ such that $\ol{V} := K^\dagger/K^+$ is an $\bF_p$-vector space that can be equipped with a symplectic form arising from composing a character $\hat \phi$ on $K^+$ that is built from the characters $\phi_i$ with the commutator on $K^\dagger$, see \S\ref{sub:construction} for precise definitions and more details. Note that $G^0_x$ acts on $K^\dagger$ and on $K^+$ via conjugation. This provides an action of $G^0_x$ on $\ol{V}$ that preserves the symplectic form. Now the representation $\kappa^{\mathrm{nt}}$ restricted to $K^\dagger$ is the unique irreducible representation whose restriction to $K^+$ is the character $\hat \phi$ (times the identity). Such a representation is obtained via the theory of Heisenberg representations. The representation  $\kappa^{\mathrm{nt}}$ restricted to $G^0_x$ is constructed by composing the map 
\begin{equation} \label{eqn:G0xtoSp} \bar\jmath :  G^0_x \ra \Sp(\ol{V}) \end{equation}
that arises from the conjugation action
with the Weil representation of $\Sp(\ol{V})$ (and twisting the resulting representation with the product of the characters $\phi_i$).

While this construction provides supercuspidal representations by \cite{Yu01, Fin21Comp}, and provides all supercuspidal representations if $p$ does not divide the order of the absolute Weyl group of $G$ (\cite{Kim07, Fin21}), this parameterization of representations in terms of supercuspidal $G$-data does not fit well into the local Langlands program and has to be twisted by a mysterious looking quadratic character introduced in \cite{FKS23} in order to obtain reduction-to-depth-zero results and calculate its Harish-Chandra character, see \S\ref{sub:motivation} for a more detailed discussion.

In this paper we present a revised approach to the construction of supercuspidal representations that resolves all these concerns. The input for our constructions is a genuine supercuspidal $G$-datum, which is again a tuple 
\[((G_i)_{1\leq i\leq n+1},x,(r_i)_{1\leq i\leq n},\rho^0_{\pm},(\phi_i)_{1\leq i\leq n}) \]
consisting of a sequence of twisted Levi subgroups and characters as above, but this time $\rho^0_{\pm}$ is not a depth-zero representation of $G^0_x$, but it is a genuine depth-zero representation of a double cover $G^0_{\pm,x}$ of $G^0_x$. This double cover is the restriction of the double cover $G^0(F)_{\pm}$ of $G^0(F)$ with respect to the ambient group $G$ as defined in \cite[\S6]{KalDC}. From such an input we can define the group $K_{\pm}:=G^0_{\pm,x}K^\dagger$, which is a double cover of $K=G^0_{x}K^\dagger$, and we can extend the genuine representation $\rho^0_{\pm}$ to a genuine representation of $K_{\pm}$ that is trivial on $K^\dagger$. We can now imitate Yu's construction by constructing a representation $\kappa_{\pm}$ of $K_{\pm}$ such that $\rho^0_{\pm} \otimes \kappa_{\pm}$ factors through a representation of $K$, which we also denote by $\rho^0_{\pm} \otimes \kappa_{\pm}$, and such that $\cind_{K}^{G(F)}(\rho^0_{\pm} \otimes \kappa_{\pm})$ is an irreducible supercuspidal representation. For the restriction of $\kappa_{\pm}$ to $K^\dagger$ we use the theory of Heisenberg representations as before, i.e., $\kappa_{\pm}|_{K^\dagger}=\kappa^{\mathrm{nt}}|_{K^\dagger}$. The challenge lies in extending the representation $\kappa_{\pm}$ to a representation of $G^0_{\pm, x}K^\dagger$. 

We achieve this by replacing the map $ G^0_x \ra \Sp(\ol{V}) $ from \eqref{eqn:G0xtoSp} by
\begin{equation} \label{eqn:G0xtoWM} \tilde\jmath :  G^0_{\pm,x} \ra \WM_\gamma(\ol{V}) ,\end{equation}
where $\WM_\gamma(\ol{V})$ is a covering group of $\Sp(\ol{V})$, more precisely, it fits into the short exact sequence
\[ 1 \ra \mu_8 \ra \WM_\gamma(\ol{V}) \ra \Sp(\ol{V}) \ra 1, \]
where $\mu_8$ denotes the eight roots of unity in the complex numbers, and is the pushout of the universal Witt-metaplectic group $\WMuniv(\ol{V})$ under the \emph{renormalized} Weil index $\tilde\gamma:\QW \ra \mu_8$ on Weissauers's Witt group $\QW$ of characters of second degree (more precisely, on a prime-to-$2$ version of it), as explained in Appendix \ref{sub:WittgroupC} in more details.

The existence and uniqueness of such a $\tilde \jmath$ is a key result of this paper. To state the result, we introduce the notation $V$ for the $F$-vector space $\mathrm{Lie}(G)/\mathrm{Lie}(G^0)$ and we equip $V$ with a symplectic form $\eta = X \circ [\phantom{X},\phantom{X}]$ for $X=\sum_{i=1}^n X_i \in \mathrm{Lie}^*(G)$, where each $X_i$ is a good element that represents $\phi_i|_{G_{i+1}(F)_{x,r_i}}$, as explained just before Theorem \ref{thm:main1b}. The relation between $X_i$ and $\phi_i$ uses the choice of a character $\psi: F \ra \bC^\times$ that is non-trivial on the ring of integers of $F$, and trivial on the maximal ideal of the ring of integers, which we fix throughout. Let $H \subset \Sp(V)$ be the twisted Levi subgroup associated to the weight decomposition of $V$ under the center $Z_{G^0}$ of $G^0$, and let $H(F)_\pm$ be its double cover associated to the set of weights.

\begin{thm}[Imprecise version of Theorem \ref{thm:main1b}]
	There exists a unique lift of \eqref{eqn:G0xtoSp} to a genuine homomorphism \eqref{eqn:G0xtoWM} such that for each semi-simple element $\dot s \in G^0_{\pm,x}$ whose image $\jmath(\dot s) \in H(F)_\pm$ is a norm (as defined in Definition  \ref{dfn:norm}, see also Theorem \ref{thm:main1b}) and has finite order prime to $p$ we have
	\[ \tilde\jmath(\dot s)=\gamma_\psi(q_{\tilde s})\boxtimes \bar\jmath(s), \] 
	where $q_{\tilde s}$ is the quadratic form on the symplectic $F$-vector space $(V,\eta)$ defined by $q_{\tilde s}(v)=\eta(\tilde sv,v)$ with $\tilde s \in \mathrm{GL}(V)$ an appropriate element whose norm is $\jmath(\dot s)$, see Theorem \ref{thm:main1b} and Lemma \ref{lem:topssgamma}(5) for details, $\gamma_\psi$ is the Weil index of quadratic forms on $F$-vector spaces recalled in \S\ref{sub:Weilindex}, and $s$ is the image of $\dot s$ in $G^0_x$. 
	
	The map $\tilde \jmath$ does not depend on the choice of $\psi$ or the choice of $X$.
	
	The value of $\tilde \jmath(\dot g)$ can be explicitly described on for all elements in $\dot g \in  G^0_{\pm,x}$.
\end{thm}

We then replace the Weil representation of $\Sp(\ol{V})$ in the construction of supercuspidal representations by the Fourier--Weil representation of $\WM_\gamma(\ol{V})$ that emerges from \cite{LionVergne}, \cite{RR93}, and \cite{Thomas2013}, and which we review in \S\ref{sub:fwh} and \S\ref{sub:WittgroupC}. Since the Fourier--Weil representation is genuine, the resulting representation $\kappa_{\pm}$ is genuine and hence the representation $\rho^0_{\pm} \otimes \kappa_{\pm}$ is non-genuine, i.e., factors through the surjection $K_{\pm} \twoheadrightarrow K$, as desired. 

 The following is a summary of some of the applications of this construction.

\begin{thm}
	\begin{enumerate}
		\item The resulting representation $\cind_K^{G(F)}(\rho^0_{\pm} \otimes \kappa_{\pm})$ constructed from a genuine supercuspidal $G$-datum is an irreducible supercuspidal representation of $G(F)$. More generally, an analogous construction provides a type from a genuine $G$-datum.
		\item When $p$ is prime to the order of the absolute Weyl group of $G$, all supercuspidal representations arise in this way, and each Bernstein block has a type attached to it that arises in this way.
\item The input for the class of supercuspidal representations from \cite{KalSLP} can also be compressed to a genuine tame $F$-non-singular elliptic pair $(S, \theta_{\pm})$ consisting of a tame elliptic maximal torus $S$ of $G$ and an appropriate genuine character $\theta_{\pm}$ of $S(F)_{\pm}$ as in Definition \ref{dfn:nonsingularpairs}. We write $\pi_{(S,\theta_{\pm})}$ for the (potentially finite direct sum of irreducible) supercuspidal representation attached to $(S,\theta_{\pm})$.
		\item The $L$-parameter of (every irreducible constituent of) the supercuspidal representation $\pi_{(S,\theta_\pm)}$ is the composition of the $L$-parameter $W_F \to {^LS}_\pm$ of $\theta_\pm$ with the canonical $L$\-embedding ${^LS}_\pm \hookrightarrow {^LG}$.
		\item The normalized Harish-Chandra character of the supercuspidal representation $\pi_{(S,\theta_{\pm})}$ at shallow elements $s \in S(F)$ is given by 
		\[ e(G)\epsilon(X^*(T)_\C-X^*(S)_\C,\psi) \sum_{w \in N(S,G)(F)/S(F)} [a_S \cdot \theta_\pm](s^w), \]
		see Theorem \ref{thm:characterformula} for the definition of the notation.
	\end{enumerate}
\end{thm}

We emphasize that the quadratic character constructed in \cite{FKS23} is not used anywhere in the construction of $\cind_K^{G(F)}(\rho^0_{\pm} \otimes \kappa_{\pm})$, and no $\chi$-data is used in the assignment of the $L$-parameter. Quite the opposite: We show in this paper that these two objects emerge naturally when one compares the Fourier--Weil representation and the pull back of the Weil representation. In other words, this paper provides a conceptual explanation for the occurrence of the quadratic character of \cite{FKS23} and of the $\chi$-data used in \cite{KalRSP} and \cite{KalSLP} that are needed when one uses Yu's original construction. They are recognized as error terms that measure the difference between two naturally occurring linearizations of the projective Weil representation of $\Sp(\oV)$. 
Note that, unlike with the case of the Weil representation, which is trivial when $\oV=\{0\}$, the group $\WM_\gamma(\ol{V})$ is non-trivial even when $\oV=\{0\}$, namely $\WM_\gamma(\{0\})=\mu_8$, and the Fourier--Weil representation, being a genuine representation, is also necessarily non-trivial. This explains the curious fact that the quadratic character of \cite{FKS23} can be non-trivial even when no Weil representation occurs in Yu's original construction.

\subsection{Motivation}\label{sub:motivation}

The construction \cite{Yu01} of supercuspidal representations of tame $p$-adic groups by J.-K.~Yu and its extension \cite{KimYu} to types has been widely used in representation theory and number theory. It was noticed by Loren Spice that some arguments in \cite{Yu01} relied on a misprint in the literature and thereby the proof as written did no longer hold and in particular the crucial intertwining statements \cite[Proposition 14.1 and Theorem 14.2]{Yu01} turned out to be false in general, see \cite[\S4]{Fin21Comp}. This prompted one of us to find alternative arguments in \cite{Fin21Comp} to prove that nonetheless the construction of \cite{Yu01} does produce irreducible supercuspidal representations. While this restored the validity of the construction with a shorter proof, we were still left with the following problems.   \begin{enumerate}
	\item The computation of the Harish-Chandra characters of the resulting representations was severely hindered by the failure of the intertwining arguments of \cite{Yu01}, and in the few cases where the character could be computed, such as in \cite{AS09}, the resulting formula contained unexpected auxiliary sign characters, see the discussion in \cite[\S4.2]{KalRSP}.
	\item The connection to the Galois side was also severly hindered because the natural notions of semi-simplicity (or non-singularity) on both sides mismatched, see the introduction of \cite{KalSLP}.
	\item Grouping together representations along stable conjugacy classes of torus-character pairs produced incorrect $L$-packets, due to the failure of stability, see \cite{DS18}.
	\item Reduction-to-depth-zero, i.e., an equivalence between Bernstein blocks attached to the type constructed by Kim and Yu from the whole $G$-datum vs just from the depth-zero part of the $G$-datum, fails, see \cite[Appendix~A.2]{AFMO2}.
\end{enumerate}

This prompted the work \cite{FKS23}, in which a subtle quadratic character was defined and used to twist the construction of Yu \cite{Yu01}. The resulting ``twisted construction'' alleviated many of the above problems: It restored the validity of Yu's intertwining claims, which allowed the computation of the Harish--Chandra character, begun in \cite{AS09}, to be completed in almost full generality in \cite{Spice18} and \cite{Spice21}. The resulting formula, recorded for example in \cite[Theorem 4.3.5]{FKS23}, did not have unexpected auxiliary sign characters any more. Grouping representations along stable classes of torus-character pairs produced correct $L$-packets, see \cite{KalRSP} and \cite{KalSLP}, for which stability and many instances of endoscopic transfer could be proved, see \cite[\S4.4]{FKS23}. Reduction-to-depth-zero became possible, see \cite{AFMO2}, \cite{DF}.

Other problems however remained. The connection to the Galois side still remained mysterious. The reason is that a (tame regular elliptic) torus character pair $(S,\theta)$ for $G$ leads to an $L$-parameter $^L\theta : W_F \to {^LS}$ for the character $\theta$, but there is no distinguished choice for an $L$-embedding ${^Lj} : {^LS} \hookrightarrow {^LG}$ with which $^L\theta$ could be composed in order to obtain an $L$-parameter for the supercuspidal representation associated to $(S,\theta)$ by the (twisted or usual) construction of Yu after obtaining a supercuspidal $G$-datum from $(S,\theta)$ via the Howe factorization algorithm of \cite[\S3.6]{KalRSP}. This problem was resolved in \cite{KalRSP} and \cite{KalSLP} by a technical work-around: The arithmetic information contained in $\theta$ was used to produce in an ad-hoc way a device called $\chi$-data, which via a procedure of Langlands--Shelstad \cite[\S2.6]{LS87} led to an $L$-embedding $^Lj$. But the ad-hoc nature of this construction was esthetically and conceptually unsatisfactory, and also led to technical difficulties.

The same ad-hoc nature was also present in the Harish-Chandra character formula, since the $\chi$-data extracted from $\theta$ was used as an ingredient in that formula. The two issues -- $L$-parameter and character formula -- were clearly related, and from the work of \cite{KalRSP} there emerged a conceptual path that leads from an $L$-parameter to a (stable) character formula, and back, but the language to formulate this precisely was not available, and neither was a construction of representations that could effectuate it.

In addition, the introduction of the quadratic character of \cite{FKS23} created its own problems. The first was more philosophical, pertaining to the rather unclear conceptual origin of this character. While an explicit construction for it was given, and it seemed to resolve many problems, it was not at all clear what the conceptual motivation of the construction is, and why it should be the ``right'' one. This philosophical problem had concrete mathematical manifestations. For example, a great deal of effort in \cite{AFMO2} was expended to prove that this quadratic character can be extended from its original domain of definition to a potentially order-four character on a somewhat larger domain, in order for the reduction to depth-zero argument to proceed. 

The problem of non-existence of a natural $L$-embedding $^Lj : {^LS} \hookrightarrow {^LG}$, and its technical manifestation as non-canonical choices of $\chi$-data, had already been elegantly resolved in the setting of real groups by work of Adams--Vogan \cite{AV92}, who showed that if one instead considers a certain topological double cover $S(\R)_\rho$ of $S(\R)$, then one can define an analog $^LS_\rho$ of an $L$-group for it, prove a local Langlands correspondence between genuine characters of $S(\R)_\rho$ and $L$-parameters valued in $^LS_\rho$, and obtain a distinguished $L$-embedding $^Lj_\rho : {^LS}_\rho \hookrightarrow {^LG}$. A regular genuine character $\theta_\rho$ of $S(\R)_\rho$ leads  to an irreducible discrete series representation $\pi_{\theta_\rho}$, whose  $L$-parameter is obtained by composing the $L$-parameter of $\theta_\rho$ with $^Lj_\rho$. We emphasize that the distinguished $L$-embedding $^Lj_\rho$ depends only on $S$ and $G$, and not on $\theta_\rho$, unlike the ad-hoc embeddings $^Lj : {^LS} \hookrightarrow {^LG}$ used in \cite{KalRSP} and \cite{KalSLP}. Moreover, the Harish-Chandra character formula of $\pi_{\theta_\rho}$ can be stated naturally in terms of $\theta_\rho$. We refer to \cite{AV16} for a modern discussion of these matters.

The idea of using double covers of tori in the $p$-adic setting was explored by Moshe Adrain, who gave in \cite{Adrian13} a clear and compelling treatment of $G=\tx{PGL}_2$ and explored it further with Lansky in \cite{AdrianLansky2015} in the setting of $\tx{PGSp}_4$.

Motivated by these ideas, one of us gave in \cite{KalDC} a uniform construction of topological double covers $S(F)_\pm$ of tori for an arbitrary local field $F$, defined their $L$-groups $^LS_\pm$, proved a local Langlands correspondence between $L$\-parameters and genuine characters, and constructed distinguished $L$\-embeddings $^Lj_\pm : {^LS}_\pm \hookrightarrow {^LG}$. Unlike the construction of $S(\R)_\rho$, which is algebraic, the construction of $S(F)_\pm$ is arithmetic, but nonetheless it produces $S(\R)_\rho$ when $F=\R$ and $S$ is elliptic. This construction provided the language to formally express the clear path between discrete semi-simple $L$-parameters and stable Harish-Chandra characters, uniformly for all local fields, in \cite[Conjecture 4.4.1]{KalDC}. 

This strongly suggested that the natural parameterization of the regular supercuspidal representations of \cite{KalRSP}, and more generally the semi-simple supercuspidal representations of \cite{KalSLP}, 
should be by pairs $(S,\theta_\pm)$ of an elliptic maximal torus $S \subseteq G$ and a genuine character $\theta_\pm : S(F)_\pm \to \C^\times$, rather than pairs $(S,\theta)$ with $\theta : S(F) \to \C^\times$. Extrapolating from this we were led to the expectation that Yu's construction will be more naturally phrased in terms of supercuspidal $G$-data that involve genuine characters of the double covers of twisted Levi subgroups introduced in \cite[\S6]{KalDC}. This is exactly what we achieve in this paper, as we outline in \S\ref{sub:introoverview}.

\subsection{\texorpdfstring{The special case of $S \subset G$ vs.\ the setting of real Lie groups}{}}

We outlined in \S\ref{sub:introoverview} that our revised input for the construction of supercuspidal representations is a supercuspidal $G$-datum 
\[((G_i)_{1\leq i\leq n+1},x,(r_i)_{1\leq i\leq n},\rho^0_{\pm},(\phi_i)_{1\leq i\leq n}) \]
and that the key technical input in our work consists in lifting the map $\bar\jmath : G^0_x \ra \Sp(\ol{V}) $ from \eqref{eqn:G0xtoSp} to a map $\tilde\jmath :  G^0_{\pm,x} \ra \WM_\gamma(\ol{V})$ in \eqref{eqn:G0xtoWM}. 
We will now explain the construction of supercuspidal representations including the construction of the map $\tilde\jmath$ in more detail in the special case of a supercuspidal $G$-datum with $n=1$ and where $G=G_1 \supsetneq G_2=G^0=:S$ is a tame elliptic maximal torus. In this case $\rho^0_{\pm}$ and $\phi_1$ are both characters of $S(F)_{\pm}$ (one of them genuine and the other a character of $S(F)$ that we view as a character of $S(F)_{\pm}$ via $S(F)_{\pm} \twoheadrightarrow S(F)$). It turns out that the resulting supercuspidal representation only depends on the product $\theta_{\pm}:=\rho^0_{\pm} \otimes \phi_1$ (see Theorem \ref{thm-screps} for such a statement for general supercuspidal $G$-data). Thus, in this special case, our starting point is a tame elliptic maximal torus $S \subset G$ and a genuine generic character $\theta_{\pm} : S(F)_{\pm} \to \C^\times$ of depth $r=r_1>0$. We denote by $x$ the point in the Bruhat--Tits building of $G$ associated to $S$. 

We first review the construction of Adler, which Yu generalizes, i.e., we first consider a generic character $\theta : S(F) \to \C^\times$ of depth $r>0$ instead of a genuine generic character of $S(F)_{\pm}$. The construction of Adler assigns to the pair $(S,\theta)$ the representation $\pi_{(S,\theta)}=\tx{c-Ind}_K^{G(F)}\tilde\rho$, where $K=S(F)G(F)_{x,r/2}$ and $\tilde\rho$ is an irreducible representation of $K$ constructed as follows. There is a straightforward way to extend $\theta|_{S(F)_{r/2+}}$ to a character $\hat\theta$ of $S(F)G(F)_{x,r/2+}$. If the latter equals $S(F)G(F)_{x,r/2}=K$, then we can take $\tilde\rho=\hat \theta$. If $G(F)_{x,r/2+} \subsetneq G(F)_{x,r/2}$, then one can easily see that if we tried to take as $\tilde\rho$ the induction of $\hat\theta$, the result would be reducible. In order to obtain an irreducible supercuspidal representation we must therefore first find a way to ``extend'' $\hat \theta$ to a representation of a larger compact subgroup before inducing.
Via the theory of Heisenberg representations, there is a unique irreducible representation $\mathrm{Heis}$ on $G(F)_{x,r/2}$ whose restriction to $S(F)_{r/2}G(F)_{x,r/2+}$ is $\hat \theta$ (times the identity). In order to extend this representation to a representation of $S(F)G(F)_{x,r/2}$, one considers the  exact sequence
\[ 1 \to S(F)_{x,r/2}G(F)_{x,r/2+} \to G(F)_{x,r/2} \to \oV \to 0 \]
and notices that $\oV$ is an $\bF_p$-vector space that carries a symplectic form obtained by composing the commutator map on $G(F)_{x,r/2}$ with the homomorphism $\hat \theta : G(F)_{x,r/2+} \to \bC^\times$. This composition factors through $\mu_p$ and by identifying $\mu_p$ with $\bF_p$ (via $\psi$), we obtain a symplectic form on $\oV$ valued in $\bF_p$. 
The action of $S(F)$ on $G(F)_{x,r/2}$ by conjugation leads to a map 
\begin{equation} \label{eq:ad1}
	\bar\jmath : S(F) \rightarrow \Sp(\ol{V}).
\end{equation}
We write $\mathrm{LWeil}$ for the Weil representation of $\Sp(\ol{V})$ that acts on the space of $\mathrm{Heis}$. Then it turns out that the representation $(\bar\jmath^*\mathrm{LWeil}\ltimes \mathrm{Heis})\otimes (\theta \ltimes 1)$ of $S(F) \ltimes G(F)_{x,r/2}$ factors through the surjection $S(F) \ltimes G(F)_{x,r/2} \twoheadrightarrow S(F)G(F)_{x,r/2}=K$, i.e., provides a representation of $K$ that we take as $\tilde \rho$.

Our modification of this construction is motivated by its structural similarity to the construction of \cite{Duflo82} of representations of arbitrary real Lie groups. If $G$ is a reductive Lie group, $S \subset G$ an elliptic maximal torus, and $\theta_\pm : S_\pm \to \C^\times$ a regular character of its double cover, Duflo pulls $\theta_\pm$ back under the distinguished splitting $S_\mathrm{sc} \to S_\pm$, where $S_\mathrm{sc}$ is the preimage of $S$ in the algebraic simply connected cover of the derived subgroup of $G$, lets $\pi_\mathrm{sc}$ be the discrete series representation of $G_\mathrm{sc}$ assigned by Harish-Chandra, and considers the exact sequence
\[ 1 \to \{\pm1\} \times S_\mathrm{sc}  \to S_\pm \ltimes G_\mathrm{sc}  \to G. \]
Duflo equips the vector space $V=\mathrm{Lie}(G)/\mathrm{Lie}(S)$ with the symplectic form that is the Lie bracket composed with the differential of $\theta_\pm$, produces a distinguished  genuine homomorphism 
\begin{equation} \label{eq:du1}
S_\pm \to \Mp(V)	
\end{equation}
to the metaplectic group $\Mp(V)$,
extends $\pi_\mathrm{sc}$ to a genuine representation $\Pi$ of $S_\pm \ltimes G_\mathrm{sc}$, and shows that $\Pi\otimes\theta_\pm$ kills the kernel of the map $S_\pm \ltimes G_\mathrm{sc}  \to G$ and hence descends to its image, and induces irreducibly to $G$.

Returning to the setting of $p$-adic groups, Duflo's construction suggests to look for a genuine homomorphism from $S(F)_\pm$ to a metaplectic group. But the following questions and observations arise.

\begin{enumerate}
	\item Should we use the  symplectic space $V=\mathrm{Lie}(G)/\mathrm{Lie}(S)$ over the $p$-adic field $F$, or the symplectic space $\oV=G(F)_{x,r/2}/(S(F)_{x,r/2}G(F)_{x,r/2+})$ over the finite residue field $\rf$ (or its restriction to $\bF_p$)?
	\item What is the meaning of ``metaplectic group'', especially in the setting of the space $\oV$ over $\rf$?
	\item If we consider the example of $G=\tx{PGL}_2$, then the construction of $\pi$ by Adler (and Yu) involves a Heisenberg--Weil representation only when $S$ is unramified, but in this case the double cover $S(F)_\pm$ splits canonically. On the other hand, $S(F)_\pm$ can often be non-split when $S$ is ramified, but in that case we have $G(F)_{x,r/2+}=G(F)_{x,r/2}$, so no Weil representation is used. At the same time, the character of \cite{FKS23} can be non-trivial for both $S$ ramified and unramified.
\end{enumerate}

These points show that the task of combining the cover $S(F)_\pm$ with the Heisenberg--Weil representation in the construction of Adler and Yu is subtle. What came as a surprise to us is that its eventual resolution, obtained after much confusion and error, is strikingly simple.

\subsection{The Witt-metaplectic group and the Fourier--Weil representation}

While the construction of supercuspidal representations of Adler and Yu requires the use of $\oV$, to connect to $S(F)_\pm$ we need to simultaneously employ $V$ and $\oV$. The version of the metaplectic group that we use is a non-trivial extension of $\Sp$ even over the finite residue field $\rf$, and can be non-trivial even when $\oV=\{0\}$. More precisely, let $k$ be any locally compact field, thus either a finite or a local field,  of characteristic not $2$. Let $\QW(k)$ denote the Witt group of equivalence classes of non-degenerate quadratic $k$-vector spaces, where equivalence is given by adding or removing hyperbolic planes. For any symplectic $k$-vector space $(W,\eta)$ a natural extension
\[ 0 \to \QW(k) \to \WMuniv(W) \to \Sp(W) \to 1 \]
and a representation $\mathrm{FWeilHeis}_\psi$ of $\WMuniv(W)$ that lifts the projective Weil representation of $\Sp(W)$ emerges from the work \cite{LionVergne} of Lion--Vergne over the real numbers, with some hints of a $p$-adic theory due to Patrice Perrin explained in the appendix of this reference, and in fuller force in \cite{RR93}, \cite{Thomas2006}, \cite{Thomas2008}, \cite{Thomas2013}.

The group $\WMuniv(W)$ is defined as the twisted product of $\QW(k)$ and $\Sp(W)$ with respect to a 2-cocycle arising from the Maslov index (in the body of the paper we will actually find a certain quotient $\QW(k)/I_k^3$ more convenient, but this is not significant). The latter is an invariant $\tau(L_1,L_2,L_3) \in \QW(k)$ assigned by Kashiwara to a triple of Lagrangian subspaces of $W$. If we fix a Lagrangian subspace $L \subset W$, then the function $z_L(g,h)=\tau(L,gL,ghL)$ is an element of $Z^2(\Sp(W),\QW(k))$ and hence determines an extension $\WM_L(W)=\QW(k)\boxtimes_{z_L}\Sp(W)$ as above, depending on the choice of $L$. Let $H(W)=k \boxtimes_{\eta/2} W$ be the  Heisenberg group. The Heisenberg representation with central character $\psi : k \to \C^\times$ can be presented as the induction $\tx{Heis}_{\psi,L}=\tx{Ind}_{k \times L}^{H(W)}(\psi\boxtimes \mathbf{1})$. For each $g \in \Sp(W)$ one can define an $H(W)$-equivariant automorphism $R_g : \tx{Heis}_{\psi,L} \to \tx{Heis}_{\psi,L}\circ g$ using a partial Fourier transform, equivalently an intertwining operator of a shape familiar from Harish-Chandra's theory of parabolic induction. A result of Lion--Vergne over $\R$, extended by Perrin to non-archimedean local fields and Rao to general locally compact fields, is that $R_{g_1}R_{g_2}=\tilde\gamma_\psi(z_L(g_1,g_2))R_{g_1g_2}$, where $\tilde\gamma_\psi(q):=\gamma_\psi(q/2)$ is the \emph{renormalized} Weil index, showing that the operators $R_g$ promote $\tx{Heis}_{\psi,L}$ to a representation $\mathrm{FWeilHeis}_{\psi,L}$ of $\WM_{L,\psi}(W) \ltimes H(W)$ that lifts the projective Heisenberg--Weil representation of $\Sp(W) \ltimes H(W)$; here $\WM_{L,\psi}(W)$ is the pushout of $\WM_L(W)$ under the renormalized Weil index $\tilde\gamma_\psi$, which is a homomorphism of $\QW(k)$ into the group of roots of unity of $\C^\times$ of a suitable order (order $8$ always works, but over non-archimedean or finite fields order $4$ or $2$ is enough depending on whether the size $q$ of the (residue) field is congruent to $3$ or $1$ modulo $4$, assuming it is odd).

The next fundamental idea was brought in by Thomas, who showed how to construct an extension $\WMuniv(W)$ of $\Sp(W)$ by $\QW(k)$ without choosing a Lagrangian subspace of $W$. Namely, for each $g \in \Sp(W)$ the graph $\Gamma_g \subset W \oplus W$ of $g : W \to W$ is a Lagrangian subspace of the doubled symplectic space $W \oplus W$ with form $(-\eta) \oplus \eta$. One can consider the element $z \in Z^2(\Sp(W),\QW(k))$ defined by $z(g,h)=\tau(\Gamma_1,\Gamma_g,\Gamma_{gh})$ and use it to construct $\WMuniv(W)$. An important observation of \cite{Thomas2013} is that one has a system of isomorphisms $\alpha_L : \WMuniv(W) \to \WM_L(W)$ and that pulling back $\mathrm{FWeilHeis}_{\psi,L}$ under $\alpha_L$ produces a representation $\mathrm{FWeilHeis}_\psi$ of $\WM_\psi(W)$ that is independent (up to isomorphism) of $L$, and its character formula becomes even cleaner and simpler than that of $\mathrm{FWeilHeis}_{\psi,L}$: for $(q,g) \in \WMuniv(W)$ with $\tx{ker}(g-1|W)=\{0\}$, the character of this pull-back is simply $\tilde\gamma_\psi(q)\cdot |\det(g-1)|^{-1/2}$.

As mentioned above, we will apply this theory to both $k=F$ and $k=\rf$ or $\bF_p$. An important fact is that there exists a reduction homomorphism $\pi_0 : \QW(F) \to \QW(\rf)$, which has the property $\gamma_\psi = \gamma_{\bar\psi}\circ\pi_0$ for any character $\psi$ of $F$ whose restriction to the ring of integers factors through a non-trivial character $\bar\psi$ of $\rf$, see \cite[Appendix A]{RR93}.

The representation $\mathrm{FWeilHeis}_\psi=\alpha_L^*\mathrm{FWeilHeis}_{\psi,L}$ is of central importance in our construction. Therefore, we spend some effort in the Appendix to present it in a clean and uniform way and to simplify and generalize some results in \cite{LionVergne}, \cite{RR93}, and \cite{Thomas2013}, especially regarding its relationship to the more classical representation $\mathrm{LWeilHeis}_\psi$ of $\Sp(W) \ltimes H(W)$ studied over finite fields by Gerardin \cite{Gerardin} and used in the construction of Adler and Yu.

\subsection{Statements of main results}

We now come to the main results of this paper in the special case discussed above. While these concern the representation theory of $p$-adic groups, they are applications of results about the interaction of the Fourier--Weil representation for symplectic groups over finite fields with the double covers of tori, or more generally twisted Levi subgroups, of a $p$-adic group. Our key technical result is purely in the setting of the symplectic group and the Fourier--Weil representation.

\begin{thm}[Imprecise version of a special case of Theorem \ref{thm:symplectic-main}]\ \\
	Let $F$ be a non-archimedean local field of residue characteristic not $2$, $(V,\eta)$ a symplectic vector space over $F$, $T \subset \Sp(V)$ a maximal torus, $\tilde T \subset \GL(V)$ its centralizer.
	\begin{enumerate}
		\item There is a surjective norm map $\tilde T(F) \twoheadrightarrow T(F)$ that factors as $\tilde T(F) \to T(F)_\pm \to T(F)$, where $T(F)_\pm$ is the double cover of $T(F)$ associated to the set $\Lambda$ of weights of $T$ acting on $V$.
		\item For $\tilde s \in \tilde T(F)$ define the quadratic form $q_{\tilde s} : V \to F$ by $q_{\tilde s}(v)=\eta(\tilde sv,v)$. There is a unique genuine function
		\[ T(F)_\pm \to \mu_8,\qquad \dot s\mapsto \gamma_\psi(q_{\dot s}), \]
		whose value is $\gamma_\psi(q_{\tilde s})$ whenever $\tilde s$ maps to $\dot s$ and
		$\tilde s_\lambda=1$ for every $\lambda \in \Lambda$ with $\lambda(s)=1$, as in Lemma \ref{lem:topssgamma}(5).
		\item Let $\mc{L}$ be a $\Lambda$-atomic lattice chain in $V$ with $\Lambda$-central quotient $(\oV,\bar\eta)$ (see Definition \ref{dfn:lambda_atomic} for the definitions). The natural homomorphism 
		$\bar\jmath : T(F)_\mc{L}/T(F)_{0+} \to \Sp(\oV)$ admits a unique genuine lift 
		\[\tilde\jmath : T(F)_{\mc{L},\pm}/T(F)_{0+} \to \WM_{\bar\psi}(\oV) \] which for topologically semi-simple elements $\dot s \in T(F)_{\mc{L},\pm}$ is given by $\tilde\jmath(\dot s) = \gamma_\psi(q_{\dot s}) \boxtimes \bar\jmath(s)$,  with $\gamma_\psi(q_{\dot s})$ as in (2).
		\item The representations $\tilde\jmath^*\mathrm{FWeil}_\psi$ and (the pull back of) $\bar\jmath^*\mathrm{LWeil}_\psi$ are twists of each other by a genuine character of $T(F)_\pm$ that is the product of the spinor norm for the action of $T(F)$ on the $\Lambda$-central quotient of an appropriate scaling 
of $\mc{L}$ and a character of $T(F)_\pm$ arising from $\chi$-data (i.e., certain Gauss sums).
	\end{enumerate}
\end{thm}

We now apply this result in the case where $V=\mathrm{Lie}(G)/\mathrm{Lie}(S)$ equipped with the symplectic form $\eta = X \circ [ \phantom{X}, \phantom{X}]$, where $X \in \mathrm{Lie}^*(S) \subset \mathrm{Lie}^*(G)$ is a good element such that $\psi \circ X$ corresponds to $\theta_{\pm}|_{S(F)_{r}}$ under the isomorphism $\mathrm{Lie}(S)_r/\mathrm{Lie}(S)_{r+} \simeq S(F)_r/S(F)_{r+}$, $T \subset \Sp(V)$ is the maximal torus associated to the weight decomposition of $V$ under $S$, and $\ol{V}$ is the space $G(F)_{x,r/2}/(S(F)_{x,r/2}G(F)_{x,r/2+})$ defined above.
The group $\WM_\psi(\oV)$ injects into $\WM_\gamma(\oV)$ and the resulting genuine homomorphism $\tilde\jmath : S(F)_\pm \to \WM_\psi(\oV) \to \WM_\gamma(\oV)$ is rather analogous to the homomorphism \eqref{eq:du1} and allows us to upgrade the homomorphism \eqref{eq:ad1} of the construction by Adler and Yu to a homomorphism
\begin{equation} \label{eq:ad2}
	\tilde\jmath : S(F)_\pm \to \WM_\gamma(\oV).
\end{equation}
The tensor product of the two genuine representation $(\tilde\jmath^*\mathrm{FWeil} \ltimes \mathrm{Heis}) \otimes (\theta_\pm \ltimes 1)$ of $S(F)_{\pm} \ltimes G(F)_{x,r/2}$ then descends to $K=S(F)G(F)_{x,r/2}$ and its compact induction is an irreducible supercuspidal representation that we denote by $\pi_{(S,\theta_\pm)}$. From the above theorem we deduce the following result.

\begin{thm}
	\begin{enumerate}
		\item Let $\chi''_S : S(F)_\pm \to \C^\times$ be the genuine character associated to the $\chi$-data defined in \cite[\S4.2]{FKS23} and let $\epsilon : S(F) \to \C^\times$ be the quadratic character of \cite[Definition 5.5.8]{FKS23}. Letting $\theta : S(F) \to \C^\times$ be the character determined by $\theta_\pm = \theta \cdot \chi''_S \cdot \epsilon$, we have 
		\[ \pi_{(S,\theta_\pm)}=\pi_{(S,\theta)}. \]
		\item The normalized Harish--Chandra character of $\pi_{(S,\theta_\pm)}$ at shallow elements of $S(F)$ is given by 	\begin{equation*} \label{eq:char}
			e(G)\epsilon(X^*(T)_\C-X^*(S)_\C,\psi) \sum_{w \in N(S,G)(F)/S(F)} [a_S \cdot \theta_\pm](s^w).		
		\end{equation*}
		\item The $L$-parameter of $\pi_{(S,\theta_\pm)}$ is the composition of the $L$-parameter of $\theta_\pm$ and the canonical $L$-embedding $^LS_\pm \hookrightarrow {^LG}$.
	\end{enumerate}
\end{thm}

We emphasize that neither the $\chi$-data $\chi''$ nor the character $\epsilon : S(F) \to \{\pm1\}$, both defined in \cite{FKS23}, are added by hand into the construction of \eqref{eq:ad2} and hence of $\pi_{(S,\theta_\pm)}$. Rather, \eqref{eq:ad2} is uniquely specified by the Ansatz
\[ \dot s \mapsto \gamma_\psi(q_{\tilde s}),\qquad q_{\tilde s}(v)=\eta(\tilde sv,v) . \]
Instead, both $\chi''$ and $\epsilon$ arise naturally when one compares the Fourier--Weil representation $\tilde\jmath^*\mathrm{FWeil}$ that we use to construct $\pi_{(S,\theta_\pm)}$ and the Weil representation $\bar\jmath^*\mathrm{LWeil}$ that is used in the construction of $\pi_{(S,\theta)}$ by Adler and Yu.

Thus, one of the contributions of this paper is to show that the approach via $\tilde\jmath$ and the Fourier--Weil representation leads to a parameterization of supercuspidal representations that is not only closer in spirit to the setting of real groups, but also considerably simplifies many applications, such as character computations and Langlands parameters, and furthermore gives a conceptual explanation for both $\chi''$ and $\epsilon$ as artifacts arising from the use of the linearization of the Weil representation over the symplectic group that exists over finite fields.

Additionally, it should be noted that there is a mistake in the published version of \cite{FKS23}, which is relevant for this paper. Namely, the sign character of \cite[Definition 5.5.7]{FKS23}, which is a factor of the sign character of \cite[Definition 5.5.8]{FKS23}, ought to be defined using the spinor norm relative to the quadratic form $2\varphi$, rather than $\varphi$, in order for the proofs of \cite[Proposition 5.5.9, Theorem 3.4]{FKS23}, to go through. This is because a factor of $2$ is mistakenly dropped in the proof of \cite[Proposition 5.5.9]{FKS23}, during the calculation of $\varphi(Y_O)$. There is a further mistake in \cite[Remark 5.5.10]{FKS23}, in which a quadratic form $\varphi_s$ is incorrectly defined and therefore fails to be rational. Both of these mistakes, and their fixes, are explained in the erratum on the website of the second author. After the fixes are preformed, the sign character of \cite[Definition 5.5.8]{FKS23} can be consistently defined using spinor norms of $2\varphi$ and $2\varphi_s$. In the above sketched theorem, as well as in Theorems \ref{thm:main1} and \ref{thm:main1b}, we have to use the corrected version of \cite[Definition 5.5.8]{FKS23}.

\subsection{The structure of this paper}

We now give a brief overview of the structure of this paper and provide a quick summary of its more technical aspects. The reader most interested in the construction of supercuspidal representations is encouraged to begin reading \S\ref{sec:supercuspidal} and refer to the other sections as the need arises.

In \S\ref{sec:symp} we prove the main technical results about the Fourier--Weil representation that underlie the difference between our construction and that of Yu. In \S\ref{sub:polar} and \S\ref{sub:lambda} we recall basic facts about twisted Levi subgroups of unitary groups over general fields and establish notation that will be used throughout the section. In \S\ref{sub:lattice} we specialize to non-archimedean local fields and introduce the concept of an \emph{atomic} lattice chain and its generalization -- a \emph{$\Lambda$-atomic} lattice chain -- in an $\epsilon$-Hermitian vector space over a non-archimedean local field. Such a lattice chain has a unique ($\Lambda$-)central quotient, which is a vector space over the residue field equipped with a pairing that is either skew-Hermitian, symplectic, or symmetric, depending on the situation. In \S\ref{sub:qs} we introduce the quadratic form $q_{\tilde s}$ that is the most fundamental object of our construction, and was already mentioned above. In \S\ref{sub:sympmain} we prove the main result of this section, Theorem \ref{thm:symplectic-main}, which establishes the existence and uniqueness of the lift $\tilde j$ in the sympletic setting, as well as the comparison of the pull-backs of the two models of the Weil representation -- the Fourier model and the linear model. The proof of this theorem is reduced by a devisage argument, provided by Proposition \ref{pro:orangebox}, to the atomic case of a twisted Siegel Levi subgroup treated in Lemma \ref{lem:tsl-gammatilde}.

This lemma contains the key calculation. To explain its statement we need to recall some basic facts about the Weil representation. The fact that the projective Weil representation linearizes over $\Sp(\ol{V})$ when $\ol{V}$ is a symplectic vector space over the finite field $\rf$ is closely linked to the fact that the universal Witt-metaplectic group, which is an extension 
\[ 0 \to \QW(\rf) \to \WM(\ol{V}) \to \Sp(\ol{V}) \to 1, \]
splits, thereby realizing $\Sp(\ol{V})$ as a subgroup of $\WM(\ol{V})$. This splitting is achieved by using orientations on Lagrangian spaces, and is realized by a 1-cochain $c : \Sp(\ol{V}) \to \QW(\rf)$, see \S\ref{sub:classmeta}. 

Returning to Lemma \ref{lem:tsl-gammatilde}, its statement compares the Weil indices of $q_{\tilde s}$ and $c(s)$ for a topologically semi-simple element $s$. Here $q_{\tilde s}$ is defined in terms of the symplectic $F$-vector space $V$, while $c(s)$ is defined in terms of a symplectic $\rf$-vector space $\ol{V}$. The two are related as follows. The twisted Levi subgroup comes from a skew-Hermitian structure on $V$ whose trace is the given symplectic structure. We are furthermore given an atomic lattice chain in this skew-Hermitian space that is preserved by $s$. Then $\ol{V}$ is the central quotient of this lattice chain. The lemma shows that the ratio $\gamma_\psi(q_{\tilde s})/\tilde\gamma_{\bar\psi}(c(s))$ equals the product of two terms. The first is the value at $\dot s$ of a genuine character of the double cover of the twisted Siegel Levi subgroup, defined in terms of $\chi$-data. The other is the spinor norm of the natural quadratic quotient of the atomic lattice chain.

In \S\ref{sec:symplectictopadic} we apply the main technical theorem \ref{thm:symplectic-main} to the setting of Yu's construction. The symplectic vector space is now $V=\tx{Lie}(G)/\tx{Lie}(G^0)$ as discussed above. The twisted Levi subgroup of $\tx{Sp}(V)$ is the one corresponding to the weight decomposition of $V$ under the center $Z_{G^0}$ of $G^0$. The $\Lambda$-atomic lattice chain comes from the Moy--Prasad filtration. We then show in Theorem \ref{thm:main1} that, in the one-step case, the $\chi$-data that occurs in the comparison result of Lemma \ref{lem:tsl-gammatilde}, hence of Theorem \ref{thm:symplectic-main}, is exactly the $\chi$-data specified in \cite[\S4.2]{FKS23}, and the spinor norm of the quadratic quotient of the lattice chain is exactly the sign character defined in \cite{FKS23}. This result is then iterated to accommodate in Theorem \ref{thm:main1b} the general multi-step set-up of Yu's construction.

In \S\ref{sec:supercuspidal} we apply the results of \S\ref{sec:symp} and \S\ref{sec:symplectictopadic} to the construction of supercuspidal representations and types, in the manner already explained earlier in the introduction. The constructions are spelled out in \S\ref{sub:construction} and \S\ref{sub:types}. In \S\ref{sub:Yucomparison} and \S\ref{sub:types} we relate our construction here to the construction of Yu (\cite{Yu01}) and Kim and Yu (\cite{KimYu}) twisted by the sign character of \cite{FKS23}. In \S\ref{sub:alternativeinput} we show that one can equally well use a slightly different version of a genuine Yu datum, namely 
\[((G_i)_{1\leq i\leq n+1},x,(r_i)_{1\leq i\leq n},\rho^0,(\phi_{i,\pm})_{1\leq i\leq n}), \]
where now $\rho^0$ is a representation of $G^0_x$ as in the original setting of Yu, rather than a genuine representation of $G^0_{x,\pm}$, but instead each $\phi_{i,\pm}$ is a genuine generic character of the double cover $G_{i+1}(F)_\pm$ relative to the ambient group $G_i$. This slight variation is sometimes more convenient, for example in the discussion of the one-step situation (where having to choose an auxiliary, necessarily non-trivial, $\rho^0_\pm$, would be awkward), or for the Howe factorization. In \S\ref{sub:nonsingular} we discuss the setting of the semi-simple supercuspidal representations introduced in \cite{KalSLP}, and show that these are in bijection with $G(F)$-conjugacy classes of pairs $(S,\theta_\pm)$ of a maximal torus $S \subset G$ and a genuine character $\theta_\pm : S(F)_\pm \to \C^\times$ subject to certain properties. In \S\ref{sub:char} we derive the character formula for non-singular supercuspidal representations at shallow elements. While this formula follows formally from the preceding discussion and the results of \cite{KalRSP}, the point we are making here is that our setting provides this formula much more simply, and in particular avoids all discussions of roots of unity that were among the most subtle parts of earlier work. In \S\ref{sub:llc} we describe the $L$-parameters of semi-simple supercuspidal representations. Here we do rely heavily on earlier work for the proof, but the point we are making is that the construction of the $L$-parameter is completely straight-forward. In \S\ref{sub:pgl2} we consider the example of the group $G=\tx{PGL}_2$ and show that our construction produces the same output as the construction of supercuspidal representations of $\tx{PGL}_2(F)$ by means of the Weil representation, which was introduced in \cite{JL79} (in the setting of $\tx{GL}_2$, but we find the case of trivial central character most illustrative), and reviewed in \cite[\S39]{BH06}.

Finally, in Appendix \S\ref{app:fw} we review the concepts of the (universal) Witt-metaplectic group and the Fourier--Weil representation. While most of the definitions and results are already present in the literature, not all of them are stated in the desired generality, so we take this opportunity to collect and present them in a form that is convenient for us. Some minor results may be new, such as the discussion of graded bilinear forms in \S\ref{sub:gradedbf} and the resulting explicit description of the splitting of the symplectic group into the universal Witt metaplectic group over unipotent elements given in Proposition \ref{pro:metasection}. These are needed in the proof of the key technical lemma \ref{lem:tsl-gammatilde}.

\subsection{Acknowledgments}

We started this project during our last visit to our mentor Dick Gross. We are very saddened by his passing and dedicate this paper to his memory, his ever-present enthusiasm for mathematics, and his generous mentorship.

\subsection{AI disclaimer}

This project began in January 2024 and was completed largely without use of AI. In the final months before completion, the second-named author used AI tools to perform some calculations and suggest simplifications of some initially longer arguments, mostly in the Appendix. AI was also used for literature search and proof-reading. No ideas (beyond simplifying already existing arguments and performing calculations) or text were generated by AI.

The authors take full responsibility for the content of the paper.

\section{Symplectic theory} \label{sec:symp}

In the first two subsections $F$ is an arbitrary field of  characteristic not $2$. Starting with \S\ref{sub:lattice}, $F$ is a  non-archimedean local field, $\cO$ its ring of integers, $\rf$ the residue field, which we assume to be of odd characteristic.

\subsection{Polarizations} \label{sub:polar}

Let $V$ be an even-dimensional $F$-vector space, $\eta : V \times V \to F$ a non-degenerate $\epsilon$-symmetric bilinear form for $\epsilon \in \{\pm1\}$, and $G:=\SU(V,\eta):=\Aut(V,\eta)^\circ$.
Thus $G=\Sp(V,\eta)$ when $\epsilon=-1$ and $G=\SO(V,\eta)$ when $\epsilon=+1$.

\begin{enumerate}
	\item A \emph{split polarization}\index[terminology]{polarization!split polarization} is an unordered pair of complementary Lagrangian subspaces of $V$. Their joint stabilizer (the subgroup preserving each member) is a Siegel Levi subgroup $M \subset G$ and every Siegel Levi subgroup of $G$ arises in this way from a unique split polarization. If we order the pair as $\{V_+,V_-\}$, then we obtain the isomorphisms
	\[ V_+ \oplus V_- \simto V,\quad V_+ \simto V_-^*,\quad \GL(V_+) \simto M. \]
	\item A \emph{twisted polarization}\index[terminology]{polarization!twisted polarization} consists of a quadratic field extension $E/F$ and a split polarization of $(V_E,\eta_E)$ that is $\Gamma_{E/F}$-stable but neither of its members are; here $V_E=V\otimes_F E$ and $\eta_E$ is the $E$-linear extension of $\eta$ to $V_E$, so that $(V_E,\eta_E)$ is an $\epsilon$-symmetric space over $E$. The joint stabilizer in $G$ of the polarization of $(V_E,\eta_E)$ 
is a non-split twisted Siegel Levi subgroup $M$ of $G$ and every non-split twisted Siegel Levi subgroup of $G$ arises this way from a unique twisted polarization. If we order the pair as $\{V_+^E,V_-^E\}$, then $V_-^E=\sigma(V_+^E)$, where $\sigma$ is the unique non-trivial automorphism of $E/F$. We obtain the $\epsilon$-Hermitian form
	\[ h_1 : V_+^E \times V_+^E \to E,\qquad h_1(v,w)=\eta_E(v,\sigma(w)).  \]
	The map
	\[ V_+^E \to V_E,\qquad \ell \mapsto \ell+\sigma(\ell) \]
	takes values in $V_E^\sigma=V$ and induces an isomorphism
	\[ V_+^E|_F \simto V \]
	that identifies $\eta$ with the trace of $h_1$ and $\U(V_+^E,h_1)$ with $M$.
	\item Returning to a split polarization with chosen order $\{V_+,V_-\}$, we can see structure similar to that of a twisted polarization by introducing the split quadratic \'etale $F$-algebra $E=F \oplus F$. Its unique non-trivial automorphism is $\sigma(x,y)=(y,x)$. We form again the $\epsilon$-symmetric space $(V_E,\eta_E)$ over $E$ by taking $V_E = V\otimes_F E=V \oplus V$, where $F$ is embedded diagonally into $F \oplus F = E$, and the $E$-linear extension $\eta_E : V_E \times V_E \to E$ of $\eta$, which one checks is given by $\eta_E((v_1,v_2),(w_1,w_2))=(\eta(v_1,w_1),\eta(v_2,w_2))$. If we let $V_+^E=V_+ \oplus V_-$ and $V_-^E=V_- \oplus V_+$, then $\{V_+^E,V_-^E\}$ is a pair of complementary Lagrangian subspaces of $V_E$ and $\sigma(V_+^E)=V_-^E$. The map
	\[ h_1 : V_+^E \times V_+^E \to E,\quad h_1(v,w) = \eta_E(v,\sigma(w))=(\eta(v_1,w_2),\eta(v_2,w_1)) \]
	is an $\epsilon$-Hermitian form with respect to the \'etale algebra $E/F$. The map
	\[ V_+^E \to V_E,\qquad \ell \mapsto \ell+\sigma(\ell) \]
	takes values in 
	\[ V_E^\sigma = \{(v,v) \in V_E|v \in V\} \cong V \]
	and induces an isomorphism
	\[ V_+^E|_F \simto V \]
	of $F$-vector spaces, which identifies $\eta$ with the trace of $h_1$, as one sees from 
	\begin{eqnarray*}
		\eta(x_1+x_2,y_1+y_2)&=&\eta(x_1,y_2)+\eta(x_2,y_1)\\
		&=&\tx{tr}_{E/F}((\eta(x_1,y_2),\eta(x_2,y_1))\\
		&=&\tx{tr}_{E/F}(h_1(x,y)),
	\end{eqnarray*}
	and identifies $\GL(V_+)=\U(V_+^E,h_1)$ with $M$.
\end{enumerate}

Let $E/F$ be a quadratic \'etale algebra. The structure map $F \to E$ identifies $F$ with $E^\sigma$. We have the norm map $\rN : E^\times \to E^\times$ sending $x \mapsto x\sigma(x)$, which takes values in $F^\times$, and we denote its kernel by $E^1$. We also have the twisted norm map $E^\times \to E^\times$ 
sending $x \mapsto x/\sigma(x)$, which takes values in $E^1$. These lead to the exact sequences
\[ 1 \to F^\times \to E^\times \to E^1 \to 1 \]
and 
\[ 1 \to E^1 \to E^\times \to F^\times \to F^\times/\rN(E^\times) \to 1. \]
If $F$ is local and $E/F$ is non-split, then $F^\times/\rN(E^\times)$ is a group of order $2$ and we denote by $\kappa : F^\times/\rN(E^\times) \to \{\pm1\}$ the unique group isomorphism. When $E/F$ is split, then $F^\times/\rN(E^\times)$ is the trivial group and we denote by $\kappa : F^\times/\rN(E^\times) \to \{\pm1\}$ the trivial character. We sometimes treat $\kappa$ as having domain $F^\times$ and kernel $\rN(E^\times)$. We will sometimes decorate $\rN$
and $\kappa$ with the subscript $E/F$.

\subsection{$\Lambda$-gradings} \label{sub:lambda}

Let $V$ be an $F$-vector space, $\eta : V \times V \to F$ a non-degenerate $\epsilon$-symmetric bilinear form for $\epsilon \in \{\pm1\}$, and $G:=\SU(V,\eta):=\Aut(V,\eta)^\circ$. Let $\Sigma=\Gamma \times \{\pm1\}$,\index[notation]{Sigma@$\Sigma$} where $\Gamma$ denotes the absolute Galois group of $F$. We denote by $\bar F$ a separable closure of $F$.

\begin{dfn} \label{dfn:lambda_basic}
	Let $\Lambda$ be a $\Sigma$-set. 
	\begin{enumerate}
		\item An element $\lambda \in \Lambda$ is called
		\begin{enumerate}
			\item \emph{symplectic}, if $-\lambda=\lambda$,
			\item \emph{polarizable}, if $-\lambda \neq \lambda$,\index[terminology]{polarizable!polarizable $\Sigma$-set}
			\item \emph{asymmetric}, if $-\lambda \notin \Gamma\lambda$ (hence $\lambda$ is automatically polarizable),
			\item \emph{symmetric}, if $-\lambda \in \Gamma\lambda$ and $\lambda$ is polarizable.			
		\end{enumerate}
		\item The set $\Lambda$ is called symmetric/asymmetric/polarizable/symplectic if all of its elements are of that type.
		\item We let $\Lambda_-$, $\Lambda^-$, $\Lambda_\tx{asym}$, and $\Lambda_\tx{sym}$, be the subsets of symplectic, polarizable, asymmetric, and symmetric  elements, respectively. 
		\item For $\lambda \in \Lambda$ let $\Gamma_\lambda$\index[notation]{__Glambda@$\Gamma_\lambda$} be the stabilizer of $\lambda$, and $\Gamma_{\pm\lambda}$\index[notation]{__Glambdapm@$\Gamma_{\pm\lambda}$} the stabilizer of the subset $\{\pm\lambda\}$, so that $\Gamma_\lambda \subset \Gamma_{\pm\lambda}$. Let and $F_\lambda$\index[notation]{__Flambda@$F_\lambda$}, and let $F_{\pm\lambda}$\index[notation]{__Flambdapm@$F_{\pm\lambda}$} be the fixed fields of $\Gamma_\lambda$ and $\Gamma_{\pm\lambda}$, respectively. Thus $F_\lambda/F_{\pm\lambda}$ is a field extension of degree $2$ when $\lambda$ is symmetric and degree $1$ otherwise. Let $E_{\pm\lambda}:=F_\lambda=F_{-\lambda}$\index[notation]{__Elambdapm@$E_{\pm\lambda}$} if $\lambda$ is symmetric and $E_{\pm\lambda}:=F_{\pm\lambda} \times F_{\pm\lambda}$ if $\lambda$ is asymmetric (if $\lambda$ is symplectic, we do not define $E_{\pm\lambda}$). Let $\sigma_\lambda \in \tx{Aut}(E_{\pm\lambda}/F_{\pm\lambda})$ be the unique non-trivial automorphism.
	\end{enumerate}
\end{dfn}

We obtain the disjoint union decompositions 
\[ \Lambda = \Lambda_- \cup \Lambda^-,\qquad \Lambda^- = \Lambda_\tx{asym} \cup \Lambda_\tx{sym}. \]

\begin{dfn} \label{dfn:lambda_grading}
	\begin{enumerate}
		\item 	Let $\Lambda$ be a $\Sigma$-set. A \emph{$\Lambda$-grading} on $V$ is a direct sum decomposition
	\[ V_{\bar F} = \bigoplus_{\lambda \in \Lambda} V_{\bar F,\lambda} \]
	\begin{enumerate}
		\item $V_{\bar F,\sigma\lambda}=\sigma(V_{\bar F,\lambda})$ for all $\lambda \in \Lambda$ and $\sigma \in \Gamma$,
		\item $V_{\bar F,\lambda} \perp V_{\bar F,\mu}$ if $\lambda \neq -\mu$,
		\item $V_{\bar F,\lambda} \neq \{0\}$.
	\end{enumerate}
	\item A \emph{$\Sigma$-grading} on $V$ is a pair consisting of a $\Sigma$-set $\Lambda$ and a $\Lambda$-grading on $V$. It is called \emph{earnest} if either $\epsilon=-1$, or $\epsilon=+1$ and there is at most one $\lambda \in \Lambda_-$ for which $\dim(V_\lambda)$ is odd and there is no $\lambda \in \Lambda_-$ for which $\dim(V_\lambda)=2$.
	\item An isomorphism betweeen two $\Sigma$-gradings is a $\Sigma$-equivariant bijection $f: \Lambda_1 \to \Lambda_2$ between their indexing sets such that $V_{\bar F, \lambda_1}=V_{\bar F, f(\lambda_1)}$ for all $\lambda_1 \in \Lambda_1$.
	\item Given a $\Sigma$-grading on $V$ we define
	\begin{enumerate}
		\item For each $\Sigma$-orbit $O \subseteq \Lambda$ the $F$-vector space $V_O:=(V_{\bar F,O})^\Gamma$ with $V_{\bar F,O}:=\bigoplus_{\lambda \in O} V_{\bar F,\lambda}$ and the form $\eta_O$ as the restriction of $\eta$ to $V_O$.
		\item For each $\lambda \in \Lambda$ the $F_\lambda$-vector space $V_\lambda:=(V_{\bar F,\lambda})^{\Gamma_\lambda}$. If $\lambda=-\lambda$, then we write $\eta_\lambda$ for the restriction of the $F_\lambda$-linear extension of $\eta$ to $V_\lambda$.
		\item For each $\{\lambda,-\lambda\} \subseteq \Lambda$ the $F_{\pm\lambda}$-vector space $V_{\pm\lambda}:=(V_{\bar F,\lambda} + V_{\bar F,-\lambda})^{\Gamma_{\pm\lambda}}$ and the  form $\eta_{\pm\lambda}$ as the restriction to $V_{\pm\lambda}$
		of the $F_{\pm\lambda}$-linear extension of $\eta$. 
\item For each $\{\lambda,-\lambda\} \subseteq \Lambda^-$ the $E_{\pm\lambda}$-module $V_{\pm\lambda}^E := V_{\pm\lambda} \otimes_{F_{\pm\lambda}} E_{\pm\lambda}$ and the form $\eta_{\pm\lambda}^E$ as the $E_{\pm\lambda}$-linear extension of $\eta_{\pm\lambda}$.
		\item For each $\lambda \in \Lambda^-$ the Lagrangian $E_{\pm\lambda}$-submodule $V_\lambda^E \subset V_{\pm\lambda}^E$ as $V_\lambda^E:=V_\lambda$ when $\lambda$ is symmetric and 
		\begin{eqnarray*}
                        V_\lambda^E&:=&(V_\lambda \oplus \{0\}) \oplus (\{0\} \oplus V_{-\lambda}) \\
                        &\subset& (V_\lambda \oplus V_{-\lambda}) \oplus (V_\lambda \oplus V_{-\lambda})\\
                        &=&(V_\lambda \oplus V_{-\lambda}) \otimes_{F_{\pm\lambda}} (F_{\pm\lambda} \times F_{\pm\lambda})
        \end{eqnarray*}
when $\lambda$ is asymmetric, and the form
		\[ h_\lambda^E : V_\lambda^E \times V_\lambda^E \to E_{\pm\lambda},\qquad h_\lambda^E(x,y)=\eta_{\pm\lambda}^E(x,\sigma_\lambda(y)). \]
	\end{enumerate}
	\end{enumerate}
\end{dfn}

\begin{fct}  \label{fct:lambda_grading}
	Assume given a $\Sigma$-grading on $V$.
	\begin{enumerate}
		\item We have the orthogonal decomposition $(V,\eta)=\bigoplus_{O \in \Lambda/\Sigma} (V_O,\eta_O)$.
		\item Given an orbit $O \in \Lambda/\Sigma$ and $\lambda \in O$ the map
		\[ V_{\pm\lambda}|_F \to V_O,\qquad v \mapsto \sum_{\sigma \in \Gamma/\Gamma_{\pm\lambda}} \sigma(v), \]
		is an isomorphism that identifies $\eta_O$ with the trace of $\eta_{\pm\lambda}$. \item Given $\lambda \in \Lambda^-$ the $E_{\pm\lambda}$-submodule $V_\lambda^E \subset V_{\pm\lambda}^E$ is Lagrangian for the form $\eta_{\pm\lambda}^E$, the form $h_\lambda^E$ is $\epsilon$-Hermitian, the structure map $F_{\pm\lambda} \to E_{\pm\lambda}$ induces an isomorphism $V_{\pm\lambda} \to (V_{\pm\lambda}^E)^{\sigma_\lambda}$, the map
		\[  V_\lambda^E \to V_{\pm\lambda}^E,\quad x \mapsto x+\sigma_\lambda(x) \]
		takes values in $V_{\pm\lambda} = (V_{\pm\lambda}^E)^{\sigma_\lambda}$ and induces an isomorphism
		\[ V_\lambda^E|_{F_{\pm\lambda}} \simto V_{\pm\lambda} \]
		which identifies $\eta_{\pm\lambda}$ with the trace of $h_\lambda^E$.
	\end{enumerate}
\end{fct}

\begin{lem}\label{lem:Levivsgrading}
	\begin{enumerate}
		\item The following maps give mutually inverse bijections between the set of connected reductive subgroups $H \subseteq G$ of maximal rank and the set isomorphism classes of earnest $\Sigma$-gradings on $V$. To $H$ we assign the weight decomposition of $V_{\bar F}$ under the center $Z(H)$ of $H$, and to a $\Sigma$-grading on $V$ we assign the identity component of the mutual stabilizer of all graded pieces.
		\item Under the above bijection we have 
\begin{eqnarray*}
			H&=&\prod_{\lambda \in \Lambda_-/\Sigma} \tx{Res}_{F_\lambda/F}\SU(V_{\lambda},\eta_{\lambda})\\
			&\times&\prod_{\lambda \in \Lambda_\tx{asym}/\Sigma} \tx{Res}_{F_\lambda/F}\GL(V_\lambda)\\
			&\times&\prod_{\lambda \in \Lambda_\tx{sym}/\Sigma} \tx{Res}_{F_{\pm\lambda}/F}\U(V_{\lambda},h_{\lambda}).
		\end{eqnarray*}
\item The subgroup $H$ is a twisted Levi subgroup if and only if $|\Lambda_-| \leq 1$.
		\item Given a $\Lambda$-grading on $V$, the identity component of the common stabilizer in $\GL(V)$ of all summands is the twisted Levi subgroup $H_{\GL} \subseteq \GL(V)$ given by 
		\[ \prod_{\lambda \in \Lambda/\Gamma} \tx{Res}_{F_\lambda/F}\GL(V_\lambda). \]
		The embedding $H \hookrightarrow H_{\GL}$ is given by the standard embeddings
		\begin{eqnarray*}
			\SU(V_{\lambda_-},\eta_{\lambda_-})&\hookrightarrow&\GL(V_{\lambda_-}),\\
			\U(V_{\lambda},h_{\lambda})&\hookrightarrow&\tx{Res}_{F_\lambda/F_{\pm\lambda}}\GL(V_{\lambda}) \hookrightarrow \GL(V_\lambda \oplus V_{-\lambda}),\\
			\GL(V_\lambda)&\hookrightarrow&\GL(V_\lambda \oplus V_{-\lambda}),\ g \mapsto (g,g^{-*}), 
		\end{eqnarray*}
		where $g^{-*}$ denotes the inverse of the adjoint of $g$.
	\end{enumerate}
\end{lem}
\begin{proof}
	It is clear that the maps in (1) are well-defined, so it is enough to show that they are inverses of each other. Given $H$, the Borel--de Siebenthal theorem (see \cite{BdS1949} for the original result for compact Lie groups and \cite[Proposition 1.1]{SPL2015} for the language we need) asserts that $H$ is the identity component of the centralizer of $Z(H)$, therefore $H$ is precisely the identity component of the common stabilizer in $G$ of the weight spaces in $V_{\bar F}$ under the action of $Z(H)$.

	In the opposite direction, give a finite $\Sigma$-set $\Lambda$ and a $\Lambda$-grading on $V$, the common stabilizer in $\GL(V)$ of all summands is easily seen to be the group $H_{\GL}$ as in (4), which is a twisted Levi subgroup. Moreover, the identity componenet of the common stabilizer in $G$ equals $(G \cap H_{\GL})^\circ$, which one computes to be the group $H$ given in (2). One now verifies immediately that the set of weights for the action of $Z(H)$ on $V$ is precisely $\Lambda$.

	We have thus proved (1), (2), (4), and it remains to show (3). Recall that a maximal rank connected reductive subgroup $H$ of $G$ is a twisted Levi subgroup if and only if it is the centralizer of $Z(H)^\circ$. From the description in (2) we see that $Z(H) = Z(H)^\circ$ if $\Lambda_- = \varnothing$. Assume $\Lambda_- \neq \varnothing$, for otherwise we are done. Then the set of weights of $Z(H)^\circ$ is $\Lambda^- \cup \{0\}$ and $V_0 := (\bigoplus_{\lambda \in \Lambda_-} V_{\bar F,\lambda})^\Gamma$. Thus the decomposition of $V$ does not change if $|\Lambda_-|=1$, but it does otherwise, and so does its centralizer, i.e., $H$ is not the centralizer of $Z(H)^\circ$ when $|\Lambda_-|>1$.
\end{proof}

\begin{rem}\label{rem:torussigmaset}
	A maximal torus of $G$, being a special case of a twisted Levi subgroup, is determined uniquely by a $\Sigma$-grading on $V$ satisfying $\dim_{F_\lambda}(V_\lambda)=1$ for all $\lambda \in \Lambda$. In particular $\Lambda_-=\varnothing$ when $\dim(V)$ is even, and $\Lambda_-=\{0\}$ with $\dim(V_0)=1$ when $\dim(V)$ is odd, and we have
		\begin{eqnarray*}
		T&=&\prod_{\lambda \in \Lambda_\tx{asym}/\Sigma} \tx{Res}_{F_\lambda/F}\mb{G}_{m,F_\lambda} \times \prod_{\lambda \in \Lambda_\tx{sym}/\Sigma} \tx{Res}_{F_{\pm\lambda}/F}\U_{1,F_\lambda/F_{\pm\lambda}},\\
		T_{\GL}&=&\prod_{\lambda \in \Lambda/\Gamma} \tx{Res}_{F_\lambda/F}\mb{G}_{m,F_\lambda}.
		\end{eqnarray*}
		We have 
		\begin{eqnarray*}
			T_{\GL}(F)&=&\{t_\lambda \in (\bar F^\times)^\Lambda|\sigma(t_\lambda)=t_{\sigma(\lambda)} \}\\
			T(F)&=&\{t_\lambda \in (\bar F^\times)^{\Lambda^-}|\sigma(t_\lambda)=t_{\sigma(\lambda)},t_{-\lambda}=t_\lambda^{-1} \}.
		\end{eqnarray*}
		From this it is clear that the centralizer of $T$ in $\GL(V)$ is $T_{\GL}$.
\end{rem}

\begin{dfn}
	\index[terminology]{polarizable!polarizbale twisted Levi subgroup}
	Let $H$ be a twisted Levi subgroup $H \subseteq G$ with associated finite $\Sigma$-set $\Lambda_H$. Then $H$ is called \emph{polarizable} if
	\begin{enumerate}
		\item $\dim(V)$ is even and $\Lambda_H$ is polarizable, i.e. $\Lambda_-=\varnothing$, or
		\item $\dim(V)$ is odd, $\Lambda_-=\{0\}$ and $\dim(V_0)=1$.
	\end{enumerate}
\end{dfn}

Note that a polarizable twisted Levi subgroup is precisely one that does not have a symplectic or orthogonal factor.

\begin{dfn} \label{dfn:achi}
	Let $\Lambda$ be a polarizable $\Sigma$-set.
	\begin{enumerate}
		\item A set of $a$-data for $\Lambda$ is a $\Sigma$-equivariant map $a : \Lambda \to \bar F^\times$, $\lambda \mapsto a_\lambda$, where $-1 \in \Sigma$ acts on $F^\times$ by multiplication by $-1 \in F^\times$. 
\item Assuming $F$ is local, a set of $\chi$-data for $\Lambda$ is a collection of extensions $\chi_\lambda : E_{\pm\lambda}^\times \to \C^\times$ of $\kappa_\lambda$ satisfying $\chi_{\sigma(\lambda)}=\chi_\lambda\circ\sigma^{-1}$ for $\sigma \in \Gamma$ and $\chi_{-\lambda}=\chi_\lambda^{-1}$.
\end{enumerate}
\end{dfn}

\begin{rem} \label{rem:adata}
	We can also think of $a$-data as a collection of $a_\lambda^E \in E_{\pm\lambda}^\times$ such that $\sigma_\lambda(a_\lambda^E)=-a_\lambda^E$ and $a_{\sigma\lambda}^E=\sigma(a_\lambda^E)$ for all $\sigma \in \Gamma$. Here $a_\lambda^E=a_\lambda \in F_\lambda^\times=E_{\pm\lambda}^\times$ when $\lambda$ is symmetric, and $a_\lambda^E=(a_\lambda,a_{-\lambda})=(a_\lambda,-a_\lambda) \in F_\lambda^\times \times F_{-\lambda}^\times = E_{\pm\lambda}^\times$ when $\lambda$ is asymmetric.
\end{rem}

\begin{lem} \label{lem:askew}
	Assume given a polarizable $\Sigma$-set $\Lambda$, a $\Lambda$-grading on $V$, and a set of $a$-data for $\Lambda$. 
	\begin{enumerate}
		\item The automorphism
	\[ T_a : V_{\bar F} \to V_{\bar F},\qquad \sum_\lambda v_\lambda \mapsto \sum_\lambda a_\lambda v_\lambda. \]
	is $\Gamma$-equivariant and skew-adjoint for $\eta$.
	\item The bilinear form $\eta_a : V \times V \to F, (v,w) \mapsto \eta(T_av,w)$ is non-degenerate and $(-\epsilon)$-symmetric.
	\index[notation]{\eta_a@$\eta_a$}
	\end{enumerate}
\end{lem}
\begin{proof}
	The $\Gamma$-equivariance of $T_a$ follows from that of $a$, while skew-adjointness follows from 
	\[ \eta(T_av,w) = \eta(\sum_\lambda a_\lambda v_\lambda,\sum_\lambda w_\lambda)  =  \sum_\lambda \eta(v_\lambda,a_\lambda w_{-\lambda}) =\eta(v,-T_aw).  \]
	This shows (1), and (2) follows immediately.
\end{proof}

\subsection{Lattice chains} \label{sub:lattice}

We assume now that $F$ is a non-archimedean local field and denote by $\cO=\cO_F$ its ring of integers, by $\mf{p}=\mf{p}_F$ the maximal ideal of $\cO$, and by $\rf=\rf_F$ the residue field, which we assume to be of odd characteristic.

Let $E$ be a trivial or quadratic \'etale $F$-algebra, $\sigma \in \tx{Aut}_F(E)$ the unique isomorphism that is non-trivial when $E/F$ is non-trivial.

If $E$ is not the split quadratic algebra $F \oplus F$, then we let $\tx{ord}_E : E^\times \to \Z$ be the normalized valuation, i.e., the one sending uniformizers to $1$, and set $\tx{ord}_E(\{0\})=\infty$. It is invariant under $\sigma$ and $\cO_E=\{x \in E|\tx{ord}_E(x)\geq 0\}$ is a local ring with maximal ideal $\mf{p}_E = \{x \in E|\tx{ord}_E(x)>0\}$. For $a \in E^\times$ with $\tx{ord}_E(a)=n$ we have $a\cO_E = \mf{p}_E^n$. Assuming $E \neq F$, there exist elements $a \in E^\times$ with $\sigma(a)=-a$. If $E/F$ is ramified, then $\tx{ord}_E(a) \in 2\Z+1$.

When $E = F \oplus F$ is the split quadratic algebra, we equip the set $E \sm \{0\}$ with the ``valuation'' $\tx{ord}_E((x_1,x_2))=\min(\tx{ord}_F(x_1),\tx{ord}_F(x_2))$ and set $\tx{ord}_E(0)=\infty$. This function is again invariant under $\sigma$, $\cO_E=\{x \in E|\tx{ord}_E(x)\geq 0\}=\cO_F \oplus \cO_F$ is a ring which is not local: it has the two maximal ideals $\mf{p}_F \oplus \cO_F$ and $\cO_F \oplus \mf{p}_F$, whose intersection $\mf{p}_E:=\mf{p}_F \oplus \mf{p}_F$ is the Jacobson radical of $\cO_E$. We still have $\mf{p}_E=\{x \in E|\tx{ord}_E(x)>0\}$. For $a \in E^\times$ with $n=\tx{ord}_E(a)$ it is generally not true that $a\cO_E = \mf{p}_E^n$, but it is true at least when $\sigma(a)=\pm a$. Such elements always exist, have the form $a=(y,\pm y)$ with $y \in F$, and $\tx{ord}_E(a)$ can be arbitrary.

Let $(V,h)$ be an $\epsilon$-Hermitian $E$-space for $\epsilon \in \{\pm1\}$. For an $\cO_E$-lattice $L \subset V$ define its $h$-dual by
\begin{equation} \label{eq:hdual}
L^h = \{v \in V|h(v,L) \subseteq \cO_E\}.
\end{equation}

\begin{lem} \label{lem:adual}
For $n \in \Z$ and an $\cO_E$-lattice $L \subset V$ we have $(\mf{p}_E^nL)^h = \mf{p}_E^{-n}L^h$.	
\end{lem}
\begin{proof}
	There exists $a \in E^\times$ with $\mf{p}_E^n=a\cO_E$ and $\sigma(a) = \pm a$. Then
	\[ (\mf{p}_E^nL)^h = (aL)^h=\{v \in V|h(v,aL) \subseteq \cO_E\} = a^{-1}L^h = \mf{p}_E^{-n}L^h. \]
\end{proof}

\begin{lem} \label{lem:different-dual}
	Let $F_0 \subseteq F$ be a subfield such that $F/F_0$ is a finite separable extension and let $E_0 \subseteq E$ be an \'etale $F_0$-subalgebra such that $E_0/F_0$ is trivial or quadratic. If $E_0/F_0$ is quadratic, assume that the multiplication map $F \otimes_{F_0} E_0 \to E$ is an isomorphism. Then $h_0=\tx{tr}_{E/E_0}\circ h$ is an $\epsilon$-Hermitian form on $V_0=V|_{E_0}$ with respect to $E_0/F_0$ and for any $\cO_E$-lattice $L \subset V$ we have
	\[ (L|_{\cO_{E_0}})^{h_0} = (\mc{D}_{E/E_0}^{-1} \cdot L^h)|_{\cO_{E_0}}, \]
	where $\mc{D}_{E/E_0}^{-1}$ is the inverse different of $E/E_0$, i.e., the dual of the $\cO_{E_0}$-lattice $\cO_E$ in the $E_0$-module $E$ with respect to the trace bilinear pairing $(x,y) \mapsto \tx{tr}_{E/E_0}(x\cdot y)$.
\end{lem}
\begin{proof}
	Let $\sigma$ and $\sigma_0$ denote the automorphisms of $E/F$ and $E_0/F_0$, respectively, taken to be non-trivial when the algebra is of dimension $2$ and trivial in dimension $1$. Our assumptions imply $\sigma_0 \circ \tx{tr}_{E/E_0} = \tx{tr}_{E/E_0} \circ \sigma$. From this the sesqui-linearity of $h_0$ follows, and we have
	\begin{eqnarray*}
		(L|_{\cO_{E_0}})^{h_0}&=&\{x \in V_0|\tx{tr}_{E/E_0}(h(x,L)) \subseteq \cO_{E_0} \}\\
		&=&\{x \in V_0|\tx{tr}_{E/E_0}(\cO_E \cdot h(x,L)) \subseteq \cO_{E_0} \}\\
		&=&\{x \in V_0|h(x,L) \subseteq \mc{D}_{E/E_0}^{-1} \}\\
		&=&(\mc{D}_{E/E_0}^{-1} \cdot \{x \in V|h(x,L) \subseteq \cO_E \})|_{\cO_{E_0}}\\
		&=&(\mc{D}_{E/E_0}^{-1} \cdot L^h)|_{\cO_{E_0}}.
	\end{eqnarray*}
	Note that, since $E/E_0$ is etale, the pairing $(a,b) \mapsto \tx{tr}_{E/E_0}(ab)$ is non-degenerate, so the non-degeneracy of $h_0$ follows from that of $h$.
\end{proof}

From now on by a lattice in $V$, we always mean an $\cO_E$-lattice unless specified otherwise.

\begin{dfn}
	Let $\mc{L}$ be a self-dual periodic lattice chain in $(V,h)$.
	\begin{enumerate}
		\item The largest lattice $L \in \mc{L}$ such that $L \subseteq L^h$ will be called the \emph{center}\index[terminology]{center} 
of $\mc{L}$.
		\item If $L \in \mc{L}$ is the center, we will call  $L/\mf{p}_EL^h$ the \emph{central quotient}.\index[terminology]{central quotient}
	\end{enumerate} 
	If we want to emphasize $h$, we will say \emph{$h$-center}\index[terminology]{center h@$h$-center} and \emph{h-central quotient}.\index[terminology]{central quotient h@$h$-central quotient}
\end{dfn}

\begin{dfn}
	A lattice chain $\mc{L}$ is called \emph{atomic}\index[terminology]{atomic} if it is generated under duality and periodicity by a single lattice.
\end{dfn}

\begin{lem} \label{lem:centequiv}
	Let $\mc{L}$ be a self-dual periodic lattice chain. If $L \in \mc{L}$ is the center, then $\mf{p}_EL^h \subseteq L \subseteq L^h$. The converse holds if $\mc{L}$ is atomic.
\end{lem}
\begin{proof}
	A lattice $L \in \mc{L}$ is the center if and only if $L \subseteq L^h$ and there is no lattice strictly between $L$ and $L^h$. If $\mf{p}_EL^h \not\subseteq L$, then $L \subsetneq \mf{p}_EL^h$ since $\mc{L}$ is totally ordered, and $\mf{p}_EL^h \subsetneq L^h$, contradicting that $L$ is the center.

	For the converse assume that $\mc{L}$ is atomic and $\mf{p}_EL^h \subseteq L \subseteq L^h$. If there existed $L \subsetneq L' \subsetneq L^h$, then by atomicity either $L'=\mf{p}_E^{-n}L$ for some $n >0$, hence $\mf{p}_E^{-1}L \subsetneq L^h$, hence $L \subsetneq \mf{p}_EL^h$, contradiction, or $L'=\mf{p}_E^nL^h$ for some $n>0$, hence again arriving at the contradiction $L \subsetneq \mf{p}_EL^h$.
\end{proof}

If $\mc{L}$ is a self-dual periodic lattice chain and $L \in \mc{L}$ is the $h$-center, then  restricting $h$ to $L$ one obtains a non-degenerate pairing $\ol{h}$ on the $h$-central quotient $\oL:=L/\mf{p}_EL^h$, which is $\epsilon$-Hermitian with respect to the trivial or quadratic \'etale algebra $\rf_E/\rf_F$, where $\rf_E:=\cO_E/\fp_E$. 

Assuming $\mc{L}$ is atomic, there is just one more quotient that is interesting, namely $L^h/L$. One can define a non-degenerate form on it by choosing a uniformizer $\varpi \in E$ such that $\sigma(\varpi)=\pm\varpi$ (the sign being $+1$ or $-1$ depending on $E/F$ being unramified or ramified, respectively) and taking the form induced by $h':=\varpi h$. Since the isometry class of this form generally depends on the choice of $\varpi$, we will prefer to not think of it  as associated to $h$, but rather as associated to $h'$. Thus, rather than having two spaces $(L/\mf{p}_EL^h,\ol{h})$ and $(L^h/L,\ol{h'})$ associated to $\mc{L}$ and $h$, we have the single space $(\oL,\ol{h})$ associated to $\mc{L}$ and $h$ but we are free to rescale $h$ by a uniformizer $\varpi \in E$ and repeat this construction. In fact, it is more convenient to allow rescaling by an arbitrary element $a \in E^\times$ with $\sigma(a)=\pm a$. We will see however that, while this does provide notational flexibility, it does not introduce any new phenomena. 

For simplicity, we set $E^\pm=\{x \in E^\times|\sigma(x)=\pm x\}$. Thus $E^\pm$ is a subgroup of $E^\times$ that contains $F^\times=(E^\times)^\sigma$.
\index[notation]{$E^\pm$@$E^\pm$}

\begin{lem} \label{lem:ssprep}
	Let $a \in E^\pm$.
	\begin{enumerate}
		\item For any $a',b \in E^\pm$ with $a=a'b\sigma(b)$ and any lattice $L \subset V$ we have  $(b^{-1}L)^{ah}=(a'b)^{-1}L^h$.
		\item There exist $a',b \in E^\pm$ with $a=a'b\sigma(b)$ and $\tx{ord}_E(a') \in \{0,1\}$.
	\end{enumerate}
\end{lem}
\begin{proof}
	(1) $(b^{-1}L)^{ah}=\{v \in V|a'b\sigma(b)h(v,b^{-1}L) \subseteq \cO_E\} = \{v \in V|a'bv \in L^h \} = (a'b)^{-1}L^h$.

	(2) Choose $\varpi \in E^\times$ with $\varpi \cO_E=\mf{p}_E$ and $\sigma(\varpi)=\pm\varpi$. Let $m=\lfloor \tx{ord}_E(a)/2 \rfloor$, $b=\varpi^m$, $a'=a/(b\sigma(b))$.
\end{proof}

\begin{lem} \label{lem:sign-switch}
	Let $\mc{L}$ be an atomic lattice chain for $(V,h)$ and let $L$ be its center.
	Let $a \in E^\pm$ and write $a=a'b\sigma(b)$ as in Lemma \ref{lem:ssprep}(2).
	\begin{enumerate}
				\item $\mc{L}$ is also atomic with respect to $ah$. 
		\item If $\tx{ord}_E(a')=0$, then the $ah$-center of $\mc{L}$ is $L_a := b^{-1}L$ and multiplication by $b$ provides an isometry $(\ol{L_a},\ol{ah}) \to (\oL,\ol{a'h})$.
		\item If $\tx{ord}_E(a')=1$, then the $ah$-center of $\mc{L}$ is $L_a := b^{-1}L^h$ and multiplication by $b$ provides an isometry $(\ol{L_a},\ol{ah}) \to (\ol{L^h},\ol{a'h})$.
	\end{enumerate}
\end{lem}
\begin{proof}
	(1) Atomicity under $ah$ means that for any fixed $L \in \mc{L}$ we have $\mc{L}=\mf{p}_E^\Z \cdot L \cup \mf{p}_E^\Z \cdot L^{ah}$. This follows from Lemma \ref{lem:ssprep}(1).

	(2),(3) From Lemma \ref{lem:centequiv} we have $\mf{p}_EL^h \subset L \subseteq L^h \subseteq \mf{p}_E^{-1}L$. If $\tx{ord}_E(a')=0$, then we multiply $\mf{p}_EL^h \subseteq L \subseteq L^h$ by $(a'b)^{-1}$ and use Lemma \ref{lem:ssprep}(1) to conclude $\mf{p}_E(b^{-1}L)^{ah} \subset b^{-1}L \subseteq (b^{-1}L)^{ah}$. If $\tx{ord}_E(a')=1$ then we multiply $L \subseteq L^h \subseteq \mf{p}_E^{-1}L$ by $b^{-1}$ to obtain $\mf{p}_E(a'b)^{-1}L \subseteq b^{-1}L^h \subseteq (a'b)^{-1}L$ and use Lemma \ref{lem:ssprep}(1) to conclude $\mf{p}_E(b^{-1}L^h)^{ah} \subseteq b^{-1}L^h \subseteq (b^{-1}L^h)^{ah}$. In both cases we apply Lemma \ref{lem:centequiv} to conclude that $L_a$ is the $ah$-center of $\mc{L}$. The isometry claim then follows from the fact that multiplication by $b$ induces an isometry $(V,ah) \to (V,a'h)$.
\end{proof}

We now revert to the notation of \S\ref{sub:lambda}. Thus $V$ is a finite-dimensional $F$-vector space and $\eta : V \times V \to F$ is a non-degenerate $\epsilon$-symmetric form. Consider given a polarizable $\Sigma$-set $\Lambda$ and a $\Lambda$-grading on $V$. Assume that for all $\lambda \in \Lambda$ the extension $F_\lambda/F$ is tame.

\begin{dfn} \label{dfn:lambda_atomic}
\begin{enumerate}
		\item A lattice chain $\mc{L}$ in $V$ is called \emph{$\Lambda$-atomic}\index[terminology]{atomicLambda@$\Lambda$-atomic} if it is periodic and
	\begin{enumerate}
		\item for each $L \in \mc{L}$ we have $L=\bigoplus_{O \in \Lambda/\Sigma} L_O$ with $L_O = L \cap V_O$,
		\item for each $O \in \Lambda/\Sigma$ and $\lambda \in O$ the isomorphism $V_\lambda^E|_F \to V_O$ of Fact \ref{fct:lambda_grading}(2,3) transports the lattice chain $\{L_O|L \in \mc{L}\}$ of $V_O$ to an atomic lattice chain of $V_\lambda^E$.
	\end{enumerate}
	\item The \emph{$\Lambda$-center}\index[terminology]{center Lambda@$\Lambda$-center} of $\mc{L}$ is the lattice $\bigoplus_{O \in \Lambda/\Sigma} L_O$,
	where $L_O \subset V_O$ is the image of the $h_\lambda^E$-center $L_\lambda^E \subset V_\lambda^E$ of $\mc{L}_\lambda^E$ under the isomorphism $V_\lambda^E|_F \to V_O$ of Fact \ref{fct:lambda_grading}(2,3).
	\item The \emph{$\Lambda$-central quotient}\index[terminology]{central quotient Lambda@$h\Lambda$-central quotient} of $\mc{L}$ is $(\ol{V},\ol{\eta})$, where $\ol{V}:=\bigoplus_{O \in \Lambda/\Sigma} L_O/L_O^+$ with $L_O$ as above and $L_O^+ \subseteq L_O$ is the image of $\mf{p}_\lambda\cdot (L_\lambda^E)^{h_\lambda^E}$ under the isomorphism $V_\lambda^E|_F \to V_O$ of Fact \ref{fct:lambda_grading}(2,3), and could equal $L_O$, and $\ol{\eta}$ is the (possibly skew)-symmetric form on $\ol{V}$ obtained by restricting $\eta_O$ to $L_O$ and taking the direct sum over $O$.
	\end{enumerate}
\end{dfn}

\begin{rem}
	\begin{enumerate}
		\item We are not requiring that a $\Lambda$-atomic lattice chain be self-dual. In particular, a $\Lambda$-atomic lattice chain need not be an atomic lattice chain. The main reasons for this are taking direct sums and using restriction of scalars.
		\item The $\Lambda$-center of a $\Lambda$-atomic lattice chain need not belong to the lattice chain.
		\item In our applications, the $\Lambda$-atomic lattice chain will come from the Moy--Prasad filtration, and will in general not be self-dual. For the symplectic form that we will work with, the $\Lambda$-center will belong to the lattice chain, but for a related symmetric form it will not.
	\end{enumerate}
\end{rem}

The following example will be of fundamental importance to us. 
\begin{exa} \label{exa:lambda_atomic}
Let $E$ be a quadratic \'etale $F$-algebra and let $\sigma$ be the non-trivial $F$-automorphism of $E$. Every $\epsilon$-Hermitian form on the 1-dimensional $E$-vector space $V:=E$ is of the form $(x,y)=ax\sigma(y)$ with some $a \in E^\times$ such that $\sigma(a)=\epsilon a$. If $E/F$ is non-split and ramified, then $\epsilon=(-1)^{\tx{ord}_E(a)}$. The lattice chain $\mc{L}=\{\mf{p}_E^n: n \in \Z\}$ is atomic for $(-,-)$; when $E$ is a field it is the unique atomic lattice chain. Let $\mc{L}|_F$ denote the lattice chain in the 2-dimensional $F$-vector space $V|_F$, thus $\mc{L}|_F=\{\mf{p}^n_E|_{\cO_F} : n \in \Z\}$. This lattice chain is periodic and self-dual with respect to the pairing $\tx{tr}_{E/F}\circ (-,-)$, but it need not be atomic. Indeed, if we denote by $[L_1:L_2]$ for $L_1,L_2 \in \mc{L}$ the set of lattices in $\mc{L}$ that contain $L_1$ and are contained in $L_2$, then we can choose $L \in \mc{L}$ with $L \subset L^\vee$ and refer to $\#[L:L^\vee]-1 \in \Z/2\Z$ as the \emph{parity} of $\mc{L}$ and to $\#[\mf{p}_E L:L]-1 \in \Z$ as the \emph{period} of $\mc{L}$ (they do not depend on $L$). Then when $E/F$ is non-split ramified and $\epsilon=-1$, the period of $\mc{L}$ is $1$ and the parity is odd,
so $\mc{L}|_F$ has period $2$ and even parity, and is therefore not atomic. It is however $\Lambda$-atomic, where $\Lambda=\tx{Hom}_F(E,\bar F)$ equipped with the $\tx{Gal}(\bar F/F)$-action coming from the action on $\bar F$, and $\{\pm1\}$-action coming from $\sigma$ acting on $E$. The $\Lambda$-center of $\mc{L}|_F$ lies in $\mc{L}|_F$.

Consider now a finite separable tame field extension $F_1/F$, a quadratic \'etale algebra $E/F_1$, a finite-dimensional $\epsilon$-Hermitian $E$-vector space $(V,(-,-))$, and an atomic $\cO_E$-lattice chain $\mc{L}$ in $V$. Then $\mc{L}|_F$ is still periodic, and is self-dual with respect to $\tx{tr}_{E/F}(-,-)$, but can fail to be atomic for any choice of $E$. It is however $\Lambda$-atomic for $\Lambda=\tx{Hom}_F(E,\bar F)$ with $\Sigma$-action defined analogously. The $\Lambda$-center of $\mc{L}|_F$ lies in $\mc{L}|_F$.

Finally, consider collections $F_i/F$, $E_i/F_i$, $(V_i,(-,-)_i)$, and $\mc{L}_i$, all indexed by some finite set $I$. Choose a strictly monotonic map $c_i : \mc{L}_i \to \R$ for each $i \in I$ subject to $c_i(\mf{p}_F \cdot L)=c_i(L) + 1$  and define the lattice $\mc{L}_{i,r} \in \mc{L}_i$ as $\sup\{L_i \in \mc{L}_i| c_i(L_i) \geq r\}$. Then
\[ \mc{L}:=\Big\{\mc{L}_r := \bigoplus_i \mc{L}_{i,r}|_F: r \in \R \Big\} \]
is a periodic lattice chain in $\bigoplus V_i|_F$ and the map and the map $c : \mc{L} \to \R$ defined by $c(L)=\sup\{r\in\R\mid\mc{L}_r=L\}$ is strictly monotonic and satisfies $c(\mf{p}_F \cdot L)=c(L) + 1$. In general $\mc{L}$ is neither self-dual with respect to $\sum_i \tx{tr}_{E_i/F}(-,-)_i$, nor atomic. It will however be $\Lambda$-atomic for the disjoint union $\Lambda = \coprod \Lambda_i$ with $\Lambda_i=\tx{Hom}_F(E_i,\bar F)$ as above. If there exists $r \in \R$ such that for each $i \in I$ the center of $\mc{L}_i$ is $\mc{L}_{i,r}$, then the $\Lambda$-center of $\mc{L}$ lies in $\mc{L}$ and equals $\mc{L}_r$. If no such $r$ exists, then the $\Lambda$-center of $\mc{L}$ does not lie in $\mc{L}$.
\end{exa}

\subsection{\texorpdfstring{The quadratic form $q_{\tilde s}$}{}} \label{sub:qs}
From now on we fix a non-trivial character $\psi: F \to \bC^\times$ that starting from \S\ref{sub:sympmain} we assume is trivial on $\fp$ and nontrivial on $\cO$.
We write $\QW(F)$\index[notation]{QW@$\QW$} for the quadratic Witt group, reviewed in Appendix \ref{sub:Wittgroup}, and $\gamma_\psi:\QW(F) \ra \mu_8$\index[notation]{gammapsi@$\gamma_\psi$} for the Weil index that is reviewed in Appendix \ref{sub:Weilindex}.

Let $(V,\eta)$ be a symplectic space over $F$. 

\begin{dfn} \label{dfn:qtildes}
	For a semi-simple element $\tilde s \in \GL(V)(F)$ define the (possibly degenerate) quadratic form\index[notation]{qstilde@$q_{\tilde s}$}
	\[ q_{\tilde s} : V \to F,\qquad v \mapsto \eta(\tilde sv,v). \]
\end{dfn}

In this subsection we will study the form $q_{\tilde s}$, compute its Weil index in different ways, and show that for a polarizable twisted Levi subgroup $H_{\Sp} \subset \Sp(V)$ the function assigning $\gamma_\psi(q_{\tilde s})$ to semi-simple $\tilde s \in H_{\GL}(F)$ descends to a double cover of $H_{\Sp}(F)$ by means of the norm map that we now define.

\begin{dfn} Let $T_{\Sp} \subset \Sp(V)$ be a maximal torus with associated $\Sigma$-set $\Lambda$ as in Remark \ref{rem:torussigmaset}. We denote $(T_{\Sp})_{\GL}$ simply by $T_{\GL}$ and define the \emph{norm}\index[terminology]{norm map} map as 
	\[ N : T_{\GL} \to T_{\Sp},\qquad (t_\lambda)_\lambda \mapsto (t_\lambda/t_{-\lambda})_\lambda. \]
\end{dfn}

\begin{lem} \label{lem:topssgamma}
	Let $T_{\Sp} \subset \Sp(V)$ be a maximal torus with associated $\Sigma$-set $\Lambda$.
	\begin{enumerate}
		\item The composition of the natural inclusion $T_{\Sp} \hookrightarrow T_{\GL}$ with the norm map $T_{\GL} \to T_{\Sp}$ is the squaring map.
		\item We denote the kernel of the norm map $T_{\GL} \to T_{\Sp}$ by $T_{\GL}^+$. It is an induced torus satisfying
		\[ T_{\GL}^+ \cong \prod_{\lambda \in \Lambda/\Sigma} \tx{Res}_{F_{\pm\lambda}/F}\mb{G}_{m,F_{\pm\lambda}}. \]
		\item Let $\kappa_\lambda : F_{\pm\lambda}^\times/N_{F_\lambda/F_{\pm\lambda}}(F_\lambda^{\times}) \to \{\pm1\}$ be the unique non-trivial character if $\lambda$ is symmetric, and the necessarily trivial character if $\lambda$ is asymmetric. Set $\kappa=\prod_{\lambda \in \Lambda/\Sigma} \tx{Res}_{F_{\pm\lambda}/F}(\kappa_\lambda)$. 
		The pushout of
		\[ \xymatrix{1 \ar[r]&T_{\GL}^+(F)\ar[r]\ar[d]^\kappa& T_{\GL}(F)\ar[r]& T_{\Sp}(F) \ar[r]&1\\
		&\{\pm1\}
		}\]
		is naturally identified with the double cover $T_{\Sp}(F)_\pm$ of $T_{\Sp}(F)$ relative to $\Lambda$.
		\item Let $(a_\lambda)_{\lambda \in \Lambda}$ be a set of $a$-data for $\Lambda$ and let $\eta_a$ be the symmetric bilinear form of Lemma \ref{lem:askew} and $q_a(v)=\eta_a(v,v)/2=\eta(T_av,v)/2$ the corresponding quadratic form. Then for  $\tilde s \in T_{\GL}(F)$ we have
		\[ \gamma_\psi(q_{\tilde s}) = \gamma_\psi(q_a|V^{s \neq 1}) \cdot \prod_{\substack{\lambda \in \Lambda_\tx{sym}/\Gamma\\ \lambda(s) \neq 1}} \kappa_\lambda\left(\frac{\tilde s_\lambda-\tilde s_{-\lambda}}{a_\lambda}\right). \]
		\item There exists a unique genuine function \[ T_{\Sp}(F)_\pm \to \mu_8,\qquad \dot s \mapsto \gamma_\psi(q_{\dot s}), \]
		such that $\gamma_\psi(q_{\dot s})=\gamma_\psi(q_{\tilde s})$ whenever $\tilde s \mapsto \dot s$ under the map $T_{\GL}(F) \to T_{\Sp}(F)_\pm$ of (3) with $\tilde s_\lambda=1$ whenever $s_\lambda = 1$. (Note that $\gamma_\psi(q_{\dot s})$ is just a symbol for an element in $\mu_4$ attached to $\dot s$, we don't define a quadratic form $q_{\dot s}$)
	\end{enumerate}
\end{lem} 
\begin{proof}
	(1) The composition sends $(t_\lambda)_\lambda \in T_{\Sp}$ to $(t_\lambda/t_{-\lambda})_\lambda \in T_{\Sp}$, but since $t_{-\lambda}=t_\lambda^{-1}$ we have $t_\lambda/t_{-\lambda}=t_\lambda^2$.

	(2) The norm map respects the decomposition of $\Lambda$ into orbits under $\Sigma$. Consider the orbit of some $\lambda \in \Lambda$. If this orbit is asymmetric, the norm map is obtained by applying the functor $\tx{Res}_{F_\lambda/F}$ to the map
	\[ \mb{G}_{m,F_\lambda} \times \mb{G}_{m,F_{\lambda}} \to \mb{G}_{m,F_\lambda},\quad (t_\lambda,t_{-\lambda}) \mapsto t_\lambda/t_{-\lambda}, \]
	whose kernel is the diagonally embedded copy of $\mb{G}_{m,F_\lambda}$; note that $F_\lambda=F_{-\lambda}=F_{\pm\lambda}$. If this orbit is symmetric, the norm map is obtained by applying the functor $\tx{Res}_{F_{\pm\lambda/F}}$ to the map
	\[ \tx{Res}_{F_\lambda/F_{\pm\lambda}}\mb{G}_{m,F_\lambda} \to \U_{1,F_\lambda/F_{\pm\lambda}},\quad t_\lambda \mapsto t_\lambda/\sigma_\lambda(t_\lambda), \]
	whose kernel is the canonically embedded copy of $\mb{G}_{m,F_{\pm\lambda}}$. Thus in both cases the kernel of the norm map is 
	\[ \tx{Res}_{F_{\pm\lambda}/F}\mb{G}_{m,F_{\pm\lambda}}, \]
	as claimed.

	(3) Follows from \cite[Remark 3.1.4]{KalDC} after observing that the tori denoted by $T^+$, $T$, and $T^-$ there coincide with the tori $T_{\GL}^+$, $T_{\GL}$, and $T_{\Sp}$ here, and moreover we have $T^-=S$ in our application.

	(4) We proceed through a series of reductions that work for both $q_{\tilde s}$ and $q_a$, so we let $q$ denote either of these quadratic forms. Decompose 
	\[ V = \bigoplus_{O \in \Lambda/\Sigma} V_O,\qquad (V_O)_{\bar F}=\bigoplus_{\lambda \in O} V_{\bar F,\lambda}. \]
	The first decomposition is defined over $F$, stable under $\tilde s$, and orthogonal with respect to $\eta$, hence orthogonal with respect to $q$. Therefore \[ \gamma_\psi(q)=\prod_O \gamma_\psi(q|_{V_O}). \]
	For $\lambda \in O$ we have from Fact \ref{fct:lambda_grading}(2) the isomorphism
	\[ (V_{\pm\lambda}|_{F},\tx{tr}_{F_{\pm\lambda}/F} \circ q_{\pm\lambda}) \to (V_O,q|_{V_O}), \]
	where we recall that $V_{\pm\lambda}=(V_{\bar F,\lambda} \oplus V_{\bar F,-\lambda})^{\Gamma_{\pm \lambda}}$ and $q_{\pm\lambda}$ is the restriction to $V_{\pm\lambda}$ of the scalar extension of $q$ to $V \otimes F_{\pm\lambda}$. According to Fact \ref{fct:gammaind} we have 
	\[ \gamma_\psi(q|_{V_O}) = \gamma_{\psi_{\pm\lambda}}(q_{\pm\lambda}), \]
	where $\psi_{\pm\lambda}=\psi\circ\tx{tr}_{F_{\pm\lambda}/F}$. If $\lambda$ is asymmetric, then $V_\lambda$ and $V_{-\lambda}$ form an $\tilde s$-stable polarization of $V_{\pm\lambda}$ for the symplectic form $\eta$, hence a polarization for the quadratic form $q_{\pm\lambda}$. Therefore, 
	\[ \gamma_{\psi_{\pm\lambda}}(q_{\pm\lambda}) = 1. \]
	If $\lambda$ is symmetric, then the subspaces $V_\lambda$ and $V_{-\lambda}$ of $V_{\pm\lambda}$ are only defined over the quadratic extension $F_\lambda/F_{\pm\lambda}$ and $V_{-\lambda}=\sigma_\lambda(V_\lambda)$. We have the isomorphism 
	\[ V_\lambda|_{F_{\pm\lambda}} \simto V_{\pm\lambda},v \mapsto v+\sigma_\lambda(v). \]
	Pulling back the quadratic form $q_{\tilde s,\pm\lambda}$ under this isomorphism produces the form
	\[ q_{\tilde s,\lambda}(v)=q_{\tilde s,\pm\lambda}(v+\sigma_\lambda(v))
		=\eta(\tilde s(v+\sigma_\lambda(v)),v+\sigma_\lambda(v))
		=(\tilde s_\lambda-\tilde s_{-\lambda})\eta(v,\sigma_\lambda(v)), \]
	while pulling back $q_{a,\pm\lambda}$ produces 
	\[ q_{a,\lambda}(v)=q_{a,\pm\lambda}(v+\sigma_\lambda(v))
		=\eta_a(v,\sigma_\lambda(v))
		=a_\lambda\eta(v,\sigma_\lambda(v)). \]

	If $\tilde s_\lambda=\tilde s_{-\lambda}$, then the quadratic form $q_{\tilde s,\pm\lambda}$ is zero, so its Weil index is trivial. Since $\lambda(s)=\tilde s_\lambda / \tilde s_{-\lambda}$, this happens if and only if  $\lambda(s)=1$, i.e., $\lambda$ is not a weight of $V^{s \neq 1}$. If that is not the case, then the above quadratic forms are the trace quadratic forms
of the non-degenerate Hermitian forms $h_{\tilde s},h_a : V_\lambda \times V_\lambda \to F_\lambda$ given by 
	\[ (v,w) \mapsto (\tilde s_\lambda-\tilde s_{-\lambda})\cdot\eta(v,\sigma_\lambda(w))\quad\tx{and}\quad (v,w) \mapsto a_\lambda\cdot\eta(v,\sigma_\lambda(w)), \]
	respectively. Here by trace quadratic form $\mathrm{tr}(V_\lambda, h)$ of a Hermitian form $(V_\lambda, h)$ we mean the quadratic form $(V_\lambda|_{F_{\pm \lambda}}, v \mapsto h(v,v))$, see also Appendix \ref{sub:Wittgroup} for a discussion of trace quadratic forms.
	From $h_{\tilde s}=a_\lambda^{-1}(\tilde s_\lambda-\tilde s_{-\lambda})h_a$ and Lemma \ref{lem:weilherm} we infer 
	\[ \gamma_{\psi_{\pm\lambda}}(q_{\tilde s,\pm\lambda})=\kappa_{\lambda}(a_\lambda^{-1}(\tilde s_\lambda-\tilde s_{-\lambda}))\gamma_{\psi_{\pm\lambda}}(q_{a,\pm\lambda}). \]

	(5) Consider $s_\lambda=\lambda(s) \neq 1$. Recall the convention $\psi_{\pm\lambda}=\psi\circ\tx{tr}_{F_{\pm\lambda}/F}$. Keeping $s_\lambda$ fixed, we can only change the pair $(\tilde s_\lambda,\tilde s_{-\lambda})$ to $(a_\lambda\tilde s_\lambda,a_\lambda\tilde s_{-\lambda})$ for some $a_\lambda \in F_{\pm\lambda}^\times$. According to Lemma \ref{lem:weilherm} we have 
	\[ \gamma_{\psi_{\pm\lambda}}(\tx{tr}(V_\lambda,h_{a\tilde s,\lambda})) = \kappa_\lambda(a_\lambda)\gamma_{\psi_{\pm\lambda}}(\tx{tr}(V_\lambda,h_{\tilde s,\lambda})).\]
	This shows that the factor in the product decomposition of (4) corresponding to $\lambda$ is genuine as a function of $\tilde s_\lambda$. Write 
	\[ T_{\GL}(F)^{\dagger}=\{(\tilde s_\lambda)|s_\lambda = 1 \Rightarrow \tilde s_\lambda \in N_{F_\lambda/F_{\pm\lambda}}(F_\lambda^\times)\}. \]
	This is a subset of $T_{\GL}(F)$ that is stable under multiplication by the subgroup $\tx{ker}(T_{\GL}(F) \to T_{\Sp}(F)_{\pm\pm})$. We just saw that the function 
	\[ T_{\GL}(F)^{\dagger} \to \mu_4, \quad \tilde s \mapsto \gamma_\psi(q_{\tilde s}) \] 
	is invariant under multiplication by $\tx{ker}(T_{\GL}(F) \to T_{\Sp}(F)_{\pm\pm})$, hence descends to $\tx{im}(T_{\GL}(F)^\dagger \to T_{\Sp}(F)_{\pm\pm})$. Since moreover $\gamma_\psi(q_{a\tilde s})=\kappa(a)\gamma_\psi(q_{\tilde s})$ whenever $a \in T_{\GL}^+(F)$ and $\tilde s,a\tilde s \in T_{\GL}(F)^\dagger$, this function further descends to a genuine function on $\tx{im}(T_{\GL}(F)^\dagger \to T_{\Sp}(F)_{\pm})$. Since 
	\[ T_{\Sp}(F)_\pm = \tx{im}(T_{\GL}(F)^\dagger \to T_{\Sp}(F)_\pm) \cdot \tx{ker}(T_{\Sp}(F)_\pm \to T_{\Sp}(F)) \]
	the claim follows.
\end{proof}

Given a polarizable twisted Levi subgroup $H_{\Sp} \subseteq \Sp(V)$ let $\Lambda \subset X^*(Z(H_{\Sp}))$ be the set of weights for the action of $Z(H_{\Sp})$ on $V$, as in Lemma \ref{lem:Levivsgrading}. Then $\Lambda^-$ is an earnest $\Sigma$-set, which injects into $X^*(H_{\Sp,\tx{ab}})$ by sending $\lambda$ to the character $\det(-|V_\lambda)$.

\begin{dfn} \label{dfn:tsl-cov}We define $H_{\Sp}(F)_\pm$ to be the pull back under $H_{\Sp}(F) \to H_{\Sp,\tx{ab}}(F)$ of the double cover of $H_{\Sp,\tx{ab}}(F)$ associated to the image of $\Lambda^-$ in $X^*(H_\tx{ab})$ by \cite[\S3.1]{KalDC}.
\end{dfn}

\begin{rem}
	We can describe this cover explicitly as follows. From Lemma \ref{lem:Levivsgrading} we see that the abelianization of $H_{\Sp}$ is realized as
	\[ 	H_{\Sp}^\tx{ab} \cong \prod_{\lambda \in \Lambda_\tx{asym}/\Sigma} \tx{Res}_{F_\lambda/F}\GL(1)
			\times\prod_{\lambda \in \Lambda_\tx{sym}/\Sigma} \tx{Res}_{F_{\pm\lambda}/F}\U_{F_\lambda/F_{\pm\lambda}}(1) \]
	where the isomorphism is given by the determinant maps $\GL(V_\lambda) \to \GL(1)$ and $\U(V_{\lambda},h_{\lambda}) \to \U_{F_\lambda/F_{\pm\lambda}}(1)$. From the double covers $F_\lambda^\times/N_{F_\lambda/F_{\pm\lambda}}(F_\lambda^\times) \to F_\lambda^1=\U_{F_\lambda/F_{\pm\lambda}}(1)(F_{\pm\lambda})$ and the split covers $F_\lambda^\times \times \{\pm1\} \to F_\lambda^\times=\GL(1)(F_\lambda)$ we obtain a double cover of $H_{\Sp}^\tx{ab}(F)$ which we pull back to $H_{\Sp}(F)$.
\end{rem}

\begin{dfn}\label{dfn:norm}
	Let $H_{\Sp} \subseteq \Sp(V)$ be a polarizable twisted Levi subgroup with associated $\Sigma$-set $\Lambda$, let $\mc{L}$ be a $\Lambda$-atomic lattice chain in $V$, and let $K_\pm \subset H_{\Sp}(F)_\pm$ be the preimage of $K=H_{\Sp}(F)_\mc{L}$.
A topologically semi-simple element $\dot s \in K_\pm$ is called a \emph{norm}\index[terminology]{norm}, if there exists a tame maximal torus $T_{\Sp} \subseteq H_{\Sp}$ that decomposes $\mc{L}$ so that $\dot s$ is in the image of $T_{\GL}(F) \to T_{\Sp}(F)_\pm \to H_{\Sp}(F)_\pm$. \end{dfn}

\begin{lem} \label{lem:qsd-tsl}
	Let 
	\begin{enumerate}
		\item $H_{\Sp} \subset \Sp(V)$ be a (twisted or usual) Siegel Levi subgroup with quadratic \'etale algebra $E/F$,
		\item $V_E = V_+^E \oplus \sigma(V_+^E)$ the unique polarization stabilized by $H_{\Sp}$,
		\item $h : V_+^E \times V_+^E \to E$ the skew-Hermitian form $h(v,w)=\eta(v,\sigma(w))$.
		\item $\dot s \in H_{\Sp}(F)_\pm$ a semi-simple element that is a norm,
\item $T_{\Sp} \subseteq H_{\Sp}$ a maximal torus with $s \in T_{\Sp}(F)$ and $\tilde s \in T_{\GL}(F)$ mapping to $\dot s$,
		\item $\tilde s^*$ be the adjoint of $\tilde s$ with respect to $\eta$, i.e., $\eta(\tilde sv,w)=\eta(v,\tilde s^*w)$.
	\end{enumerate}
	Assume $\ker(s-1)=\{0\}$. Then
	\begin{enumerate}
		\item 
		\[ \det(\tilde s-\tilde s^*|V_+^E) = \det(\tilde s|V_+^E) \cdot \det(1-s^{-1}|V_+^E), \]
		in particular the image in $E^\times/N_{E/F}(E^\times)$ of this element depends only on $\dot s$ and not on $\tilde s$ or $T_{\Sp}$.
	\item Writing $E^0=\mathrm{ker}(\mathrm{tr}_{E/F} : E \to F)$, we have
	\[ \det(\tilde s-\tilde s^*|V_+^E), \det(h) \in \begin{cases}
		E^0/N_{E/F}(E^\times),&2\nmid\dim_E(V_+^E),\\
		F^\times/N_{E/F}(E^\times),&2\mid \dim_E(V_+^E).
	\end{cases}\]
	\item We have 
	\[ \gamma_\psi(q_{\tilde s})=\kappa_{E/F}(\det(\tilde s-\tilde s^*|V_+^E) \cdot \det(h)) \cdot \epsilon_L(1/2,\kappa_{E/F},\psi)^{\dim_E(V_+^E)}.  \]
	\item Let $a \in E^\times$ be an element of trace zero and let $\chi : E^\times \to \C^\times$ be a character extending $\kappa : F^\times \to \{\pm1\}$ such that $\chi(a)=\epsilon_L(1/2,\kappa_{E/F},\psi)$. Then \[ \gamma_\psi(q_{\tilde s}) = \chi(\det(\tilde s|V_+^E))\cdot \chi(\det(1-s^{-1}|V_+^E))\cdot \kappa_{E/F}(\det(ah)).\]
	\end{enumerate}
	In particular, $\gamma_\psi(q_{\tilde s})$ depends on $\dot s$ alone, and may thus be denoted by $\gamma_\psi(q_{\dot s})$.
\end{lem}
\begin{proof}
	(1) Let $\Lambda$ correspond to $T$ and $\Lambda_0$ to $H$. Writing $\tilde s = (\tilde s_\lambda)_\lambda$ we infer from $V_\lambda \perp V_\mu$ for $\lambda \neq -\mu$ that $\tilde s^*_\lambda=\tilde s_{-\lambda}$. Since $s=N(\tilde s)$ we have $s_\lambda=\tilde s_\lambda/\tilde s_{-\lambda}$, hence $s=\tilde s/\tilde s^*$. This implies 
	\[ \det(\tilde s-\tilde s^*|V_+^E) = \det(\tilde s|V_+^E) \cdot \det(1-s^{-1}|V_+^E), \]
	and by construction of $H_{\Sp}(F)_\pm$ the term $\det(\tilde s|V_+^E)$ depends only on $\dot s$ modulo $N_{E/F}(E^\times)$.

	(2) By construction of $h$ we see that $\tilde s^*|_{V_+^E}$ is the adjoint to $\tilde s|_{V_+^E}$ with respect to $h$. Therefore $(\tilde s - \tilde s^*)|_{V_+^E}$ is skew-adjoint for $h$. The claim about $\det(\tilde s - \tilde s^*|V_+^E)$ follows. The claim about $\det(h)$ is immediate from the fact that $h$ is skew-Hermitian.

	(3) If $E/F$ is split the claim is immediate: the form $q_{\tilde s}$ is hyperbolic, while $\kappa_{E/F}$ is trivial. We may thus assume $E/F$ is non-split. The isomorphism 
	\[ V_+^E|_F \to V,\qquad \ell \mapsto \ell+\sigma(\ell) \]
	identifies the trace of $h$ with $\eta$, hence $q_{\tilde s}$ with the trace quadratic form of the Hermitian form 
	\[ h_{\tilde s} : V_+^E \times V_+^E \to E, \quad (v,w) \mapsto h((\tilde s-\tilde s^*)v,w). \]
	It follows that 
	\begin{eqnarray*}
		\gamma_\psi(q_{\tilde s})&=&\kappa_{E/F}(\det(h_{\tilde s}))\cdot \gamma_\psi(N_{E/F})^{\dim V_+^E}\\
		&=&\kappa_{E/F}(\det(\tilde s-\tilde s^*|V_+^E) \cdot \det(h)) \cdot \gamma_\psi(N_{E/F})^{\dim V_+^E}.
	\end{eqnarray*}
	According to Lemma \ref{lem:weilherm} we have $\gamma_\psi(N_{E/F})=\epsilon_L(1/2,\kappa_{E/F},\psi)$.

	(4) We have
	\begin{eqnarray*}
	\gamma_\psi(q_{\dot s})&=&\kappa_{E/F}(\det(\tilde s-\tilde s^*|V_+^E) \cdot \det(h)) \cdot \epsilon_L(1/2,\kappa_{E/F})^{\dim_E(V_+^E)}\\
	&=&\kappa_{E/F}(\det(a^{-1}(\tilde s-\tilde s^*)|V_+^E))\cdot \kappa_{E/F}(\det(ah)) \cdot \epsilon_L(1/2,\kappa_{E/F})^{\dim_E(V_+^E)}\\
	&=&\chi(\det(\tilde s-\tilde s^*|V_+^E))\cdot \kappa_{E/F}(\det(ah))\\
	&=&\chi(\det(\tilde s|V_+^E))\cdot \chi(\det(1-s^{-1}|V_+^E))\cdot \kappa_{E/F}(\det(ah)),
	\end{eqnarray*}
	where the first equality is due to (3), the second because $a^{-1}(\tilde s-\tilde s^*)$ is self-adjoint (with respect to both $h$ and $h_a$) so its determinant lies in $F^\times$, and the determinant of $h_a$ lies in $F^\times/N_{E/F}(E^\times)$, the third because by construction $\chi|_{F^\times}=\kappa_{E/F}$ and $\chi(a)=\epsilon_L(1/2,\kappa_{E/F})$, and the fourth from the formula in (1).
\end{proof}

\begin{cor} \label{cor:topsgamma-1}
	Let $H \subset \Sp(V)$ be a polarizable twisted Levi and $\dot s \in H(F)_\pm$ a semi-simple element with image $s \in H(F)$. Choose a maximal torus $T \subset H$ with $s \in T(F)$. Then the element $\gamma_\psi(q_{\dot s}) \in \C^\times$ of Lemma \ref{lem:topssgamma}(5) is independent of the choice of $T$, and depends only on $\dot s$.
\end{cor}
\begin{proof}
	According to Lemma \ref{lem:Levivsgrading} there is a $\Sigma$-decomposition of $V$ associated to $H$. By passing to each $\Sigma$-orbit of the indexing set $\Lambda$ and applying restriction of scalars (using Fact \ref{fct:gammaind}) we reduce to the case that $H$ is a (split or twisted) Siegel Levi subgroup, where we can apply Lemma \ref{lem:qsd-tsl}(1)--(3).
\end{proof}

\subsection{Main results} \label{sub:sympmain}

Let $(V,\eta)$ be symplectic vector space over $F$, $G=\Sp(V,\eta)$, $H$ a polarizable twisted Levi subgroup that splits over a tame extension, and $\Lambda$ the corresponding $\Sigma$-set with associated $\Lambda$-grading of $V$ as in Definition \ref{dfn:lambda_grading}. Let $\psi : F \to \C^\times$ be a character non-trivial on $\cO_F$ but trivial on $\mf{p}_F$. Its restriction to the ring of integers factors factors through the non-trivial character $\bar\psi$ of $\rf$.
Let $\mc{L}$ be a $\Lambda$-atomic lattice chain in $V$ and let $(\oVS,\ol{\eta})$ be its $\Lambda$-central quotient.  
We have the natural map $\bar\jmath : H_\mc{L}(\rf) \to \Sp(\oVS,\ol{\eta})$.

Recall from \S\ref{sub:wittmeta} the universal Witt-metaplectic group $\WMuniv(\oVS,\ol{\eta})$ and from \S\ref{sub:classmeta} the classical metaplectic group, which in this case specializes to $\Sp(\oVS,\ol{\eta})$ but is equipped with a section into $\WMuniv(\oVS,\ol{\eta})$ given by a map $c : \Sp(\oVS,\ol{\eta}) \to \QW(\rf)$,\index[notation]{c@$c$} see Corollary \ref{cor:unisec}.

In this section we will compare the Weil indices of $(q_{\dot s}|V)$ and $(c(s)|\oVS)$ for topologically semi-simple elements $\dot s \in H(F)_{\mc{L},\pm}$ that are norms. The calculation will be done in two steps: first in the case when $H$ is a twisted Siegel Levi subgroup (Lemma \ref{lem:tsl-gammatilde}), and then in the general case (Proposition \ref{pro:orangebox}).

From this we will derive the key theorem (Theorem \ref{thm:symplectic-main}), which is the existence and uniqueness of a genuine homomorphism $\tilde\jmath : H_\mc{L}(\rf)_\pm \to \WMuniv_\psi(\oVS,\ol{\eta})$ that lifts $\bar\jmath$ and satisfies $\tilde\jmath(\dot s) = \gamma_\psi(q_{\dot s})\boxtimes \bar\jmath(s)$ for all topologically semi-simple elements $\dot s \in H(F)_{\mc{L},\pm}$ that are norms.

\begin{lem} \label{lem:tsl-gammatilde}
	Assume that $H$ is a (potentially split) twisted Siegel Levi subgroup, $E/F$ the associated quadratic \'etale algebra, $V_E = V_+^E \oplus \sigma(V_+^E)$ the unique $H$-stable polarization. Let $\mc{L}^E$ be the atomic  lattice chain in the Hermitian space $(V_+^E,h)$ whose restriction to $\cO_F$ is $\mc{L}$. Let $0 \neq a \in E^0$ and $\chi : E^\times \to \C^\times$ be determined as follows:
	\begin{enumerate}[label=(\alph*)]
		\item $E/F$ is split: $\mathrm{ord}_E(a)=0$, $\chi=1$, \item $E/F$ non-split unramified: $\mathrm{ord}_E(a)=0$, $\chi(x)=(-1)^{\mathrm{ord}_E(x)}$,
		\item $E/F$ is non-split ramified: $a$ is unconstrained, $\chi|_{F^\times}=\kappa_{E/F}$, $\chi(1+\mf{p}_E)=1$, $\chi(a)=\epsilon_L(1/2,\kappa_{E/F},\psi)$.
	\end{enumerate}
	Let $\eta_a$ be the symmetric bilinear form on $V$ as in Lemma \ref{lem:askew}(2) and
let $(\oVO,\ol{\eta_a})$ be the symmetric bilinear $\rf$-space arising as the $(\eta_a,\Lambda)$-central quotient of $\mc{L}$.

	For any topologically semi-simple $\dot s \in H(F)_{\mc{L},\pm}$ that is a norm and $\ker(s-1) = \{0\}$ we have 
	\[ \gamma_\psi(q_{\dot s}) = \tilde\gamma_{\bar\psi}(c(s)|\oVS) \cdot \chi(\det(\tilde s|V_+^E)) \cdot  \mathrm{sgn}_\rf(\mathrm{sn}(s|\ol{q_a})),  \]
	where $\ol{q_a}(v)=\ol{\eta_a}(v,v)$ on $\oVO$ and $\mathrm{sn}$ denotes the spinor norm.
\end{lem}
\begin{proof}
	We apply Lemma \ref{lem:qsd-tsl}(4) but with $a$ replaced by $a'$ defined as $a'=a$ unless $E/F$ non-split unramified, in which case $a'=a\pi$ for a uniformizer $\pi \in F^\times$; then we have $\chi(a')=\epsilon(1/2,\kappa_{E/F},\psi)$ as required by that lemma. It, together with Proposition \ref{pro:metasection}(2), and Lemma \ref{lem:weilindexsubqout}(4) to reduce to showing that $\tx{sgn}_\rf(\tx{sn}(s|ol{q_a}))$ equals
	\[ \chi(\det(1-s^{-1}|V_+^E))\cdot \kappa_{E/F}(\det(a'h)) \cdot \tx{sgn}_\rf\big(\det(s-1|\oVS)\cdot(-1)^{\dim(\oVS)/2}\big)^{-1}.\]
	To that end, we first compute $\tx{sn}(s|\ol{q_a})$. We have
	\[ \tx{sn}(s|\ol{q_a})=\tx{sn}(-1|\ol{q_a})\cdot \tx{sn}(-s|\ol{q_a}) = \det(\ol{q_a}) \cdot \tx{sn}(-s|\ol{q_a}). \]
	By assumption $s$ induces an orthogonal transformation of $\oVO$ that has no eigenvalue $1$, hence $-s$ induces an orthogonal transformation of $\oVO$ that has no eigenvalue $-1$ and to which therefore Zassenhaus' formula applies (see \cite[\S2]{Zass1962} for the original result, and \cite[Theorem 3.5]{CT2020} for more modern language). This shows
	\begin{equation} \label{eq:sn1}
		\tx{sn}(s|\ol{q_a}) = \det(\ol{q_a})\det\left(\frac{1-s}{2}\Big| \oVO \right) = \det(\ol{q_a}/2)\det(1-s|\oVO).	
	\end{equation}
	Next, we claim that 
	\begin{equation} \label{eq:sn2}
		\kappa_{E/F}(\det(a'h)) = \tx{sgn}_\rf(\det(\ol{q_a}/2) \cdot (-1)^{\dim(\oVS)/2}).	
	\end{equation}
	Postponing the proof for now, \eqref{eq:sn1} and \eqref{eq:sn2} reduce the proof to
	\[ \chi(\det(1-s^{-1}|V_+^E)) = \tx{sgn}_\rf\Big(\det(1-s|\oVO)\cdot\det(s-1|\oVS)\Big). \]

	The eigenvalues of $s$ on $V_+^E$, hence also of $s^{-1}$ on $L$, are roots of unity of order prime to $p$, and never equal to $1$ by assumption, so $\det(1-s^{-1}|V_+^E) \in \cO_E^\times$. 
	
	We now consider two cases. If $E/F$ is split or unamified then $\chi$ is trivial on $\cO_E^\times$ by assumption and $\tx{ord}_E(a)=0$. Lemma \ref{lem:sign-switch} implies that the vector spaces $\oVS$ and $\oVO$ are equal (and even-dimensional since they carry a symplectic form), so the argument of $\tx{sgn}_\rf$ is a square. Thus both sides of the above identity are equal to $1$.

	If $E/F$ is non-split and ramified then the restriction of $\chi$ to $\cO_E^\times$ factors through $\tx{sgn}_\rf$, so 
	\[ \chi(\det(1-s^{-1}|V_+^E))=\tx{sgn}_\rf(\ol{\det(1-s^{-1}|V_+^E)}) = \tx{sgn}_\rf(\ol{\det(1-s^{-1}|L/\mf{p}_EL)}), \]
	where the overline denotes the projection map $\cO_E^\times \to \rf^\times$ and $L \subset V_+^E$ is any $s$-stable lattice. We choose $L$ to be the $(a'h)$-center of $\mc{L}^E$.  Lemma \ref{lem:sign-switch} produces the exact sequence 
	\[ 0 \to \oVS \to L/\mf{p}_E L \to \oVO \to 0, \]
	which shows
	\[ \ol{\det(1-s^{-1}|V_+^E)} = \det(1-s^{-1}|L/\mf{p}_E L) = \det(1-s^{-1}|\oVS) \cdot \det(1-s^{-1}|\oVO). \]
	We have
	\begin{eqnarray*}
	\det(1-s^{-1}|\oVS)&=&\det(s^{-1}|\oVS)\det(s-1|\oVS)=\det(s-1|\oVS)\\
	\det(1-s^{-1}|\oVO)&=&\det(s^{-1}|\oVO)\det(s-1|\oVO)=\det(1-s|\oVO),
	\end{eqnarray*}
	where in the first line we have used that $s|\oVS$ is a symplectic transformation, so its determinant is $1$, and in the second line we have used that $s|\oVO$ is an orthogonal transformation without eigenvalue $1$, which implies $\det(s^{-1}|\oVO)=(-1)^{\dim(\oVO)}$, since the eigenvalues of $s|\oVO$ are either $-1$ or pairs of reciprocial elements of $\bar \rf^\times$.

	We now come to the proof of \eqref{eq:sn2}. We first treat the case that $E/F$ is split. Then $\kappa_{E/F}$ is trivial. The split polarization $V_+^E=V_+\oplus V_-$ induces a polarization of $\oVO=\oVS$ for the non-degenerate quadratic form $\ol{q_a}/2$. Thus this form is hyperbolic, and $\det(\ol{q_a}/2)=(-1)^{\dim(\oVS)/2}$ in $\rf^\times/\rf^{\times,2}$,
	which proves \eqref{eq:sn2} in the split case.

	Next we come to the case that $E/F$ is non-split. Choose a uniformizer $\varpi \in E$ with $\sigma(\varpi) = \pm \varpi$ with sign $+1$ or $-1$ depending on whether $E/F$ is unramified or ramified, respectively. Let $L \in \mc{L}$ be the $(a'h)$-center. According to \cite[Proposition 4.3]{Jacobowitz1962} we can write $L = L_1 \oplus L_2$ such that $V_+^E=V_{+,1}^E \oplus V_{+,2}^E$ is an $a'h$-orthogonal decomposition, $L_i$ is a lattice in $V_{+,i}^E$, and $L^\vee=L_1 \oplus \varpi^{-1}L_2$. The Gram matrix of $a'h$ with respect to an adapted $\cO_E$-basis then has the form 
	\[ \begin{bmatrix}
		A&0\\ 0&\sigma(\varpi) B
	\end{bmatrix}, \]
	where $A$ reduces mod $\mf{p}_E$ to the matrix of the form $\ol{a'h}$ on $\ol{L}$ while $B$ reduces mod $\mf{p}_E$ to the matrix of the form $\ol{\varpi a'h}$ on $\ol{L^{a'h}}$. Both of these matrices have determinants in $\cO_E^\times$ and we see
	\[ 	\kappa_{E/F}(\det(a'h))=\kappa_{E/F}(\det(A) \cdot \varpi^{\dim_{\rf_E}(\ol{L^{a'h}})}\cdot\det(B)). \]
	If $E/F$ is unramified then $a'=a\varpi$, which implies $\det(A)\cdot\det(B) \in \cO_F^\times$. Since $\kappa_{E/F}$ is the unramified quadratic character on $F^\times$ we see
	\[ \kappa_{E/F}(\det(a'h)) = (-1)^{\dim_{\rf_E}(\ol{L^{a'h}})}. \]
	By Lemma \ref{lem:sign-switch} $\varpi L^{a'h}$ is the center of $\mc{L}$ with respect to both $h$ and $(ah)$ and multiplication by $\varpi$ identifies $\ol{L^{a'h}}=L^{a'h}/L$ with the $ah$-central quotient, on which we obtain a symmetric bilinear form $\ol{\eta_a}$ by taking the trace of $ah$. Thus $(-1)^{\dim_{\rf_E}(\ol{L^{a'h}})}=(-1)^{\dim_{\rf}(\oVO)/2}=(-1)^{\dim_{\rf}(\oVS)/2}$. It is therefore enough to show 
	\[ \tx{sgn}_\rf(\det(\ol{q_a}/2)) = (-\tx{sgn}_\rf(-1))^{\dim_{\rf}(\oV_O)/2}. \] 
	This is a general statement about the trace quadratic form of a Hermitian form, and after diagonalization reduces to the 1-dimensional Hermitian form $cx\bar y$ on $\rf_E$ for $c \in \rf^\times$. Its determinant is $-c^2a^2$, which up to $c^2 \in \rf^{\times,2}$ equals the norm of $a$. Now $\tx{sgn}_{\rf}(-a^2)=\tx{sgn}_{\rf_E}(a)=-\tx{sgn}_\rf(-1)$.

	If $E/F$ is ramified then we have $a'=a$, $\ol{L}=\oVO$, and 
	\[ \ol{q_a}(v)=\ol{\eta_a}(v,v)=\overline{ah(v,v)+\sigma(ah(v,v))} = 2\overline{ah(v,v)}, \]
	hence $A$ reduces to the matrix of $\ol{q_a}/2$. On the other hand, Lemma \ref{lem:sign-switch} shows that $B$ reduces to the matrix of $\bar\eta$ on $\oVS$, and in particular has determinant equal to $1$ in $k^\times/k^{\times,2}$. Since $\kappa_{E/F}(\varpi^2)=\kappa_{E/F}(-1)=\tx{sgn}_\rf(-1)$, we conclude
	\begin{align*}
	\kappa_{E/F}(\det(ah))&=\kappa_{E/F}(\det(A) \cdot \varpi^{\dim\ol{V}_S}\det(B)) \\
	&=\tx{sgn}_\rf(\det(A) \cdot \det(B) \cdot (-1)^{\dim\oVS/2})\\
	&=\tx{sgn}_\rf(\det(\ol{q_a}/2) \cdot (-1)^{\dim\oVS/2}).\hfill\qedhere
	\end{align*}
\end{proof}

We now return to a general polarizable twisted Levi subgroup $H \subseteq G$. It is associated to a finite polarizable $\Sigma$-set $\Lambda$ and hence has the associated double cover $H(F)_\pm$ of Definition \ref{dfn:tsl-cov}. This cover splits uniquely over $H(F)_{\mc{L},0+}$ and hence induces a double cover of $H_\mc{L}(\rf)=H(F)_\mc{L}/H(F)_{\mc{L},0+}$. We write $H_\mc{L}(\rf)_\pm$ for it. 

Let $(a_\lambda)_\lambda$ be a set of $a$-data for $\Lambda$. For each $\lambda$ we have the skew-Hermitian form $h_\lambda^E : V_\lambda^E \times V_\lambda^E \to E_{\pm\lambda}$ of Definition \ref{dfn:lambda_grading} and obtain the Hermitian form $a_\lambda^E h_\lambda^E : V_\lambda^E \times V_\lambda^E \to E_{\pm\lambda}$. Its trace with respect to $E_{\pm\lambda}/F_{\pm\lambda}$ provides a symmetric bilinear form on $V_\lambda^E|_{F_{\pm\lambda}}$.

\begin{pro} \label{pro:orangebox}
	 Let $(a_\lambda)$ and $(\chi_\lambda)_\lambda$ be sets of $a$-data and $\chi$-data for $\Lambda$ (Definition \ref{dfn:achi}) such that
	\begin{enumerate}
		\item $a_\lambda \in \cO_\lambda^\times$ and $\chi_\lambda=1$ if $\lambda \in \Lambda_\mathrm{asym}$,
		\item $a_\lambda \in \cO_\lambda^\times$ and $\chi_\lambda(x)=(-1)^{\mathrm{ord}_{F_\lambda}(x)}$ if $\lambda \in \Lambda_\mathrm{sym,unram}$,
		\item $a_\lambda$ any, $\chi_\lambda(1+\mf{p}_E)=1$, and $\chi_\lambda(a_\lambda)=(-1)^{[\rf_{\pm\lambda}:\rf]+1}\gamma_{\bar\psi\circ \mathrm{tr}_{\rf_{\pm\lambda}/\rf}}([1])$, if $\lambda \in \Lambda_\mathrm{sym,ram}$.
	\end{enumerate}
	Let $\chi_H : H(F)_\pm \to \C^\times$ be the genuine character associated to this $\chi$-data as in \cite[\S3.2]{KalDC}. Let $\eta_a : V \times V \to F$ be the symmetric bilinear form of Lemma \ref{lem:askew}(2) and let $(\ol{V_a},\ol{\eta_a})$ be the symmetric bilinear $\rf$-space obtained from the $(\Lambda,\eta_a)$-central quotient of $\mc{L}$. Let $\ol{q_a}(v)=\ol{\eta_a}(v,v)$.
	
	Then for every topologically semi-simple elements $\dot s \in H(F)_{\mc{L},\pm}$ that is a norm we have
	\begin{equation}\label{eqn:orangebox} \gamma_\psi(q_{\dot s}|V^{s \neq 1}) = \tilde\gamma_{\bar\psi}(c(s)|\ol{V}^{s \neq 1}) \cdot \mathrm{sgn}_\rf(\mathrm{sn}(s|\ol{q_a})) \cdot \chi_H(\dot s). 
		\end{equation}
\end{pro}
\begin{proof}
	Since both sides of the identity are genuine functions in $\dot s$, we may always assume in what follows that $\dot s$ is a norm.

	\ul{Step 1: Reduce to the case $\tx{ker}(s-1)=\{0\}$.}\\
	Let $V^{s=1}=\tx{ker}(s-1)$ and $V^{s \neq 1}=\tx{im}(s-1)$. We have the orthogonal decomposition $V=V^{s=1} \oplus V^{s \neq 1}$. The $\Sigma$-grading of $V$ by $\Lambda$ given by the $Z(H)$-decomposition of $V$ into weights is compatible with this decomposition, i.e. $V_\lambda=(V_\lambda)^{s=1} \oplus (V_\lambda)^{s \neq 1}$ and $(V_\lambda)^{s=1}=(V^{s=1})_\lambda$, $(V_\lambda)^{s \neq 1}=(V^{s \neq 1})_\lambda$. It is possible that some of these weight spaces become zero, so we let $\Lambda^{s=1}, \Lambda^{s \neq 1} \subset \Lambda$ consist of those $\lambda$ for which $(V_\lambda)^{s=1} \neq \{0\}$ resp. $(V_\lambda)^{s \neq 1} \neq \{0\}$. These decompositions determine polarized twisted Levi subgroups $H^{s=1} \subset G^{s=1}=\Sp(V^{s=1})$ and $H^{s \neq 1} \subset G^{s \neq 1}=\Sp(V^{s \neq 1})$. We have $G^{s=1} \times G^{s \neq 1} \subset G$ and $H^{s=1} \times H^{s \neq 1} = H \cap (G^{s=1} \times G^{s \neq 1})$.
	
	Each lattice $L \in \mc{L}$ also decomposes as $L=L^{s=1} \oplus L^{s \neq 1}$ and we let $\mc{L}^{s=1}$ and $\mc{L}^{s \neq 1}$ be the lattice chains generated by the respective summands, as $L \in \mc{L}$ varies. Then $\mc{L}^{s=1}$ is $\Lambda^{s=1}$-atomic and $\mc{L}^{s \neq 1}$ is $\Lambda^{s \neq 1}$-atomic.

	The terms $\gamma_\psi(q_{\dot s}|V^{s \neq 1})$ and $\tilde\gamma_{\bar\psi}(c(s)|\ol{V}^{s \neq 1})$ are unchanged if we pass from $V$ to $V^{s\neq 1}$, once we note that $\ol{V^{s\neq 1}}=\ol{V}^{s \neq 1}$. We must check that the terms $\tx{sn}(s|\ol{q_a})$ and $\chi(\dot s)$ are also unchanged. For the term $\tx{sn}(s|\ol{q_a})$ we note that we have the $s$-stable orthogonal decomposition $\ol{V_a} = \ol{(V^{s=1})_a} \oplus \ol{(V^{s \neq 1})_a}$, hence $\tx{sn}(s|\ol{q_a})=\tx{sn}(s|\ol{q_a^{s \neq 1}}) \cdot \tx{sn}(s|\ol{q_a^{s=1}})$ and moreover $\tx{sn}(s|\ol{q_a^{s=1}})=1$ since $s=1$ on $(V^{s=1})_a$. For the term $\chi(\dot s)$ we note that the inclusion $H^{s \neq 1} \to H$  lifts naturally to an inclusion $H^{s \neq 1}(F)_\pm \to H(F)_\pm$. The element $\dot s$ lies in the image of this lift, since $s=1$ in $H^{s=1}(F)$. Pulling back $\chi_H$ under this lift produces the character $\chi_{H^{s \neq 1}}$, and we conclude $\chi_H(\dot s)=\chi_{H^{s \neq 1}}(\dot s)$. \\

	\ul{Step 2: Reduce to the case of a Siegel Levi subgroup.}\\
	We have the orthogonal decomposition $(V,\eta)=\bigoplus_{O \in \Lambda/\Sigma} (V_O,\eta_O)$ of Fact \ref{fct:lambda_grading}. The symplectic space $(V_O,\eta_O)$ is endowed with an $O$-grading, which leads to the twisted Levi subgroup $H_O \subset G_O=\Sp(V_O,\eta_O)$. The twisted Levi subgroup $H$ is the product $H=\prod_O H_O$ and $H(F)_\pm$ is the Baer sum of the covers $H_O(F)_\pm$. Choosing a presentation $\dot s = (\dot s_O)_O$ we thus obtain
	\[ \gamma_\psi(q_{\dot s}|V) = \prod_O \gamma_\psi(q_{\dot s_O}|V_O),\qquad \chi_H(\dot s) = \prod_O \chi_{H_O}(\dot s_O). \]
	
	For each $L \in \mc{L}$ define $L_O = L \cap V_O$ and let $\mc{L}_O=\{L_O|L \in \mc{L}\}$. Then $\mc{L}_O$ is $O$-atomic, the $(\Lambda,\eta)$-center of $\mc{L}$ is the direct sum of the $(O,\eta_O)$-centers of $\mc{L}_O$, and the $(\Lambda,\eta_a)$-center of $\mc{L}$ is the direct sum of the $(O,\eta_{O,a})$-centers of $\mc{L}_O$. This implies 
	\[ \tilde\gamma_{\bar\psi}(c(s)|\ol{V}) = \prod_O \tilde\gamma_{\bar\psi}(c(s)|\ol{V}_O),\qquad \tx{sn}(s|\ol{q_a}) = \prod_O \tx{sn}(s|\ol{q_{O,a}}). \]
	This reduces to the case that $\Lambda$ is a transitive $\Sigma$-set. In that case, choosing $\lambda \in \Lambda$ we obtain the isomorphism
	\[ V_{\pm\lambda}|_F \to V,\qquad v \mapsto \sum_{\sigma \in \Gamma/\Gamma_{\pm\lambda}}\sigma(v) \]
	which identifies $\eta$ with the trace of $\eta_{\pm\lambda}$, $\eta_a$ with the trace of $\eta_{a,\pm\lambda}$, and $H$ with $\tx{Res}_{F_{\pm\lambda}/F}\U(V_{\lambda}^E,h_\lambda^E)$. We have the $\{\pm\lambda\}$-atomic lattice chain $\mc{L}_{\pm\lambda}$ in $V_{\pm\lambda}$ and the isomorphism
	\[ \ol{V_{\pm\lambda}} \to \ol{V},\qquad \ol{V} \mapsto \sum_{\sigma \in \tx{Gal}(\rf_{\pm\lambda}/\rf)} \sigma(\ol{V}) \]
	which identifies $c(s)|\ol{V^{s\neq 1}}$ with the trace of $c(s)|\ol{V_{\pm\lambda}^{s \neq 1}}$. We have by Fact \ref{fct:gammaind}
	\[ \gamma_\psi(q_{\dot s}|V) = \gamma_{\psi\circ\tx{tr}_{F_{\pm\lambda}/F}}(q_{\dot s}|V_{\pm\lambda}),\qquad \tilde\gamma_{\bar\psi}(c(s)|\ol{V}) = \tilde\gamma_{\bar\psi\circ\tx{tr}_{\rf_{\pm\lambda}/\rf}}(c(s)|\ol{V}). \]
	We will use the character $\psi_{\pm\lambda}=\psi\circ\tx{tr}_{F_{\pm\lambda}/F}$ on the base field $F_{\pm\lambda}$. Then we have $\ol{\psi_{\pm\lambda}}(x) = \bar\psi\circ\tx{tr}_{\rf_{\pm\lambda}/\rf}(e_{\pm\lambda}\cdot x)$ for all $x \in \rf_{\pm\lambda}$. At the moment the discrepancy between this character and $\ol{\psi_{\pm\lambda}}$ does not matter: the Weil index of $c(s)$ is independent of the character being used, since the vector space is even-dimensional. But in a moment this discrepancy will make a difference.

	Turning to the spinor norm, we have the isomorphism
	\[ \ol{V_{\pm\lambda,a}}|_{\rf} \to \ol{V_a},\qquad \ol{V} \mapsto  \sum_{\sigma \in \tx{Gal}(\rf_{\pm\lambda}/\rf)} \sigma(\ol{V}) \]
	which identifies $\ol{q_a}$ with the trace of $e_{\pm\lambda}\cdot\ol{q_{\pm\lambda,a}}$. By \cite[Chapter 9, \S3, Example 3.7]{Sch85} we have
	\[ \tx{sn}(s|\ol{q_a}) = N_{\rf_{\pm\lambda}/\rf}(\tx{sn}(s|e_{\pm\lambda}\cdot \ol{q_{\pm\lambda,a}})) \cdot \Delta_{\pm\lambda}^{\epsilon(s)}, \]
	where $\Delta_{\pm\lambda} \in \rf^\times$ is the discriminant of the trace pairing for $\rf_{\pm\lambda}/\rf$ and $\epsilon(s) \in \Z/2\Z$ is determined by $\det(s|\ol{V_{\pm\lambda,a}})=(-1)^{\epsilon(s)}$. It is known that $\Delta_{\pm\lambda}$ is non-square precisely when $[\rf_{\pm\lambda}:\rf]$ is even, from which we conclude
	\[ \tx{sgn}_\rf(\tx{sn}(s|\ol{q_a})) = \tx{sgn}_{\rf_{\pm\lambda}}(\tx{sn}(s|\ol{q_{\pm\lambda,a}})) \cdot \Big(\tx{sgn}_{\rf_{\pm\lambda}}(e_{\pm\lambda})\cdot (-1)^{[\rf_{\pm\lambda}:\rf]+1}\Big)^{\epsilon(s)}. \]
	Note that, when $\lambda$ is not ramified symmetric, then $\ol{V_{\pm\lambda,a}}$ is even-dimensional, and the assumption $\tx{ker}(s-1)=0$ that is in place due to Step 1 implies $\det(s|\ol{V_{\pm\lambda,a}})=+1$, and hence $\epsilon(s)=0$, since all eigenvalues of $s$ that are not equal to $\pm 1$ come in reciprocal pairs.

	We turn to $\chi_H$. We have the isomorphism $\U(V_\lambda^E,h_\lambda^E)(F_{\pm\lambda}) \to H(F)$ which by construction of $H(F)_\pm$ gives an isomorphism $\U(V_\lambda^E,h_\lambda^E)(F_{\pm\lambda})_\pm \to H(F)_\pm$ which identifies $\chi_H$ with the chararcter $\chi_\lambda$ that is part of the $\chi$-data giving rise to $\chi_H$. When $\lambda$ is not ramified symmetric we see that switching from $V_O$ to $V_{\pm\lambda}$ has no effect on the definition of $\chi_\lambda$. But when $\lambda$ is symmetric ramified we obtain the following discrepancy: For $V_{\pm\lambda}$ the base field is $F_{\pm\lambda}$ rather than $F$, and from the point of view of $F_{\pm\lambda}$ we have $\rf_{\pm\lambda}=\rf$, which removes the factor $(-1)^{[\rf_{\pm\lambda}:\rf]}$ from the definition of $\chi_\lambda$, and moreover we would be using the character $\ol{\psi_{\pm\lambda}}$ in place of $\bar\psi\circ\tx{tr}_{\rf_{\pm\lambda}/\rf}$ in that same definition. Let us denote the resulting character of $F_\lambda^\times$ by $\chi_\lambda^{F_\lambda}$. Now both $\chi_\lambda$ and $\chi_\lambda^{F_\lambda}$ are tamely ramified characters of $F_\lambda^\times$ that extend $\kappa_\lambda$, hence their quotient is either trivial or the unique quadratic sign character. Since $a_\lambda$ has odd valuation, this character is determined by its value at $a_\lambda$. We have 
	\begin{eqnarray*}
		\chi_\lambda(a)&=&(-1)^{[\rf_{\pm\lambda}:\rf]+1}\cdot\gamma_{\bar\psi\circ \tx{tr}_{\rf_\lambda/\rf}}([1])\\
		&=&(-1)^{[\rf_{\pm\lambda}:\rf]+1}\cdot\gamma_{e_{\pm\lambda}\ol{\psi_{\pm\lambda}}}([1])\\
		&=&(-1)^{[\rf_{\pm\lambda}:\rf]+1}\tx{sgn}_{\rf_{\pm\lambda}}(e_{\pm\lambda})\cdot\gamma_{\ol{\psi_{\pm\lambda}}}([1])\\
		&=&(-1)^{[\rf_{\pm\lambda}:\rf]+1}\tx{sgn}_{\rf_{\pm\lambda}}(e_{\pm\lambda})\chi_\lambda^{F_\lambda}(a),
	\end{eqnarray*}
	hence
	\[ \Big(\chi_\lambda^{F_\lambda}/\chi_\lambda\Big)(x) = \Big((-1)^{[\rf_{\pm\lambda}:\rf]+1}\tx{sgn}_{\rf_{\pm\lambda}}(e_{\pm\lambda})\Big)^{\tx{ord}_{F_\lambda}(x)},\qquad x \in F_\lambda^\times. \]
	We now claim that $\epsilon(s)=\tx{ord}_{F_\lambda}(z_\lambda)$ in $\Z/2\Z$, where $z_\lambda \in F_\lambda$ is any element with $z_\lambda/\sigma_\lambda(z_\lambda)=\det_{F_\lambda}(s|V_\lambda)$. We have 
	\[ 
		\det_{\rf_{\pm\lambda}}(s|\ol{V_{\pm\lambda,a}})
		=\ol{\det_{F_\lambda}(s|V_{\lambda})}.
	\]
	We have $\det_{F_\lambda}(s|V_{\lambda}) \in F_\lambda^1$ and its projection to $\rf_\lambda^\times=\rf_{\pm\lambda}^\times$ is either $+1$ or $-1$. It is $-1$ precisely when $\tx{ord}_{F_\lambda}(z_\lambda)$ is odd. We conclude that the quantity $\tx{sgn}_\rf(\tx{sn}(s|\ol{q_a})) \cdot \chi_H(\dot s)$ remains unchanged when we switch from the symplectic $F$-space $(V_O,\eta_O)$ to the symplectic $F_{\pm\lambda}$-space $(V_{\pm\lambda},\eta_{\pm\lambda})$ and from the character $\psi$ to the character $\psi_{\pm\lambda}$.
	\ \\

	\ul{Step 3: The case of a Siegel Levi subgroup}.\\
	This is exactly Lemma \ref{lem:tsl-gammatilde}, using the fact that $\gamma_{\ol{\psi}}([1])=\epsilon(1/2,\kappa_{E/F},\psi)$ when $E/F$ is ramified.
\end{proof}

We now come to the main result of this section. Let $\WM_{\bar\psi}(\oVS,\bar\eta)$\index[notation]{WMpsibar@$\WM_{\bar\psi}$} denote the pushout under $[q] \mapsto \tilde\gamma_{\bar\psi}([q]) : \QW(\rf) \to \mu_8$ 
of the universal Witt-metaplectic group $\WMuniv(\oVS,\bar\eta)$ of \S\ref{sub:wittmeta}.

\begin{thm}\ \\[-20pt] \label{thm:symplectic-main}
	\begin{enumerate}
		\item The natural projection $\bar\jmath : H_\mc{L}(\rf) \to \Sp(\oVS,\bar\eta)(\rf)$ lifts uniquely to a genuine homomorphism $\tilde\jmath : H_\mc{L}(\rf)_\pm \to \WM_{\bar\psi}(\oVS,\bar\eta)$ with the property that for any semi-simple element $\dot s \in H_\mc{L}(\rf)_\pm$ we have
	\[ \tilde\jmath(\dot s) = \gamma_\psi(q_{\dot s}) \boxtimes \bar\jmath(s). \]
		\item Choose $a$-data and $\chi$-data for $\Lambda$ as in Proposition \ref{pro:orangebox} and let $(\ol{V_a},\ol{\eta_a})$ and $\chi_H$ be the resulting symmetric bilinear $\rf$-space and genuine character. Set again $\ol{q_a}(v)=\ol{\eta_a}(v,v)$. The following identity of isomorphism classes of representations holds, where $\mathrm{FWeil}$ and $\mathrm{LWeil}$ are the Fourier--Weil and Weil representations defined in Appendix \ref{sub:fwh}\[ \tilde\jmath^*\mathrm{FWeil}_{\bar\psi} = \bar\jmath^*\mathrm{LWeil}_{\bar\psi} \otimes \chi_H \otimes \mathrm{sgn}_\rf(\mathrm{sn}(-|\ol{q_a})). \]
\end{enumerate}
	
\end{thm}
\begin{proof}
	(1) \underline{Proof of uniqueness:} Assume $\tilde\jmath$ as given. Another such homomorphism is necessarily of the form $f \cdot \tilde\jmath$ for a homomorphism $f : H_\mc{L}(\rf) \to \mu_8$. Since $p \neq 2$ we see that $f$ kills all elements of $p$-power order, i.e. all unipotent elements of $H_\mc{L}(\rf)$. On othe other hand, every semi-simple element of $H_\mc{L}(\rf)$ is represented by a topologically semi-simple element of $H(F)_\mc{L}$. But every topologically semi-simple element of $H(F)_\mc{L}$ is a norm and the specification of $\tilde\jmath(\dot s)$ for topologically semi-simple norms $\dot s$ implies that $f(s)$ is trivial. This shows that $f$ is trivial.

	\underline{Proof of existence:} Every element of $H(F)_{\mc{L},\pm}$ has a topological Jordan decomposition $\dot su$, with $\dot s \in H(F)_{\mc{L},\pm}$ topologically semi-simple and $u \in H(F)_{\mc{L},\pm}$ topologically unipotent. Let $s \in H(F)_\mc{L}$ be the image of $\dot s$. For any $h \in H(F)_{\mc{L}}$ we define $V^{h=1}=\tx{ker}(h-1|V)$ and $V^{h \neq 1}=V/V^{h=1}$, which maps under $h-1$ isomorphically onto $\tx{im}(h-1|V)$. The same definitions apply to $\oVS$.
	Define
	\begin{equation} \label{eq:tildejansatz}
		\tilde \jmath(\dot su)=\gamma_\psi(q_{\dot s}|V^{s \neq 1}) \cdot \tilde\gamma_{\bar\psi}(q_u|\oVS^{s=1,u \neq 1}) \boxtimes \bar\jmath(su).
	\end{equation}
	This is evidently a set-theoretic map that descends to $H_\mc{L}(\rf)_\pm$ and lifts $\bar\jmath$. Its multiplicativity is equivalent to the multiplicativity of the map 
	\[ \dot f : H_\mc{L}(\rf)_\pm \to \mu_4,\qquad \dot su \mapsto \frac{\gamma_\psi(q_{\dot s}|V^{s \neq 1})\tilde\gamma_{\bar\psi}(c(u)|\oVS^{s=1,u \neq 1})}{\tilde\gamma_{\bar\psi}(c(su)|\oVS^{su \neq 1}])}. \]
	According to Proposition \ref{pro:metasection} we have
	\[ \tilde\gamma_{\bar\psi}(c(su)|\oVS^{su \neq 1}) = \tilde\gamma_{\bar\psi}(c(s)|\oVS^{s \neq 1}) \cdot \tilde\gamma_{\bar\psi}(c(u)|\oVS^{s=1,u \neq 1}) \] 
	and conclude $\dot f(\dot su)=\gamma_\psi(q_{\dot s}|V^{s \neq 1})/\tilde\gamma_{\bar\psi}(c(s)|\oVS^{s \neq 1})$. Proposition \ref{pro:orangebox} shows that 
	\[ \dot f(\dot su) = \epsilon(s)\chi_H(\dot s) = \epsilon(su)\chi_H(\dot su), \]
	with $\epsilon(h)=\tx{sgn}_\rf(\tx{sn}(h|\ol{q_a}))$, 
	where the second equality follows from the fact that both $\epsilon$ and $\chi_H$ are multiplicative and trivial on topologically unipotent elements. This shows the multiplicativity of $\dot f$. 

	(2) Let $\iota : \Sp(\oV) \to \WM_{\bar\psi}(\oV)$ denote the inclusion discussed in \S\ref{sub:classmeta}. According to Lemma \ref{lem:fweilheisrest} the right hand side is 
	\[ (\bar\jmath^*\iota^*\mathrm{FWeil}_{\bar\psi}) \otimes \chi_H \otimes \epsilon. \]
	But Proposition \ref{pro:orangebox} shows $\tilde\jmath = (\iota \circ \bar\jmath) \otimes \chi_H \otimes \epsilon$.
\end{proof}

\section{Symplectic spaces associated to good elements}
\label{sec:symplectictopadic}

In this subsection we will see how Theorem \ref{thm:symplectic-main} applies in the setting of an arbitrary connected reductive group and a good element of its dual Lie algebra. 

\subsection{The setting}
We keep the notation for the local field $F$ the same as for \S\ref{sec:symp}. Thus, $F$ is a non-archimedean local field, $\cO$ is the ring of integers, $\mf{p}$ the maximal ideal, and $\rf=\cO/\mf{p}$ the residue field, whose characteristic we assume odd. We denote by $\mathrm{ord}: F \twoheadrightarrow \bZ \cup \{\infty\}$ the normalized valuation on $F$ and also its extension to algebraic field extensions of $F$.

Unlike \S\ref{sec:symp}, we now take $G$ to be an arbitrary connected reductive $F$-group and $H \subseteq G$ a tame twisted Levi subgroup. We write $R(T,G)$ and $R(T,H)$ for the absolute root systems of $G$ and $H$ with respect to a maximal torus $T$ of the respective group, and $R(T,G/H)=R(T,G) \smallsetminus R(T,H)$. If $G_\mathrm{sc}$ is the algebraic simply connected cover of the derived subgroup of $G$ and $T_\tx{sc}$ is the preimage of $T$ in $G_\tx{sc}$, the natural map $T_\tx{sc} \to T$ gives an identification $R(T_\tx{sc},G_\tx{sc})=R(T,G)$. Therefore, for $\alpha \in R(T,G)$ we have $H_\alpha=d\alpha^\vee(1) \in \tx{Lie}(T_\tx{sc})$.

Write $H_\mathrm{sc}$ for the inverse image of $H$ in $G_\tx{sc}$. Note that $H_\mathrm{sc}$ is in general not itself simply connected, but it does have a simply connected derived subgroup. We let $H_\mathrm{sc,ab}$ be the maximal abelian quotient of $H_\tx{sc}$ and assume given a good element $X \in \mathrm{Lie}^*(H_\mathrm{sc,ab})(F)$ of depth $-r$ for a fixed $r>0$, following the notation of \cite[\S3]{FKS23}. This means that for one, hence any, tame maximal torus $T \subseteq H$ we have $\tx{ord}(\<X,H_\alpha\>)=-r$ for all $\alpha \in R(T,G/H)$.

The $F$-vector space $V=\tx{Lie}(G)/\tx{Lie}(H)$ is endowed with the $H$-invariant symplectic form $\eta = X \circ \tx{pr}_H \circ [-,-]'$, where $X \in \mathrm{Lie}^*(H_\mathrm{sc,ab})(F) \subset \mathrm{Lie}^*(H_\mathrm{sc})(F)$, $\tx{pr}_H : \mathrm{Lie}(G_\mathrm{sc}) \to \mathrm{Lie}(H_\mathrm{sc})$ is the projection onto the $Z(H_\tx{sc})^\circ$-invariant subspace, and $[-,-]'$ is the canonical lift of the Lie bracket specified by the following  well-known lemma, for which we have not been able to locate a convenient reference.
\begin{lem}
	Let $k$ be an arbitrary field, $G_\mathrm{sc}$ a semi-simple and simply connected $k$-group, and $G_\mathrm{ad}$ the adjoint quotient of $G_\mathrm{sc}$. We consider triples $(G,q,p)$, where $G$ is a connected reductive group, $p : G \to G_\mathrm{ad}$ is a homomorphism identifying the adjoint quotient of $G$ with $G_\mathrm{ad}$, $q : G_\mathrm{sc} \to G$ is a homomorphism idenifying the simply connected cover of the derived subgroup of $G$ with $G_\mathrm{sc}$, and $r=p\circ q$ is the canonical central isogeny $G_\mathrm{sc} \to G_\mathrm{ad}$.
	\begin{enumerate}
		\item There exists a unique morphism $c' : G_\mathrm{ad} \times G_\mathrm{ad} \to G_\mathrm{sc}$ of $k$-schemes  such that for all triples $(G,p,q)$ as above
		the commutator map $c : G \times G \to G$, $c(x,y)=xyx^{-1}y^{-1}$ factors as $c(a,b)=q(c'(p(a),p(b)))$.
		\item There exists a unique map $[-,-]' : \mathrm{Lie}(G_\mathrm{ad}) \times \mathrm{Lie}(G_\mathrm{ad}) \to \mathrm{Lie}(G_\mathrm{sc})$ that is alternating $k$-bilinear such that for all triples $(G,p,q)$ as above the Lie bracket $[-,-] : \mathrm{Lie}(G) \times \mathrm{Lie}(G) \to \mathrm{Lie}(G)$ factors as $[X,Y]=dq([dp(X),dp(Y)]')$.
	\end{enumerate}
\end{lem}
\begin{proof}
	(1) The case of a field $k$ of characteristic zero is essentially trivial: letting $\bar k$ be an algebraic closure of $k$ we have $x',y' \in G_\mathrm{sc}(\bar k)$ such $p(x)=r(x')$ and $p(y)=r(y')$. Then $x'y'(x')^{-1}(y')^{-1} \in G_\mathrm{sc}$ depends only on $p(x),p(y)$. This defines $c' : G_\tx{ad}(\bar k) \times G_\tx{ad}(\bar k) \to G_\mathrm{sc}(\bar k)$, which is clearly algebraic any by uniqueness and Galois descent also $k$-rational. 
	
	For general fields $k$ we cannot use Galois (i.e. \'etale) descent, but we can use faithfully flat descent. Indeed, the map $r:=p \circ q : G_\mathrm{sc} \to G_\tx{ad}$ is a quotient map, thus faithfully flat. Let $c_\mathrm{sc} : G_\mathrm{sc} \times G_\mathrm{sc} \to G_\mathrm{sc}$ be the commutator map. In order to show that it descends to $G_\tx{ad} \times G_\tx{ad}$ it is enough to check that the pull-backs of this map under the two projections 
	\[ (G_\mathrm{sc} \times G_\mathrm{sc}) \times_{(G_\tx{ad} \times G_\tx{ad})} (G_\mathrm{sc} \times G_\mathrm{sc}) \rightrightarrows (G_\mathrm{sc} \times G_\mathrm{sc}), \]
	where the fiber product is taken with respect to $r \times r$ on both sides, agree. But this is immediate from the fact that for any $k$-algebra $R$ we have $c_\mathrm{sc}(xz,yw)=c_\mathrm{sc}(x,y)$ for $x,y \in G_\mathrm{sc}(R)$ and $z,w \in Z(G_\mathrm{sc})(R)$. 
	
	(2) This is again trivial if the isogeny $G_\mathrm{sc} \to G_\tx{ad}$ is \'etale, for then also $G_\mathrm{sc} \times Z(G)^\circ \to G$ is \'etale, hence both induce isomorphisms on Lie algebras.
	
	In the general case, if we write $\mathrm{Lie}(G_\mathrm{sc}) = \tx{ker}(G_\mathrm{sc}(k[\epsilon]/(\epsilon^2) \to G_\mathrm{sc}(k)))$ and have the analogous description for $\mathrm{Lie}(G_\tx{ad})$, then we can define $[-,-]' : \mathrm{Lie}(G_\tx{ad}) \times \mathrm{Lie}(G_\tx{ad}) \to \mathrm{Lie}(G_\mathrm{sc})$ via $c'(1+\epsilon_1 X,1+\epsilon_2 Y)=1+\epsilon_1\epsilon_2[X,Y]'$ over $k[\epsilon_1,\epsilon_2]/(\epsilon_1^2,\epsilon_2^2)$, where $c'$ is the map from (1).
	
	Slightly more formally, let $I_\mathrm{sc} \subset k[G_\mathrm{sc}]$ be the kernel of the counit $e^* : k[G_\mathrm{sc}] \to k$ and let $I_\tx{ad} \subset k[G_\tx{ad}]$ be the analogous kernel for $G_\tx{ad}$. Thus $\mathrm{Lie}(G_\mathrm{sc})=(I_\mathrm{sc}/I_\mathrm{sc}^2)^*$ and $\mathrm{Lie}(G_\tx{ad})=(I_\tx{ad}/I_\tx{ad}^2)^*$. From $c'(a,1)=c'(1,a)=1$ for $a \in G_\tx{ad}$ we see $c'^*(I_\mathrm{sc}) \subset I_\tx{ad} \otimes I_\tx{ad}$. Since $c'^*$ is a $k$-algebra homomorphism this implies $c'^*(I_\mathrm{sc}^2) \subset (I_\tx{ad} \otimes I_\tx{ad})^2= I_\tx{ad}^2 \otimes I_\tx{ad}^2$, and thus we obtain from $c'^*$ a $k$-linear map $I_\mathrm{sc}/I_\mathrm{sc}^2 \to (I_\tx{ad}/I_\tx{ad}^2) \otimes (I_\tx{ad}/I_\tx{ad}^2)$ whose dual is the desired pairing $[-,-]'$.

	To verify alternation and uniqueness even when $r$ is inseparable, choose a $k$-torus $S$ and an embedding $Z(G_{\mathrm{sc}})\hookrightarrow S$. Form $G=(G_{\mathrm{sc}}\times S)/Z(G_{\mathrm{sc}})$ using the antidiagonal embedding. Then $q:G_{\mathrm{sc}}\to G$ is a closed immersion and $p:G\to G_{\mathrm{ad}}$ is smooth and surjective, with kernel $S$. Consequently $dq$ is injective and $dp$ is surjective. Differentiating the commutator identity from (1) twice gives
	 \[ dq([dp(X),dp(Y)]')=[X,Y] \]
	for $X,Y\in\operatorname{Lie}(G)$. Since the bracket on the right is alternating, $dq$ injective, and $dp$ surjective, we see that $[-,-]'$ is alternating and unique. Differentiating (1) for any of the given triples proves its required factorization property.
\end{proof}

  \begin{rem}
	Note that if $X \in \mathrm{Lie}^*(H)^H(F)$ is a $(G,H)$-generic element of depth $-r$ in the sense of \cite[Definition~4.1.1]{FS25}, we can use it to endow $V$ with the symplectic form $X \circ [-,-]$. Using the natural map $\mathrm{Lie}^*(H) \to \mathrm{Lie}^*(H_\mathrm{sc})$ we obtain from $X$ an element of $\mathrm{Lie}^*(H_{\mathrm{sc}})^{H_{\mathrm{sc}}}$, the latter being isomorphic to $\mathrm{Lie}^*(H_\mathrm{sc,ab})$ by \cite[Lemma~4.1.5]{FKS23}. The symplectic form defined above then coincides with $X \circ [-,-]$. While in the applications of \S\ref{sec:supercuspidal} we will only use these more special good elements, our set up here allows the more general setting of good elements in $\mathrm{Lie}^*(H_\mathrm{sc,ab})$ in view of potential future applications.
  \end{rem}

	The adjoint action induces a map $\jmath : H \to \Sp(V)$ of algebraic groups. Let $\Lambda=R(Z(H),G)$ be the set of (absolute) weights for the action of the center $Z(H)$ (also sometimes denote by $Z_H$)
of $H$ on $V$. Note that this action factors through the quotient $Z(H)/Z(G)$, which is a torus. The set $\Lambda$ is a polarizable $\Sigma$-set in the sense of Definition \ref{dfn:lambda_basic}, and we have the tautological $\Lambda$-grading on $V$ in the sense of Definition \ref{dfn:lambda_grading}. Via Lemma \ref{lem:Levivsgrading} this determines a polarizable Levi subgroup $H_{\Sp} \subseteq \Sp(V)$ and $\jmath$ factors through a map $H \to H_{\Sp}$, which we shall also denote by $\jmath$. In an analogous manner, we obtain from any maximal torus $T \subseteq H$ a maximal torus $T_{\Sp} \subseteq \Sp(V)$ and $\jmath$ restricts to a map $\jmath : T \to T_{\Sp}$ of algebraic groups.

The earnest $\Sigma$-set $\Lambda \subset X^*(Z(H))$ leads via \cite[\S6]{KalDC} to a double cover $H(F)_\pm$, while we also have the double cover $H_{\Sp}(F)_\pm$ of Definition \ref{dfn:tsl-cov}. The map $\jmath : H(F) \to H_{\Sp}(F)$ lifts tautologically to a genuine homomorphism $\jmath : H(F)_\pm \to H_{\Sp}(F)_\pm$.

Let $\mc{B}(G)$ be the (reduced) Bruhat--Tits building of $G$ and let $\mc{B}(G)_H = \bigcup_T \mc{A}(T)^\Gamma$ as $T$ runs over all tame maximal tori of $H$. We have the natural projection map $\mc{B}(G)_H \to \mc{B}(H)$\index[notation]{BG@$\cB(G)$}\index[notation]{BGH@$\cB(G)_H$} which is an isomorphism when $Z(H)/Z(G)$ is anisotropic. 

We assume given a point $x \in \mc{B}(G)_H$ and write $H_x$\index[notation]{Gx@$G_x$} for the stabilizer of $x$ in $H(F)$ and denote by $G_{x,s}$\index[notation]{Gxs@$G_{x,s}$} the Moy--Prasad filtration subgroup of $G(F)$ at depth $s \in \bR_{\geq 0}$. The quotient 
\[ \oV=G_{x,r/2}/(H_{x,r/2}G_{x,r/2+}) \]
is an $\rf$-vector space and the restriction of $\eta$ induces a $H_x$-invariant symplectic form $\ol{\eta}$ on it, which leads to the map $\bar\jmath : H_x/H_{x,0+} \to \Sp(\oV)$.

Let $H_{x,\pm} \subset H(F)_\pm$ denote the preimage of $H_x$. The subgroup $H_{x,0+}$ splits uniquely into $H_{x,\pm}$. Our goal is to prove the existence and uniqueness of a genuine homomorphism $\tilde\jmath$ lifting $\bar\jmath$ as follows \begin{equation}
	\xymatrix{
	H_{x,\pm}/H_{x,0+}\ar[d]\ar@{..>}[r]^{\tilde\jmath}&\WM_{\bar\psi}(\oV)\ar[d]\\
	H_x/H_{x,0+}\ar[r]^{\bar\jmath}&\Sp(\oV),}
\end{equation}
where $\WM_{\bar\psi}(\oV)$\index[notation]{WMpsibar@$\WM_{\bar\psi}$} is the pushout of the universal Witt-metaplectic group $\WMuniv(\oV)$ of \S\ref{sub:wittmeta} under the \emph{renormalized} Weil index $\tilde\gamma_{\bar\psi} : \QW(\rf) \to \mu_8$. It will be of the form
\[ \tilde\jmath(\dot g)=\tilde\gamma_{\bar\psi}([q_{\dot g}]) \boxtimes \bar\jmath(g) \] 
for some Witt class $[q_{\dot g}] \in \QW(\rf)$ whose Weil index is uniquely determined by $\dot g$.

\subsection{Lattice chains from Moy--Prasad filtrations}	
	
\begin{rem}
	If $H$ is a maximal torus that only has symmetric roots, then the Moy--Prasad filtration $t \mapsto V_{x,t}:=\mathrm{Lie}(G)(F)_{x,t}/\mathrm{Lie}(H)(F)_{x,t}$ is $\Lambda$-atomic in the sense of Definition \ref{dfn:lambda_atomic}. When $H$ has asymmetric roots, or is a general twisted Levi subgroup, this is no longer the case, but the next lemma shows that the lattice $V_{x,r/2}$ generates a good $\Lambda$-atomic filtration.
\end{rem}

On the $E_{\pm\lambda}$-module $V_\lambda^E$ we have the graded lattice chain $t \mapsto (V_\lambda^E)_{x,t}$ given by the Moy--Prasad filtration of the $\pm\lambda$-weight spaces, where we use the normalization given by the valuation of $F^\times$ having image $\Z$ and the valuation on any finite extension of $F$ extending the valuation on $F$. We recall from \eqref{eq:hdual} the notation $(V_\lambda^E)_{x,t}^{h_\lambda^E}$ for the lattice in $V_\lambda^E$ that is dual to $(V_\lambda^E)_{x,t}$ with respect to the skew-Hermitian form $h_\lambda^E$. We write $\cO_\lambda^E$ for the maximal order of $E_{\pm\lambda}$, $\mf p_\lambda^E$ for its Jacobson radical, and $\mathfrak e_{\pm\lambda}=\cO_\lambda^E/\mf p_\lambda^E$ for its residue algebra. 

\begin{lem} \label{lem:bt-align1}
	\leavevmode
	\begin{enumerate}
		\item For each $\lambda \in \Lambda$ and each $t \in \R$ we have $(V_\lambda^E)_{x,t}^{h_\lambda^E}=\mf{p}_\lambda^{-1}(V_\lambda^E)_{x,(r-t)+}$, using the notation of \eqref{eq:hdual}. 
\item For each $\lambda \in \Lambda$ set $a_\lambda=\ell_{p'}(\alpha^\vee)/\<X,H_\alpha\>$, where $\ell_{p'}(\alpha^\vee)$ is the prime-to-$p$ part of the normalized square length of $\alpha^\vee$, for an arbitrary tame maximal torus $T \subseteq H$ and $\alpha \in R(T,G)$ restricting to $\lambda$. Then $a_\lambda$ depends only on $\lambda$ and $X$, not $T$ and $\alpha$, and the symmetric bilinear form $\eta_{\pm\lambda,a}$ on $V_{\pm\lambda}$ is the polar form of the quadratic form $\varphi$ defined in \cite[Definitions 5.4.5, 5.4.8]{FKS23}.
	\end{enumerate}
\end{lem}
\begin{proof}
	(1) Let $\lambda \in \Lambda$. Then
	\begin{enumerate}
		\item $h_\lambda^E((V_\lambda^E)_{x,t},(V_\lambda^E)_{x,u}) \subset \mf{p}_\lambda^E$ if $t+u>r$,
		\item The pairing $(V_\lambda^E)_{x,t}/(V_\lambda^E)_{x,t+} \times (V_\lambda^E)_{x,(r-t)}/(V_\lambda^E)_{x,(r-t)+} \to  \mathfrak e_{\pm\lambda}$ induced by $h_\lambda^E$ is perfect.
	\end{enumerate}
	The first property already implies $(V_\lambda^E)_{x,t}^{h_\lambda^E} \supseteq (\mf{p}_\lambda^E)^{-1}(V_\lambda^E)_{x,(r-t)+}$. To see the opposite inclusion, take $v \in (V_\lambda^E) \sm (V_\lambda^E)_{x,(r-t)+}$ and let $u$ be the unique real number such that $v \in (V_\lambda^E)_{x,u} \sm (V_\lambda^E)_{x,u+}$. By assumption $u \leq (r-t)$ and by the second property there exists $v' \in (V_\lambda^E)_{x,r-u} \sm (V_\lambda^E)_{x,(r-u)+}$ such that $h_\lambda^E(v,v') \in \cO_\lambda^E \sm \mf{p}_\lambda^E$. Hence $v \notin \mf{p}_\lambda^E\cdot (V_\lambda^E)_{x,r-u}^{h_\lambda^E} \supseteq \mf{p}_\lambda^E\cdot (V_\lambda^E)_{x,t}^{h_\lambda^E}$.

	(2) The independence of $a_\lambda$ of $T$ and $\alpha$ is recorded in \cite[Lemma 5.4.9]{FKS23}. We now compute the form $\eta_{\pm\lambda,a}$ on $V_{\pm\lambda}$. In our setting $V=\mathrm{Lie}(G)/\mathrm{Lie}(H)$ with symplectic form $\eta(v,w)=X([v,w]')$, where $X \in \mathrm{Lie}^*(H_\mathrm{sc,ab})(F)$ is good of depth $-r$ with $r>0$. On the other hand, $\varphi$ is defined by
\[ v \mapsto \sum_{\pm\alpha} \ell_{p'}(\alpha^\vee)\frac{[v_\alpha,v_{-\alpha}]'}{H_\alpha}, \]
where we have chosen a maximal torus $T \subset H$ and the sum is running over a set of representatives for each pair $\{\pm\alpha\}$ of absolute roots of $T$ outside $H$ that restrict to $\pm\lambda$ on $Z(H)^\circ$. We show in \cite[Proposition 5.4.6]{FKS23} that $\varphi$ is $H$-invariant, in particular independent of the choce of $T$, as well as Galois-equivariant and non-degenerate.

The polar form of $\varphi$ is computed as 
\begin{eqnarray*}
	\beta_\varphi(v,w)&=&\varphi(v+w)-\varphi(v)-\varphi(w)\\
	&=&\sum_{\pm\alpha}\ell_{p'}(\alpha^\vee)\frac{[v_\alpha,w_{-\alpha}]'+[w_\alpha,v_{-\alpha}]'}{H_\alpha}\\
	&=&\sum_{\pm\alpha}\ell_{p'}(\alpha^\vee)\left(\frac{[v_\alpha,w_{-\alpha}]'}{H_\alpha}+\frac{[w_\alpha,v_{-\alpha}]'}{{H_\alpha}}\right)\\
	&=&\sum_{\pm\alpha}\ell_{p'}(\alpha^\vee)\left(\frac{[v_\alpha,w_{-\alpha}]'}{H_\alpha}+\frac{[v_{-\alpha},w_\alpha]'}{{H_{-\alpha}}}\right)\\
	&=&\sum_{\alpha}\ell_{p'}(\alpha^\vee)\frac{[v_\alpha,w_{-\alpha}]'}{H_\alpha}.
\end{eqnarray*}
where the final sum now runs over all roots of $T$ outside of $H$.

To compare $\eta$ and $\beta_\varphi$ we compute
\begin{eqnarray*}
	\beta_\varphi(v,w)&=&\sum_{\alpha}\ell_{p'}(\alpha^\vee)\frac{[v_\alpha,w_{-\alpha}]'}{H_\alpha}\\
	&=&\sum_{\alpha}\ell_{p'}(\alpha^\vee)\frac{[v_\alpha,w_{-\alpha}]'}{H_\alpha} \cdot \frac{X(H_\alpha)}{X(H_\alpha)}\\
	&=&\sum_{\alpha}\ell_{p'}(\alpha^\vee)\frac{X([v_\alpha,w_{-\alpha}]')}{X(H_\alpha)}.
\end{eqnarray*}
Recall from \cite[Lemma 5.4.9]{FKS23} and the paragraph after \cite[(4.2.5)]{FKS23} that the quantity $a_\lambda:=\ell_{p'}(\alpha^\vee)/X(H_\alpha)$ depends only on $\lambda$ and not on $\alpha$. In particular it lies in $F_\lambda^\times$, and if moreover $\lambda$ is symmetric then $\sigma_\lambda(a_\lambda)=-a_\lambda$. With this, we see that on $V_{\pm\lambda}$ we have
\[ \beta_\varphi(v,w)= X([a_\lambda v_{\lambda},w_{-\lambda}]'+[a_{-\lambda} v_{-\lambda},w_{\lambda}]').\]
If we thus define $T_a$ to be the operator on $V_{\pm\lambda}$ that multiplies $V_\lambda$ by $a_\lambda$ and $V_{-\lambda}$ by $a_{-\lambda}$, then $T_a$ is a semi-simple automorphism of the $F_{\pm\lambda}$-vector space $V_{\pm\lambda}$ and 
\[ \beta_\varphi(v,w)=\eta(T_av,w). \]
Under the isomorphism
\[ V_\lambda^E|_{F_{\pm\lambda}} \to V_{\pm\lambda},\qquad v \mapsto v+\sigma_{\lambda}(v) \]
both $\beta_\varphi$ and $\eta$ pull back to forms that are traces: the pull-back of $\eta$ is the trace of the skew-Hermitian form $h(v,w)=\eta(v,\sigma_{\lambda}(w))$, while $\beta_\varphi$ pulls back to the trace of the Hermitian form $a_\lambda^E \cdot h_\lambda^E$.
\end{proof}

Fix $O \in \Lambda/\Sigma$ and choose $\lambda \in O$. Let $\mc{L}_\lambda^E$ be the subchain of the Moy--Prasad lattice chain in $V_\lambda^E$ that is generated by the lattice $V_{\lambda,x,r/2}^E$ under periodicity by $\mf{p}_\lambda^E$ and duality by $h_\lambda^E$. By Lemma \ref{lem:bt-align1}(1) this is the chain
\[ \mf{p}_\lambda^\Z \cdot V_{\lambda,x,r/2}^E \cup \mf{p}_\lambda^\Z \cdot V_{\lambda,x,r/2+}^E. \]
We keep the Moy--Prasad grading on it, i.e., for $t \in \R$ we let $\mc{L}_{\lambda,t}^E$ be the smallest lattice in the chain $\mc{L}_{\lambda}$ that contains the Moy--Prasad lattice $V_{x,\lambda,t}^E$. The transport $\mc{L}_O$ of this chain under the isomorphism $V_\lambda^E|_F \to V_O$ of Fact \ref{fct:lambda_grading}(2,3) does not depend on the choice of $\lambda \in O$ and is $O$-atomic. Let $\mc{L}$ be the lattice chain in $V$ consisting of the lattices
\[ \mc{L}_t = \{\bigoplus_{O \in \Lambda/\Sigma} \mc{L}_{O,t} \} \]
for all $t \in \R$.
Then $\mc{L}$ is $\Lambda$-atomic. We think of $\mc{L}$ as being the $\Lambda$-atomic lattice chain in $V$ generated by $V_{x,r/2}$.

\begin{lem} \label{lem:bt-align}
		 \leavevmode \begin{enumerate} 
		\item The $\Lambda$-center of $\mc{L}$ is $V_{x,r/2}$ and the $\Lambda$-central quotient is $V_{x,r/2}/V_{x,r/2+}$.
		\item Let $O \in \Lambda/\Sigma$ be ramified symmetric and choose $\lambda \in O$. Let $a=(a_\lambda)_{\lambda \in O}$ be as in Lemma \ref{lem:bt-align1}(2) and let $\beta_O$ be the symmetric bilinear form corresponding to  $\mathrm{tr}_{F_\lambda/F}\circ a_\lambda h_\lambda$ under the isomorphism $V_\lambda|_F \to V_O$ of Fact \ref{fct:lambda_grading}(2,3). The $(O,\beta_O)$-center  of $\mc{L}_O$ is $V_{O,x,(-1/2e_\lambda)+}$, the $(O,\beta_O)$-central quotient is $V_{O,x,(-1/2e_\lambda)+}/V_{O,x,1/2e_\lambda}$, and the symmetric bilinear form on this quotient induced by $\beta_O$ is the polar form $\beta_\varphi$ of the quadratic form $\varphi$ of \cite[Definitions 5.4.5, 5.4.8]{FKS23}. 
\item Let $O$ be not symmetric ramified and choose $\lambda \in \Lambda$. Let $a=(a_\lambda)_{\lambda \in O}$ be $a$-data comprised of $a_\lambda \in \cO_\lambda^\times$ and let $\beta_O$ be the symmetric bilinear form corresponding to  $\mathrm{tr}_{E_{\pm\lambda}/F}\circ a_\lambda^E h_\lambda^E$ under the isomorphism $V_\lambda^E|_F \to V_O$ under the isomorphism of Fact \ref{fct:lambda_grading}(2,3). The $(O,\beta_O)$-center and central quotient of $\mc{L}_O$ are the same as those for $(O,\eta_O)$. The symmetric bilinear form on this quotient equals the polar form of the quadratic form $\varphi_s$ of \cite[Remark 5.5.10]{FKS23}.
	\end{enumerate}
\end{lem}
\begin{proof}
	(1) This follows from Lemma \ref{lem:bt-align1} and Lemma \ref{lem:centequiv}.
	
	(2) Using Lemma \ref{lem:bt-align1} and Lemma \ref{lem:ssprep}(1) we see
	\[ (V_\lambda^E)_{x,r/2}^{a_\lambda^E h_\lambda^E} = a_\lambda^{-1} (V_\lambda^E)_{x,r/2}^{h_\lambda^E} = a_\lambda^{-1} \mf{p}_\lambda^{-1}(V_\lambda^E)_{x,r/2+} = (V_\lambda^E)_{x,(-r/2-e_\lambda^{-1})+}, \]
	which shows
	\begin{eqnarray*}
	((V_\lambda^E)_{x,-1/2e_\lambda})^{a_\lambda^E h_\lambda^E} &=& (\mf{p}_\lambda^{-(1+e_\lambda r)/2}(V_\lambda^E)_{x,r/2})^{a_\lambda^E h_\lambda^E}\\
	& =& \mf{p}_\lambda^{(1+e_\lambda r)/2} (V_\lambda^E)_{x,(-r/2-e_\lambda^{-1})+}	\\
	&=&(V_\lambda^E)_{x,(-1/2e_\lambda)+}.
	\end{eqnarray*}
	From this we infer 
	\[ \begin{array}{ccccc}
	 \mf{p}_\lambda (V_\lambda^E)_{x,-1/2e_\lambda} &\subset& (V_\lambda^E)_{x,(-1/2e_\lambda)+}& \subset& (V_\lambda^E)_{x,-1/2e_\lambda}\\
	 \Vert&&&&\Vert\\
	 \mf{p}_\lambda (V_\lambda^E)_{x,(-1/2e_\lambda)+}^{a_\lambda^E h_\lambda^E} &&&& (V_\lambda^E)_{x,(-1/2e_\lambda)+}^{a_\lambda^E h_\lambda^E},
	\end{array} \]
	and Lemma \ref{lem:centequiv} shows that $(V_\lambda^E)_{x,(-1/2e_\lambda)+}$ is the $(a_\lambda h_\lambda)$-center and, since $(V_\lambda^E)_{x,1/2e_\lambda}=\mf{p}_\lambda (V_\lambda^E)_{x,-1/2e_\lambda}$, that
	\[ (V_\lambda^E)_{x,(-1/(2e_\lambda))+}/(V_\lambda^E)_{x,1/(2e_\lambda)} \]
	is the $(a_\lambda h_\lambda)$-central quotient. The claim about the form follows from Lemma \ref{lem:bt-align1}(2).

	(3) Since $a_\lambda$ is a unit, switching from $h_\lambda$ to $a_\lambda h_\lambda$ does not change the duality operation on $\mc{L}_\lambda$, hence the center and central quotient remain unchanged. When $\lambda$ is symmetric unramified the image of $a_\lambda$ in $\rf_\lambda$ has zero trace in $\rf_{\pm\lambda}$, and when $\lambda$ is asymmetric $a_\lambda \in \rf_\lambda^\times=\rf_{\pm\lambda}^\times$. 
	We have the polarization $V_{\pm\lambda}^E = V_\lambda^E \oplus \sigma_\lambda(V_\lambda^E)$ and the isomorphism $V_\lambda^E|_{F_{\pm\lambda}} \to V_{\pm\lambda}$ given by $x \mapsto x + \sigma_\lambda(x)$, see Definition \ref{dfn:lambda_grading} and Fact \ref{fct:lambda_grading}. Let $L \subset V_\lambda^E$ be the $h_\lambda^E$-center, i.e. such that $\mf{p}_\lambda L^{h_\lambda^E} \subset L \subset L^{h_\lambda^E}$. We have the central quotient $\ol{V}_\lambda^E=L/\mf{p}_\lambda L^{h_\lambda^E}$ with skew-Hermitian form induced by restricting $h_\lambda^E$. Its trace is the sympectic space $\ol{V}_{\pm\lambda} := V_{\pm\lambda,x,r/2}/V_{\pm\lambda,x,(r/2)+}$ with symplectic form $\bar\eta_{\pm\lambda}$ given by restricting the $F_{\pm\lambda}$-linear extension of $\eta$ according to (2).  The form $\varphi_s$ is defined on this quotient as the descent of 
\[ \ol{V}_{\pm\lambda}\otimes_{\rf_{\pm\lambda}} \mathfrak{e}_{\pm\lambda} = \ol{V}_\lambda^E \oplus \sigma_\lambda(\ol{V}_\lambda^E) \to \mathfrak{e}_{\pm\lambda},\qquad (Y^++Y^-) \mapsto a\bar\eta_{\pm\lambda}^E(Y^+,Y^-), \]
where $a \in \mathfrak{e}_{\pm\lambda}$ is an arbitrary non-zero element whose trace in $\mathfrak{f}_{\pm\lambda}$ is zero, as discussed in the erratum to \cite{FKS23}. This element is needed to ensure that $\varphi_s$ is $\sigma_\lambda$-equivariant so descends to a quadratic form $\ol{V}_{\pm\lambda} \to \rf_{\pm\lambda}$.

For the element $Y^++Y^-$ to be $\sigma_\lambda$-stable we must have $Y^-=\sigma_\lambda(Y^+)$. In other words, it is the image of $Y=Y^+ \in \ol{V}_\lambda^E$ under the isomorphism $\ol{V}_\lambda^E|_{\rf_{\pm\lambda}} \to \ol{V}_{\pm\lambda}$ sending $Y$ to $Y+\sigma_\lambda(Y)$. Pulling back the quadratic form $\varphi_s$ under this isomorphism gives the form
\[ Y \mapsto \varphi_s(Y+\sigma_\lambda(Y)) = a\bar\eta_{\pm\lambda}^E(Y,\sigma(Y))=a\bar h_\lambda^E(Y,Y). \]
Thus $\varphi_s$ is the quadratic form associated to the Hermitian form $a\bar h_\lambda^E$. In contrast, the diagonal quadratic form of its trace bilinear form is
\[ Y \mapsto \tx{tr}_{\mf{e}_{\pm\lambda}/\rf_{\pm\lambda}}\big(a\bar h_\lambda^E(Y,Y)\big)=2a\bar h_\lambda^E(Y,Y)=2\varphi_s(Y+\sigma_\lambda(Y)). \]
We take $a$ to be the reduction of $a_\lambda^E$.
\end{proof}

\subsection{The genuine homomorphism replacing the symplectic action}

\begin{thm}
	\label{thm:main1}
	 \leavevmode
	\begin{enumerate}
		\item \label{item:thm:main1} There is a unique genuine homomorphism 
		\[ \tilde \jmath : H_{x,\pm}/H_{x,0+} \to \WM_{\bar\psi}(\oV,\ol{\eta}) \]
		that lifts $\bar\jmath$ and has the following property. Given an element $\dot s \in H(F)_{x,\pm}$ whose image in $H_{\Sp}(F)_\pm$ is a norm and has order prime to $p$, we have $\tilde\jmath(\dot s) = \gamma_\psi(q_{\jmath(\dot s)}) \boxtimes \bar\jmath(s)$. See Remark \ref{rem:jvalues} below for general values of $\tilde\jmath$.
		
		\item \label{item:FWeilvsLWeil} The homomorphism $\tilde \jmath$ satisfies the following identity, where we refer to Appendix \ref{sub:fwh} for the definition of the Fourier--Weil and Weil representation:
\[ \tilde\jmath^*\mathrm{FWeil}_{\bar\psi} = \bar\jmath^*\mathrm{LWeil}_{\bar\psi} \otimes \chi_H^{-1} \otimes \epsilon_x^{G/H},\] 
		where $\chi_H$ is the restriction to $H_{x,\pm}$ of the genuine character of $H(F)_\pm$ produced  as in \cite[\S3.2]{KalDC} from the $\chi$-data for $R(Z_H,G)$ 
defined in \cite[\S4.2]{FKS23} and denoted by $\chi''$ there, and $\epsilon_x^{G/H}$ is the character of \cite[Lemma~4.1.2]{FKS23}. 

		\item \label{item:independenceofX} The homomorphism $\tilde \jmath$ only depends on the coset $X + \mathrm{Lie}^*(H_\mathrm{sc,ab})_{(-r)+}$ and not on the choice of good element $X$ within this coset.
	\end{enumerate}
\end{thm}

We note that the inverse on $\chi_H$ here, which does not appear in the equation of Theorem \ref{thm:symplectic-main}(2), is intentional. We note further that in \cite[\S4.2]{FKS23} it was assumed that $p$ does not divide any bond of the Dynkin diagram, i.e. that in addition to the assumption $p \neq 2$ made in this paper also the assumption $p \neq 3$ holds when the absolute root system of $G$ has a component of type $G_2$. We do not make this assumption here, and interpret the left-hand side of \cite[(4.2.5)]{FKS23} as $\chi''_{\alpha_0}(\ell_{p'}(\alpha^\vee)a_\alpha)$.

\begin{proof}
	\ref{item:thm:main1} \underline{Proof of uniqueness:} The argument is the same as for Theorem \ref{thm:symplectic-main}.

	\underline{Proof of existence:} The map $\jmath : H \to \Sp(V)$ factors through the inclusion $H_{\Sp} \to \Sp(V)$; we shall denote the resulting map $H(F) \to H_{\Sp}(F)$ again by $\jmath$. It lifts canonically to $\jmath : H(F)_\pm \to H_{\Sp}(F)_\pm$ and induces a map $H(F)_{x,\pm} \to H_{\Sp}(F)_{\mc{L},\pm}$,
	where $\mc{L}$ is the lattice chain of Lemma \ref{lem:bt-align}. The composition of this map with the genuine homomorphism produced by Theorem \ref{thm:symplectic-main}(1) provides a genuine homomorphism $\tilde\jmath : H(F)_{x,\pm} \to \WM_{\bar\psi}(\oV)$. The composition of this homomorphism with the unique splitting $H(F)_{x,0+} \to H(F)_{x,\pm}$ projects trivially to $\Sp(\oV)$, hence lands in the central $\mu_8 \subset \WM_{\bar\psi}(\oV)$. But $H(F)_{x,0+}$ is pro-$p$ and $p$ is odd, so we conclude that $\tilde\jmath$ factors through a genuine homomorphism $H(F)_{x,\pm}/ H(F)_{x,0+} \to \WM_{\bar\psi}(\oV)$.

	\ref{item:FWeilvsLWeil} We apply Theorem \ref{thm:symplectic-main}(2), where we take as $a$-data the one specified in Lemma \ref{lem:bt-align}(2) when $\lambda$ is ramified symmetric, and in Lemma \ref{lem:bt-align}(3) otherwise, and obtain the very similar looking identity 
	\[ \tilde\jmath^*\mathrm{FWeil}_{\bar\psi} = \bar\jmath^*\mathrm{LWeil}_{\bar\psi} \otimes \chi_{H_{\Sp}} \otimes \epsilon,\] 
	where now $\bar\jmath$ and $\tilde\jmath$ are with respect to a polarizable twisted Levi subgroup $H_{\Sp} \subset \Sp(V)$, $\chi_{H_{\Sp}} : H_{\Sp}(F)_\pm \to \C^\times$ is the genuine character determined by the $\chi$-data for $\Lambda$ specified in Proposition \ref{pro:orangebox}, and $\epsilon(h)=\tx{sgn}_\rf(\tx{sn}(h|\ol{q_a}))$. It remains to show that the pull-back of $\chi_{H_{\Sp}}$ to $H(F)_\pm$ agrees with $\chi_H^{-1}$ and the pull-back of $\epsilon$ to $H(F)_x$ agrees with $\epsilon_x^{G/H}$.

	We begin with $\chi_H$. The $\chi$-data for $\chi_H$ is specified to be the same as for $\chi_{H_{\Sp}}$ when $\lambda$ is not ramified symmetric, in which case it is of order $1$ or $2$ and hence invariant under inversion. When $\lambda$ is ramified symmetric let us write $\chi_\lambda''$ for the definition of \cite[(4.2.5)]{FKS23} and $\chi_\lambda$ for the one in Proposition \ref{pro:orangebox}. Then we have
	\[ \chi_\lambda''(\ell_{p'}(\alpha^\vee)\<X,H_\alpha\>)=\chi_\lambda(\ell_{p'}(\alpha^\vee)/\<X,H_\alpha\>). \]
	Since $\ell_{p'}(\alpha^\vee) \in \Z$ we have $N_{F_\lambda/F_{\pm\lambda}}(\ell_{p'}(\alpha^\vee))=\ell_{p'}(\alpha^\vee)^2$. Since $\chi$-data are trivial on norms we have $\chi_\lambda''=\chi_\lambda^{-1}$ and conclude $\chi_H=\chi_{H_{\Sp}}^{-1}$, as required.

	We now turn to $\epsilon_x^{G/H}$. Recall that according to \cite[Definition 5.5.8]{FKS23} it is the product of (the pullback of) three separately defined characters $_H\epsilon_\tx{sym,ram}$, $_H\epsilon_{r/2}^\tx{sym,ram}$, and $_H\epsilon_0$. While the definitions of these characters are all of different nature, \cite[Remark 5.5.10]{FKS23} reinterprets all of them as instances of the spinor norm. Since each of the three characters is defined as a product over the indexing set $\Lambda/\Sigma$, we may focus on the contribution of a single $\Sigma$-orbit $O$. 

	If $O$ is not symmetric ramified, the characters $_H\epsilon_\tx{sym,ram}$ and $_H\epsilon_0$ are trivial according to \cite[Definitions 5.5.1, 5.5.7]{FKS23}, while the character $_H\epsilon_{r/2}^\tx{sym,ram}$ is reinterpreted by \cite[Remark 5.5.10]{FKS23} as the spinor norm of the quadratic form $2\varphi_s$ on the $\rf$-vector space $V_{x,O,r/2}/V_{x,O,r/2+}$. According to Lemma \ref{lem:bt-align}(3) this quadratic space is of the $O$-central quotient of the lattice chain $\mc{L}_O$ for the symmetric bilinear form $\eta_a$.

	If $O$ is symmetric ramified, the character $_H\epsilon_{r/2}^\tx{sym,ram}$ is trivial by \cite[Definition 5.5.5]{FKS23}, while the characters $_H\epsilon_\tx{sym,ram}$ and $_H\epsilon_0$ are reinterpreted by \cite[Remark 5.5.10]{FKS23} as the spinor norm of the quadratic form $2\varphi$ on the $\rf$-vector spaces $\bigoplus (V_{x,O,t}/V_{x,O,t+})$ and $V_{x,O,0}/V_{x,O,0+}$ respectively, where the direct sum runs over $t \in (-1/(2e_\lambda),1/(2e_\lambda)) \sm \{0\}$, where $\lambda \in O$ arbitrary. Therefore the product of these two characters is given by the spinor norm of $2\varphi$ on $\bigoplus (V_{x,O,t}/V_{x,O,t+})$, where now the direct sum runs over $t \in (-1/(2e_\lambda),1/(2e_\lambda))$. This is the associated graded of $V_{x,O,(-1/(2e_\lambda))+}/V_{x,O,1/(2e_\lambda)}$ with respect to the Moy--Prasad filtration. According to Lemma \ref{lem:bt-align}(2), the latter is the $\eta_a$-central quotient of $\mc{L}_O$, and the Moy--Prasad filtration is adapted to the polar form of $\varphi$. We conclude from Lemma \ref{lem:graded-spinor} that the spinor norm on this associated graded is equal to the spinor norm of the diagonal quadratic form on the $(\eta_a,O)$-central quotient of the lattice chain $\mc{L}_O$.
	
	We have thus obtained, for all $O$, that the character $\epsilon_x^{G/H}$ of \cite{FKS23} equals the spinor norm, composed with the sign character on $\rf$, of the diagonal quadratic form of the symmetric bilinear $\rf$-space that is the $(\eta_a,O)$-central quotient of the lattice chain $\mc{L}_O$. Taking the direct sum over all $O$, this is precisely the character $\epsilon$ of Theorem \ref{thm:symplectic-main}(2).
	
	\ref{item:independenceofX} By the uniqueness claim in \ref{item:thm:main1} it is enough to show that $\gamma_\psi(q_{\dot s})$ for $\dot s$ as in \ref{item:thm:main1} only depends on this coset. We consider the right-hand side of the formula provided by Proposition \ref{pro:orangebox}, in the set-up of Lemma \ref{lem:bt-align}(2). Adding an element of $\mathrm{Lie}^*(H_\mathrm{sc,ab})_{(-r)+}$ to $X$ does not change the symplectic or quadratic spaces over $\rf$ obtained as central quotients, so we focus on the character $\chi_H$, built from $\chi$-data. Here the $a$-data changes by the addition of an element of valuation $r+$. Thus, the corresponding mod-$a$-data, in the sense of \cite[Definition 4.6.8]{KalRSP} does not change, and the paragraph following that definition shows that the $\chi$-data, being tamely ramified, also does not change.
\end{proof}

\begin{rem} \label{rem:jvalues}
	We can describe the values of $\tilde\jmath$ at an arbitrary element $\dot g \in K_\pm$ as follows. The group $H_{x,\pm}/Z(G)$ is locally-pro-p, hence there is a unique topological Jordan decomposition $\dot g=\dot s \cdot u$ modulo $Z(G)$ according to \cite[Propositions 1.7, 1.8]{Spice08}. If $\epsilon$ is the non-trivial element in the kernel of $H_{x,\pm} \to H_x$, then $\epsilon \dot g = (\epsilon\dot s) \cdot u$ is also a topological Jordan decomposition. There exists a tame maximal torus $T \subset G$ with $x \in \mc{B}(G)_T$ and $s \in T(F)$, and $n \in \Z/2\Z$ such that $\epsilon^n \dot s$ is the norm of an element $\tilde s \in T_{\GL}(F)^\dagger$ as in the proof of Lemma \ref{lem:topssgamma}(5). Then 
	\begin{equation}
		\tilde\jmath(\dot g)=(-1)^n\gamma_\psi(q_{\epsilon^n\dot s}|V)\tilde\gamma_{\bar\psi}(\overline{q_u}|\oV^{s=1}) \boxtimes \bar\jmath(g),
	\end{equation}
	where $g \in H_x$ is the image of $\dot g$, $\gamma_\psi(q_{\epsilon^n\dot s}|V)$ is the Weil index of the quadratic form $q_{\tilde s}$ of Definition \ref{dfn:qtildes}, and $\overline{q_u}$ is the quadratic form on $\oV^{s=1}$ defined by $\ol{q_u}(v)=\ol{\eta}(uv,v)$, where $\oV = V_{x,r/2}/V_{x,(r/2)+}$ and $s \in K$ is the image of $\dot s$. 

	The validity of this formula follows from \eqref{eq:tildejansatz}, where we note that the subspaces $V^{s=1} \subseteq V$, resp., $\oV^{s=1,u=1} \subseteq \oV^{s=1}$ belong to the radicals of the form $q_{\tilde s}$ resp. $\ol{q_u}$, and since we allow Witt classes to be represented by degenerate quadratic forms the result is the same.

\end{rem}

\section{Supercuspidal representations and types} \label{sec:supercuspidal}

In this section we will apply the results of \S\ref{sec:symp} and \ref{sec:symplectictopadic} in order to give an alternative construction of supercuspidal representations and types for $p$-adic groups and demonstrate how it simplifies many applications.

We keep the notation for the local field $F$ the same as for \S\ref{sec:symplectictopadic}. Thus, $F$ is a non-archimedean local field, $\cO$ is the ring of integers, $\mf{p}$ the maximal ideal, and $\rf=\cO/\mf{p}$ the residue field, whose characteristic we assume odd.

Let $G$ be a connected reductive $F$-group that splits over a tamely ramified extension of $F$. For a tame twisted Levi subgroup $H \subseteq G$ we write $R(Z_H,G)$ for the set of weights of the action of the center of $H$ on $\mathrm{Lie}(G)/\mathrm{Lie}(H)$. The topological double cover $H(F)_\pm$ of $H(F)$ relative to the earnest $\Sigma$-set $R(Z_H,G)$ was introduced in \cite[\S6.1]{KalDC}.

\subsection{The input for the construction}\label{sub:input}

Let $H\subseteq G$ be a twisted Levi subgroup and let $x \in \cB(H,F)$. Then write $H_{\pm,x}$, resp., $H_{\pm,x,0}$ for the preimages of $H_x$, resp., $H_{x,0}$ in $H(F)_\pm$. Since $p$ is odd, any pro-$p$ subgroup of $H_{x,0}$ splits uniquely into $H_{\pm,x,0}$ and we use that splitting to identify pro-$p$ subgroups of $H_{x,0}$ with isomorphic ones in $H_{\pm, x,0}$. In particular, we can form the quotient $H_{\pm, x}/H_{x,0+}$, which is a double cover of $H_x/H_{x,0+}$ and contains $H_{\pm,x,0}/H_{x,0+}$. Likewise, the unipotent radical of any $\rf$-parabolic subgroup of $H_{x,0}/H_{x,0+}$ splits uniquely into $H_{\pm, x,0}/H_{x,0+}$ and we identify it with the image of the splitting.

\begin{dfn} \label{def:cuspidal} \ \\[-18pt]
	\begin{enumerate}
		\item We call a representation of a double cover $K_{\pm}$ of a group $K$ \emph{genuine}\index[terminology]{genuine} if the non-trivial element of the kernel of the cover acts as multiplication by $-1$, and we call it \emph{non-genuine}\index[terminology]{non-genuine} if the non-trivial element of the kernel of the cover acts as multiplication by $1$, i.e., if the representation factors through the group $K$, in which case we may view the representation also as a representation of $K$. \item We say that a representation of $H_{\pm,x}$ has \emph{depth zero}, if its kernel contains $H_{x,0+}$, and we call it \emph{essentially depth-zero},\index[terminology]{essentially depth-zero representation} if there exists a character of $H_{\pm,x}$ whose tensor product with the representation has depth zero. 
		\item For a subgroup $H' \subseteq H_{\pm, x}$ that contains $H_{\pm, x,0}$, we call a representation of $H'/H_{x,0+}$ and its inflation to $H'$ \emph{cuspidal} if no non-zero element of the representation space is fixed by the (unique splitting of the) unipotent radical of a proper parabolic subgroup of $H_{x,0}/H_{x,0+}$. 
\end{enumerate}	 
\end{dfn}

The input for the construction of irreducible supercuspidal representations is going to be a genuine irreducible, essentially depth-zero representation attached to a supercuspidal $G$-datum defined as follows.

\begin{dfn}\label{dfn:supercuspidalGdatum} 
	We call a tuple 
\[\Upsilon= ((G_i)_{1\leq i\leq n+1},x,(r_i)_{1\leq i\leq n},\rho^0,(\phi_i)_{1\leq i\leq n}) \]
for some non-negative integer $n$, where  
\begin{enumerate}[(a)]
	\item $G=G_1 \supseteq G_2 \supsetneq G_3 \supsetneq \hdots \supsetneq G_{n+1}=:G^0$ are twisted Levi subgroups of $G$ that split over a tamely ramified field extension of $F$,
	\item  $x \in \cB(G,F)_{G^0} \subseteq \cB(G,F)$, 
	\item $r_1 > r_2 > \hdots > r_n >0$ are real numbers,
	\item \label{item:rho} $\rho^0$ is an irreducible depth-zero representation of $G^0_{\pm, x}$, where $G^0(F)_{\pm}$ is the double cover attached to $R(Z_{G^0}, G)$,
	\item $\phi_i$, for $1 \leq i \leq n$, is a character of $G_{i+1}(F)$ of depth $r_i$,
\end{enumerate}
that satisfies the following conditions 
\begin{enumerate}[(i)]
	\item  the quotient of the centers $Z(G^0)/Z(G)$ is anisotropic, 
	\item  $x$ is a vertex in $\cB(G^0,F)$,
	\item  $\rho^0$ is cuspidal,
	\item $\phi_i$ is $(G_i, G_{i+1})$-generic relative to $x$ of depth $r_i$ for all $1 \leq i \leq n$.
\end{enumerate}
a \emph{supercuspidal $G$-datum}.\index[terminology]{supercuspidal $G$-datum}

We call the supercuspidal $G$-datum \emph{genuine} if $\rho^0$ is genuine, and we call it \emph{non-genuine} if $\rho^0$ is non-genuine.\index[terminology]{supercuspidal $G$-datum!genuine}\index[terminology]{supercuspidal $G$-datum!non-genuine}

We call a tuple \(\Upsilon= ((G_i)_{1\leq i\leq n+1},x,(r_i)_{1\leq i\leq n},\rho^0,(\phi_i)_{1\leq i\leq n}) \)
as above except where we replace \eqref{item:rho} by
\begin{enumerate}[(a')]
	\setcounter{enumi}{3}
	\item $\rho^0$ is a finite dimensional semi-simple genuine or non-genuine
	depth-zero representation of $G^0_{\pm, x}$
\end{enumerate}
a \emph{potentially reducible supercuspidal $G$-datum}.\index[terminology]{supercuspidal $G$-datum!potentially reducible}

Given a potentially reducible supercuspidal $G$-datum \(\Upsilon= ((G_i)_{1\leq i\leq n+1},x,$ $(r_i)_{1\leq i\leq n},\rho^0,(\phi_i)_{1\leq i\leq n}) \), we may view $\phi_i|_{G^0_x}$ also as a representation of $G^0_{\pm, x}$ and we refer to the representation $\rho:=\rho^0 \otimes \otimes_{i=1}^n \phi_i|_{G_x^0}$ of $G^0_{\pm,x}$ as \emph{the essentially depth-zero representation attached to $\Upsilon$}.\index[terminology]{essentially depth-zero representation}
\end{dfn}

\begin{rem}	
	Here we use the notion of generic characters as defined in \cite[Definition~3.8]{Fintzen-IHES}. Alternatively one can work with the weaker notion as defined in \cite[Definition~4.1.1]{FS25}
	and then the below construction works exactly the same way, but yields in general a finite direct sum of supercuspidal representations instead of an irreducible supercuspidal representation, as in \cite{FS25}. Then considering $\rho^0$ in a non-genuine supercuspidal $G$-datum as a representation of $G^0_x$ recovers the notion of a supercuspidal $G$-datum in \cite[\S4.1]{FS25}.
\end{rem}

\begin{dfn} Let $G^0 \subseteq G$ be a tamely ramified twisted Levi subgroup of $G$ and let $x \in \cB(G,F)_{G^0}$. We say that an essentially depth-zero representation $\rho$ of $G^0_{\pm, x}$ \emph{admits a supercupsidal $G$-datum} if there exists a potentially reducible supercuspidal $G$-datum  \(\Upsilon= ((G_i)_{1\leq i\leq n+1},x,(r_i)_{1\leq i\leq n},\rho^0,(\phi_i)_{1\leq i\leq n}) \) with $G_{n+1}=G^0$ such that $\rho$ is the essentially depth-zero representation attached to $\Upsilon$.	
\end{dfn}

 \begin{rem}
 	The input for our construction of (irreducible) supercuspidal representations will be an (irreducible) finite-dimensional essentially depth-zero representation of $G^0_{\pm, x}$ that admits a supercuspidal $G$-datum. Our definition of a genuine supercuspidal $G$-datum is analogous to the input for Yu's construction of supercuspidal representations \cite[\S3]{Yu01} except that we work with the double cover $G^0_{\pm, x}$ of $G^0_{x}$ instead of with $G^0_{x}$. Alternatively one may also define a genuine supercuspidal $G$-datum by replacing the characters $\phi_i$ of $G_{i+1}(F)$ by characters of appropriate double covers of $G_{i+1}(F)$. We discuss this in more detail in \S\ref{sub:alternativeinput}.
 \end{rem}

 \subsection{Construction of supercuspidal representations} \label{sub:construction}
In order to construct supercuspidal representations, let $G^0 \subseteq G$ be a tamely ramified twisted Levi subgroup of $G$ and let $x \in \cB(G,F)_{G^0}$. Let $\rho$ be a genuine essentially depth-zero representation of $G^0_{\pm, x}$ that admits a supercupsidal $G$-datum.

To explain the construction of a supercuspidal representation from this input, we let \(\Upsilon= ((G_i)_{1\leq i\leq n+1},x,(r_i)_{1\leq i\leq n},\rho^0,(\phi_i)_{1\leq i\leq n}) \) be a potentially reducible supercuspidal $G$-datum whose attached essentially depth-zero representation is $\rho$. The resulting supercuspidal representation will not depend on the choice of $\Upsilon$. We define the following compact-mod-center, open subgroups of $G(F)$: \begin{eqnarray*}
	K & =& G_1(F)_{x,\frac{r_1}{2}}G_2(F)_{x,\frac{r_2}{2}}\hdots G_n(F)_{x,\frac{r_n}{2}}G_{n+1}(F)_{x} \\
	K^\dagger & =& G_1(F)_{x,\frac{r_1}{2}}G_2(F)_{x,\frac{r_2}{2}}\hdots G_n(F)_{x,\frac{r_n}{2}}G_{n+1}(F)_{x,0+}\\
	K^+ & =& G_1(F)_{x,\frac{r_1}{2}+}G_2(F)_{x,\frac{r_2}{2}+}\hdots G_n(F)_{x,\frac{r_n}{2}+}G_{n+1}(F)_{x,0+}
\end{eqnarray*}

We set 
\begin{eqnarray*}
	K_{\pm} & =& G_1(F)_{x,\frac{r_1}{2}}G_2(F)_{x,\frac{r_2}{2}}\hdots G_n(F)_{x,\frac{r_n}{2}}(G_{n+1})_{\pm,x} 
\end{eqnarray*}
and also view the pro-$p$ groups $K^\dagger$ and $K^+$ as subgroups of $K_{\pm}$.

We will construct a genuine representation $\kappa_{\pm}$ of $K_{\pm}$ by the analogue of Yu's construction replacing the use of the Heisenberg--Weil representation by the Fourier--Heisenberg--Weil representation, which is recalled in Appendix \ref{sub:fwh}. Then we will show that $\rho^0 \otimes \kappa_{\pm}$ (where $\rho^0$ also denotes the inflation to $K_\pm$ that is trivial on $K^\dagger$) descends to a representation of $K$ whose induction to $G(F)$ is a supercuspidal representation, which is irreducible if $\rho^0$ is irreducible, see Theorem \ref{thm-screps}.

Following \cite[\S4]{Yu01} and \cite[Lemma~4.2.1]{FS25} we denote by $\hat \phi_i$ the character of $(G_{i+1})_{x}G_{x,r_i/2+}$ that agrees with $\phi_i$ on $(G_{i+1})_{x}$ and is trivial on $(G_{i+1},G)_{x,r_i+,r_i/2+}$, and we set 
\[\hat \phi =\prod_{i=1}^n \hat \phi_i|_{G^0(F)_xK^+} .\] We also view $\hat \phi$ as a character of $G^0_{\pm, x}$ by inflation via $G^0_{\pm, x} \twoheadrightarrow G^0_{x}$ (and restriction).
We write $\ol{V}=K^\dagger/K^+$. Then the map $(\bar k_1, \bar k_2) \mapsto \hat \phi(k_1k_2k_1^{-1}k_2^{-1})$ (where $k_i$ is a lift of $\bar k_i$ to $K^\dagger$) yields a well-defined non-degenerate symplectic form valued in $C:=\mu_p \cdot \hat \phi (K^+) \subset \bC^\times$ (cf.\ \cite[2.5.13]{DF}), which we denote by $\zeta$.\index[notation]{zeta@$\zeta$} In line with \cite[3.6.1]{DF} we denote by $C \widetilde{\boxtimes} \ol{V}=C \boxtimes_{\sqrt{\zeta}} \ol{V}$
 the group with underlying set $C \times \ol{V}$ and group law given by $(c_1,v_1)\cdot (c_2,v_2)=(c_1c_2\sqrt{\zeta(v_1,v_2)}, v_1+v_2)$. Since $\zeta$ is valued in $\mu_p$, the group $C \widetilde{\boxtimes} \ol{V}$ contains as a subgroup the Heisenberg group $\mu_p \boxtimes_{\sqrt{\zeta}} \ol{V}$, and $C \widetilde{\boxtimes} \{0\}$ is the center of  $C \widetilde{\boxtimes} \ol{V}$. Note that the subgroup $\Sp(\ol{V},\zeta)$ of linear automorphisms of $\ol{V}$ that preserve $\zeta$ acts on  $C \widetilde{\boxtimes} \ol{V}$ by acting on $\ol{V}$ (and trivially on $C$).
 We denote by $\WM_\gamma(\ol{V}, \zeta)$ the covering group of $\Sp(\ol{V}, \zeta)$ that fits into the short exact sequence
 \[ 1 \ra \mu_8 \ra \WM_\gamma(\ol{V}, \zeta) \ra \Sp(\ol{V}, \zeta) \ra 1 \]
 and is defined in Appendix \ref{sub:WittgroupC}. It is the pushout of the universal Witt-metaplectic group $\WMuniv(\ol{V}, \zeta)$ under the renormalized Weil index $\tilde\gamma:\QW_{2'} \ra \mu_8$\index[notation]{gamma@$\gamma$} on the prime-to-2 part of Weissauers's Witt group of characters of second degree $\QW$, see  Appendix \ref{sub:WittgroupC} for details.
Via the projection $\WM_\gamma(\ol{V}, \zeta) \twoheadrightarrow \Sp(\ol{V},\zeta)$ we then obtain an action of $\WM_\gamma(\ol{V}, \zeta)$ on $C \widetilde{\boxtimes} \ol{V}$ extending its action on $\mu_p \boxtimes_{\sqrt{\zeta}} \ol{V}$.
 
\begin{dfn}
	We denote by $\mathrm{FWeilHeis}$\index[notation]{FWeilHeis} the (up to isomorphism) 
unique representation of $\WM_\gamma(\ol{V}, \zeta) \ltimes (C \widetilde{\boxtimes} \ol{V})$ whose restriction to $\WM_\gamma(\ol{V}, \zeta) \ltimes (\mu_p \boxtimes_{\sqrt{\zeta}} \ol{V})$ is the Fourier--Heisenberg--Weil representation as defined in Appendix \ref{sub:WittgroupC}  and whose restriction to $C \widetilde{\boxtimes} \{0\}$ is multiplication by $C$. We denote its restriction to $\WM_\gamma(\ol{V}, \zeta)$ by $\mathrm{FWeil}$\index[notation]{FWeil} and its restriction to $C \widetilde{\boxtimes} \ol{V}$ by $\tx{Heis}$. \end{dfn}

In order to compare our construction with the traditional construction of Yu, we also introduce notation for the traditional Heisenberg--Weil representation.

\begin{dfn}
	We denote by $\mathrm{LWeilHeis}$\index[notation]{LWeilHeis} the (up to isomorphism) unique representation of $\Sp(\ol{V}, \zeta) \ltimes (C \widetilde{\boxtimes} \ol{V})$ whose restriction to $\Sp(\ol{V}, \zeta) \ltimes (\mu_p \boxtimes_{\sqrt{\zeta}} \ol{V})$ is the Heisenberg--Weil representation and whose restriction to $C \widetilde{\boxtimes} \{0\}$ is multiplication by $C$, and we denote its restriction to $\Sp(\ol{V}, \zeta)$ by $\mathrm{LWeil}$\index[notation]{LWeil}.
\end{dfn}

By \cite[Lemma~3.6.5]{DF}\footnote{Strictly speaking the reference treats the case $(G^0)^1_x$, but everything works equally well for $G^0_x$ in place of $(G^0)^1_x$, and the existence of such a homomorphism is essentially due to Yu (\cite{Yu01}) who uses the notion of ``special isomorphism'' and ``symplectic action'' and carefully chooses one that he prefers to work with.} there exists a homomorphism 
\[\bar \jmath \ltimes f: G^0(F)_x \ltimes K^\dagger \ra \Sp(\ol{V}, \zeta) \ltimes(C \widetilde{\boxtimes} \ol{V})\]
whose restriction to $G^0(F)_x$ is $\bar \jmath$ and whose restriction $f$ to $K^\dagger$ fits in the following commutative diagram
\[\xymatrix{
	K^+ \ar[d]_{\hat \phi} \ar@{^{(}->}[r] & K^\dagger \ar@{->>}[r] \ar[d]^{f} & \ol{V} \ar@{=}[d]\\
	C \ar@{^{(}->}[r] & C \widetilde{\boxtimes} \ol{V} \ar@{->>}[r] & \ol{V} .}\]
Moreover, by \cite[Lemma~3.6.5]{DF} a homomorphism $f: K^\dagger \ra C \widetilde{\boxtimes} \ol{V}$ with the above properties is unique up to conjugation by an element in $\{1 \} \widetilde{\boxtimes} \ol{V}^{\bar\jmath (G^0_x)}$.

In order to lift the homomorphism $\bar \jmath \ltimes f$ to a group homomorphism of covering groups
\[\tilde\jmath \ltimes f: G^0_{\pm, x} \ltimes K^\dagger \ra \WM_\gamma(\ol{V}, \zeta) \ltimes(C \widetilde{\boxtimes} \ol{V}),\]
we adjust Theorem \ref{thm:main1} from the setting of one twisted Levi subgroup to our setting with potentially more than one twisted Levi subgroup as follows. To do so, let $X_i \in \mathrm{Lie}^*(G_{i+1})^{G_{i+1}}(F)\subset \mathrm{Lie}^*(G)$ be $(G_i, G_{i+1})$-generic of depth $-r_i$ such that $\phi_i|_{(G_{i+1})_{x,r_i}}$ is given on $(G_{i+1})_{x,r_i}/(G_{i+1})_{x,r_i+} \simeq (\mathrm{Lie}(G_{i+1}))_{x,r_i}/(\mathrm{Lie}(G_{i+1}))_{x,r_i+}$ by $\psi \circ X_i$. 
Let $X=\sum_{i=1}^n X_i \in \mathrm{Lie}^*(G)$. Then $\eta:=X \circ [\phantom{X},\phantom{X}]$ defines a non-degenerate symplectic form on $V:=\mathrm{Lie}(G)/\mathrm{Lie}(G^0)$ and we denote by $\Sp(V)$ the resulting symplectic group over $F$.
\begin{thm}\label{thm:main1b}
	 \leavevmode \begin{enumerate}
		\item \label{item:existenceoftildej} There is a unique genuine homomorphism
		\[ \tilde \jmath : G^0_{\pm, x}/G^0_{x,0+} \to \WM_\gamma(\ol{V}, \zeta) \]
		that lifts $\bar\jmath$ and has the following property. Given an element $\dot s \in G^0_{\pm,x}$ whose image $\jmath(\dot s)$ has finite prime-to-$p$ order and is a norm in the sense of Definition \ref{dfn:norm}, we have \[ \tilde\jmath(\dot s) = \gamma_\psi(q_{\jmath(\dot s)}) \boxtimes \bar\jmath(s), \]
		 where $q_{\jmath(\dot s)}$ is the quadratic form on $V$ given by $q_{\jmath(\dot s)}(v)=\eta(\tilde s v, v)$ for an appropriate choice of $\tilde s$ whose norm is $\jmath(\dot s)$, see Lemma \ref{lem:topssgamma}(5) for details, and $s$ is the image of $\dot s$ in $G_{x}^0$.
		
		\item\label{item:independenceofXi} The homomorphism $\tilde \jmath$ does not depend on the choice of $\psi$ or the choices $X_i$ of generic elements representing $\phi_i|_{(G_{i+1})_{x,r_i}}$.
		
		\item \label{item:description} Let $\dot g \in G^0_{\pm, x}$. Then 
			\begin{equation}
			\tilde\jmath(\dot g)=(-1)^n\gamma_\psi(q_{\epsilon^n\dot s}|V)\tilde\gamma(\overline{q_u}|\ol{V}^{s=1}) \boxtimes \bar\jmath(g),
		\end{equation}
		where the objects are defined as follows, see also Remark \ref{rem:jvalues}. $\epsilon$ is the non-trivial element in the kernel of $G^0(F)_{\pm} \twoheadrightarrow G^0(F)$. $\dot g = \dot s \cdot u$ with $\dot s$ topologically semisimple modulo $Z(G)(F)$ and $u$ topologically unipotent is the topological Jordan decomposition of $\dot g$ modulo $Z(G)(F)$. $n \in \bZ/2\bZ$ is such that $\jmath(\eps^n\dot s)$ is the norm of an appropriate element $\tilde s$ as in Lemma \ref{lem:topssgamma}(5). $g$ and $s$ are the images of $\dot g$ and $\dot s$ in $G^0_{x}/G^0_{x,0+}$, respectively. $\ol{V}^{s=1}$ is the subspace of $\ol{V}$ on which $s$ acts trivially, $\gamma_\psi(q_{\epsilon^n\dot s}|V)$ is the Weil index of the quadratic form $q_{\tilde s}(v):=\eta(\tilde s v, v)$, and $\tilde\gamma(\overline{q_u}|\ol{V}^{s=1})$ is the renormalized Weil index of the  character of second degree on the space $\ol{V}^{s=1}$ defined by $\ol{q_u}(v)=\zeta(uv,v)$.
		
		\item\label{item:organbox} The homomorphism $\tilde \jmath$ satisfies the following identity:
		\[ \tilde\jmath^*\mathrm{FWeil} = \bar\jmath^*\mathrm{LWeil} \otimes \chi_{G^0}^{-1} \otimes \epsilon,\] 
		where $\chi_{G^0}$ is the restriction to $G^0_{\pm, x}$ of the genuine character of $G^0(F)_\pm$ produced from the $\chi$-data for $R(Z_{G^0},G)$ that is obtained by composing the $\chi$-data for $R(Z_{G_{i+1}}, G_i)$ defined in \cite[\S4.2]{FKS23} and denoted by $\chi''$ there with the appropriate norm maps, and $\epsilon:=\prod_{i=1}^n\epsilon_x^{G_i/G_{i+1}}$, i.e., the quadratic character as defined on page 2259 of \cite{FKS23} (just before \cite[Theorem~4.1.13]{FKS23}).
	\end{enumerate}
\end{thm}
\begin{proof}
	\ref{item:existenceoftildej} 
	\underline{Proof of uniqueness:} The argument is the same as for Theorem \ref{thm:main1}, i.e., as for Theorem \ref{thm:symplectic-main}. 
	
	\underline{Proof of existence:} 	
	Let $G^0_{\pm^i}$ be the pull back of the double cover $G_{i+1}(F)_{\pm^i}$ of $G_{i+1}$ corresponding to the  $\Sigma$-set $R(Z_{G_{i+1}}, G_i)$ for $1 \leq i \leq n$. Then $G^0_{\pm}$ is the Baer sum of $G^0_{\pm^i}$ over $1 \leq i \leq n$. Thus in order to define a genuine homomorphism $\tilde \jmath : G^0_{\pm, x}/G^0_{x,0+} \to \WM_\gamma(\ol{V}, \zeta)$ it suffices to provide genuine homomorphisms $\tilde \jmath_i: G^0_{\pm^i, x}/G^0_{x,0+} \to \WM_\gamma(\ol{V}, \zeta)$ for $1 \leq i \leq n$. Let 
	\begin{eqnarray*}
		\ol{V}_i &=& (G_i)_{x,r_i/2}/((G_{i+1})_{x,r_i/2}(G_i)_{x,r_i/2+})\\
		&=&((K^\dagger\cap G_{i+1})(K^\dagger\cap G_{i}))/((K^\dagger\cap G_{i+1})(K^+\cap G_{i}))
	\end{eqnarray*}
	and denote by $\zeta_i$ the restriction of $\zeta$ to $\ol{V}_i$. Then $\zeta_i = \bar\psi \circ \bar \eta_i$, where $\bar \eta_i :=X_i \circ [\phantom{X},\phantom{X}] \mod \fp$ when identifying $\ol{V}_i$ with  $\mathrm{Lie}(G_i)_{x,r_i/2}/(\mathrm{Lie}(G_{i+1})_{x,r_i/2}+\mathrm{Lie}(G_i)_{x,r_i/2+})$. We write $\Sp(\ol{V}_i, \bar \eta_i)$ for the corresponding symplectic group.
	Then $\ol{V} = K^\dagger/K^+=\bigoplus_{i=1}^n \ol{V}_i$ and $\zeta = \prod_i \zeta_i$ and hence we have embeddings
	\[\Sp(\ol{V}_j, \bar \eta_j) \hookrightarrow \textstyle\prod_{i=1}^n \Sp(\ol{V}_i,\zeta_i) \hookrightarrow \Sp(\ol{V}, \zeta)\]
	 for $1 \leq j \leq n$, which lift to embeddings $\iota_j: \WM_{\bar\psi}(\ol{V}_j) \hookrightarrow \WM_{\gamma}(\ol{V}, \zeta)$. For $1 \leq i \leq n$, we set $\tilde \jmath_i$ to be the pull back to $G^0_{\pm^i}$ of the genuine homomorphism $G_{i+1}(F)_{\pm^i,x}/G_{i+1}(F)_{x,0+} \ra \WM_{\bar\psi}(\ol{V}_i)$ from Theorem \ref{thm:main1} composed with $\iota_i$. Then $\tilde \jmath_i: G^0_{\pm^i, x}/G^0_{x,0+} \to \WM_\gamma(\ol{V}, \zeta)$ is genuine and hence the collection of $\tilde \jmath_i$ leads to a genuine homomorphism $\tilde \jmath : G^0_{\pm, x}/G^0_{x,0+} \to \WM_\gamma(\ol{V}, \zeta)$.
	
	Note that if the image of $\dot s \in G^0_{\pm, x}$ in the double cover $G^0_{\Sp}(F)_\pm$ of $G^0_{\Sp}(F) \subseteq \Sp(V, \eta)$ is of finite prime-to-$p$ order and is a norm, then it can be represented by elements $(\dot s_i \in G^0_{\pm^i,x})_{1 \leq i \leq n}$ with the property that the image of $\dot s_i$ in the double cover $(G_{i+1})_{\Sp}(F)_{\pm^i}$ of $(G_{i+1})_{\Sp}(F)\subseteq \Sp(V_i, \eta_i)$ is of finite prime-to-$p$ order and is a norm. Here we write $V_i:=\mathrm{Lie}(G_i)/\mathrm{Lie}(G_{i+1})$, which is equipped with the symplectic form $\eta_i$ such that  $V=\oplus_{i=1}^n V_i$ and $\eta = \sum_{i=1}^n \eta'_i$. Note that this means that  $\eta'_i:= X_i' \circ [\phantom{X},\phantom{X}]$ for some element $X_i' \in X_i + \mathrm{Lie}^*(G_{i+1})^{G^0}_{x,(-r_i)+}$. Let $\eta_i:= X_i \circ [\phantom{X},\phantom{X}]$, and let ${q^i}'_{j(\dot s_i)}$ and $q^{i}_{j(\dot s_i)}$ be the quadratic forms on $V_i$ analogues to $q_{j(\dot s)}$ for $\eta_i'$ and $\eta_i$, respectively.  We claim that  $\gamma_\psi({q^i}'_{j(\dot s_i)})=\gamma_\psi(q^{i}_{j(\dot s_i)})$. Granting this claim, using that the Weil index of the direct sum of two quadratic spaces is the product of their Weil indices, we see that the desired property $\tilde\jmath(\dot s) = \gamma_\psi(q_{\jmath(\dot s)}) \boxtimes \bar\jmath(s)$ follows from the analogous properties for $\dot s_i$ in Theorem \ref{thm:main1}\ref{item:thm:main1}.
	
	It remains to prove that $\gamma_\psi({q^i}'_{j(\dot s_i)})=\gamma_\psi(q^{i}_{j(\dot s_i)})$. Note that these Weil indices can be computed using \eqref{eqn:orangebox} of Proposition \ref{pro:orangebox} for $\cL$ the $\Lambda$-atomic lattice chain in $V_i$ generated by $(V_i)_{x,r_i/2}$ as defined just before Lemma \ref{lem:bt-align}. We can use the same $a$-data and $\chi$-data in both cases, which means that $\chi_H(\dot s_i)$ agrees in both cases. Moreover, the first two factors of the right hand side of \eqref{eqn:orangebox} only depend on the coset $X_i+ \mathrm{Lie}^*(G_{i+1})^{G^0}_{x,(-r_i)+}$. Hence $\gamma_\psi({q^i}'_{j(\dot s_i)})=\gamma_\psi(q^{i}_{j(\dot s_i)})$.

	\ref{item:independenceofXi} The independence on the choices of $X_i$ follows from Theorem \ref{thm:main1}\ref {item:independenceofX}. Let $\psi': F \ra \bC^\times$ be another non-trivial character, then there is $a \in F^\times$ such that $\psi'(x)=\psi(ax)$. We need to show that $\gamma_{\psi'}(q'_{j(\dot s)})=\gamma_{\psi}(q_{j(\dot s)})$, where $q'$ corresponds to the quadratic form arising from $\eta' = X' \circ [\phantom{X}, \phantom{X}]$ with $\psi' \circ X'= \psi \circ X$. We may choose $aX'=X$. Then $\gamma_{\psi'}(q'_{j(\dot s)})=\gamma_{\psi}(aq'_{j(\dot s)})=\gamma_{\psi}(q_{j(\dot s)})$.

	\ref{item:description} This follows from Remark \ref{rem:jvalues} and the construction of $\tilde \jmath$.
	
	\ref{item:organbox} For $1 \leq i \leq n$, we write $\mathrm{LWeil_i}$ and $\mathrm{FWeil_i}$ for the Weil and Fourier--Weil representations of $\Sp(\ol{V}_i, \bar\eta_i)$ and $\WM_{\bar \psi}(\ol{V}_i, \bar \eta_i)$, respectively.
	Then the restriction of $\mathrm{LWeil}$ from $\Sp(\ol{V}, \zeta)$ to (the image of) $\prod_{i=1}^n\Sp(\ol{V}_i, \bar \eta_i)$ is $\otimes_{i=1}^n\mathrm{LWeil_i}$ by \cite[Corollary~2.5~and~Corollary~2.6]{Gerardin}.  The restriction of $\mathrm{FWeil}$ from $\WM_\gamma(\ol{V}, \zeta)$ to (the image of) $\prod_{i=1}^n\WM_{\bar\psi}(\ol{V}_i, \bar \eta_i)$ is $\otimes_{i=1}^n\mathrm{FWeil_i}$ by Lemmas \ref{lem:fweilrest} and \ref{lem:fweilproduct} and Definition \ref{dfn:FWeilHeis}.
	Note that $\eps=\prod_{i=1}^n \eps_x^{G_i/G_{i+1}}$ and $\chi_{G^0}=\prod_{i=1}^n \chi_{G_{i+1},G_{i}}$, where $\chi_{G_{i+1},G_{i}}$ denotes the pull back to $G^0_{\pm^i}$ of the genuine character of  $G_{i+1}(F)_{\pm^i}$ produced from the $\chi$-data for $R(Z_{G_{i+1}}, G_i)$ defined in \cite[\S4.2]{FKS23}. Hence the identity
	\( \tilde\jmath^*\mathrm{FWeil} = \bar\jmath^*\mathrm{LWeil} \otimes \chi_{G^0}^{-1} \otimes \epsilon\)
	follows from the construction of $\tilde \jmath$ and Theorem \ref{thm:main1}\ref{item:FWeilvsLWeil}.

\end{proof}
\begin{rem}
	Note that only Part \ref{item:existenceoftildej} of Theorem \ref{thm:main1b} is needed for the construction of supercuspidal representations (except for simplicity we use Part \ref{item:organbox} in one of the proofs below in the case where the residue field has only three elements). Otherwise Part\ref{item:organbox} is only used to compare our construction with the construction of Yu \cite{Yu01} twisted by the quadratic character of \cite{FKS23}. The detailed description in Part \ref{item:description} is only for the convenience of the reader.
\end{rem}

\begin{dfn}
	We define the representation $\kappa_{\pm}$ of $K_{\pm}$ to be the unique representation of $K_{\pm}$ such that the composition of $G^0_{\pm, x} \ltimes K^\dagger \twoheadrightarrow K_{\pm}$ with $\kappa_{\pm}$ yields $(\hat \phi \ltimes 1)\otimes (\mathrm{FWeilHeis} \circ (\tilde\jmath \ltimes f))$.
\end{dfn}

\begin{lem}
	The (isomorphism class of the) representation $\kappa_{\pm}$ does not depend on the choice of $f$.  Moreover, $\kappa_{\pm}$ is a genuine representation of $K_{\pm}$, i.e.,  the central subgroup $\{ \pm 1\} \subset K_{\pm}$ acts by multiplication by the sign character.
\end{lem}
\begin{proof}
	The proof of the first part is a variant of \cite[3.6.6]{DF} adjusted to our setting. Let $f'$ be another choice of homomorphism $K^\dagger \ra (C \widetilde{\boxtimes}\ol{V})$ with the above properties and let $\kappa_{\pm}'$ denote the attached representation of $K_{\pm}$. Then by \cite[Lemma~3.6.5]{DF} there exists $(1, \ol{V}) \in C \widetilde{\boxtimes}\ol{V}$ such that $f'=\Inn((1,\ol{V})) \circ f$, where $\Inn$ denotes the conjugation by $(1 , \ol{V})$. Hence the action of $\mathrm{FWeilHeis}(1 \ltimes (1,\ol{V} ))$ on the space of $\tx{Heis}$ provides an isomorphism between $\tx{Heis}\circ f$ and $\tx{Heis}\circ f'$. Since $\ol{V} \in \ol{V}^{\bar\jmath(G_x^0)}$, the action of $\mathrm{FWeilHeis}(1 \ltimes (1,\ol{V} ))$ commutes with the action of $\mathrm{FWeil}\circ \tilde \jmath (G^0_{\pm,x})$, and hence provides an isomorphism between $\mathrm{FWeilHeis} \circ (\tilde\jmath \ltimes f)$ and $\mathrm{FWeilHeis} \circ (\tilde\jmath \ltimes f')$ and thus also between $\kappa_{\pm}$ and $\kappa_{\pm}'$.
	
	That $\kappa_{\pm}$ is genuine follows from $\tilde \jmath$ and $\mathrm{FWeil}$ being genuine, see Definitions \ref{dfn:weiluniv} and \ref{dfn:FWeilHeis}.
\end{proof}

\begin{thm} \label{thm-screps}
	The representation $\rho^0 \otimes \kappa_{\pm}$ of $K_{\pm}$ factors through a representation of $K$, which we also denote by $\rho^0 \otimes \kappa_{\pm}$. The compactly induced smooth representation $\pi_{\Upsilon}:=\cind_K^{G(F)}(\rho^0 \otimes \kappa_{\pm})$ is a finite direct sum of irreducible supercuspidal representations. If $\rho$, equivalently, $\rho^0$, is irreducible, then $\pi_\Upsilon$ is irreducible.
	
	The representation $\pi_{\Upsilon}$ depends only on the genuine essentially depth-zero representation $\rho=\rho^0 \otimes_{i=1}^n \phi_i|_{G^0_x}$ of $G^0_{\pm, x}$ and not choice of $\Upsilon$. Hence we may also write $\pi_\rho:=\pi_{\Upsilon}$.
\end{thm}
\begin{proof}
	Since $\rho^0$ and $\kappa_{\pm}$ are both genuine representations, the central $\{ \pm 1\} \subset K_{\pm}$ acts trivially on $\rho^0 \otimes \kappa_{\pm}$, so the representation factors through $K_{\pm} \twoheadrightarrow K$.
	
	Since $\rho^0$ is a finite direct sum of irreducible representations and compact induction preserves finite direct sums, it suffices to consider the case that $\rho^0$ is irreducible, which we assume from now on.
	
	To show that $\cind_K^{G(F)}(\rho^0 \otimes \kappa_{\pm})$ is irreducible and supercuspidal it suffices to show that $(K, \rho^0 \otimes \kappa_{\pm})$ also appears as inducing data in Yu's construction and then cite \cite{Fin21Comp, Yu01} to conclude that the resulting induced representation is irreducible and supercuspidal.
	Here we will only sketch the proof under the assumption that the residue field $\rf$ of $F$ has more than three elements. For a precise proof in full generality we refer to Proposition \ref{pro:agreeswithYu}.
	 We write $\kappa^{\mathrm{nt}}$ for the representation of $K$ constructed via the theory of Heisenberg--Weil representation in Yu's construction, see \S\ref{sub:Yucomparison} for details. Then the inducing data for Yu's construction is of the form $\rho^0_{Yu} \otimes \kappa^{\mathrm{nt}}$ for an irreducible representation $\rho^0_{Yu}$ of $K/K^\dagger$ whose restriction to $G^0_{x,0}/G^0_{x,0+}$ is cuspidal.
	Let $\chi$ be a genuine character of $K_{\pm}/K^\dagger$. Then $\chi \otimes \kappa_{\pm}$ extends the irreducible representation $\kappa_{\pm}|_{K^\dagger}=\kappa^{\mathrm{nt}}|_{K^\dagger}$ and is a linearization of the projective representation of $K$ acting on $\kappa^{\mathrm{nt}}|_{K^\dagger}$ arising from the theory of Heisenberg--Weil representations. Hence  $\chi \otimes \kappa_{\pm}$ and $\kappa^{\mathrm{nt}}$ differ by a depth-zero character. Replacing $\chi$ by the product of $\chi$ with this character, we may assume that $\chi \otimes \kappa_{\pm}$ and $\kappa^{\mathrm{nt}}$ agree. It remains to show that the restriction of the irreducible representation $\chi^{-1} \otimes \rho^0$ to  $G^0_{x,0}/G^0_{x,0+}$ is cuspidal. This follows from the restriction of $\rho^0$ being cuspidal and all characters of $K/K^\dagger$ and hence of $K_{\pm}/K^\dagger$ (which is non-canonically isomorphic to $\{\pm 1\} \times K/K^\dagger$
) being trivial on the unipotent radicals of parabolic subgroups of $G^0_{x,0}/G^0_{x,0+}$ (cf.\ \cite[Theorem~C.1]{FS25}).
	
	To observe that $\pi_{\Upsilon}$ only depends on the representation $\rho^0 \otimes_{i=1}^n \phi_i|_{G^0_x}$ of $G^0_{\pm, x}$ (which encodes the information about the twisted Levi sequence $G_1, \hdots, G_{n+1}$ and the depths $r_1> \hdots >r_n$), let $\Upsilon'= ((G_i)_{1\leq i\leq n+1},x,(r_i)_{1\leq i\leq n},{\rho^0}',(\phi_i')_{1\leq i\leq n})$ be a different supercuspidal $G$-datum with $\rho^0 \otimes_{i=1}^n \phi_i|_{G^0_x}={\rho^0}' \otimes_{i=1}^n \phi_i'|_{G^0_x}$. Since $\rho^0$ and ${\rho^0}'$ have depth zero, the character $\varphi:=\otimes_{i=1}^n \phi_i|_{G^0_x}(\otimes_{i=1}^n\phi'_i|_{G^0_x})^{-1}$ has depth-zero, and $\hat \phi = \varphi \hat \phi'$ and $\varphi \otimes \rho^0 = {\rho^0}'$, where $\varphi$ also denotes the extension of $\varphi$ that is trivial on $K^+$. Hence by construction of $\kappa_{\pm}$ and $\kappa'_{\pm}$ from $\Upsilon$ and $\Upsilon'$, respectively, we obtain that $\kappa_{\pm}=\varphi \otimes \kappa'_{\pm}$. Hence ${\rho^0}'\otimes \kappa'_{\pm} \simeq \varphi \otimes \rho^0 \otimes \varphi^{-1} \otimes \kappa= \rho^0 \otimes \kappa$, and therefore their compact induction $\pi_{\Upsilon}$ and $\pi_{\Upsilon'}$ are isomorphic.
\end{proof}

\begin{rem}\label{rem:finitesumofsc}
	If one uses the weaker notion of generic characters as defined in \cite[Definition~4.1.1]{FS25}, then one can still perform all the above constructions with the only difference that $\cind_K^{G(F)}(\rho^0 \otimes \kappa_{\pm})$ is no longer irreducible in general even if $\rho^0$ is irreducible. Nevertheless $\cind_K^{G(F)}(\rho^0 \otimes \kappa_{\pm})$ is a finite direct sum of irreducible supercuspidal representations by using \cite[Theorem~4.6.9.(b)]{FS25} instead of \cite{Fin21Comp, Yu01}.
\end{rem}

\subsection{Comparison with Yu's construction}\label{sub:Yucomparison}
Let \(\Upsilon= ((G_i)_{1\leq i\leq n+1},x,(r_i)_{1\leq i\leq n},\rho^0,(\phi_i)_{1\leq i\leq n}) \) be a genuine supercuspidal $G$-datum as defined in \S\ref{sub:input}, and let $\chi_{G^0}$ be the restriction to $K_{\pm}$ of the genuine character of $G^0(F)_{\pm}$ attached to $\chi$-data as in Theorem \ref{thm:main1b}.

\begin{lem}
	The representation $\rho^0 \otimes \chi_{G^0}^{-1}$  of $K_{\pm}$ factors through $K/K^\dagger$ and its restriction to  $G^0_{x}/G^0_{x,0+}$ is irreducible and cuspidal.
\end{lem}
\begin{proof}
	This follows from $\chi_{G^0}^{-1}$ being genuine and trivial on the unipotent radicals of the parabolic subgroups of $G^0_{x,0}/G^0_{x,0+}$
\end{proof}

As a consequence of the lemma, we may also write $\rho^0 \otimes \chi_{G^0}^{-1}$ for the representation of $K$ through which it factors. Then the tuple 
\[\Upsilon_\chi=((G_i)_{1\leq i\leq n+1},x,(r_i)_{1\leq i\leq n},(\rho^0 \otimes \chi_{G^0}^{-1})|_{G^0_x},(\phi_i)_{1\leq i\leq n}) \] 
is a non-genuine supercuspidal $G$-datum and hence the input to which Yu (\cite{Yu01}) attaches the irreducible supercuspidal representation 
\[\pi_{\Upsilon_\chi}^{\mathrm{nt}}:=\cind_K^{G(F)}( \rho^0 \otimes \chi_{G^0}^{-1} \otimes \kappa^{\mathrm{nt}}) ,\]
 where $\kappa^{\mathrm{nt}}$ is defined as follows.
\begin{dfn}\label{dfn:kappant}
	The representation $\kappa^{\mathrm{nt}}$ is the unique representation of $K$ such that the composition of $G^0_{x} \ltimes K^\dagger \twoheadrightarrow K$ with $\kappa^{\mathrm{nt}}$ yields $(\hat \phi \ltimes 1)\otimes  (\mathrm{LWeilHeis} \circ (\bar\jmath \ltimes f))$. \end{dfn}

We also recall the following twisted version of Yu's construction from \cite[Theorem~4.1.13]{FKS23}:
\begin{dfn}\label{dfn:twistedYu}
	We define $\kappa:=\kappa^{\mathrm{nt}} \otimes \eps$, where $\eps=\prod_{i=1}^n\eps_{x}^{G_i/G_{i+1}}$ is as defined in \cite[Theorem~4.1.13]{FKS23}, and set $\pi_{\Upsilon_\chi}:=\cind_K^{G(F)}(\rho^0 \otimes \chi_{G^0}^{-1} \otimes \kappa)$.
\end{dfn}

\begin{pro}\label{pro:agreeswithYu}
	We have $\pi_{\Upsilon} \simeq \pi_{\Upsilon_\chi}$.
\end{pro}
\begin{proof}
	By construction it suffices to show that $\rho^0 \otimes \kappa_{\pm} \simeq \rho^0 \otimes \chi_{G^0}^{-1} \otimes \kappa^{\mathrm{nt}} \otimes \eps$. This follows from the construction of $\kappa_{\pm}$ and $\kappa^{\mathrm{nt}}$ and the relation $\tilde\jmath^*\mathrm{FWeilHeis} = \bar\jmath^*\mathrm{LWeilHeis} \otimes \chi_{G^0}^{-1} \otimes \epsilon$ from Theorem \ref{thm:main1b}.
\end{proof}

\begin{cor}
	If $p$ does not divide the order of the absolute Weyl group of $G$, then all supercuspidal representations arise from the construction in \S\ref{sub:construction}.
\end{cor}
\begin{proof}
	This follows by observing that every input for Yu's construction is of the form $\Upsilon_\chi$ for some $\Upsilon$ as above (since $\chi_{G^0}$ is trivial on unipotent radicals of parabolic subgroups of $G^0_{x,0}/G^0_{x,0+}$), and then combining Proposition \ref{pro:agreeswithYu} with \cite[Theorem~8.1]{Fin21}.
\end{proof}

\subsection{More general notion of supercuspidal $G$-datum} \label{sub:alternativeinput}
In \S\ref{sub:construction} we have attached a supercuspidal representation to an essentially depth-zero representation of $G^0_{\pm, x}$ that admits a supercuspidal $G$-datum. In the construction we chose and used a potentially reducible supercuspidal $G$-datum whose attached essentially depth-zero representation is $\rho$, but we observed that the resulting supercuspidal representation did not depend on this choice. For applications and to offer a different viewpoint on the input, one may therefore also generalize the notion of supercuspidal $G$-datum as follows.
\begin{dfn}
	We call a tuple 
	\[\Upsilon= ((G_i)_{1\leq i\leq n+1},x,(r_i)_{1\leq i\leq n},\rho^0,(\phi_i)_{1\leq i\leq n}) \]
	for some non-negative integer $n$, where $((G_i)_{1\leq i\leq n+1},x,(r_i)_{1\leq i\leq n},\rho^0)$ are as in (a)-(d) of Definition \ref{dfn:supercuspidalGdatum} satisfying Properties (i)-(iii) of Definition \ref{dfn:supercuspidalGdatum} and 
	\begin{enumerate}[(a)]
		\setcounter{enumi}{4}
		\item $\phi_i$, for $1 \leq i \leq n$, is a depth-$r_i$ character of the double cover $G_{i+1}(F)_{\pm^i}$ of $G_{i+1}(F)$ associated to the $\Sigma$-set $R(Z_{G_{i+1}},G_i)$, 
	\end{enumerate}
	that satisfy the following conditions 
	\begin{enumerate}[(i)]
		\setcounter{enumi}{3}
		\item $\phi_i$ is $(G_i, G_{i+1})$-generic relative to $x$ of depth $r_i$ for all $1 \leq i \leq n$, i.e., $\phi_i$ is trivial on $(G_{i+1})_{x,r_i+}$ and the restriction of $\phi_i$ to $(G_{i+1})_{x,r_i}/(G_{i+1})_{x,r_i+}$ is given by a $(G_i, G_{i+1})$-generic element $X_i \in \mathrm{Lie}^*(G_{i+1})(F)$ of depth $-r_i$ as in \cite[Definition~3.8]{Fintzen-IHES} \item either all $\phi_i$ for $1 \leq i \leq n$ are genuine or all $\phi_i$ for $1 \leq i \leq n$ are non-genuine 
	\end{enumerate}
	a \emph{supercuspidal $G$-datum}.\index[terminology]{supercuspidal $G$-datum} 
	
	Note that $G^0(F)_{\pm}$ is the Baer sum  over $1 \leq i \leq n$ of the double covers of $G^0(F)$ pulled back from the double covers $G_{i+1}(F)_{\pm^i}$ of $G_{i+1}(F)$. Since all the characters $\phi_i$ are either all genuine or all non-genuine, the product of their pull backs defines a character $\otimes_{i=1}^n \phi_i$ on $G^0(F)_{\pm}$. 
	We call a supercuspidal $G$-datum \emph{genuine} if $\rho^0 \otimes \otimes_{i=1}^n \phi_i|_{G^0_{\pm, x}}$ is genuine.\index[terminology]{supercuspidal $G$-datum!genuine} 

	We call a tuple \(\Upsilon= ((G_i)_{1\leq i\leq n+1},x,(r_i)_{1\leq i\leq n},\rho^0,(\phi_i)_{1\leq i\leq n}) \)
	as above where we replace the requirement that $\rho^0$ is irreducible by the requirement that it is finite dimensional semi-simple genuine or non-genuine a \emph{potentially reducible supercuspidal $G$-datum}.\index[terminology]{supercuspidal $G$-datum!potentially reducible}
	
	Given a potentially reducible supercuspidal $G$-datum \(\Upsilon= ((G_i)_{1\leq i\leq n+1},x,$ $(r_i)_{1\leq i\leq n},\rho^0,(\phi_i)_{1\leq i\leq n}) \), we refer to the representation $\rho:=\rho^0 \otimes \otimes_{i=1}^n \phi_i|_{G_{\pm,x}^{0}}$ of $G^0_{\pm,x}$ as \emph{the essentially depth-zero representation attached to $\Upsilon$}.\index[terminology]{essentially depth-zero representation}
\end{dfn}
One can perform the construction in \S\ref{sub:construction} directly from such a more general (potentially reducible) genuine supercuspidal $G$-datum, but in this case $\hat \phi$ has to be replaced by a potentially genuine character, which means that $\kappa_{\pm}$ might potentially no longer be genuine, but the representation $\rho^0 \otimes \kappa_{\pm}$ of $K_{\pm}$ still factors through a representation of $K$ and the main result, Theorem \ref{thm-screps}, holds true.

\subsection{{Supercuspidal representations from genuine tame $F$-non-singular elliptic pairs}{}}
\label{sub:nonsingular}
A large subset of the supercuspidal representations constructed by Yu has been parameterized in terms of elliptic maximal tori and characters thereof. In this section we will parameterize these representations using characters of double covers, which makes them fit more canonical into the local Langlands correspondence. Starting after Definition \ref{dfn:nonsingularpairs} we will assume that $p$ is not a bad prime of the absolute root system of $G$ and does not divide the orders of $\pi_1(G_{\mathrm{der}})$ and $\pi_0(Z(G))$.

\begin{dfn} \label{dfn:splitunram}
	Let $S \subset G$ be a maximal torus. Let $S' \subset S$ be any subtorus to which every symmetric ramified member of $R(S,G)$ has zero restriction. The \emph{canonical} splitting of $S(F)_\pm$ over $S'(F)$ is the map $S'(F) \to S(F)_\pm$ constructed as follows.
\end{dfn}
	Recall from \cite[\S3.1]{KalDC} that $S(F)_\pm$ is defined as the pushout of the big cover $S(F)_{\pm\pm}$, so it is enough to produce a map $S'(F) \to S(F)_{\pm\pm}$. In turn, the big cover is defined via pull-back of the covers $J_O(F)_\pm \to J_O(F)$ for each $\Sigma$-orbit $O \subseteq R(S,G)$. It is therefore enough to construct for each $\alpha \in R(S,G)$ a lift $S'(F) \to J_\alpha(F)_\pm$ of the restriction of $\alpha : S(F) \to J_\alpha(F)$ to $S'(F)$.

	If $\alpha$ is ramified, then its restriction to $S'$ vanishes by definition. If $\alpha$ is unramified or asymmetric, the cover $J_\alpha(F)_\pm \to J_\alpha(F)$ splits canonically: in the asymmetric case by definition, and in the symmetric unrammified case by \cite[Remark 3.6.1(1)]{KalDC}.

\begin{rem}
	An example of $S'$ in the above definition is the maximal unramified subtorus of $S$, since by definition $X^*(S')$ is the torsion-free quotient of $X^*(S)_I$.
\end{rem}

\begin{rem} \label{rem:spllevi}
	The double cover of a twisted Levi subgroup $G' \subset G$ is equipped with a splitting over $G'_\tx{sc}(F)$ given as follows. According to \cite[Definition~6.4.1]{KalDC}, the double cover $G'(F)$ is pulled back from a double cover of $[G'/G'_\tx{der}](F)$. Since the map $G'_\tx{sc}(F) \to G'(F)$ has trivial composition with the map $G'(F) \to [G'/G'_\tx{der}](F)$, the pull-back to $G'_\tx{sc}(F)$ of the double cover is canonically split. Note that if $G'$ is isotropic then $G'_\tx{sc}(F)$ has no non-trivial characters, so this splitting is unique.
\end{rem}

\begin{dfn} \label{dfn:nonsingularpairs}
	We call a pair $(S, \theta)$ consisting of a maximal torus $S \subset G$ and a character $\theta: S(F)_{\pm} \ra \bC^\times$, where $S(F)_{\pm}$ denotes the double cover associated to $R(S,G)$, a \emph{tame $F$-non-singular elliptic pair}\index[terminology]{tame $F$-non-singular elliptic pair} if the following conditions are satisfied:
	\begin{enumerate}
		\item $S$ is elliptic and splits over a tamely ramified field extension $E/F$.
		\item \label{cond:G0} The action of the inertia subgroup of the absolute Galois group of $F$ on the root system
		\[ R_{0+}= \{ \alpha \in R(S,G) \, | \, \theta|_{S(F)_{0+}}\circ \rN_{E/F}\circ \alpha^\vee(1+\mf{p}_E)=1 \} \]
		preserves a subset of positive roots of $R_{0+}$, where $\theta|_{S(F)_{0+}}$ is the composition of $\theta$ with the unique splitting of $S(F)_\pm$ over $S(F)_{0+}$.
		\item Let $S'$ be the maximal unramified subtorus of $S$, and let $F'/F$ be a finite, unramified field extension over which  $S'$ is split. We write $R_{0+}(S')$ for the restriction of the roots in $R_{0+}$ to $S'$. Note that $R_{0+}$ corresponds to a reductive subgroup $G^0$ of $G$, and $R_{0+}(S')$ is the relative root system of $G^0 \times_FF^{\mathrm{ur}}$. We ask that for every $\alpha' \in R_{0+}(S')$, the character 
		\[\theta|_{S'(F)} \circ \rN_{F'/F} \circ \alpha'^\vee :F'^\times \ra \bC^\times \]
		has non-trivial restriction to $\cO_{F'}^\times$, where $\theta|_{S'(F)}$ is the composition of $\theta$ with the splitting of Definition \ref{dfn:splitunram}
	\end{enumerate} 
	We call such a pair $(S,\theta)$ \emph{genuine}\index[terminology]{genuine} if $\theta$ is a genuine character of $S(F)_{\pm}$, and \emph{non-genuine}\index[terminology]{non-genuine} otherwise.
\end{dfn}
Note that if $(S,\theta)$ is a non-genuine tame $F$-non-singular elliptic pair, then viewing $\theta$ as a character of $S(F)$ recovers the notion of tame $F$-non-singular elliptic pairs introduced in \cite[Definition~3.4.1]{KalSLP}.

Let $(S,\theta)$ be a genuine tame $F$-non-singular elliptic pair. We will now attach a finite direct sum of supercuspidal representations to this pair via the construction in \S\ref{sub:construction}, for which we first exhibit a genuine potentially reducible supercuspidal $G$-datum $\Upsilon$ as in \S\ref{sub:alternativeinput}. 

For this, we employ the Howe factorization algorithm of \cite[\S3.6]{KalRSP}. While this algorithm was phrased without covers, it extends immediately to the case of covers. More precisely, let $(S,\theta)$ be a pair consisting of a tame maximal torus $S \subset G$ and a genuine character $\theta : S(F)_\pm \to \C^\times$. Using that $S(F)_\pm$ splits uniquely over $S(F)_{0+}$ the filtration \cite[(3.6.1)]{KalRSP} of the root system $R(S,G)$ is defined and provides a tower of twisted Levi subgroups $G^0 \subsetneq \dots \subsetneq G^d=G$ of $G$, which we write as $G_{n+1} \subsetneq \dots \subsetneq G_2 \subseteq G_1=G$, where we set $G_1=G_2$ if and only if $r_d>r_{d-1}$ in the notation of \cite[\S3.6]{KalRSP} (which could also be phrased more directly in the process of performing the Howe factorization). We let $G_{n+2}:=S$ and for $i=1,\dots,n+1$ write $G_{i+1}(F)_{\pm}$ for the double cover associated to the earnest $\Sigma$-set $R(Z_{G_{i+1}},G_i)$. Note that, when $(S,\theta)$ is a tame $F$-non-singular elliptic pair, the double cover of $G_{n+2}(F)$ carries the canonical splitting of Definition \ref{dfn:splitunram}, because $S$ is maximally unramified in $G_{n+1}$, and note that $S(F)_{\pm}$ is in general not the same as $G_{n+2}(F)_\pm$. 

\begin{dfn} \label{dfn:howe} A \emph{Howe factorization} of $(S,\theta)$ is a sequence of genuine characters $\phi_i : G_{i+1}(F)_\pm \to \C^\times$ for $i=1,\dots,n+1$ with the following properties. 
\begin{enumerate}
	\item \begin{equation} \label{eq:thetaprod} \theta = \prod_{i=1}^{n+1}\phi_i, \end{equation}
	where the product is to be interpreted as first restricting $\phi_i$ to the preimage of $S(F)$ in $G_{i+1}(F)_{\pm}$, and then taking the product of these restrictions as a character of the Baer sum of these preimages, which is canonically identified with $S(F)_\pm$.
	\item For all $1 \leq i \leq n$ the character $\phi_i$ is trivial on the canonical splitting of $(G_{i+1})_{\mathrm{sc}}(F)$ into $G_{i+1}(F)_\pm$ given by Remark \ref{rem:spllevi}. 
	\item For all $1 \leq i \leq n$, $\phi_i$ is $(G_i,G_{i+1})$-generic of depth $r_i$. 
For $i=n+1$, $\phi_{n+1}$ is trivial on the unique splitting of $G_{n+2}(F)_{0+}$ into $G_{n+2}(F)_\pm$.
\end{enumerate}
\end{dfn}

\begin{rem}
	 This is the analogue of the Howe factorization of non-genuine characters of $S(F)$ defined in \cite[\S3.6]{KalRSP}, except that we have adjusted to the convention of notation used in \cite{Fin21, Fin21Comp, FS25} rather then in \cite{Yu01} and \cite{KalRSP}, as the former is the one we used in the definition of a potentially reducible supercuspidal $G$-datum to which we want to connect in this section.  See \cite[Remark~2.4]{Fin21Comp} for how to translate between these two conventions to obtain the notation used in \cite{KalRSP}. Using this convention of notation and moreover allowing $\phi_{n+1}$, whose non-genuine analogue is denoted by $\phi_{-1}$ in \cite{KalRSP}, to be an arbitrary genuine depth-zero character even if $G^0=S$, allows us to present the conditions for a Howe factorization more uniformly rather than having to distinguish different cases to formulate the conditions on the characters $\phi_i$.
\end{rem}

Note the notation of refactorization of \cite[Definition 4.19]{HM08} also makes sense in this setting, in fact verbatim, because the characters $\chi_i$ defined in terms of two genuine Howe factorizations are themselves non-genuine. 

\begin{pro}
	Every pair $(S,\theta)$ consisting of a tame maximal torus and a genuine character $\theta : S(F)_\pm \to \C^\times$ has a Howe factorization. Moreover, any two Howe factorizations are refactorizations of each other.
\end{pro}
\begin{proof}
	The proofs of Lemma 3.6.6 and Proposition 3.6.7 of \cite{KalRSP} go through essentially verbatim, using the fact that $S(F)_\pm$ splits uniquely over $S(F)_{0+}$, $H(F)_\pm$ splits uniquely over $H(F)_{x,0+}$, and at the end of the proof of Lemma 3.6.9 one extends via Pontryagin duality from $D(F)_r$ to $D(F)_\pm$.
\end{proof}

Consider now a genuine tame $F$-non-singular pair $(S,\theta)$. Apply Howe factorization to obtain the sequence of characters $(\phi_1,\dots,\phi_n, \phi_{n+1})$. 
We interpret $\phi_{n+1}$ as a character of $G_{n+2}(F)=S(F)$ via the canonical splitting (Definition \ref{dfn:splitunram}) of $G_{n+2}(F)_{\pm}$. Let $\rho^0:=\kappa_{S,\phi_{n+1}}^{G^0}$ be the cuspidal representation of $G^0(F)_x/G^0(F)_{x,0+}$ of \cite[Definition 2.6.8]{KalSLP}, inflated to  $G^0(F)_x$, where $x \in \sB(G)_{G^0}$ denotes the unique point corresponding to $S$.
Then
\[ \Upsilon_{(S,\theta)}= ((G_i)_{1\leq i\leq n+1},x,(r_i)_{1\leq i\leq n},\rho^0,(\phi_i)_{1\leq i\leq n})\]
is a potentially reducible supercuspidal $G$-datum in the sense of \S\ref{sub:alternativeinput}, genuine if $n>0$. When $n=0$ the datum just consists of the cuspidal $\rho^0$ of $H_x/H_{x,0+}$ inflated to $H_x$, and compact induction to $G(F)$ produces a depth-zero supercuspidal representation according to \cite{MP96}. When $n>0$ we have the following result.

\begin{lem}\label{lem:Sthetarep}
	The representation $\pi_{(S,\theta)}$\index[notation]{piStheta@$\pi_{(S,\theta)}$} constructed from $\Upsilon_{(S,\theta)}$ as in \S\ref{sub:construction}, which is a finite direct sum of supercuspidal representations by Theorem \ref{thm-screps}, only depends on $(S,\theta)$ and not on the choice of Howe factorization. 
\end{lem}
\begin{proof}
	This follows by combining Theorem \ref{thm-screps} and \cite[Lemma~2.6.10]{KalSLP}.
\end{proof}

\subsection{Character formula at shallow elements} \label{sub:char}

\begin{thm}\label{thm:characterformula}
	Let $(S,\theta)$ be a genuine tame $F$-non-singular elliptic pair.
	Let $\gamma \in G(F)$ be regular elliptic and topologically semi-simple modulo center. The character of $\pi_{(S,\theta)}$ vanishes at $\gamma$ unless $\gamma$ is conjugate to an element of $S(F)$. If $\gamma \in S(F)$, then the normalized character of $\pi_{(S,\theta)}$ at $\gamma$ is given by \cite[(4.3.2)]{KalDC}, i.e.,
	\[ e(G)\epsilon(1/2,X^*(T)_\C-X^*(S)_\C,\psi)\sum_{w \in N(S,G)(F)/S(F)}[a_S\cdot\theta](\gamma^w), \]
	where $a_S : S(F)_\pm \to \{\pm1\}$ is the genuine function associated to the $a$-data given by $a_\alpha=\<X,H_\alpha\>$ for $\alpha \in R(S,G/G^0)$ and $a_\alpha \in \cO_\alpha^\times$ for $\alpha \in R(S,G^0)$, $e(G)$ is the Kottwitz sign of $G$, and $T$ is the universal maximal torus of $G$.
\end{thm} \begin{proof}
	We have $\pi_{(S,\theta)}=\tx{c-Ind}_K^{G(F)}\sigma$, where $\sigma=\rho^0 \otimes \kappa_{\pm}$ is constructed as above from a genuine supercuspidal $G$-datum $\Upsilon_{(S,\theta)}$. We follow the arguments of \cite[\S4.4,\S4.10]{KalRSP}, which are based on \cite{AS09}, especially the proof of Theorem 6.4 there. We arrive at 
\[ \Theta_\pi(\gamma) = (-1)^{\tx{rk}_F(G^0)-\tx{rk}_F(S)}\sum_{\substack{g \in S(F) \lmod G(F)/S(F)\\ {^g\gamma} \in S(F)}} \tx{tr}\,(\tilde\jmath^*\mathrm{FWeil} \otimes\theta)({^g\gamma}), \]
	which shows that the sum vanishes unless $\gamma$ is $G(F)$-conjugate to an element of $S(F)$. Assuming $\gamma \in S(F)$ the above becomes
	\[ \Theta_\pi(\gamma) = (-1)^{\tx{rk}_F(G^0)-\tx{rk}_F(S)}\sum_{w \in N(S,G)(F)/S(F)} \tx{tr}\,(\tilde\jmath^*\mathrm{FWeil} \otimes\theta)({^w\gamma}). \]
	Choose a lift $\dot\gamma\in S(F)_\pm$ of $\gamma$ whose image under $\jmath$ is a norm; such a lift exists by Lemma \ref{lem:topssgamma}(5). Then $\tilde \jmath(\dot \gamma) = \gamma_\psi(q_{\jmath(\dot \gamma)}) \boxtimes \bar\jmath(\gamma)$ and using \cite[Theorem B]{Thomas2013} we see
	\[ \tx{tr}\,(\tilde\jmath^*\mathrm{FWeil} \otimes\theta)({^w\gamma}) = \tx{tr}\,\mathrm{FWeil}(\tilde\jmath({^w}\dot\gamma)) \cdot \theta({^w\dot\gamma}) = \gamma_\psi(q_{\jmath(^w\dot\gamma)})\cdot\theta({^w}\dot\gamma). \]
	To evaluate the right-hand side we lighten notation and replace ${}^w\dot\gamma$ by $\dot\gamma$. To compute $\gamma_\psi(q_{\jmath(\dot\gamma)})$ we fix a lift $\tilde\gamma \in S_{\GL}(F)$ of $\jmath(\dot\gamma)$. Let $(a'_\alpha)$ denote the $a$-data from Lemma \ref{lem:bt-align1}(2). Then the quadratic form $q_{a'}$ induced by $\eta_{\lambda,a'}$ is the form of \cite[Definition 5.4.5]{FKS23}. Apply Lemma \ref{lem:topssgamma}(4) to the torus $T_{\Sp}=S_{\Sp}$ and observe that $\Lambda=R(S,G/G^0)$ and that $\lambda(\gamma) \neq 1$ for all $\lambda \in \Lambda$ by the regularity of $\gamma$ to obtain
	\[ \gamma_\psi(q_{\tilde \gamma}) = \gamma_\psi(q_{a'}|V) \cdot \prod_{\alpha \in R(S,G/G^0)_\tx{sym}/\Gamma} \kappa_\alpha\left(\frac{\tilde \gamma_\alpha-\tilde \gamma_{-\alpha}}{a'_\alpha}\right). \]
	In \cite[\S1.3.3]{KotWeil} Kottwitz defines a vector space $V=\tx{Lie}(G)/\tx{Lie}(S)$ and a quadratic from $Q_V$ on it, which we shall henceforth call $Q_G$. From the definition it is clear that this form is the orthogonal direct sum $Q_{G^0} \oplus Q_{G/G^0}$, where $Q_{G^0}$ is the Kottwitz form on $V^0=\tx{Lie}(G^0)/\tx{Lie}(S)$, and $Q_{G/G^0}$ is a quadratic form on $V_0=\tx{Lie}(G)/\tx{Lie}(G^0)$.The form $q_{a'}$ differs from $Q_{G/G^0}$
	by the factor $\ell_{p'}(\alpha^\vee)$ on $V_{\pm\alpha}$ for each $\alpha \in R(S,G/G^0)$. In fact, we will need to use $-Q_{G/G^0}$ in place of $Q_{G/G^0}$, which also introduces the factor $-1$ on each $V_{\pm\alpha}$. This implies using Fact \ref{fct:gammaind} and Lemma \ref{lem:weilherm} that
	\[ \gamma_\psi(q_{a'}) = \gamma_\psi(-Q_{G/G^0})\cdot \prod_{\alpha \in R(S,G/G^0)_\tx{sym}/\Gamma} \kappa_\alpha(-\ell_{p'}(\alpha^\vee)). \]
	Applying \cite[Theorem 1.1]{KotWeil} to $(S,G^0)$ and $(S,G)$ we obtain
	\[ \gamma_\psi(Q_{G/G^0})=e(G)e(G^0)\epsilon(1/2,X^*(T^{G^0})_\C-X^*(T^G)_\C,\psi), \]
	where $T^{G^0}$ and $T^G$ are the universal maximal tori of $G^0$ and $G$. Recalling that $a'_\alpha=\ell_{p'}(\alpha^\vee)/\<X,H_\alpha\>$, we obtain
	\[ \gamma_\psi(q_{\tilde\gamma}) = e(G)e(G^0)\epsilon(1/2,X^*(T)_\C-X^*(T^{G^0})_\C,\psi)\cdot \prod_{\alpha \in R(S,G/G^0)_\tx{sym}/\Gamma}\kappa_\alpha\left(\frac{\tilde \gamma_\alpha-\tilde \gamma_{-\alpha}}{-\<X,H_\alpha\>^{-1}}\right). \]
	Since $\<X,H_\alpha\> \in F_\alpha^\times$ has zero trace to in $F_{\pm\alpha}$ we see that $-\<X,H_\alpha\>^{-1}$ equals $\<X,H_\alpha\>$ multiplied by the norm of $\<X,H_\alpha\>^{-1}$. Therefore the product over $\alpha$ is precisely $a_S(\dot \gamma)$. On the other hand, according to \cite[Lemma 4.9.1]{KalRSP} we have $e(G^0)\epsilon(1/2,X^*(T^{G^0})_\C-X^*(S)_\C,\psi)=(-1)^{\tx{rk}_F(G^0)-\tx{rk}_F(S)}$, while the last paragraph of the proof of \cite[Proposition 4.9.2]{KalRSP} shows that $\kappa_\alpha((\tilde\gamma_\alpha-\tilde\gamma_{-\alpha})/a_\alpha)=1$ when $\alpha \in R(S,G^0)_\tx{sym}$, hence 
	\[ (-1)^{\tx{rk}_F(G^0)-\tx{rk}_F(S)} \gamma_\psi(q_{\tilde\gamma}) = e(G)\epsilon(1/2,X^*(T)_\C-X^*(S)_\C,\psi)a_S(\dot\gamma). \]
\end{proof}

\subsection{The local Langlands correspondence for supercuspidal parameters} \label{sub:llc}

Let $(S,\theta)$ be a genuine tame $F$-non-singular elliptic pair and let $\pi:=\pi_{(S,\theta)}$ be the finite direct sum of supercuspidal representations attached to it.
Let $\varphi_S : W_F \to {^LS}_\pm$ be  the parameter of $\theta : S(F)_\pm \to \C^\times$ in the sense of \cite[Theorem 3.4.3]{KalDC}. Let $^Lj : {^LS}_\pm \to {^LG}$ be the canonical $L$-embedding  of \cite[\S4.1]{KalDC}.

\begin{thm} \label{thm:llc}
	The composition of $\varphi_S$ and $^Lj$ equals the $L$-parameter of any irreducible summand of $\pi$ as constructed in \cite[\S4]{KalSLP}.
\end{thm}
\begin{proof}
	We recall the $L$-parameter  $\varphi_\pi$ associated to $\pi$ by the construction of \cite[\S4.1,\S4.2]{KalSLP}. We represent this schematically as 
	\[ \pi \mapsto (S,\hat\jmath,\chi_0,\theta',(G,\xi,z),j) \mapsto \varphi_\pi, \] 
	where the middle datum is called a semi-simple Deligne--Lusztig packet datum in \cite{KalSLP}. We describe first the contents of this datum: $(\xi,z) : G_0 \to G$ is a rigid inner twist relating $G$ to its quasi-split inner form, $j : S \to G$ is an embedding, $\hat\jmath : \hat S \to \hat G$ is an embedding dual to $j$ as discussed in \cite[\S5.1]{KalRSP}, $\theta' : S(F) \to \C^\times$ is an $F$-non-singular character, and $\chi_0$ is a set of $\chi$-data for $R(Z_{G^0},G)$, where $Z_{G^0}$ is the center of the smallest twisted Levi subgroup in a potentially reducible supercuspidal $G$-datum from which $\pi$ can be constructed.

	The association between this datum and $\varphi_\pi$ is that $\varphi_\pi = {^Lj}_{\chi_0}\circ{^L\theta'}$, where $^L\theta' : W_F \to {^LS}$ is the $L$-parameter of $\theta'$ and $^Lj_{\chi_0}$ is the $L$-embedding $^LS \to {^LG}$ produced by the recipe of \cite[\S2.6]{LS87} from the $\chi$-data for $R(S,G)$ obtained by inflating $\chi_0$ to $\chi$-data for $R(S,G/G^0)$ in the sense of \cite[Definition 5.3.5]{KalDC} and adding to this the canonical unramified $\chi$-data for $R(S,G^0)$. 
	
	The association between this datum and $\pi$ is that $\pi$ should be the (generally reducible) supercuspidal representation $\pi_\Upsilon$ obtained from the twist of Yu's construction in \cite{FKS23}, as in Definition \ref{dfn:twistedYu} above, applied to the non-genuine supercuspidal $G$-datum $\Upsilon$ that arises from applying the Howe factorization algorithm of \cite[\S3.6]{KalRSP} to the pair $(jS,j\theta')$. Since so far this association does not depend on $\chi_0$, we must link $\chi_0$ to the remaining pieces of the datum to ensure that $\pi \mapsto \varphi_\pi$ is well-defined. This is done in \cite[\S4.2]{KalSLP}, which refers to \cite[\S3.5]{KalSLP}, which specifies $\chi_0$ to be given as in \cite[\S4.2]{FKS23} in terms of the character $\theta'$. 

	We now come to the actual proof. We have
	\[ \varphi_\pi = {^Lj_{\chi_0}}\circ{^L\theta'} = {^Lj}\circ {^L(\theta'\cdot\chi_S)}, \]
	where $^Lj : {^LS}_\pm \hookrightarrow {^LG}$ is the canonical $L$-embedding and $\chi_S$ is the genuine character of $S(F)_\pm$ obtained by the $\chi$-data for $R(S,G)$ used to construct ${^{L}}j_{\chi_0}$. The first identity is by construction and the second is immediate from the constructions of the objects $^Lj_{\chi_0}$, $^Lj$, and $^L(\theta'\cdot\chi_S)$. Under the Baer-sum decomposition of $S(F)_\pm$ into the covers associated to $R(S,G/G^0)$ and $R(S,G^0)$, the character $\chi_S$ is the product of the pullback of $\chi_{G^0}$ and the genuine character attached to the canonical unramified $\chi$-data for $R(S,G^0)$. The latter is trivial on the canonical splitting. Applying the comparison of Proposition \ref{pro:agreeswithYu} factor by factor therefore gives $\theta=\theta'\chi_S$, which proves the claim.
\end{proof}

\subsection{The example of $\tx{PGL}_2$} \label{sub:pgl2}

A classical construction of the irreducible representations of $\tx{GL}_2(F)$ when $p \neq 2$ proceeds via the Weil representation over $F$. An early account of this appears in \cite{JL79} and a more recent exposition in \cite[Chapter 9]{BH06}. We give here a short summary of the construction and compare it with the one of this paper for the group $\tx{PGL}_2$.

Working through \cite[\S39]{BH06} one can recognize the classical construction as an application of the similitude theta correspondence, for a description of which we refer to \cite[\S20.4.2]{GanKudlaTakedaTheta}. The general setting is that of a symplectic space $W$ and an even-dimensional orthogonal space $V$, both over $F$. Then $V \otimes W$ is again a symplectic space and acting on each tensor factor separately leads to the map $\tx{GO}(V) \times \tx{GSp}(W) \to \tx{GSp}(V \otimes W)$. Under this map the similitude character of $\tx{GSp}(V \otimes W)$ restricts to the product of the similitude characters $\mathrm{sim}$ of $\tx{GO}(V)$ and $\tx{GSp}(W)$, hence the map restricts to a map $\mc{G} \to \tx{Sp}(V \otimes W)$, where 
\[ \mc{G} = \{(g,h) \in \tx{GO}(V) \times \tx{GSp}(W)|\tx{sim}(g)\cdot \tx{sim}(h)=1\}. \]
Brooks Roberts \cite{Roberts96} constructs a lifting of this map to a map $\mc{G} \to \tx{Mp}(V \otimes W)$, via which the Weil representation can be restricted and decomposed, yielding a theta correspondence between irreducible representations of $\tx{GO}(V)$ (we note that the projection $\tx{GO}(V) \times \tx{GSp}(W) \to \tx{GO}(V)$ remains surjective on $\mc{G}$) and irreducible representations of the image $\tx{GSp}(W)'$ of the restriction to $\mc{G}$ of the projection $\tx{GO}(V) \times \tx{GSp}(W) \to \tx{GSp}(W)$.

To construct representations of $\tx{GL}_2(F)$ this general theory is specialized to the setting of $W=F \oplus F$ equipped with the standard symplectic form $x_1y_2-x_2y_1$, and $V=E$ for a separable quadratic extension $E/F$ equipped with the norm map $E \to F$. The action of $E^\times$ on $E$ by multiplication gives an embedding $E^\times \to \tx{GO}(V)$ of index $2$, the missing element being supplied by the non-trivial automorphism $\sigma$ of $E/F$. Thus every character $\phi : E^\times \to \C^\times$ with $\phi\circ \sigma \neq \phi$ has irreducible induction to $\tx{GO}(V)$ and via the similitude theta correspondence provides an irreducible representation of $\tx{GSp}(W)'$. The latter group is the subgroup of $\tx{GL}_2(F)$ of elements whose determinant is a norm from $E$. This is again a subgroup of index $2$ in $\tx{GL}_2(F)$, and induction provides an irreducible supercuspidal representation $\pi_\phi$ of $\tx{GL}_2(F)$.

We shall now compare the assignment $\phi \mapsto \pi_\phi$ with the construction of this paper. Let $G=\tx{PGL}_2$ and let $S \subset G$ be the anisotropic maximal torus associated to $E/F$. 

\begin{lem}
	The representation $\pi_\phi$ descends to $\mathrm{PGL}_2(F)$ if and only if $\phi$ induces a genuine character of the double cover $S(F)_\pm$.
\end{lem}
\begin{proof}
To descend from $\tx{GL}_2(F)$ to $\tx{PGL}_2(F)$ the representation $\pi_\phi$ must have trivial central character. It is shown in \cite[Proposition 39.2(3)]{BH06} that the central character of $\pi_\phi$ equals $\phi \cdot \kappa_{E/F}$, hence the condition is that $\phi$ is trivial on $N_{E/F}(E^\times) \subset E^\times$ but non-trivial on $F^\times \subset E^\times$, which means that $\phi$ induces a genuine character of the double cover $S(F)_{\pm} \simeq E^\times/N_{E/F}(E^\times)$.
\end{proof}

We now let $\phi: E^\times \ra \bC^\times$ be a character that induces a genuine character of the double cover $S(F)_{\pm} \simeq E^\times/N_{E/F}(E^\times)$, and we also view $\phi$ as a character of $S(F)_{\pm}$.

\begin{lem}
	Let $(S,\phi)$ be a tame $F$-non-singular elliptic pair in the sense of Definition \ref{dfn:nonsingularpairs}. The $L$-parameter of $\pi_\phi$ constructed in \cite[Theorem 40.1]{BH06} coincides with the composition of the $L$-parameter $W_F \to {^LS}_\pm$ of $\phi$ with the canonical $L$-embedding $^Lj : {^LS}_\pm \hookrightarrow {^LG}$.
\end{lem}
\begin{proof}
To compare the $L$-parameters, recall that Bushnell--Henniart define the $L$-parameter of $\pi_\phi$ as $\tx{Ind}_{E^\times}^{W_{E/F}}\phi$. Choosing a lift $w \in W_{E/F}$ of $\sigma \in \Gamma_{E/F}$ we can explicitly present this induced representation as 
	\[ E^\times \ni e \mapsto \begin{bmatrix}
		\phi(e)\\&\phi(e)^{-1}
	\end{bmatrix},\qquad w \mapsto \begin{bmatrix}
		0&1\\-1&0
	\end{bmatrix}.\]
	According to \cite[Appendix B.2]{KalDC} we have $^LS_\pm = \C^\times_{(-1)}\boxtimes \Gamma_{E/F}$, the $L$-parameter of $\phi$ is $e \mapsto \phi(e)\boxtimes 1$ and $w \mapsto 1 \boxtimes \sigma$, and the canonical $L$-embedding $^Lj$ is given by 
	\[ z \boxtimes 1 \mapsto \begin{bmatrix}
		z\\&z^{-1}
	\end{bmatrix},\qquad 1 \boxtimes \sigma \mapsto \begin{bmatrix}
		0&1\\-1&0
	\end{bmatrix}. \qedhere\]
\end{proof}

\begin{cor}
	The representation $\pi_\phi$ coincides with $\pi_{(S,\phi)}$ constructed in \S\ref{sub:nonsingular}.
\end{cor}
\begin{proof}
	It was shown by \cite{OT21} and \cite{Tok2023} that the $L$-parameter constructed in \cite{KalRSP} coincides with the one of \cite{BH06}. The claim now follows from the above lemmas and Theorem \ref{thm:llc}.
\end{proof}

\begin{rem}
One can make the above analysis for the group $\tx{GL}_2(F)$, not just $\tx{PGL}_2(F)$, but the coordinates become a bit more complicated.
		
\end{rem}

\subsection{Construction of types and reduction to essentially depth-zero} \label{sub:types}
Kim and Yu (\cite{KimYu}) modified the non-genuine supercuspidal datum slightly, mostly relaxing some assumptions, in order to obtain a construction of types following Yu's construction, which provides types for every Bernstein block if $p$ does not divide the order of the absolute Weyl group of $G$ (\cite{Kim07, Fin21}). We will adjust the input analogously to the definition of a genuine supercuspidal $G$-datum to also obtain a more ``natural'' construction of types. In order to do so, we first introduce some notation. 
We denote by $G(F)^1$ those elements of $G(F)$ that get send to torsion elements under the Kottwitz homomorphism, 
and for any subgroup $K\subseteq G(F)$ we write $K^1$ for $K \cap G(F)^1$.
Given a sequence $G=G_1 \supseteq G_2 \supsetneq G_3 \supsetneq \hdots \supsetneq G_{n+1}=:G^0$ of tame, twisted Levi subgroups and a point $x \in \cB(G,F)_{G^0} \subseteq \cB(G,F)$ we define Levi subgroups $M_i \subseteq G_i$ ($1 \leq i \leq n+1$) as follows. We let $M_{n+1}=M^0 \subseteq G_{n+1}=G^0$ be the Levi subgroup as defined in \cite[\S6.3]{MP96} that satisfies that the inclusion $(M^0)_{x,0}/(M^0)_{x,0+}\simto (G^0)_{x,0}/(G^0)_{x,0+}$ is an isomorphism and such that $x$ is contained in $\cB(G)_{M^0}$ and its image in $\cB(M^0)$ is a vertex. We let $Z_{\mathrm{split}}(M^0)$ be the maximal split torus in the center of $M^0$ and we define $M_i$ to be the centralizer of $Z_{\mathrm{split}}(M^0)$ in $G_i$.

\begin{dfn}\label{dfn:Gdatum} 
	We call a tuple 
	\[\Upsilon= ((G_i)_{1\leq i\leq n+1},x,(r_i)_{1\leq i\leq n},(K^0_{\pm},\rho^0),(\phi_i)_{1\leq i\leq n}) \]
	for some non-negative integer $n$, where  
	\begin{enumerate}[(a)]
		\item $G=G_1 \supseteq G_2 \supsetneq G_3 \supsetneq \hdots \supsetneq G_{n+1}=:G^0$ are twisted Levi subgroups of $G$ that split over a tamely ramified field extension of $F$,
		\item  $x \in \cB(G,F)_{G^0} \subseteq \cB(G,F)$, 
		\item $r_1 > r_2 > \hdots > r_n >0$ are real numbers,
		\item \label{item:rhotypes} $\rho^0$ is an irreducible depth-zero representation of a compact subgroup $K^0_{\pm, x}$ with $(G^0)_{\pm, x, 0} \subseteq K^0_{\pm, x} \subseteq (G^0)_{\pm, x, 0}(M^0)^1_{\pm, x}$, where $G^0(F)_{\pm}$ is the double cover attached to $R(Z_{G^0}, G)$ and $(M^0)^1_{\pm,x}$ is its pullback to $M^0(F)^1_x \subseteq G^0(F)$,
		\item	$\phi_i$, for $1 \leq i \leq n$, is a depth-$r_i$ character of the double cover $G_{i+1}(F)_{\pm^i}$ of $G_{i+1}(F)$ associated to the $\Sigma$-set $R(Z_{G_{i+1}},G_i)$, 
	\end{enumerate}
	that satisfies the following conditions 
	\begin{enumerate}[(i)]
		\setcounter{enumi}{1}
\item the point $x$ is chosen so that $(M_i)_{x,r_i/2}/(M_i)_{x,r_i/2+}\simto (G_i)_{x,r_i/2}/(G_i)_{x,r_i/2+}$ for $1 \leq i \leq n$
		\item  $\rho^0$ is cuspidal (as in Definition \ref{def:cuspidal}),
		\item $\phi_i$ is $(G_i, G_{i+1})$-generic relative to $x$ of depth $r_i$ for all $1 \leq i \leq n$, i.e., $\phi_i$ is trivial on $(G_{i+1})_{x,r_i+}$ and the restriction of $\phi_i$ to $(G_{i+1})_{x,r_i}/(G_{i+1})_{x,r_i+}$ is given by a $(G_i, G_{i+1})$-generic element $X_i \in \mathrm{Lie}^*(G_{i+1})(F)$ of depth $-r_i$ as in \cite[Definition~3.8]{Fintzen-IHES} \item either all $\phi_i$ for $1 \leq i \leq n$ are genuine or all $\phi_i$ for $1 \leq i \leq n$ are non-genuine 
	\end{enumerate}
	a \emph{$G$-datum}.\index[terminology]{G-datum@$G$-datum}
	
	We write $K^0_x$ for the image of $K^0_{\pm, x}$ in $G^0(F)$. 
	Since all the characters $\phi_i$ are either all genuine or all non-genuine, the product of their pull backs defines a character $\otimes_{i=1}^n \phi_i$ on $K^0_{\pm,x}$, which is the Baer sum  over $1 \leq i \leq n$ of the double covers of $K^0_x$ pulled back from the double covers $G_{i+1}(F)_{\pm^i}$ of $G_{i+1}(F)$.
	We call a $G$-datum \emph{genuine} if the $K^0_{\pm, x}$-representation $\rho:=\rho^0 \otimes_{i=1}^n \phi_i$ is genuine.\index[terminology]{G-datumgenuine@$G$-datum!genuine} 
	We refer to the representation $\rho$ of $K^0_{\pm,x}$ as \emph{the essentially depth-zero representation attached to $\Upsilon$}.\index[terminology]{essentially depth-zero representation}
	
	If the representations $\rho^0$ and $\phi_i$ ($1 \leq i \leq n$) are all non-genuine, then we call the $G$-datum \emph{truly non-genuine}.\index[terminology]{G-datumgenuine@$G$-datum!truly non-genuine}  
	\end{dfn}

\begin{rem}
	Considering the non-genuine representations $\rho^0$ and $\phi_i$ ($1 \leq i \leq n$) in a truly non-genuine $G$-datum as representations of $K^0_x$ and $G_{i+1}(F)$ (through which they factor), we recover the notion of a ``datum'' in \cite[\S7.2]{KimYu}, called ``$G$-datum'' in \cite[\S4.1]{AFMO2}. We present the datum slightly differently here to better match the above definition of a supercuspidal $G$-datum. 
\end{rem}

Let \(\Upsilon= ((G_i)_{1\leq i\leq n+1},x,(r_i)_{1\leq i\leq n},\rho^0,(\phi_i)_{1\leq i\leq n}) \) be a genuine $G$-datum. To align directly with \S\ref{sub:input}--\S\ref{sub:Yucomparison}, we assume that $\phi_i$ is non-genuine. But the construction can also be performed without this assumption leading to a genuine $\hat \phi$ and a non-genuine $\kappa_{\pm}$ and non-genuine $\rho^0$ (rather than to a non-genuine $\phi$ with genuine $\kappa_{\pm}$ and genuine $\rho^0$) as discussed in \S\ref{sub:alternativeinput}.

We define 
\begin{eqnarray*}
	K_{\pm} & =& G_1(F)_{x,\frac{r_1}{2}}G_2(F)_{x,\frac{r_2}{2}}\hdots G_n(F)_{x,\frac{r_n}{2}}K^0_{\pm,x}  \\
	K & =& G_1(F)_{x,\frac{r_1}{2}}G_2(F)_{x,\frac{r_2}{2}}\hdots G_n(F)_{x,\frac{r_n}{2}}K^0_{x} \\
	K^\dagger & =& G_1(F)_{x,\frac{r_1}{2}}G_2(F)_{x,\frac{r_2}{2}}\hdots G_n(F)_{x,\frac{r_n}{2}}G_{n+1}(F)_{x,0+}\\
	K^+ & =& G_1(F)_{x,\frac{r_1}{2}+}G_2(F)_{x,\frac{r_2}{2}+}\hdots G_n(F)_{x,\frac{r_n}{2}+}G_{n+1}(F)_{x,0+} 
\end{eqnarray*}

Now we can follow the construction in \S\ref{sub:construction} verbatim with this notation except for replacing $G^0_{\pm, x}$ by its subgroup $K^0_{\pm, x}$ to obtain a genuine representation $\kappa_{\pm}$ of $K_{\pm}$. 

\begin{thm} \label{thm-types}
	The representation $\rho^0 \otimes \kappa_{\pm}$ of $K_{\pm}$ factors through a representation of $K$, which we also denote by $\rho^0 \otimes \kappa_{\pm}=:\sigma_\Upsilon$. The representation $\sigma_\Upsilon$ only depends on the point $x$ and the genuine essentially depth-zero representation $\rho=\rho^0 \otimes_{i=1}^n \phi_i|_{G^0_x}$ of $K^0_{\pm, x}$ and not on the choice of factorization into depth-zero representation and characters in $\Upsilon$, so we may also denote it by $\sigma_\rho$.
	
	The pair $(K, \sigma_\rho)$ is a type for a finite union of Bernstein blocks. If $K^0_{\pm, x}= (G^0)_{\pm, x, 0}(M^0)^1_{\pm, x}$, then $(K, \sigma_\rho)$ is a type for a single Bernstein block.
	
	If $p$ does not divide the order of the absolute Weyl group of $G$, then each Bernstein block has a type of this form.
\end{thm}
\begin{proof}
	The claim that the representation $\rho^0 \otimes \kappa_{\pm}$ of $K_{\pm}$ factors through a representation of $K$ and that it only depends on the the point $x$ and the genuine essentially depth-zero representation $\rho=\rho^0 \otimes \otimes_{i=1}^n \phi_i|_{G^0_x}$ of $K^0_{\pm, x}$ and not the choice of $\Upsilon$ is proved the same way as in Theorem \ref{thm-screps}.
	
	The remaining claims follow from Proposition \ref{pro:comparisonKimYu} (and \cite{KimYu, Fin21Comp, Kim07, Fin21}), which compares $(K, \rho^0 \otimes \kappa_{\pm})$  with the types constructed by Kim and Yu.
\end{proof}

Let $\chi_{G^0}$ be the genuine character of $G^0(F)_{\pm}$ produced from the $\chi$-data for $R(Z_{G^0},G)$ as in Theorem \ref{thm:main1b}\eqref{item:organbox}, and let $\eps:=\prod_{i=1}^n \eps_x^{G_i/G_{i+1}}$ be the quadratic character as defined just before \cite[Theorem~4.1.13]{FKS23}, see also Theorem \ref{thm:main1b}\eqref{item:organbox}. As in \S\ref{sub:Yucomparison}, the representation $\rho^0 \otimes \chi_{G^0}^{-1}$ is non-genuine and we obtain the truly non-genuine $G$-datum 
\[\Upsilon_\chi=((G_i)_{1\leq i\leq n+1},x,(r_i)_{1\leq i\leq n},(K^0_x, \rho^0 \otimes \chi_{G^0}^{-1} \otimes \eps)),(\phi_i)_{1\leq i\leq n}) .\] 
To such a datum, Kim and Yu (\cite{KimYu}) attach the type $(K,  \rho^0 \otimes \chi_{G^0}^{-1} \otimes \eps \otimes \kappa^{\mathrm{nt}})$, where $\kappa^{\mathrm{nt}}$ is as in Definition \ref{dfn:kappant}. We write $\sigma_{\Upsilon_\chi}:= \rho^0 \otimes \chi_{G^0}^{-1} \otimes \eps \otimes  \kappa^{\mathrm{nt}}$.

\begin{pro}\label{pro:comparisonKimYu}
	We have $\sigma_\Upsilon = \sigma_{\Upsilon_\chi}$.
\end{pro}
\begin{proof}
	This follows from the relation $\tilde\jmath^*\mathrm{FWeilHeis} = \bar\jmath^*\mathrm{LWeilHeis} \otimes \chi_{G^0}^{-1} \otimes \epsilon$ from Theorem \ref{thm:main1b} and the construction of $\kappa_{\pm}$ and $\kappa^{\mathrm{nt}}$.
\end{proof}

We denote by $\mathrm{Rep}(G)^{(K, \sigma_\rho)}$ the subcategory of the category of a smooth complex representations of $G(F)$ described by the type $(K, \sigma_\rho)$, which is a finite union of Bernstein blocks. Likewise the pair $(K^0_{\pm, x}, \rho)$ defines a Serre subcategory  $\mathrm{Rep}(G^0_{\pm})^{(K^0_{\pm, x}, \rho)}$ of the category of smooth complex representations of $\mathrm{Rep}(G^0(F)_{\pm})$, which contains those irreducible representations whose restriction to $K^0_{\pm, x}$ contain $\rho$.

One can now rephrase the recent reduction-to-depth-zero results by one of us in their joint work with Adler, Mishra and Ohara (\cite{AFMO2}), and in their joint work with Dat (\cite{DF}), as follows. \begin{pro}
	Suppose that $p$ is not a torsion prime of $G$ and not a torsion prime of its Langlands dual group $\hat G$, or suppose that $K^0_x$ is either $G^0_{x,0}$ or $G^0_{x,0}(M^0_x)^1$. 
	Then we have an equivalence of categories 
	\[\mathrm{Rep}(G)^{(K, \sigma_\rho)} \simeq \mathrm{Rep}(G^0_{\pm})^{(K^0_{\pm,x}, \rho)} ,\] 
	and an isomorphism of (Hecke) algebras 
	\[\mathrm{End}_G(\cind_K^{G(F)} \sigma_\rho) \simeq \mathrm{End}_{G^0_{\pm}}(\cind_{K^0_{\pm, x}}^{G^0(F)_{\pm}} \rho) .\]
\end{pro}
\begin{proof}
	Tensoring with $(\chi_{G^0}\otimes \otimes_{i=1}^n \phi_i|_{G^0(F)})^{-1}$ and viewing a non-genuine representation of $G^0(F)_{\pm}$ as a representation of $G^0(F)$ provides an equivalence of categories 
	\begin{equation} \label{eq:reduction1}
		\mathrm{Rep}(G^0_{\pm})^{(K^0_{\pm, x}, \rho)} \simeq  \mathrm{Rep}(G^0)^{(K^0_{x}, \rho^0 \otimes \chi_{G^0}^{-1})} .
		\end{equation}
	Moreover, we have 
	\[\cind_{K^0_{\pm, x}}^{G^0(F)_\pm} \rho=\cind_{K^0_{\pm, x}}^{G^0(F)_\pm}(\rho^0 \otimes \otimes_{i=1}^n \phi_i|_{G^0(F)})=\cind_{K^0_{\pm, x}}^{G^0(F)_\pm}(\rho^0) \otimes \otimes_{i=1}^n \phi_i|_{G^0(F)}, \]
	 \[\cind_{K^0_{\pm,x}}^{G^0(F)_{\pm}} (\rho^0) \otimes \chi_{G^0}^{-1} = \cind_{K^0_{\pm,x}}^{G^0(F)_{\pm}} (\rho^0 \otimes \chi_{G^0}^{-1}) = \cind_{K^0_{x}}^{G^0(F)}(\rho^0 \otimes \chi_{G^0}^{-1}) ,\]
	  and hence we obtain
	\begin{equation}\label{eq:reduction2}
		\mathrm{End}_{G^0_{\pm}}(\cind_{K^0_{\pm, x}}^{G^0(F)_\pm} \rho) \simeq \mathrm{End}_{G^0}(\cind_{K^0_{x}}^{G^0(F)}(\rho^0 \otimes \chi_{G^0}^{-1})) .
		 \end{equation}
	On the other hand, by \cite[Theorem~4.3.11]{AFMO2} and \cite[Corollary~4.1.2]{DF} we have an isomorphism
	\begin{eqnarray*}
		\mathrm{End}_G(\cind_K^{G(F)} \sigma_\rho)&=&\mathrm{End}_G(\cind_K^{G(F)}(\rho^0 \otimes \chi_{G^0}^{-1} \otimes \eps \otimes \kappa^{\mathrm{nt}}))\\
		&\simeq&\mathrm{End}_{G^0}(\cind_{K^0_{x}}^{G^0(F)} (\rho^0 \otimes \chi_{G^0}^{-1}) ),
	\end{eqnarray*}
	from which one obtains (see also \cite[Corallary~4.5.1]{AFMO2}) 
	\[\mathrm{Rep}(G)^{(K, \sigma_\rho)} \simeq \mathrm{Rep}(G^0)^{(K^0_{x}, \rho^0 \otimes \chi_{G^0}^{-1})} .\] 
	Combining the last two equations with \eqref{eq:reduction1} and \eqref{eq:reduction2} provides the desired results.
\end{proof}

\begin{appendices}
	 
	\section{The Fourier--Weil representation} \label{app:fw}

\subsection{Review of the Witt groups} \label{sub:Wittgroup}

We review some material from \cite{EKM2008} and \cite[Appendix 2]{MH73}. Let $k$ be a field of characteristic not $2$. The quadratic Witt group $\QW(k)$\index[notation]{QW@$\QW$} is the Grothendieck group on the semi-group whose underlying set is the set of equivalence classes of non-singular quadratic spaces over $k$, i.e., pairs $(V,q)$ consisting of a finite-dimensional $k$-vector space $V$ and a quadratic form $q : V \to k$ whose polar form $b_q(x,y)=q(x+y)-q(x)-q(y)$ is non-degenerate, and where equivalence is defined as $q_1 \sim q_2$ if and only if $q_1 \perp H^{\perp r}$ is isometric to $q_2 \perp H^{\perp s}$ for some natural numbers $r,s$, where $H : k \oplus k \to k$ is the hyperbolic form $H(x,y)=xy$. The semi-group operation is given by orthogonal direct sum.

Similarly, the bilinear Witt
ring $\BW(k)$\index[notation]{BW@$\BW$} is the Grothendieck ring associated to the semi-ring whose underlying set is the set of equivalence classes of non-singular symmetric bilinear spaces over $k$, i.e., pairs $(V,b)$ consisting of a finite-dimensional $k$-vector space $V$ and a non-degenerate symmetric bilinear form $b : V \times V \to k$, and where equivalence is defined as $b_1 \sim b_2$ if and only if $b_1 \perp M_1$ is isometric to $b_2 \perp M_2$ for metabolic forms $M_1,M_2$. Recall that a metabolic form is a form on a vector space of even dimension $2n$ that possesses a totally isotropic subspace of dimension $n$. The semi-ring operations are given by orthogonal direct sum and tensor product. The group $\QW(k)$ is a module over $\BW(k)$ by tensor product. 

Finally, if we are given a subfield $k_0 \subset k$ such that $k/k_0$ has degree $2$, we can define the Hermitian Witt
ring $\HW(k)$\index[notation]{HW@$\HW$} in exact analogy with $\BW(k)$. 

The assignments $q_b(v)=b(v,v)$ and $b_q(v,w)=q(v+w)-q(v)-q(w)$ give group isomorphisms
\[ \BW(k) \simto \QW(k)\qquad \tx{and} \qquad \QW(k) \simto \BW(k),\]
whose compositions in either order are rescaling the form by multiplication by $2$. In the Hermitian setting of $k/k_0$ we have the homomorphism\index[notation]{tr}
\[ \tx{tr} : \HW(k) \to \QW(k_0) \]
that sends $[(V,h)]$ to $[(V|_{k_0},q_h)]$, with $q_h(v)=h(v,v)$, which we call the \emph{trace quadratic space} of $(V,h)$ and the quadratic form the \emph{trace quadratic form}.\index[terminology]{trace quadratic form} Note that $b_{q_h}=\tx{tr}_{k/k_0}\circ h$.

If $\ell/k$ is a finite separable field extension we have the group homomorphisms
\[ \tx{tr} : \BW(\ell) \to \BW(k),\qquad \tx{tr} : \QW(\ell) \to \QW(k), \]
which send $(V,b)$, resp., $(V,q)$, over $\ell$ to $(V|_k,\tx{tr}_{\ell/k}\circ b)$, resp., $(V|_k,\tx{tr}_{\ell/k}\circ q)$, which we call the \emph{trace bilinear space, resp., trace quadratic space}. We also have the homomorphisms
\[ \tx{bc} : \BW(k) \to \BW(\ell),\qquad \tx{bc} : \QW(k) \to \QW(\ell) \]
given by extension of scalars. We have the exact sequence of $\QW(k_0)$-modules
\[ 0 \to \HW(k) \xrightarrow{\mathrm{tr}} \QW(k_0) \xrightarrow{\mathrm{bc}} \QW(k). \]

For $a \in k^\times$ we denote by $[a]_Q \in \QW(k)$\index[notation]{aQ@$[a]_Q$} the class of the quadratic form $k \to k$, $x \mapsto ax^2$, and by $[a]_B \in \BW(k)$\index[notation]{aB@$[a]_B$} the class of the symmetric bilinear form $k \times k \to k$, $(x,y) \mapsto axy$. Then $[1]_B$ is the unit of the ring $\BW(k)$. For $a \in k_0^\times$ we denote by $[a]_H \in \HW(k)$\index[notation]{aH@$[a]_H$} the Hermitian form $k \times k \to k$, $(x,y) \mapsto ax\sigma(y)$, where $\sigma$ is the non-trivial automorphism of $k/k_0$. Then $[1]_H$ is the unit of the ring $\HW(k)$.

Let $(V,q)$ be a singular quadratic space, and let $N \subseteq V$ be the radical of $b_q$. Since $q=q_{b_q}/2$, the form $q$ descends to a non-singular form $\bar q$ on $\ol{V}:=V/N$. In this way, we obtain a class $[(V,q)] \in \QW(k)$. In the same way, we obtain a class $[(V,b)] \in \BW(k)$ from a degenerate symmetric bilinear space.

The inverse of $[(V,q)]$ is $[(V,-q)]$ and the additive inverse of $[(V,b)]$ is $[(V,-b)]$.

\begin{lem}[Witt reduction] \label{lem:wittred}
	Let $(V,q)$ be a (possibly degenerate) quadratic space and $W \subseteq V$ isotropic for $b_q$. Then $q$ induces a quadratic form on $W^\perp/W$, and we have $[(V,q)]=[(W^\perp/W,q)]$ in $\QW(k)$.
\end{lem}
\begin{proof}
	If $(V,q)$ is non-singular this is \cite[Lemma 8.10]{EKM2008}. In general pass first to $(\ol{V},\bar q)$.
\end{proof}

Let $I_k \subset \BW(k)$ be the subset generated by even-dimensional forms. Then $I_k$ is an ideal in $\BW(k)$. We have the following group isomorphisms ([\cite{EKM2008}])
\begin{enumerate}
	\item $\dim : \BW(k)/I_k \simto \Z/2\Z$,
	\item $s\det : I_k/I_k^2 \simto k^\times/k^{\times,2}$ sending $[(V,b)]$ to $(-1)^{n(n-1)/2}\det(b)$ with $n=\dim(V)$, \item $I_k^2/I_k^3 \simto H^2(\Gamma,\mu_2)$ by the Clifford invariant, see \cite[\S14, Fact~16.2, and Chapter~VIII]{EKM2008}.
\end{enumerate}

\subsection{Review of the Maslov index}

Let $k$ be an arbitrary field of characteristic not $2$, and let $(V,\eta)$ be a symplectic $k$-vector space.

The following useful description of Lagrangian subspaces in a polarized symplectic vector space is well-known, but we have not found a convenient reference in the given generality.

\begin{lem} \label{lem:lagpol}
	Assume $V=X \oplus X^*$ with $\eta((x,\alpha),(y,\beta))=\beta(x)-\alpha(y)$. 
	\begin{enumerate}
		\item There is a bijection between the set of Lagrangian subspaces $L \subset V$ and the set of pairs $(E,b)$, where $E \subset X$ is a subspace and $b$ is a symmetric bilinear form on $E$, given as follows. To a pair $(E,b)$ assign $L=\{(x,\alpha) \, | \, x \in E, \alpha \in X^*, \alpha|_E=b(x,-)\}$. To a Lagrangian subspace $L \subset V$ assign $E=p_X(L)$, where $p_X : V \to X$ is the projection with kernel $X^*$, and $b(x,y)=\alpha(y)$ for any $\alpha \in X^*$ s.t.\ $(x,\alpha) \in L$.
		\item If $L$ corresponds to $(E,b)$ under the above bijection, then
		\begin{enumerate}
			\item $L \cap X=\mathrm{rad}(b) \subset E$,
			\item $L \cap X^*=E^\perp$,
			\item $L$ is transversal to $X$ if and only if  $b$ is non-degenerate,
			\item $L$ is transversal to $X^*$ if and only if $E=X$. In that case $b$ induces a linear map $X \to X^*$ and $L$ is its graph.
		\end{enumerate}
	\end{enumerate}
	
\end{lem}
\begin{proof} 
	(1) We first verify that the maps are well-defined. Given a Lagrangian subspace $L$ set $E=p_X(L)$. Given $x,y \in E$ let $\alpha_1,\alpha_2 \in X^*$ be such that $(x,\alpha_1),(x,\alpha_2) \in L$. Then $L \ni (x,\alpha_1) - (x,\alpha_2) = (0,\alpha_1-\alpha_2)$. Choosing $\beta \in X^*$ such that $(y,\beta) \in L$ we see $0=\eta((y,\beta),(0,\alpha_1-\alpha_2))=(\alpha_1-\alpha_2)(y)$, hence $\alpha_1(y)=\alpha_2(y)$, showing that $b$ is well-defined. To see that $b$ is symmetric, we let $(x,\alpha),(y,\beta) \in L$ and check
	\[ 0=\eta((x,\alpha),(y,\beta))=\beta(x)-\alpha(y) = b(y,x)-b(x,y). \]

	Conversely, given a pair $(E,b)$, we define $L$ as in the statement and note that the last displayed computation converts the symmetry of $b$ into the isotropy of $L$. Furthermore, we have the Cartesian square
	\[ \xymatrix 	{
		L \ar[r]^{p_X}\ar[d]_{p_{X^*}}&E\ar[d]^-b\\
		X^*\ar[r]&E^* }
	\] 
	which gives the exact sequences 
	\begin{equation} \label{eq:lagseq}
		0 \to E^\perp \to L \to E \to 0\quad\text{and}\quad 0 \to \mathrm{rad}(b) \to L \to X^* ,
 	\end{equation}
	where $E^\perp \subset X^*$ is the annihilator of $E \subset X$. From the first sequence we infer $\dim(L)=\dim(E)+\dim(E^\perp)=\dim(X)=\dim(V)/2$. Therefore $L$ is an isotropic subspace of a maximal dimension, hence a Lagrangian subspace.

	Next we check that the maps are mutually inverse. Starting with $L$ define $(E,b)$ and let $L'$ be associated to $(E,b)$. It is immediate that $L \subset L'$, and since both are Lagrangian, we conclude $L=L'$. Conversely, starting with $(E,b)$ let $L$ be associated to it and let $(E',b')$ be associated to $L$. Since $L \to E$ is surjective by the above Cartesian diagram, we see $E'=E$. Moveover, for $(x,\alpha) \in L$ we have $b'(x,-)=\alpha|_E$, but $\alpha|_E=b(x,-)$, so $b'(x,-)=b(x,-)$ for all $x \in E$, hence $b'=b$.

	(2) Using the first sequence in \eqref{eq:lagseq} we see $L \cap (\{0\} \oplus X^*)=\tx{ker}(p_X|_L)=E^\perp$, which proves (b) and (d) follows. Using the second sequence in \eqref{eq:lagseq} we see $L \cap (X \oplus \{0\})=\tx{ker}(p_{X^*}|_L)=\tx{ker}(b : E \to E^*)=\tx{rad}(b)$, which proves (a), and (c) follows.
\end{proof}

\begin{dfn}
	Let $(L_1,L_2,L_3)$ be a triple of Lagrangian subspaces of $V$. The \emph{Kashiwara--Maslov index} $\tau[(V,\eta)](L_1,L_2,L_3) \in \QW(k)$ is the class of the (possibly degenerate) quadratic space
		\[ q_{L_1,L_2,L_3} : L_1 \oplus L_2 \oplus L_3 \to k,\quad (\ell_1,\ell_2,\ell_3) \mapsto \eta(\ell_1,\ell_2)+\eta(\ell_2,\ell_3)+\eta(\ell_3,\ell_1). \]
		We may also simply write $\tau(L_1,L_2,L_3)$ for $\tau[(V,\eta)](L_1,L_2,L_3)$.
\end{dfn}

\begin{lem} \label{lem:maslov_basic1}
	Let $L_1,L_2,L_3$ be three Lagrangian subspaces of $V$ and consider $W=L_2 \cap (L_1+L_3)$.
	\begin{enumerate}
		\item The bilinear form $b : W \times W \to k$ defined by $b(x,y)=\eta(x_1,y_3)$ for $x=x_1+x_3$ and $y=y_1+y_3$ with $x_1, y_1 \in L_1$ and $x_3, y_3 \in L_3$ is well-defined and symmetric and its radical is $(L_1 \cap L_2) + (L_3 \cap L_2)$.
		\item The quadratic form $q'(x)=b(x,x)$ on $W$ represents the Witt class $\tau(L_1,L_2,L_3)$.
	\end{enumerate}
	
\end{lem}
\begin{proof}
	(1) That the form is well-defined follows from $L_1$ and $L_3$ being isotropic. Using that each $L_i$ is Lagrangian and that $x_1+x_3,y_1+y_3 \in L_2$ we have $\eta(x_1,y_3)=\eta(x_1+x_3,y_3)=\eta(x_1+x_3,-y_1)=\eta(y_1,x_1+x_3)=\eta(y_1,x_3)$, showing the symmetry of $b$. To compute its radical, note that the set of $y_3 \in L_3$ for which there exists $y_1 \in L_1$ such that $y_1+y_3 \in L_2 \cap (L_1+L_3)$ is precisely $L_3 \cap (L_1+L_2)$. Therefore, the radical of $b$ is precisely 
	\begin{eqnarray*}
	(L_2 \cap (L_1+L_3)) \cap (L_3 \cap (L_1+L_2))^\perp&=& (L_2 \cap (L_1+L_3)) \cap (L_3 + (L_1 \cap L_2))\\
	&=&(L_1 \cap L_2) + (L_3 \cap L_2).	
	\end{eqnarray*}

	(2) We compute the polar form of $q$ and obtain
	\begin{eqnarray*}
	b_q((x_1,x_2,x_3),(y_1,y_2,y_3))&=&(\eta(x_1,y_2)+\eta(y_1,x_2))\\
	&+&(\eta(x_2,y_3)+\eta(y_2,x_3))\\
	&+&(\eta(x_3,y_1)+\eta(y_3,x_1)).\end{eqnarray*}
	This shows that $L_1 \oplus \{0\} \oplus \{0\}$ is an isotropic subspace and 
	\[ L_1^\perp=\{(\ell_1,\ell_2,\ell_3)|\ell_2-\ell_3 \in L_1\}. \]
	Lemma \ref{lem:wittred} shows that $\tau(L_1,L_2,L_3)$ is represented by $L_1^\perp/L_1$ with quadratic form induced by $q$. The projection $L_1 \oplus L_2 \oplus L_3 \to L_2$ induces an surjective map $L_1^\perp \to L_2 \cap (L_1+L_3)$ whose kernel is $N:=L_1 \oplus\{0\} \oplus\{0\} + \{0\} \oplus \{0\} \oplus (L_1 \cap L_3)$. We note that $b_q(N,L_1^\perp)=0$, therefore $q$ descends via this surjective map to a form on $L_2 \cap (L_1 + L_3)$, representing the same Witt class. A preimage of $\ell_1+\ell_3\in L_2 \cap (L_1+L_3)$ in $L_1^\perp$ is given by $(\ell_1,\ell_1+\ell_3,\ell_3)$, so the descended quadratic form sends $\ell_1+\ell_3$ to
	\begin{eqnarray*}
		q(\ell_1,\ell_1+\ell_3,\ell_3)
		&=&\eta(\ell_1,\ell_1+\ell_3)+\eta(\ell_1+\ell_3,\ell_3)+\eta(\ell_3,\ell_1)\\
		&=&\eta(\ell_1,\ell_3)+\eta(\ell_1,\ell_3)+\eta(\ell_3,\ell_1)\\
		&=&q'(\ell_1+\ell_3),
	\end{eqnarray*}
	as claimed. 
\end{proof}

\begin{lem} \label{lem:maslov_basic2}
	Fix a polarization $V=X \oplus Y$. Let $L_1,L_2,L_3 \subset V$ be Lagrangian subspaces with associated pairs $(E_i,\beta_i)$ as in Lemma \ref{lem:lagpol}. Let $q_i(x)=\beta_i(x,x)$. 
	\[ \tau(L_1,L_2,L_3) = [q_2-q_1|_{E_2 \cap E_1}] + [q_3-q_2|_{E_3 \cap E_2}] + [q_1-q_3|_{E_1 \cap E_3}]. \]
\end{lem}
\begin{proof}
	Let again $\beta_q$ be the polar form of $q_{L_1,L_2,L_3}$. One checks that $N:=(L_1 \cap Y)\oplus (L_2 \cap Y) \oplus (L_3 \cap Y) \subset L_1 \oplus L_2 \oplus L_3$ is totally isotropic for $\beta_q$. Lemma \ref{lem:wittred} shows that $\tau(L_1,L_2,L_3)$ is represented by the quadratic form on $N^\perp/N$ induced by $q$. The projection $V \to X$ along $Y$ restricts to a surjective map $p_i : L_i \to E_i$. The map $(p_1,p_2,p_3)$ induces an isomorphism
	\[ N^\perp/N \to E = \{(x_1,x_2,x_3) \in E_1 \oplus E_2 \oplus E_3|x_3-x_2 \in E_1, x_1-x_3 \in E_2, x_2-x_1 \in E_3\}. \]
	For any $(x_i,y_i) \in L_i$ one has
	\[ q_{L_1,L_2,L_3}((x_1,y_1),(x_2,y_2),(x_3,y_3)) = \eta(x_3-x_2,y_1)+\eta(x_1-x_3,y_2)+\eta(x_2-x_1,y_3)\]
	and from this one infers that the form on $N^\perp/N$ induced by $q$ transports to $E$ to the form 
	\[ (x_1,x_2,x_3) \mapsto \beta_1(x_1,x_3-x_2)+\beta_2(x_2,x_1-x_3)+\beta_3(x_3,x_2-x_1). \]
	We have the further map
	\[  (E_1 \cap E_2) \oplus (E_2 \cap E_3) \oplus (E_3 \cap E_1) \to E,\qquad (u,v,w) \mapsto (u+w,u+v,v+w), \]
	which is an isomorphism since the characteristic of $k$ is not $2$. Pulling the above quadratic form under this isomorphism produces $(q_2-q_1) \oplus (q_3-q_2) \oplus (q_1-q_3)$.
\end{proof}

\begin{pro} \label{pro:maslov_basic}
	Let $L_1,L_2,L_3 \subset V$ be Lagrangian subspaces.
	\begin{enumerate}
		\item For $\sigma \in S_3$ we have $\tau(L_{\sigma(1)},L_{\sigma(2)},L_{\sigma(3)})=\epsilon(\sigma)\tau(L_1,L_2,L_3)$, where $\epsilon(\sigma)$ is the sign of the permutation $\sigma$ interpreted using $[-1]_B \in \BW(k)$.
		\item For $g \in \Sp(V)(k)$ we have $\tau(gL_1,gL_2,gL_3)=\tau(L_1,L_2,L_3)$.
		\item The function $\tau$ is a cocycle, i.e. given a fourth Lagrangian $L_4 \subset V$ we have
		\[ \tau(L_2,L_3,L_4)-\tau(L_1,L_3,L_4)+\tau(L_1,L_2,L_4)-\tau(L_1,L_2,L_3)=0.\]
		\item If any two of $L_1,L_2,L_3$ are equal, then $\tau(L_1,L_2,L_3)=0$.
	\end{enumerate}
\end{pro}
\begin{proof}
	(1) Consider the swapping of $L_1,L_2$. The isomorphism
	\[ L_1 \oplus L_2 \oplus L_3 \to L_2 \oplus L_1 \oplus L_3,\qquad (\ell_1,\ell_2,\ell_3) \mapsto (\ell_2,\ell_1,\ell_3) \]
	identifies the form $q_{L_1,L_2,L_3}$ with the form $-q_{L_2,L_1,L_3}$. The analogous holds when $L_2,L_3$ are swapped. Since these two transpositions generate $S_3$ the claim is proved.

	(2) This is immediate from the definition of $\tau$.

	(3,4) These follows from Lemma \ref{lem:maslov_basic2}. Note that (4) does not formally follow from (1) or (3) ince $\QW(k)$ often has $2$-torsion.	
\end{proof}

\begin{lem}
	Let $U \subset V$ be an isotropic space, i.e. $U \subset U^\perp$. For any subspace $W \subset V$ define $W^U=(W \cap U^\perp)+U$. 
	\begin{enumerate}
		\item We have $(W^U)^\perp = (W^\perp)^U$.
		\item If $L$ is a Lagrangian subspace of $V$, then so is $L^U$, and $L^U/U$ is a Lagrangian subspace of $U^\perp/U$.
	\end{enumerate}
	
\end{lem}
\begin{proof}
	(1) We have 
	\[ ((W\cap U^\perp)+U)^\perp=(W \cap U^\perp)^\perp \cap U^\perp=(W^\perp+U)\cap U^\perp = (W^\perp \cap U^\perp) + U. \]
	(2) That $L^U$ is a Lagrangian subspace follows from (1). It is furthermore contained in $U^\perp$ and its image modulo $U$ is a totally isotropic subspace of dimension $\dim(V)/2-\dim(U)=\dim(U^\perp/U)/2$, which is the maximal possible dimension in $U^\perp/U$.
\end{proof}

\begin{pro} Let $L_1,L_2,L_3$ be Lagrangian subspaces of $V$. Any $U \subset L_1 \cap L_2 + L_2 \cap L_3 + L_3 \cap L_1$ is isotropic and we have 
	\begin{eqnarray*}
		\tau[(V,\eta)](L_1,L_2,L_3)&=&\tau[(V,\eta)](L_1^U,L_2^U,L_3^U)\\
		&=&\tau[(U^\perp/U,\eta_U)](L_1^U/U,L_2^U/U,L_3^U/U).
	\end{eqnarray*}
\end{pro}
\begin{proof}
	We have $L_1^U=(L_1 \cap L_1^U)+(L_2 \cap L_1^U)=(L_1 \cap L_1^U) + (L_2^U \cap L_1^U)$. From Lemma \ref{lem:maslov_basic1} we conclude $\tau(L_1,L_1^U,L_2)=0=\tau(L_1,L_1^U,L_2^U)$. Applying Proposition \ref{pro:maslov_basic}(3) to combinations of the six Lagrangian subspaces $L_1,L_2,L_3,L_1^U,L_2^U,L_3^U$ we conclude $\tau[(V,\eta)](L_1,L_2,L_3)=\tau[(V,\eta)](L_1^U,L_2^U,L_3^U)$. The second identity follows from Lemma \ref{lem:maslov_basic1}, since now $U \subset (L_1^U \cap L_2^U) + (L_3^U \cap L_2^U)$.
\end{proof}

\begin{dfn}
	For $n>3$ and an $n$-tuple $(L_1,\dots,L_n)$ of Lagrangian subspaces of $V$ define the \emph{Kashiwara--Maslov} index recursively on $n$ by
	\[ \tau[(V,\eta)](L_1,\dots,L_n) :=  \tau[(V,\eta)](L_1,L_2,L_3)+\tau[(V,\eta)](L_1,L_3,L_4,\dots,L_n). \]
\end{dfn}

\subsection{The Witt-metaplectic groups} \label{sub:wittmeta}

We continue with an arbitrary field of characteristic not $2$ and a symplectic $k$-vector space $(V,\eta)$. Let $H=\mb{G}_a \boxtimes_{\eta/2} V$ be the Heisenberg group.

\begin{lem} \label{maslov_cocycle}
	Let $L \subset V$ be a Lagrangian subspace. Then 
	\[ z_L(g,h):=\tau[(V,\eta)](L,gL,ghL) \]
	defines an element $z_L \in Z^2(\Sp(V),\QW(k))$.
\end{lem}
\begin{proof}
	This follows at once from Proposition \ref{pro:maslov_basic}(3), since the homogeneous form of $z_L$ sends $(g,h,k)$ to $\tau(gL,hL,kL)$.
\end{proof}

\begin{dfn}
	Fix a Lagrangian subspace $L \subset V$. 
The \emph{Witt-metaplectic group relative to $L$} is the group\index[notation]{WML@$\WM_L$} 
		\[ \WM_L(V) := \QW(k)/I_k^3 \boxtimes_{z_L} \Sp(V)(k). \]

\end{dfn}

\begin{lem}
	\begin{enumerate}
		\item Let $g \in \Sp(V)(k)$. The graph $\Gamma_g \subset V \oplus V$ of the endomorphism $g : V \to V$ is a Lagrangian subspace with respect to the symplectic form $\eta \oplus (-\eta)$.
		\item The function $z(g,h):=\tau[(V\oplus V,(-\eta) \oplus \eta)](\Gamma_1,\Gamma_g,\Gamma_{gh})$ is an element of $Z^2(\Sp(V)(k),\QW(k))$.
	\end{enumerate}
\end{lem}
\begin{proof}
	(1) is an immediate computation. For (2), letting $\Sp(V,\eta)$ act on $V \oplus V$ by the identity on the first factor and the natural action on the second, we obtain an embedding $\Sp(V,\eta) \to \Sp(V \oplus V,(-\eta) \oplus \eta)$. Under this embedding we have for $g \in \Sp(V,\eta)$ the identity $\Gamma_g=g\Gamma_1$ and see that $z(g,h)$ is the pull-back under this embedding of the Maslov cocycle for the Lagrangian subspace $\Gamma_1 \subset V \oplus V$. 
\end{proof}

\begin{dfn} \label{dfn:meta_univ} The \emph{universal Witt-metaplectic group} is $\WMuniv(V) := \QW(k)/I_k^3 \boxtimes_z \Sp(V)(k)$.\index[notation]{WM@$\WMuniv$}
\end{dfn}

\subsection{Orientations and a reduction of the Maslov index}

We continue with an arbitrary field $k$ of characteristic not $2$. In this subsection we recall orientations and prove a formula for the Masolv index in terms of them, Proposition \ref{pro:maslov_break}. This result appears in \cite[Proposition 2.1]{Parimala2000}, but our proof is entirely different and less computational. See also \cite[Theorem 1.7.6]{LionVergne} and \cite[Theorem A.15]{LionVergne} for conceptual proofs of weaker results.

While orientations in \cite{LionVergne}, \cite{Parimala2000}, and \cite{Thomas2013}, are taken is taken up to squares, but we'll set things up more strictly in order to allow geometric usage of the term as well.

\begin{dfn}
	\begin{enumerate}
		\item An \emph{orientation} on a $d$-dimensional $k$-vector space $W$ is a non-zero element $o \in \Wedge^dW$.
		\item An \emph{oriented vector space} is a pair $\tilde W=(W,o)$ of a finite-dimensional $k$-vector space $W$ and an orientation $o$ on $W$.
		\item The \emph{dual} of an orientation $o$ on $W$ is the orientation $o^*$ of $W^*$ such that $\<o,o^*\>=1$, where $\<-,-\> : \Wedge^d(W) \otimes \Wedge^d(W^*) \to k$ is the pairing defined by $\<w_1\wedge \dots \wedge w_d,w'_1 \wedge \dots \wedge w'_d\>=\det((\<w_i,w'_j\>)_{i,j})$. Given $\tilde W=(W,o)$ define $\tilde W^*=(W^*,o^*)$.
		\item Given two oriented vector spaces $\tilde W_i$ of the same dimension, an isomorphism $\tilde f : \tilde W_1 \to \tilde W_2$ is an isomorphism $f : W_1 \to W_2$, and $\det(\tilde f):=\Wedge^d f(o_1)/o_2 \in k^\times$.
	\end{enumerate}
\end{dfn}

\begin{fct}
	\begin{enumerate}
		\item For an automorphism $f : \tilde W \to\tilde W$ we have $\det(\tilde f)=\det(f)$.
		\item For automorphisms $\tilde f : \tilde W_1 \to \tilde W_2$ and $\tilde g : \tilde W_2 \to \tilde W_3$ we have $\det(\tilde g\circ \tilde f)=\det(\tilde f)\cdot \det(\tilde g)$.
	\end{enumerate}
\end{fct}

Let $(V,\eta)$ be a symplectic $k$-vector space. Given two Lagrangian subspaces $L_1,L_2 \subset V$ that are complementary consider the isomorphism $g_{L_2,L_1} : L_1 \to L_2^*$ defined by $g_{L_2,L_1}(l_1)(l_2)=\eta(l_1,l_2)$. If $L_1,L_2$ are not assumed complementary, the restriction of $\eta$ to $L_1+L_2$ factors through a symplectic form on $U=(L_1+L_2)/(L_1 \cap L_2)$ and the images $\bar L_1,\bar L_2 \subset U$ are complementary Lagrangian subspaces in $U$, hence we have the isomorphism $g_{L_2,L_1} : \bar L_1 \to \bar L_2^*$.

\begin{dfn}
	\begin{enumerate}
		\item An \emph{oriented Lagrangian subspace} is a pair $(L,o)$ consisting of a Lagrangian subspace $L \subset V$ and an orientation $o$ on $L$.
		\item Given oriented Lagrangian subspaces $\tilde L_1=(L_1,o_1)$ and $\tilde L_2=(L_2,o_2)$ of $V$ we choose arbitrarily an orientation $o_{12}$ of $L_1 \cap L_2$ and endow $\bar L_i=L_i/(L_1 \cap L_2)$ with the orientation $\bar o_i := o_i/o_{12}$. Define 
\[ m(\tilde L_1,\tilde L_2) \in \QW(k)/I_k^2 \] 
in the generic case $L_1 \neq L_2$ to be the class of quadratic forms on $L_1/(L_1 \cap L_2)$ whose determinant equals that of $\tilde g_{L_2,L_1}$; in the degenerate case $L_1=L_2$ there are no non-trivial forms on $L_1/(L_1 \cap L_2)$, but we can cover this case by formally letting $n_{12}$ be any positive integer congruent modulo $4$ to $\dim(L_1/(L_1 \cap L_2))$ and taking the class of forms with signed determinant $(-1)^{n_{12}(n_{12}-1)/2}\det(\tilde g_{L_2,L_1}) \in k^\times/k^{\times,2}$.	
	\end{enumerate}

\end{dfn}

\begin{rem}
		The class $m(\tilde L_1,\tilde L_2)$ does not depend on the choice of $o_{12}$: while $\det(\tilde g_{L_2,L_1}) \in k^\times$ does depend on $o_{12}$, its square class, i.e. its image in $k^\times/k^{\times,2}$, does not.
\end{rem}

\begin{pro} \label{pro:maslov_break}
	Let $(V,\eta)$ be a symplectic $k$-vector space. For three oriented Lagrangian subspaces $\tilde L_1,\tilde L_2,\tilde L_3 \subset V$. In $W(k)/I_k^2$ we have the identity
	\[ \tau(L_1,L_2,L_3) = m(\tilde L_1,\tilde L_2) + m(\tilde L_2,\tilde L_3) + m(\tilde L_3,\tilde L_1). \] 
\end{pro}
\begin{proof}
	\ul{Step 1:} The right hand side is independent of the choices of orientations. Indeed, changing multiplying $o_i$ by $a \in k^\times$ multiplies the determinants of the two terms involving $\tilde L_i$ by $a$. To fix notation let us assume $i=2$. In $\QW(k)/I_k^2$ we have $m(\tilde L_1,\tilde L_2)=[D_{12}]+(n_{12}-1)[1]$ where $D_{12}=\det(\tilde g_{L_2,L_1})$. Multiplying $o_2$ by $a$ replaced $D_{12}$ by $aD_{12}$ and therefore adds $[aD_{12}]-[D_{12}]$ to $m(\tilde L_1,\tilde L_2)$. Now $[aD_{12}]-[D_{12}] = [aD_{12},-D_{12}] \in I_k$ and modulo $I_k^2$ this element only depends on its signed determinant, which equals $aD_{12}^2 \equiv a$ and is thus independent of $D_{12}$. Therefore we see that $m(\tilde L_1,\tilde L_2)$ and $m(\tilde L_2,\tilde L_3)$ modify the same way, each being added the element $[a]-[1]$. Their sum therefore changes by adding $2([a]-[1])$. But the abelian group $I_k/I_k^2$ has exponent $2$, so $2([a]-[1])=0$.

	\ul{Step 2:} Both sides are compatible with orthogonal direct sums. This is evident.

	\ul{Step 3:} There exists an orthogonal decomposition $V=\bigoplus_i V_i$ such that $\dim(V_i)=2$ and $\dim(V_i \cap L_j)=1$ for all $i,j$. To see this, we induce on the dimension of $V$. If $\dim(V)=2$ there is nothing to show. For $\dim(V)>2$ we consider three cases. 
	
	If $L_1 \cap L_2 \cap L_3 \neq \{0\}$ choose $0 \neq e \in L_1 \cap L_2 \cap L_3$ and $f \in V$ such that $\eta(e,f) = 1$ and set $V_1=ke \oplus kf$ and $V'=V_1^\perp$. Then $V=V_1 \oplus V'$. Any $x \in L_i$ has the form $ae+bf+v'$ with $a,b \in k$ and $v' \in V'$. Since $x,e \in L_i$ we have $0=\eta(x,e)=-b$, showing that $L_i=ke \oplus L_i'$ with $L_i' = L_i \cap V'$. We can now apply inductive hypothesis to the three Lagrangian subspaces $L_1',L_2',L_3'$ of $V'$.

	If $L_1 \cap L_2 \cap L_3 = \{0\}$ but, up to reindexing, $L_1 \cap L_2 \neq \{0\}$, choose $0 \neq e \in L_1 \cap L_2$. We claim there exists $f \in L_3$ such that $\eta(e,f) = 1$. If not, then $\eta(e,L_3)=0$, hence $e \in L_3^\perp = L_3$, contradicting $L_1 \cap L_2 \cap L_3=\{0\}$. As in the previous paragraph, setting $V_1=ke \oplus kf$, $V'=V_1^\perp$, and $L_i' = V' \cap L_i$, we have $V=V_1 \oplus V'$, $L_1=ke \oplus L_1'$, $L_2=ke \oplus L_2'$, and $L_3=kf \oplus L_3'$. Apply inductive hypothesis to to the three Lagrangian subspaces $L_1',L_2',L_3'$ of $V'$.

	If $L_1$, $L_2$, $L_3$ are pairwise transverse then $V=L_1 \oplus L_2$ and there exists an isomorphism $T : L_1 \to L_2$ whose graph is $L_3$. The bilinear form $\beta(x,y)=\eta(x,Ty)$ on $L_1$ is then non-degenerate and symmetric. Let $e_1,\dots,e_n \in L_1$ be an orthogonal basis for $\beta$. Set $f_i=T(e_i) \in L_2$ and $V_i=ke_i \oplus kf_i$. Then $V=V_1 \oplus \dots \oplus V_n$, $L_1 \cap V_i=ke_i$, $L_2 \cap V_i=kf_i$, and $L_3 \cap V_i=k(e_i+f_i)$.

	\ul{Step 4:} Completion of the argument. We may assume $\dim(V)=2$ by Steps 2+3. Then $\dim(L_i)=1$ and there are three possibilities. If $L_1 \cap L_2 \cap L_3 \neq \{0\}$ then $L_1=L_2=L_3$ and choosing the same orientation on all $L_i$ we see that both sides of the identity vanish. If $L_1 \cap L_2 \cap L_3 = \{0\}$ but, up to reindexing, $L_1 \cap L_2 \neq \{0\}$, then $L_1=L_2 \neq L_3$ and we can find $e,f \in V$ with $\eta(e,f)=1$ such that $L_1=L_2=ke$ and $L_3=kf$. By Lemma \ref{lem:maslov_basic1}(1) the left-hand side of the identity vanishes, while using $e$ and $f$ as orientations we see that the right hand side is $[0]+[1]+[-1]=0$. Finally, if $L_1,L_2,L_3$ are pairwise transverse, then choosing $0 \neq e \in L_1$ and $0 \neq f \in L_2$ with $\eta(e,f)=1$ there exists $0 \neq a \in k$ such that $e+af \in L_3$. Using these vectors as orientations we compute the right-hand side as $[1]+[-1]+[-a]=[1,-1,-a]=[-a]$, while the left-hand side is the quadratic form $q(xe,yf,z(e+af))=xy-yz-axz=(x-z)(y-az)-az^2$. The subspace $L_1 \oplus L_2 \oplus \{0\}$ is a hyperbolic plane, while its orthocomplement is the line spanned by $(e,af,(e+af))$ and is isometric to the quadratic form $[-a]$.
\end{proof}

\subsection{Graded bilinear forms} \label{sub:gradedbf}

We continue with arbitrary field $k$, which we allow to have characteristic $2$ until Lemma \ref{lem:grbetasym}.
In this subsection we will discuss some basic linear algebra that will help us understand one aspect of the relationship between $\mathrm{FWeilHeis}_\psi$ and $\mathrm{LWeilHeis}_\psi$.

Let $V$ be a finite-dimensional $k$ vector space and $\beta$ a (not necessarily symmetric or non-degenerate) bilinear form on $V$.

\begin{dfn} \label{dfn:grbeta}
	\begin{enumerate}
		\item A an ascending filtration $V_0 \subset V_1 \subset \dots \subset V_n=V$ is called \emph{adapted} to $\beta$ if $\beta(V_i,V_{n-i})=0$ for all $i=0,\dots,n$.
		
		\item Given an adapted filtration, the \emph{associated graded} bilinear form $\mathrm{gr}(\beta)$ on $\mathrm{gr}(V)=\bigoplus_{i=1}^n V_i/V_{i-1}$ is defined as the direct sum of the pairings
		\[ V_i/V_{i-1} \otimes V_{n-i+1}/V_{n-i} \to k \]
		induced by restricting $\beta$ to $V_i \otimes V_{n-i+1}$. Note that we are excluding the summand $V_0$ from $\mathrm{gr}(V)$.
	\end{enumerate}
\end{dfn}

\begin{rem}
	\begin{enumerate}
		\item Part (1) of the definition implies $\beta(V_{n-i},V_i)=0$ by reindexing. Taking $i=0$ we conclude that $V_0$ is contained in both the left kernel and the right kernel of $\beta$. Thus, if $\beta$ is non-degenerate (it satisfies the equivalent conditions that its left resp. right kernel vanishes), then $V_0=\{0\}$ and our definition of $\mathrm{gr}(V)$ agrees with the usual one.
		
		\item The form $\mathrm{gr}(\beta)$ is non-degenerate if and only if $V_{n-i}$ equals both the left and the right annihilator of $V_i$ under $\beta$ for all $i$; in particular $V_0$ is equal to the left and the right kernel of $\beta$. In that case, $\mathrm{gr}(\beta)$  equals $\mathrm{gr}(\bar\beta)$, where $\bar\beta$ is the non-degenerate form on $V/V_0$ induced by $\beta$.
	\end{enumerate}
\end{rem}

\begin{lem} \label{lem:grbetadet}
	Assume $\beta$ and $\mathrm{gr}(\beta)$ are non-degenerate. In $k^\times/k^{\times,2}$ we have the identity $\det(\beta) = \det(\mathrm{gr}(\beta))$.
\end{lem}
\begin{proof}
	Choose a basis $(v_1,\dots,v_m)$ of $V$ adapted to the filtration, i.e. each $V_i$ has as basis those elements of the basis of $V$ that lie in $V_i$. Projecting those elements that lie in $V_i$ but not in $V_{i-1}$ modulo $V_{i-1}$ produces a basis of $V_i/V_{i-1}$. In this way, the basis produces a basis of $\mathrm{gr}(V)$. With respect to those bases the matrix of $\beta$ resp. $\mathrm{gr}(\beta)$ is given by 
	
\[
\setlength{\arraycolsep}{0pt}
\left(
\begin{array}{c|c|c|c|c}
  \diagblock{6mm}{0} & 0 & 0 & 0 & \star \\ \hline
  0 & \diagblock{8mm}{0} & 0 & \star & \star \\ \hline
  0 & 0 &
  \diagblock{10mm}{\star}
  & \star & \star \\ \hline
  0 & \star & \star & \diagblock{8mm}{\star} & \star \\ \hline
  \star & \star & \star & \star & \diagblock{6mm}{\star}
\end{array}
\right)
\quad\tx{resp.}\quad
\left(
\begin{array}{c|c|c|c|c}
  \diagblock{6mm}{0} & 0 & 0 & 0 & \star \\ \hline
  0 & \diagblock{8mm}{0} & 0 & \star & 0 \\ \hline
  0 & 0 &
  \diagblock{10mm}{\star}
  & 0 & 0 \\ \hline
  0 & \star & 0 & \diagblock{8mm}{0} & 0 \\ \hline
  \star & 0 & 0 & 0 & \diagblock{6mm}{0}
\end{array}
\right).
\]
The determinants of these two matrices are equal.
\end{proof}

\begin{lem} \label{lem:grbetasym}
	Assume the characteristic of $k$ is not $2$, and write $\beta^t(x,y)=\beta(y,x)$ and let $\beta_s=\tfrac{1}{2}(\beta+\beta^t)$ and $\beta_a=\tfrac{1}{2}(\beta-\beta^t)$, so that $\beta_s$ is symmetric, $\beta_a$ skew-symmetric, and $\beta=\beta_s+\beta_a$. A filtration adapted to $\beta$ is also adapted to $\beta_s$ and $\beta_a$ and we have
	\[ (\mathrm{gr}\beta)_s=\mathrm{gr}(\beta_s)\qquad\text{and}\qquad (\mathrm{gr}\beta)_a = \mathrm{gr}(\beta_a). \]
\end{lem}
\begin{proof}
	That the filtration is adapted to $\beta_s$ and $\beta_a$ is immediate. The formation of $\tx{gr}$ is linear, hence $\tx{gr}(\beta)=\tx{gr}(\beta_s)+\tx{gr}(\beta_a)$. It is immediate that the first summand is symmetric and the second skew-symmetric, hence the claim.
\end{proof}

\begin{lem} \label{lem:q-u-s} Assume the characteristic of $k$ is not $2$. If $\mathrm{gr}(\beta)$ is symmetric and non-degenerate then $\beta$ is non-degenerate on $V/V_0$ and the determinants of the quadratic form $q(v)=\beta(v,v)$ and the (possibly non-symmetric) bilinear form $\beta$, both on $V/V_0$, are equal in $k^\times/k^{\times,2}$.
\end{lem}
\begin{proof}
	We pass from $V$ to $V/V_0$. The non-degeneracy of the associated graded implies that of both $\beta$ and $\beta_s$. By definition we have $\det(q)=\det(\beta_s)$. Using Lemmas \ref{lem:grbetadet} and \ref{lem:grbetasym} we see
	\[ \det(q)=\det(\beta_s)=\det(\tx{gr}(\beta_s))=\det(\tx{gr}(\beta)_s) \stackrel{(*)}{=}\det(\tx{gr}(\beta))=\det(\beta), \]
	where the identity marked with $(*)$ is due to the assumption that $\tx{gr}(\beta)$ is symmetric.
\end{proof}

For the remainder of this subsection we shall assume that the characteristic of $k$ is not $2$. 

\begin{lem} \label{lem:graded-spinor}
	Assume the characteristic of $k$ is not $2$, and $\beta$ and $\mathrm{gr}(\beta)$ are non-degenerate symmetric. 
	\begin{enumerate}
		\item There exists a filtration-preserving isometry $(\mathrm{gr}(V),\mathrm{gr}(\beta)) \to (V,\beta)$ that induces the identity on each graded piece.
		\item If $s \in O(V,\beta)$ preserves the filtration, then the spinor norm of $s$ with respect to $q(v)=\beta(v,v)$ equals the spinor norm of $\mathrm{gr}(s) \in O(\mathrm{gr}(V),\mathrm{gr}(\beta))$ with respect to $\mathrm{gr}(q)(v) = \mathrm{gr}(\beta)(v,v)$.
	\end{enumerate}
\end{lem}
\begin{proof}
	(1) Let $s=\lfloor n/2\rfloor$ and $r=\dim(V_s)$. Choose an adapted basis $e_1,\dots,e_r$ of $V_s$ with respect to the filtration $V_0 \subset \dots \subset V_s$ and let $\bar f_1,\dots,\bar f_r$ be the $\beta$-dual basis on $V/V_{n-s}$. If $n$ is odd then take a orthogonal basis of $V_{s+1}/V_s$ and let $z_1,\dots,z_m \in V_{s+1}$ be arbitrary lifts of it. There exist lifts $f_1,\dots,f_r$ of $\bar f_1,\dots,\bar f_r$ in their respective filtration steps which are still dual to $e_1,\dots,e_r$ and orthogonal to $z_1,\dots,z_m$. Replacing the $f_i$ simultaneously by $f_i-\tfrac{1}{2}\sum_{j=1}^r\beta(f_i,f_j)e_j$
	preserves these properties and their images in the associated graded, and ensures that in addition $\beta(f_i,f_j)=0$ for all $i,j$. Putting all these together we obtain a Witt basis of $V$ with respect to $\beta$. The unique isomorphism $V \to \mathrm{gr}(V)$ sending each basis vector of $V$ to its image in the appropriate graded piece is a filtration-preserving isometry that induces the identity on each graded piece.

	(2) Choose an isometry as in (1) and let $t$ be the image of $\mathrm{gr}(s)$ in $O(V,\beta)$ under this isometry. Both $t$ and $s$ are filtration-preserving automorphisms that induce the same action on each graded piece, hence $u:=st^{-1}$ is a unipotent orthogonal transformation of $(V,\beta)$. Write $u=1+n$ for a nilpotent endomorphism $n$ of $V$. Zassenhaus' formula (see \cite[\S2]{Zass1962} for the original result, and \cite[Theorem 3.5]{CT2020} for more modern language) shows that the spinor norm of $u$ equals the determinent of $(u+1)/2=(n+2)/2=n/2+1$. This is still a unipotent element, so its determinant is $1$. The multiplicativity of the spinor norm now implies that the spinor norms of $s$ and $t$ are equal.
\end{proof}

The following lemma appears essentially as \cite[Proposition 2.15]{Mezer2023}, formulated in the setting of local fields. Since we need a slightly different formulation (in particular the case of finite fields) and the proof is short and works for arbitrary fields, we reproduce it here for the convenience of the reader. Recall that given an endomorphism $X$ of a $k$-vector space $W$, an $X$-cyclic vector in $W$ is a vector $w \in W$ such that $w,Xw,\dots,X^{n-1}w$ generates $W$, and hence is a basis of $W$, where $n=\dim(W)$.

\begin{lem} \label{lem:spnormal}
	Assume the characteristic of $k$ is not $2$. Let $\eta$ be a symplectic form on $V$ and let $X \in \mf{sp}(V,\eta)$ be nilpotent. There exists a (non-canonical) $X$-stable orthogonal decomposition $V=\bigoplus V_i$ such that for each $i$ one of the following mutually exclusive conditions holds
	\begin{enumerate}
		\item $V_i$ has an $X$-cyclic vector.
		\item There exists an $X$-stable polarization $V_i=L_i \oplus L_i'$ such that $L_i$ and $L_i'$ have $X$-cyclic vectors.
	\end{enumerate}
\end{lem}
\begin{proof}
	We induct on $\dim(V)$, beginning with the trivial case $\dim(V)=0$. Let $n$ be the least positive integer such that $X^n=0$ and consider the bilinear form $b(v,w)=\eta(v,X^{n-1}w)$. Then $b$ is a non-zero bilinear form and satisfies $b(w,v)=(-1)^n b(v,w)$. 
	
	If $n$ is even then $b$ is a non-zero symmetric bilinear form on $V$ so there exists $v \in V$ such that $b(v,v) \neq 0$. The sequence $v,Xv,\dots,X^{n-1}v$ is then linearly independent and the subspace $V_1$ that it generates is $X$-stable. It is also non-degenerate, becase the matrix of $\eta|_{V_1}$ in terms of this basis has its $(i,j)$ entry equal to zero when $i+j>n-1$, and non-zero when $i+j=n-1$, where $0 \leq i,j \leq n-1$, hence its determinant is non-zero. Therefore $V=V_1 \oplus V_1^\perp$ and we can apply inductive hypthesis to $V_1^\perp$.

	If $n$ is odd then $b$ is non-zero skew-symmetric, so there exist $v,w \in V$ such that $b(v,w)=1$. We will replace $v,w$ by $v',w'$ whose cyclic subspaces
	are isotropic. We first construct $v'\in V$ and set $\ell_i=X^iv'$. From $\eta(\ell_i,\ell_j)=\eta(X^iv',X^jv')=(-1)^{i}\eta(v',X^{i+j}v')$ we see that the condition on $v'$ we need is $\eta(v',X^iv')=0$ for all $0 \leq i \leq n-1$. Note that, if $i$ is even, then we have $\eta(v',X^iv')=\eta(X^iv',v')=-\eta(v',X^iv')$. Thus $\eta(v',X^iv')=0$ is automatic for even $i$ and we only need to force it for odd $i$. We manufacture such $v'$ out of $v$ recursively as follows. Set $v_0'=v$ and for each $0<r<n/2$ set $v_r'=v_{r-1}'-(\eta(v_{r-1}',X^{n-2r}v_{r-1}')/2)X^{2r-1}w$. We claim that $v':=v'_{(n-1)/2}$ satisfies $\eta(v',X^iv')=0$ for all odd $i$, hence for all $i$. Indeed, for any $r>0$ we have $\eta(v_r',X^{n-1}w)=\eta(v_{r-1}',X^{n-1}w)$ and by downward recursion on $r$ we obtain
	\begin{equation} \label{eq:vr'}
		b(v'_r,w)=\eta(v_r',X^{n-1}w)=\eta(v,X^{n-1}w)=b(v,w)=1.
	\end{equation}
	Furthermore, for $r>0$ and odd $i$ we have 
	\begin{eqnarray*}
	\eta(v_r',X^iv_r')&=&\eta(v_{r-1}'-(\eta(v_{r-1}',X^{n-2r}v_{r-1}')/2)X^{2r-1}w,X^iv_{r-1}'-(\eta(v_{r-1}',X^{n-2r}v_{r-1}')/2)X^{2r+i-1}w)\\
	&=&\eta(v_{r-1}',X^iv_{r-1}')\\
	&-&(\eta(v_{r-1}',X^{n-2r}v_{r-1}')/2)\cdot(\eta(v_{r-1}',X^{2r-1+i}w)+\eta(X^{2r-1}w,X^iv_{r-1}'))\\
	&+&(\eta(v_{r-1}',X^{n-2r}v_{r-1}')/2)^2\cdot\eta(X^{2r-1}w,X^{2r+i-1}w).\\	
	(1)&=&\eta(v_{r-1}',X^iv_{r-1}')\\
	(2)&-&\eta(v_{r-1}',X^{n-2r}v_{r-1}')\cdot\eta(v_{r-1}',X^{2r-1+i}w)\\
	(3)&-&(\eta(v_{r-1}',X^{n-2r}v_{r-1}')/2)^2\cdot\eta(w,X^{4r+i-2}w).\\	
	\end{eqnarray*}
	If $i\geq n-2r$ then $4r+i-2 \geq n+2r-2 \geq n$ and the term marked (3) vanishes.  If $i>n-2r$ then $2r-1+i>n-1$ and the term marked (2) vanishes, while if $i=n-2r$ then $2r-1+i=n-1$ and using \eqref{eq:vr'} we see that the term marked (2) cancels the term marked (1). We thus obtain for $r>0$ and odd $i$
	\begin{equation}
			\eta(v'_r,X^iv'_r)=
			\begin{cases}
				0,&i=n-2r\\
				\eta(v'_{r-1},X^iv'_{r-1}),&i>n-2r.
			\end{cases}
	\end{equation}
	Recursion on $r$ shows $\eta(v',X^iv')=0$ for all odd $i$, hence for all $i$, as desired. The construction of $v'$ is thus complete.

	To obtain $w'$ we apply the same construction to the pair $(w,-v')$, for which $b(w,-v')=1$. This gives $w' \in V$ such that $\eta(X^iw',X^jw')=0$ for all $i,j$ and $b(w',-v')=1$. Now define the subspaces $L=k[X]v'$ and $L'=k[X]w'$ of $V$; both are isotropic.

	It remains to show that the sum $L + L'$ is direct and the resulting vector subspace $V_1 \subset V$ is non-degenerate. For this we take the ordered list of vectors $(v',Xv',\dots,X^{n-1}v',w',Xw',\dots,X^{n-1}w')$ with $v'$ and $w'$ just constructed. Since $b(v',w')=b(v',w)=b(v,w)=1$, and $\eta(X^iv',X^jw')=(-1)^i\eta(v',X^{i+j}w')=0$ when $i+j \geq n$, we see that the matrix of $\eta$ with respect to this list of vectors consists of four $n \times n$ blocks, of which the top left and bottom right are zero, and the top right and bottom left have zeroes below the anti-diagonal and $\pm 1$ at the anti-diagonal. This matrix has determinant $1$ and is thus non-degenerate. Therefore the list of vectors is linearly independent and $V_1=L \oplus L'$ is non-degenerate. Apply inductive hypothesis to $V_1^\perp$.
\end{proof}

\begin{cor} \label{cor:spnormal}
Assume the characteristic of $k$ is not $2$. Let $\eta$ be a symplectic form on $V$ and let $u \in \Sp(V)(k)$ be unipotent. The determinant of the quadratic form $q(v)=\eta(v,uv)$ on $V/V^u$ equals the determinant of the (generally non-symmetric) bilinear form $\beta(v,w)=\eta(v,(u-1)w)$ on $V/V^u$.
\end{cor}
\begin{proof}
	Recall the Cayley transform $c(x)=(1+x)^{-1}(1-x)$. It is a self-inverse bijection on the set of endomorphisms of $V$ that do not have eigenvalue $-1$. It induces a bijection between the subset of nilpotent elements of $\mf{sp}(V)(k)$ and the subset of unipotent elements of $\Sp(V)(k)$. Letting $x=c(u)$ we have $u=c(x)$ and since 
	\[ u-1 = c(x)-1=\frac{1-x}{1+x}-1 = x\cdot \frac{-2}{1+x}\]
	and $(-2)(1+x)^{-1}$ is invertible, we see that $x^i$ and $(u-1)^i$ have the same kernel.

	We apply Lemma \ref{lem:spnormal} to $x=c(u)$ and obtain an $\eta$-orthogonal decomposition $V=\bigoplus V_j$ that is stable under $x$, hence also under $u=c(x)$, and is therefore orthogonal for both $q$ and $\beta$. We may thus reduce to considering each individual $V_j$, i.e. we may assume that $V$ is either $x$-cyclic or has a polarization $V=L \oplus L'$ with $x$-cyclic Lagrangian subspaces $L,L'$.

	Define $\delta=1$ if $V$ is $x$-cyclic and $\delta=2$ if $V=L \oplus L'$ with $x$-cyclic $L,L'$. Let $n=\dim(V)/\delta-1$ and define $V_i=\tx{ker}(x^{i+1})=\tx{ker}((u-1)^{i+1})$ for $i=-1,\dots,n$ and $\ol{V}_i=V_i/V_0$ for $i=0,\dots,n$. Then $\{0\}=V_{-1}\subset V_0 \subset \dots \subset V_n=V$ is a filtration of $V$ and $\{0\}=\ol{V}_0 \subset \dots \subset \ol{V}_n=\ol{V}$ is a filtration in $\ol{V}=V/V^u$, both with successive subquotients of dimension $\delta$; note that $V_0=V^u$. We have $(u-1)V_{i+1}=V_i=xV_{i+1}$ and hence $\beta(V_i,V_{n-i})=\eta(V_i,(u-1)V_{n-i})=\eta(V_i,V_{n-i-1})=\eta(x^{n-i}V,x^{i+1}V)=(-1)^{n-i}\eta(V,x^{n+1}V)=0$, showing that the filtration $V_i$ is adapted to $\beta$. Next we claim that $\tx{gr}(\beta)$ is symmetric. This is equivalent to $\beta(x,y)=\beta(y,x)$ for $x \in V_i$ and $y \in V_{n-i+1}$. We compute
	\begin{eqnarray*}
		\beta(x,y)-\beta(y,x)&=&\eta(x,(u-1)y)-\eta(y,(u-1)x)\\
		&=&\eta(x,(u-1)y)-\eta((u^{-1}-1)y,x)\\
		&=&\eta(x,(u+u^{-1}-2)y)\\
		&=&\eta(x,u^{-1}(u-1)^2y)\\
		&\in&\eta(V_i,V_{n-i-1})\\
		&=&0.
	\end{eqnarray*}
	The left and right radicals of $\beta$ on $V$ are both $V_0=V^u$. On $\bar V=V/V_0$, we have $\dim\bar V_i=i\delta$ and $\dim\bar V=n\delta$. Therefore  $\bar V_{n-i}$ is both the left and right annihilator of $\bar V_i$, and $\tx{gr}(\beta)$ is non-degenerate. The claim now follows from Lemma \ref{lem:q-u-s}.	
\end{proof}

\subsection{The classical metaplectic groups} \label{sub:classmeta}

We continue with an arbitrary field $k$ of characteristic not $2$ and a symplectic vector space $(V,\eta)$. Following \cite{Thomas2013}, we introduce a subgroup $\Mp(V) \subset \WMuniv(V)$ that is an extension 
\[ 1 \to I_k^2/I_k^3 \to \Mp(V) \to \Sp(V)(k) \to 1, \]
as well as a variant $\Mp_L(V) \subset \WM_L(V)$.

\begin{cor}
	Let $L \subset V$ be a Lagrangian subspace. Choose an orientation $o$ on $L$ and set $\tilde L=(L,o)$. For $g \in \Sp(V)(k)$ let $g\tilde L=(gL,go)$. Then $c_L(g) := m(g\tilde L,\tilde L) \in \QW(k)/I_k^2$ does not depend on the choice of $o$ and the image of $z_L \in Z^2(\Sp(V)(k),\QW(k))$ in $Z^2(\Sp(V)(k),\QW(k)/I_k^2)$ equals the negative of the coboundary of $c_L \in C^1(\Sp(V)(k),\QW(k)/I_k^2)$.
\end{cor}
\begin{proof}
	If we replace $o$ by $ao$ for $a \in k^\times$ then the determinant of the quadratic form representing $m(g\tilde L,\tilde L)$ rescales by $a^2 \equiv 1 \in k^\times/k^{\times,2}$, showing that $m(g\tilde L,\tilde L)$ is unchanged in $\QW(k)/I_k^2$. The second statement follows from Proposition \ref{pro:maslov_break} upon noting that $m(gh\tilde L,g\tilde L)=m(h\tilde L,\tilde L)$.
\end{proof}

\begin{cor} \label{cor:unisec}
	Fix an orientation on $V$ and for every $g \in \Sp(V)(k)$ use the isomorphism $\Gamma_g \to V$ given by composing the inclusion $\Gamma_g \subset V \oplus V$ with the first projection. Write $\tilde \Gamma_g$ for the corresponding oriented Lagrangian subspace of the space $U=V \oplus V$ with symplectic form $\eta_U=(-\eta) \oplus \eta$. Then $c(g):=m(\tilde\Gamma_g,\tilde\Gamma_1)$\index[notation]{c@$c$} does not depend on the choice of $o$ and the image of $z \in Z^2(\Sp(V)(k),\QW(k))$ in $Z^2(\Sp(V)(k),\QW(k)/I_k^2)$ equals the negative of the coboundary of $c \in C^1(\Sp(V)(k),\QW(k)/I_k^2)$.
\end{cor}
\begin{proof}
	If we replace $o$ by $ao$ with $a \in k^\times$ then the orientations on $\tilde \Gamma_1$ and $\tilde \Gamma_g$ both multiply by $a$, hence the the determinant of the quadratic form representing $m(\tilde \Gamma_g,\tilde\Gamma_1)$ rescales by $a^{-2} \equiv 1 \in k^\times/k^{\times,2}$, showing that $m(\tilde \Gamma_g,\tilde\Gamma_1)$ is unchanged in $\BW(k)/I_k^2$. The second statement follows from Proposition \ref{pro:maslov_break} upon noting that $m(\tilde \Gamma_{gh},\tilde\Gamma_{g})=m(\tilde \Gamma_h,\tilde\Gamma_1)$. \end{proof}

\begin{dfn}
\begin{enumerate}
	\item Define the \emph{universal metaplectic group} as 
	\[ \Mp(V) = \{w \boxtimes g \in \WMuniv(V)| w-c(g) \in I_k^2\}. \]
	\item For a fixed Lagrangian subspace $L \subset V$ define the \emph{metaplectic group relative to $L$} as
	\[ \Mp_L(V) = \{w \boxtimes g \in \WMuniv_L(V)| w-c_L(g) \in I_k^2 \}. \]
\end{enumerate}
\end{dfn}

\begin{rem}
	It is clear that we obtain the exact sequences
	\[ 0 \to I_k^2/I_k^3 \to \Mp(V) \to \Sp(V)(k) \to 1 \]
	and 
	\[ 0 \to I_k^2/I_k^3 \to \Mp_L(V) \to \Sp(V)(k) \to 1. \]
	The Clifford invariant provides the isomorphism $I_k^2/I_k^3 \to H^2(\Gamma,\mu_2)=\tx{Br}_k[2]$. If $k$ is a local field that is not $\C$ then pusing out by the invariant map $\tx{Br}_k[2] \to \tfrac{1}{2}\Z/\Z$ of local class field theory we obtain the usual metaplectic double cover of $\Sp(V)(k)$. If $k$ is a finite field or $\C$ then $\tx{Br}_k=\{0\}$ and we obtain the isomorphisms $\Mp(V) = \Sp(V)(k) = \Mp_L(V)$.
\end{rem}

\begin{pro}
	For $g \in \Sp(V)(k)$ the element $c(g) \in \QW(k)/I_k^2$ is represented by the class of quadratic forms on $(1-g)V$ whose determinant is equal to that of the usually non-symmetric bilinear form 
	\[ (1-g)V \times (1-g)V \to k,\qquad ((1-g)x,(1-g)y) \mapsto \eta(x,(1-g)y). \]
\end{pro}
\begin{proof}
	Choose a basis $e_1,\dots,e_{2n}$ of $V$ whose first $d$-many vectors are a basis of $V^g$. Then $(e_1,e_1),\dots,(e_{2n},e_{2n})$ is a basis of $\Gamma_1$ and $(e_1,ge_1),\dots,(e_{2n},ge_{2n})$ is a basis of $\Gamma_g$, each containing the basis $(e_1,e_1),\dots,(e_{d},e_{d})$ of $\Gamma_1 \cap \Gamma_g$. The orientations corresponding to these bases provide under the first projection the same orientation on $V$, so can be used to compute $c(g)$.The vectors in the two bases that do not belong to the basis of $\Gamma_1 \cap \Gamma_g$ project to bases of $\Gamma_1/(\Gamma_1 \cap \Gamma_g)$ and $\Gamma_g/(\Gamma_1 \cap \Gamma_g)$.
	with their induced orientations.
	
	The determinant of the class of forms representing $c(g)$ is computed as the determinant of the matrix with $i,j$-th entry
	\[ \eta_U((e_j,ge_j),(e_i,e_i)) =\eta(e_i,e_j)-\eta(e_i,ge_j)=\eta(e_i,(1-g)e_j). \]
	for $i,j$ ranging over the indices of the last
$2n-d$ basis vectors. The vectors $(1-g)e_i$
with these indices form a basis of $(1-g)V$,
and the above scalars are the entries of the matrix of the bilinear form in the statement (here the determinant
of the empty matrix is understood to be $1$.
\end{proof}

\begin{rem}
	\begin{enumerate}
		\item The endomorphism $1-g : V \to V$ induces an isomorphism $V/V^g \to (1-g)V$ and we can alternatively describe $c(g) \in \QW(k)/I_k^2$ as the class of quadratic forms on $V/V^g$ whose determinant equals that of the (generally non-symmetric) bilinear form $(\bar x,\bar y) \mapsto \eta(x,(1-g)y)$.
		\item If $s$ is the semi-simple part of $g$ then the restriction of $\eta$ to $V^s$ is a symplectic form. If the unipotent part $u$ of $g$ belongs to $\Sp(V)(k)$ then it induces a symplectic automorphism of $(V^s,\eta)$ and we can consider $c(u|V^s)$.
	\end{enumerate}
\end{rem}

\begin{pro} \label{pro:metasection}
	Assume that $g \in \Sp(V)(k)$ has Jordan decomposition $g=su$ with $s,u \in \Sp(V)(k)$. 
	\begin{enumerate}
		\item In $\QW(k)/I_k^2$ we have the identity $c(g|V) = c(s|V) + c(u|V^s)$.
		\item We have $c(s|V) \in I_k/I_k^2$ and $s\det(c(s|V))=\det(s-1|V_s)\cdot (-1)^{\dim(V_s)/2}$.
		\item The class $c(u|V^s)$ is represented by the quadratic form $v \mapsto \eta(uv,v)$ on $(V^s)/V^{su}$.
	\end{enumerate}
\end{pro}
\begin{proof}
	(1) We have the $g$-stable orthogonal direct sum $V=V^s \oplus V_s$, where $V_s$ is the direct sum of all $s$-eigenspaces in $V$ with eigenvalue not equal to $1$. We have $V^g=(V^s)^u$ and obtain the direct sum $V/V^g = V_s \oplus V^s/(V^s)^u$ that is orthogonal with respect to the form $(\bar x,\bar y)=\eta(x,(1-g)y)$. We thus see $c(g|V)=c(g|V_s)+c(g|V^s)$. Since $g|_{V^s}=u|_{V^s}$ it remains to show $c(g|V_s)=c(s|V_s)$, noting the isomorphism $V_s \to V/V^s$ induced by the inclusion $V_s \to V$. This is equivalent proving $s\det(c(g|V_s))=s\det(c(s)|V_s))$, which will be done in (2) below.

	(2) By definition $c(s|V)=c(s|V_s)$, so we may assume $V=V_s$, i.e. $\tx{ker}(s-1)=0$. Since $V$ is symplectic, it is even-dimensional, hence $c(s|V) \in I_k/I_k^2$. We now compute its signed determinant. In order to complete the proof of (1) we consider not just $s$ but also $g=su$ and using the fact that $V$ is symplectic, so its dimension is even, we let $2m=\dim_k(V)$ and compute
		\begin{eqnarray*}
		s\det(c(g))&=&(-1)^{2m(2m-1)/2}\det(c(g)|V)\\
		&=&(-1)^m\det(g-1|V)\cdot\det(\eta)\\
		&=&(-1)^m\det(g-1|V),
	\end{eqnarray*}
	where the last line comes from the fact that $\eta$ is symplectic, so its Gram determinant is $1$. Taking $g=s$ this completes the proof of (2). To complete the proof of (1) we must confirm the outstanding identity $s\det(c(g|V))=s\det(c(s)|V))$. Using that $s-1$ is invertible we have
	\begin{eqnarray*}
		\det(g-1|V)&=&\det((s-1)(1+(s-1)^{-1}s(u-1))|V)\\
		&=&\det(s-1|V),
	\end{eqnarray*}
	the last line coming from the fact that $(s-1)^{-1}s(u-1)$ is nilpotent. 

	(3) follows from Corollary \ref{cor:spnormal} applied to $u \in \Sp(V^s,\eta|_{V^s})$.
\end{proof}

\begin{exa}
	Both $c(s),c(u) \in W(k)/I_k^2$ are in general non-trivial. Consider $V=k \oplus k$ with the standard symplectic form $x_1y_2-x_2y_1$. 
	
	For $s \in \mathrm{Sp}(V)=\mathrm{SL}_2(k)$ a diagonal matrix with entries $x,x^{-1}$ assumed not equal to $1$ we have $s\det(c(s))=-(x-1)(x^{-1}-1)=x^{-1}(x-1)^2\equiv x^{-1} \equiv x$ in $k^\times/k^{\times,2}$. 
	
	For $u \neq 1$ an upper triangular unipotent matrix with top right entry $a$, we have $V^u=k \oplus \{0\}$, so $\{0\} \oplus k$ maps isomorphically to $V/V^u$, and the quadratic form representing $c(u)$ sends $(0,y) \in \{0\} \oplus k$ to $ay^2$, i.e. $c(u)=[a]_Q$.
\end{exa}

\subsection{Review of the Weil index}\label{sub:Weilindex} 

In this subsection $k$ is a locally compact field, thus a finite extension of one of the following fields: $\R$, $\Q_p$, $\F_p((t))$, $\F_p$. We assume that the characteristic of $k$ is not $2$. Let  $\psi : k \to \C^\times$ be a non-trivial unitary character. The Weil index is a unitary character\index[notation]{gammapsi@$\gamma_\psi$}
\[ \gamma_\psi = \gamma_\psi^k : \QW(k) \to \mu_8 \]
defined in \cite{Weil64}, see also \cite{Cartier64}. 
Useful information about $\gamma_\psi$ is collected in \cite[Appendix A]{RR93}.

	The following is immediate from Weil's definition, and is also implicitly used in \cite[\S30]{Weil64}.

	\begin{fct} \label{fct:gammaind}
    Let $\ell/k$ be a finite separable extension, $\psi : k \to \C^\times$ non-trivial. We have
    \[ \gamma^\ell_{\psi\circ \mathrm{tr}_{\ell/k}} = \gamma^k_\psi \circ \mathrm{tr}_{\ell/k}.\]
\end{fct}

\begin{lem} \label{lem:weilherm}
	Let $k/k_0$ be a quadratic extension, $\kappa : k_0^\times/N_{k/k_0}(k^\times) \to \{\pm1\}$ the associated character (trivial if $k$ finite), $\psi : k_0 \to \C^\times$ a non-trivial character, and $(V,h)$ a Hermitian space with respect to $k/k_0$.
	\begin{enumerate}
		\item $\gamma_\psi^{k_0}(\mathrm{tr}(V,h))=\kappa(\det(h))\cdot \gamma_\psi(N_{k/k_0})^{\dim_k(V)}$, where $N_{k/k_0} : k \to k_0$ is the norm map considered as a 2-dimensional quadratic space over $k_0$.
		\item $\gamma_\psi(N_{k/k_0})$ equals $\epsilon_L(1/2,\kappa,\psi)$ when $k$ is local and $-1$ when $k$ is finite.
	\end{enumerate}
\end{lem}
\begin{proof}
	(1) Diagonalizing $h$ we immediately reduce to considering only the elements $[a]_H \in \HW(k)$ for $a \in k_0^\times$. Then $\mathrm{tr}([a]_H)=aN_{k/k_0}$. Choosing a trace-zero element $0 \neq b \in k$ we have the $k_0$-basis $\{1,b\}$ of $k$ which is diagonal for $aN_{k/k_0}$ and we see $aN_{k/k_0}=[a]_Q \oplus [-ab^2]_Q$. We must therefore show
	    \[ \gamma_\psi([a]_Q)\gamma_\psi([-ab^2]_Q)\gamma_\psi([1]_Q)^{-1}\gamma_\psi([-b^2]_Q)^{-1}=\kappa_{k/k_0}(a). \]
    Using the notation $\gamma(a,\psi)=\gamma_\psi([a]_Q)/\gamma_\psi([1]_Q)$ the left-hand side above becomes
    \[ \gamma(a,\psi)\gamma(a,-b^2\psi), \]
    which using \cite[Corollary A.5(1)]{RR93} equals 
    \[ \gamma(a,\psi)\cdot (a,-b^2)_{k_0} \cdot \gamma(a,\psi) \]
    the symmetry and multiplicativity of the Hilbert symbol, together with $k=k_0(b)$, imply $(a,-b^2)_{k_0}=(-1,a)_{k_0}\cdot (a,b^2)_{k_0}=(-1,a)_{k_0} \cdot \kappa(a)$, which together with \cite[Corollary A.5(3)]{RR93} completes the proof.

	(2) For local $k$ this is \cite[Lemma 1.2]{JL79}. Suppose that $k$ is finite and put $q=|k_0|$. The norm has one preimage of $0$ and $q+1$ preimages of every element of $k_0^\times$. Therefore the finite Gauss-sum formula gives
	\[ \gamma_\psi(N_{k/k_0})=q^{-1}\sum_{x\in k}\psi(N_{k/k_0}(x))=q^{-1}\left(1+(q+1)\sum_{t\in k_0^\times}\psi(t)\right)=-1.\qedhere\]
\end{proof}

Let $I_k \subset \BW(k)$ be the fundamental ideal of even-dimensional forms, which we also consider as a subgroup of $\QW(k)$ via the isomorphism $\BW(k) \to \QW(k)$. 

\begin{lem} \label{lem:weilindexsubqout}
\begin{enumerate}
	\item $I_k^3=\{0\}$ unless $k=\R$, in which case $\QW(k)=\Z$ via the signature map and $I_k^r=2^r\Z$.
	\item The restriction of $\gamma_\psi$ to $I_k^3$ is trivial.
	\item The restriction of $\gamma_\psi$ to $I_k^2$ induces an isomorphism $I_k^2/I_k^3 \to Z_k$, where $Z_k \subset \mu_8$ is the image of quadratic Hilbert symbol, thus $Z_k=\{1\}$ for $k$ finite or complex, and $Z_k=\{\pm1\}$ for the other locally compact fields.
	\item Assume $k$ is finite. The restriction of $\gamma_\psi$ to $I_k$ equals the composition
		\[ \xymatrix{ I_k\ar[r]^-\cong &I_k/I_k^2 \ar[r]^-\cong_-{s\det}& k^\times/k^{\times,2}\ar[r]^-\cong_{\mathrm{sgn}}& \{\pm1\}}.\] 
\end{enumerate}	
\end{lem}
\begin{proof}
	(1) is well-known, see e.g. \cite[Chapters 2,6]{Lam2005} and the proof of \cite[Theorem A.2]{Thomas2013}, while 
(2) and (3) are proved in \cite[A.4.1(d)]{Thomas2013}. For (4) we first note that $I_k^2=I_k^3=\{0\}$ by (3), hence the isomorphism $I_k \to I_k/I_k^2$. The second isomorphism is \cite[Chapter II, Proposition 2.1]{Lam2005}, while the third is clear from the fact that $k^\times$ is a cyclic group of even order. To see that the composition of the three equals the Weil index (for any $\psi$), we use \cite[A.4.1(b,c)]{Thomas2013} and the fact that the Hilbert symbol and Hasse invariants are trivial for finite $k$ to compute that for $[q] \in I_k$ of dimension $n$ we have 
\[ \gamma_\psi(q)=\gamma_\psi(1)^n\mathrm{sgn}_k(\det(q)) = \mathrm{sgn}_k((-1)^{n/2}\cdot\det(q))=\mathrm{sgn}_k(s\det(q)), \] 
using that $(-1)^{n/2}=(-1)^{n(n-1)/2}$ since $n$ is even.
\end{proof}

Consider a non-archimedean local field $F$ with ring of integers $\cO_F$ and residue field $\rf$. Assume the characteristic of $\rf$ is not $2$. Choose a uniformizer $\pi \in F$ and write each $a \in F^\times$ as $\pi^{\mathrm{ord}_F(a)}\dot a$ with $\dot a \in \cO_F^\times$. Then
\[ \pi_0 : \QW(F) \to \QW(\rf),\qquad \pi_0([a]_Q)=\begin{cases}
	[\dot a]_Q,&2 \mid \mathrm{ord}_F(a)\\
	0,&2 \nmid \mathrm{ord}_F(a)
\end{cases}\]
is a well-defined group homomorphism, independent of the choice of $\pi$. Let $\psi : F \to \C^\times$ be such that its restriction to the ring of integers of $F$ factors through a non-trivial character $\bar\psi : \rf \to \C^\times$. It is shown in \cite[Proposition A.11]{RR93} that 
\[ \gamma_\psi = \gamma_{\bar\psi} \circ \pi_0. \]

\begin{dfn}
	The \emph{renormalized} Weil index is $\tilde\gamma_\psi(q)=\gamma_\psi(q/2)$.
	\index[notation]{gammapsitilde@$\tilde\gamma_\psi$}
\end{dfn}

\subsection{The Fourier--Weil representation} \label{sub:fwh}

In this section we will review some basic theory of the Weil representation, in a form suitable for our purposes. For a general and thorough introduction to the Weil representation and surrounding topics we refer to \cite{GanKudlaTakedaTheta}. Our focus here will be more in line with \cite{LionVergne} and \cite{Thomas2013}, in that we will use a more arithmetically motivated version of the metaplectic group, where the kernel of the cover is the Witt group of the base field. The reference \cite{LionVergne} is written mostly for the real numbers, with some crucial results extended to arbitrary local fields of characteristic different from $2$ in Appendix A, but with proofs often just sketched. The reference \cite{Thomas2013} works in the desired generality of any locally compact field of characteristic not $2$, but omits some arguments.
 
In this section $k$ is any locally compact field of characteristic different from $2$. Let $(V,\eta)$ be a finite-dimensional symplectic $k$-vector space. Let  $H(V)=\mb{G}_a \boxtimes_{\eta/2} V$ be the Heisenberg group.

Let $\psi : k \to \C^\times$ be a non-trivial continuous unitary character. The Stone--von-Neumann theorem states that there exists a unique irreducible unitary representation $\tx{Heis}_\psi$ of $H(V)(k)$ with central character $\psi$, called the \emph{Heisenberg} representation. Its uniqueness implies that its isomorphism class is fixed by the action of $\Sp(V)(k)$, hence that it extends to a projective representation $\mathrm{PWeilHeis}_\psi : \Sp(V)(k) \ltimes H(V)(k) \to \tx{PGL}(\tx{Heis}_\psi)$, where the semi-direct product corresponds to the action of $\Sp(V)$ on $H(V)$ via the natural action on $V$ and the trivial action on $\mb{G}_a$. It is known that this projective representation lifts to a linear representation  when $k$ is complex or finite, and that this lift is unique unless $k=\F_3$, where one can specify a preferred lift, see Lemma \ref{lem:weilheisuniq} below. For the other local $k$ it is known that $\mathrm{PWeilHeis}_\psi$ lifts uniquely to $\Mp(V) \ltimes H(V)(k)$, where $\Mp(V)$ is the unique 2-sheeted cover of $\Sp(V)(k)$. We will call this lift\index[notation]{LWeilHeispsi@$\mathrm{LWeilHeis}_\psi$} $\mathrm{LWeilHeis}_\psi : \Mp(V) \ltimes  H(V)(k) \to \tx{GL}(\tx{Heis}_\psi)$, where $\Mp(V)$ will be taken to be $\Sp(V)(k)$ whenever $k$ is complex or finite, and the unique $2$-sheeted cover of $\Sp(V)(k)$ otherwise. We denote its restriction to $\Mp(V)$ by $\mathrm{LWeil}_\psi$.\index[notation]{LWeilpsi@$\mathrm{LWeil}_\psi$}

However, there is a different natural lift of $\mathrm{PWeilHeis}_\psi$, involving a different natural cover of $\Sp(k)$, that will be relevant to us. Recall that a choice of a Lagrangian subspace $L \subset V$ provides the following presentation of $\tx{Heis}_\psi$:
\[ \tx{Heis}_{\psi,L} = \tx{Ind}_{k \times L}^{H(V)(k)}(\psi \boxtimes \mathbf{1}_{L}), \]
where the underlying vector space is the Hilbert space consisting of measurable functions $f : H(V)(k) \to \C$ such that $f((x,l)h)=\psi(x)f(h)$ for all $x \in k$, $l \in L$, $h \in H(V)(k)$, and $\int_{(k \times L) \lmod H(V)(k)}|\varphi(\ol{V})|^2d\ol{V} < \infty$ for some, hence any, Haar mesure $d\bar h$ on $(k \times L) \lmod H(V)(k) = L \lmod V$. To fix an $H(V)(k)$-invariant scalar product on this space one fixed $d\ol{V}$ and takes
\[ \<f_1,f_2\> := \int_{(k \times L) \lmod H(V)(k)} f_1(\ol{V})\cdot \ol{f_2(\ol{V})}d\ol{V}. \]

The uniqueness statement and Schur's lemma imply that for two Lagrangian subspaces $L_1,L_2 \subset V$ there exists an isomorphism $\tx{Heis}_{\psi,L_1} \to \tx{Heis}_{\psi,L_2}$ that is unique up to multiplication by a non-zero complex number. Once scalar products on both sides are fixed and the isomorphism is required to preserve them, the scalar ambiguity reduces to a complex number of absolute value $1$. 

It turns out that there exists an explicit choice 
\[ \mc{F}_{L_2,L_1} : \tx{Heis}_{\psi,L_1} \to \tx{Heis}_{\psi,L_2} \]
for this isomorphism, in the form of an intertwining operator, or partial Fourier transform, as follows.

Choose Haar measures $d\ol{V}_1$ on $V/L_1$ and $d\ol{V}_2$ on $V/L_2$. The operator $\mc{F}_{L_2,L_1}$ will depend on these choices. Choose further Haar measures $dy$ on $L_1 \cap L_2$ and $dv$ on $V$. The operator$\mc{F}_{L_2,L_1}$ will be independent of these choices. There are unique measures $d\bar \ell_1$ on $\bar L_1$ and $d\bar\ell_2$ on $\bar L_2$ such that the filtration $L_1 \cap L_2 \subset L_i \subset V$ leads to $dv = dy \times d\bar\ell_i \times d\ol{V}_i$. For $f : H(V)(k) \to \C$ lying in $\tx{Heis}_{\psi,L_1}$ define
\begin{equation} \label{eq:iop}
	(\mc{F}_{L_2,L_1}f)(x) = \int_{\bar L_2}f(\bar\ell_2 x)\Big|\frac{(d\bar\ell_2)^*}{(g_{L_2,L_1})_*(d\bar\ell_1)}\Big|^{1/2}d\bar\ell_2.
\end{equation}
Here $(d\bar\ell_2)^*$ is the measure on $\bar L_2^*$ that is $\psi$-dual to $d\bar\ell_2$.

\begin{lem}
	Assuming the integral \eqref{eq:iop} converges, it is independent of the choices of $dy$ and $dv$. Replacing $d\ol{V}_i$ by $a_id\ol{V}_i$ for some $a_1,a_2>0$ multiplies $\mc{F}_{L_2,L_1}$ by $(a_1/a_2)^{1/2}$.
\end{lem}
\begin{proof}
	If we replace $dv$ by $adv$ and $dy$ by $bdy$ with $a,b \in \R_{>0}$ and keep $d\ol{V}_i$ fixed, then $d\bar \ell_i$ multiplies by $a/b$, so $(d\bar\ell_2)^*$ multiplies by $b/a$, so the term in absolute values multiplies by $(b/a)^2$, and the square root of this change cancels the change of $d\bar\ell_2$.

	Keeping $dv$ and $dy$ fixed and multiplying $d\ol{V}_i$ by $a_i \in \R_{>0}$ the measure $d\bar\ell_i$ multiplies by $a_i^{-1}$, the term in absolute values multiplies by $a_1a_2$, so the integral by $(a_1/a_2)^{1/2}$. 
\end{proof}

\begin{thm}[Lion--Vergne A.1]
	\begin{enumerate}
		\item The integral \eqref{eq:iop} converges for all $f$ in a dense subspace of $\mathrm{Heis}_{\psi,L_1}$ and $\mc{F}_{L_2,L_1}$ extends to an isomorphism of unitary $H(V)(k)$-representations $\mathrm{Heis}_{\psi,L_1} \to \mathrm{Heis}_{\psi,L_2}$, where the scalar product on $\mathrm{Heis}_{\psi,L_i}$ is normalized by the Haar measure $d\ol{V}_i$.
		\item We have $\mc{F}_{L_2,L_1}^{-1}=\mc{F}_{L_1,L_2}$ provided the same measures $d\ol{V}_1$ and $d\ol{V}_2$ are used for both operators. 
	\end{enumerate}
\end{thm}

Thus, once the scalar product structures on $\mathrm{Heis}_{\psi,L_1}$ and $\mathrm{Heis}_{\psi,L_2}$ are fixed, the above formula determines a unique unitary isomorphism $\mc{F}_{L_1, L_2}$.

These isomorphisms are not just useful for explicating the independence of the isomorphism class of $\mathrm{Heis}_{\psi,L}$ of the parameter $L$. They can also be used to endow $\mathrm{Heis}_{\psi,L}$ with an action of $\Sp(V)(k)$, as follows. Choose a Haar measure $d \ol{V}$ on $V/L$. An element $g \in \Sp(V)(k)$ induces an isomorphism $V/L \to V/gL$ and on $V/gL$ we take the measure $g_*(d\ol{V})$. With respect to the scalar products on $\mathrm{Heis}_{\psi,L}$ and $\mathrm{Heis}_{\psi,gL}$ induced by these measures we obtain the unitary operator 
\[ A_g : \mathrm{Heis}_{\psi,L} \to \mathrm{Heis}_{\psi,gL},\qquad f \mapsto gf, \]
where $gf(v)=f(g^{-1}(v))$ for any function $f : V \to \C$. Then $R_g:=\mc{F}_{L,gL}\circ A_g$ becomes an isomorphism of unitary representations
\[ R_g : \mathrm{Heis}_{\psi,L} \to \mathrm{Heis}_{\psi,L}\circ g. \]
However, for $g_1,g_2 \in \Sp(V)(k)$ it turns out that $R_{g_1}\circ R_{g_2}$ is not exactly equal to $R_{g_1g_2}$.

The following is a result of Perrin, stated as Lion--Vergne Theorem A.8, in the setting of a local field $k$ of characteristic different from $2$. The proof carries over verbatim to any locally compact field of characteristic different from $2$, i.e., it applies to finite fields as well.

\begin{thm}[{Lion--Vergne, Perrin, \cite[A.8]{LionVergne}}]
	For three Lagrangian subspaces $L_1,L_2,L_3$ we have
		\[ \mc{F}_{L_1, L_2} \circ \mc{F}_{L_2,L_3} \circ \mc{F}_{L_3,L_1} = \tilde\gamma_\psi(\tau(L_1,L_2,L_3)). \]
\end{thm}

\begin{dfn}	
		 The \emph{Fourier--Heisenberg--Weil}\index[terminology]{Fourier--Heisenberg--Weil representation} representation relative to $L$ is the homomorphism\index[notation]{FWeilHeispsiL@$\mathrm{FWeilHeis}_{\psi,L}$} $\mathrm{FWeilHeis}_{\psi,L} : \WM_{L}(V) \ltimes H(V)(k) \to \mathrm{GL}(\mathrm{Heis}_\psi)$ that sends $(w,g, h)$ to $\tilde\gamma_\psi(w) \cdot  R_g \circ \mathrm{Heis}_\psi(h) $, where the semi-direct product corresponds to $\WM_{L}(V)$ acting on $H(V)(k)$ through its quotient $\Sp(V)(k)$.
\end{dfn}

\begin{cor} 
	$\mathrm{FWeilHeis}_{\psi,L}$ is a unitary representation of $\WM_{L}(V) \ltimes H(V)(k)$ that extends the unitary representation $\mathrm{Heis}_{\psi,L}$ of $H(V)(k)$ and whose projectivization equals $\mathrm{PWeilHeis}_\psi$.
\end{cor}

In \cite{Thomas2013}, Thomas observed that the theory becomes cleaner and more balanced if one removes the dependence on the Lagrangian subspace $L$ by considering the universal Witt-metaplectic group $\WM(V)$ of Definition \ref{dfn:meta_univ}. We will also consider the pushout $\WM_\psi(V)$\index[notation]{WMpsi@$\WM_\psi$} of $\WM(V)$ under $\tilde\gamma_\psi : \QW(k)/I_k^3 \to \mu_8^\times$, where $\mu_8$ denotes the group of complex roots of unity of order $8$.

\begin{lem} \label{lem:weil_univ}
	\begin{enumerate}
		\item  \label{eqn:weilunivalphaL} Let $L \subset V$ be a Lagrangian subspace. The map
		\[ \alpha_L : \WMuniv(V) \simto \WM_L(V),\qquad (w,g) \mapsto (w+b_L(g),g),  \]
		where $b_L(g)=\tau[U,\eta_U](L \oplus L,\Gamma_1,\Gamma_g,L \oplus gL)$, is an isomorphism of extensions of $\Sp(V)(k)$, in fact the unique such isomorphism unless $\dim(V)=2$ and $k=\F_3$.
		\item \label{eqn:weiluniv} The isomorphism class of the pull back of $\mathrm{FWeilHeis}_{\psi,L}$ under the isomorphism of (1) is a representation of $\WMuniv(V)$ that is independent of $L$. This representation factors through the pushout $\WM_\psi(V)$ of $\WM(V)$ under $q\mapsto \tilde\gamma_\psi(q)$ and extends to a genuine representation of $\WM_\psi(V)$ where the kernel $\mu_8$ of $\WM_\psi(V) \twoheadrightarrow \Sp(V)(k)$ acts via multiplication by $\mu_8$.
	\end{enumerate}
\end{lem}
\begin{proof}
	\eqref{eqn:weilunivalphaL} This is \cite[Proposition 3.4]{Thomas2013}.

	\eqref{eqn:weiluniv} One can either appeal to Thomas' formula \cite[Theorem B]{Thomas2013} for the character of the pull back of $\mathrm{FWeilHeis}_{\psi,L}$ to $\WMuniv(V)$, which evidently does not depend on $L$, or give the following more direct argument. Given two Lagrangian subspaces $L_1,L_2 \subset V$ choose $g \in \Sp(V)$ such that $L_2=gL_1$. Choose a Haar measure $d\ol{V}_1$ on $V/L_1$ to endow $\tx{Heis}_{\psi,L_q}$ with a unitary structure and take $d\ol{V}_2:=g_*(d\ol{V}_1)$ on $V/L_2$. As above we have the unitary isomorphism $A_g : \tx{Heis}_{\psi,L_1} \to \tx{Heis}_{\psi,L_2}$. Letting $\pi_{L_i}=\mathrm{FWeilHeis}_{\psi,L_i}\circ\alpha_{L_i}$ we see from the definition that for any $j \in \WMuniv(V) \ltimes H(V)(k) $ we have
	\[ A_g \circ \pi_{L_1}(j) \circ A_g^{-1} = \pi_{L_2}(B_g(j)), \]
	where $B_g((h,q,s))=(gh,q,gsg^{-1})$. Taking any lift $\tilde g=(q,g) \in \WMuniv(V)$ of $g$, we see that $B_g(j)=\tilde g \cdot j \cdot \tilde g^{-1}$, so the above identity shows $\pi_{L_1} \cong \pi_{L_2}$.
	
	Since $\tilde\gamma_\psi(c(1))=1$, the kernel $\QW(k)/I_k^3$ of $\WMuniv(V) \twoheadrightarrow \Sp(V)(k)$ acts via the character $\tilde\gamma_\psi:\QW(k)/I_k^3 \ra \mu_8$, and hence the claimed descent and extend of the representation $\mathrm{FWeilHeis}_{\psi,L}\circ\alpha_{L}$ to $\WM_{\psi}(V)$ follows.
\end{proof}

\begin{dfn} \label{dfn:weiluniv} The \emph{Fourier--Heisenberg--Weil representation}\index[terminology]{Fourier--Heisenberg--Weil representation} is the (isomorphism class of the) representation\index[notation]{FWeilHeispsi@$\mathrm{FWeilHeis}_\psi$} $\mathrm{FWeilHeis}_\psi$ of $\WM_\psi(V) \ltimes H(V)(k)$ that is the genuine extension of the descent of the pull-back of $\mathrm{FWeilHeis}_{\psi,L}$ under $\alpha_L$ for any Lagrangian subspace $L \subset V$, i.e., the representation from Lemma \ref{lem:weil_univ}\eqref{eqn:weiluniv}. We denote its restriction to $\WM_\psi(V)$ by $\mathrm{FWeil}_\psi$\index[notation]{FWeilpsi@$\mathrm{FWeil}_\psi$} and call it the \emph{Fourier--Weil representation}.\index[terminology]{Fourier--Weil representation}
\end{dfn}
 
\begin{lem} \label{lem:fweilrest} 
Let $k_0 \subset k$ be a subfield such that $k/k_0$ is a finite separable field extension. Let $V_0=V|_{k_0}$ and $\eta_0=\mathrm{tr}_{k/k_0}\circ\eta$. Let $\psi_0 : k_0 \to \C^\times$ be a non-trivial character and set  $\psi=\psi_0 \circ \mathrm{tr}_{k/k_0}$. The embedding $\Sp(V,\eta) \to \Sp(V_0,\eta_0)$ and the identity on $\mu_8$ splice to an embedding $\WM_\psi(V,\eta) \to \WM_{\psi_0}(V_0,\eta_0)$ and the pull-back of $\mathrm{FWeil}_{\psi_0}$ under this embedding is isomorphic to $\mathrm{FWeil}_\psi$.
\end{lem}

\begin{lem} \label{lem:fweilproduct} 
	Let $(V, \eta)=(\oplus_{i=1}^n V_i, \sum_{i=1}^n\eta_i)$ be the direct sum of (pairwise orthogonal) symplectic vector spaces $(V_i,\eta_i)$, and denote by $\mathrm{FHeisWeil}_{\psi,i}$ the Fourier--Heisenberg--Weil representation of $\WM_\psi(V_i) \ltimes H(V_i)(k)$. We have an inclusion $\prod_{i=1}^n\Sp(V_i) \hookrightarrow \Sp(V)$ and a surjection $\prod_{i=1}^n H(V_i)(k) \twoheadrightarrow H(V)(k)$ that splice to a homomorphism 
	\[\prod_{i=1}^n (\WM_\psi(V_i) \ltimes H(V_i)(k)) \ra \WM_\psi(V) \ltimes H(V)(k). \]
	The pull back of $\mathrm{FHeisWeil}_\psi$ under the above homomorphism is $\otimes_{i=1}^n\mathrm{FHeisWeil}_{\psi,i}$.
\end{lem}

We close this section by comparing the restriction of $\mathrm{FWeilHeis}_\psi$ to $\Mp(V)$ with $\mathrm{LWeilHeis}_\psi$. To do this, we must first address the issue that $\mathrm{LWeilHeis}_\psi$ is not unique when $\dim(V)=2$ and $k=\F_3$, because then $\Sp(V)=\mathrm{SL}_2(\F_3)$ is not a perfect group.

\begin{lem} \label{lem:weilheisuniq}
	The following characterizations of $\mathrm{LWeilHeis}_\psi$ are equivalent and specify it uniquely among the possible linearizations of $\mathrm{PWeilHeis}_\psi$ when $k$ is finite.
	\begin{enumerate}
		\item If $\dim(V)=2$ and $k=\F_3$, then $\mathrm{LWeil}_\psi$ is disjoint from its complex conjugate.
		\item Let $L \subset V$ be a Lagrangian subspace and $U \subset \Sp(V)$ the pointwise stabilizer of $L$. The restriction of $\mathrm{LWeilHeis}_{\psi,L}$ to $U(k) \ltimes H(V)(k)$ equals $\mathrm{Ind}_{U(k) \ltimes (k \times L) }^{U(k) \ltimes H(V)(k)}(\mathbf{1} \boxtimes \psi \boxtimes  \mathbf{1})$.
		\item The action of $U(k)$ on the $L$-fixed part of $\mathrm{LWeilHeis}_{\psi,L}$ is trivial.
	\end{enumerate}
\end{lem}
\begin{proof}
	Property (2) is proved in \cite[Theorem 2.4(b)]{Gerardin} for any finite $k$. We now check $(2) \Rightarrow (3)$. It follows from the Mackey formula that the $L$-fixed part of $\mathrm{LWeilHeis}_{\psi,L}$ is 1-dimensional: The restriction of $\mathrm{LWeilHeis}_{\psi,L}$ to $U(k) \times L \subset U(k) \ltimes H(V)(k)$ is a direct sum over $(U(k) \times k \times L) \lmod U(k) \ltimes H(V)(k)/(U(k)\times L)$ of twists of $(k \times L,\psi\boxtimes \mathbf{1})$. Choosing a complementary Lagrangian $L' \subset V$ can represent any element of this double coset space by $v \in L'$.
	But $k \times L$ is normal in $H(V)(k)$, while for each $0 \neq v \in L'$ the $v$-twist of $\psi$ has non-trivial restriction to $L \subset k \times L$. Therefore, the $L$-fixed points of $\mathrm{LWeilHeis}_{\psi,L}$ are precisely the summand indexed by the trivial double coset, i.e., $(U(k) \times k \times L,\mathbf{1} \boxtimes \psi \boxtimes \mathbf{1})$. On this summand, $U(k)$ acts evidently trivially.

	The lemma is thus proved unless $\dim(V)=2$ and $k = \F_3$, which we henceforth assume. Then $(1) \Rightarrow (2)$ is \cite[Theorem 2.4(a,b)]{Gerardin}, while $(2) \Rightarrow (3)$ was just established. To see $(3) \Rightarrow (1)$ note that any other linearization of $\mathrm{PWeilHeis}_\psi$ is of the form $\mathrm{LWeilHeis}_\psi\otimes\chi$ for some non-trivial character $\chi$ of $\mathrm{SL}_2(\F_3)$. Then $U(k)$ acts by $\chi|_{U(k)}$ on the $L$-fixed points of this linearization, but the restriction of $\chi$ to $U(k)$ is non-trivial. Therefore (3) can only hold for $\mathrm{LWeilHeis}_\psi$, which by definition implies (1).
\end{proof}

\begin{lem} \label{lem:fweilheisrest}
	The restriction of $\mathrm{FWeilHeis}_\psi$ to $\Mp(V)$ equals $\mathrm{LWeilHeis}_\psi$.
\end{lem}
\begin{proof}
	Except in the case of $\dim(V)=2$ and $k=\F_3$ the claim is automatic, since the restriction of $\mathrm{FWeilHeis}_\psi$ to $\Mp(V)$ is a linear representation that lifts $\mathrm{PWeilHeis}_\psi$, but there is a unique such.

	The outstanding case can be treated by the following argument that is valid for any finite field $k$ and any $\dim(V)$. By Lemma \ref{lem:weilheisuniq} it is enough to fix a Lagrangian $L \subset V$ and check that the unipotent radical $U(k) \subset \Sp(V)(k)=\Mp(V)(k) \subset \WM_L(k)$ associated to $L$ acts trivally on $(\mathrm{FWeilHeis}_{\psi,L})^L$.

	Let $f$ be the function on $H(V)(k)$ that is supported on $k \times L$ and equals $\psi \boxtimes \mathbf{1}$. Then $(\mathrm{FWeilHeis}_{\psi,L})^L=\C \cdot f$. The action of $u \in U(k)$ on this function is by $(uf)(x)=\tilde\gamma_\psi(c_L(u))\cdot f(u(x))$. Since $u$ acts trivially on $k \times L$ this becomes $(uf)(x)=\tilde\gamma_\psi(c_L(u))\cdot f(x)$. On the other hand, $c_L(u)=m(\tilde L,u\tilde L)$, where $\tilde L$ is $L$ equipped with an arbitrary choice of orientation. Since $uL=L$ the orientation on $u\tilde L$ is the orientation on $L$ multiplied by $\det(u|L)$. But $u$ acts trivially on $L$, so $u\tilde L=\tilde L$, hence $m(\tilde L,u\tilde L)=0$, hence $uf=f$.
\end{proof}

\subsection{A different viewpoint}  \label{sub:WittgroupC} 

Let $A$ be a locally compact abelian group (written multiplicatively) and let $\zeta : A \times A \to \mb{S}^1 \subset \C^\times$ be a bicharacter in the sense of Weil, i.e., a continuous map such that for every $b \in A$, the maps $a \mapsto \zeta(a, b)$ and $a \mapsto \zeta(b, a)$ are characters of $A$. Here $\mb{S}^1$ denotes the complex numbers of absolute value 1. We require that this bicharacter is alternating, i.e., $\zeta(a,a)=1$ for $a \in A$, and non-degenerate, i.e., the map $b \mapsto (a \mapsto \zeta(a,b))$ from $A$ to the group of continuous group homomorphism of $A$ into $\mb{S}^1 \subset \bC^\times$ is an isomorphism of topological groups. We define the \emph{symplectic group} 
\[ \Sp(A,\zeta)=\{g \in \Aut(A)|\zeta(g(a),g(b))=\zeta(a,b)\ \forall a,b \in A\}. \] 
A Lagrangian subgroup $L \subset A$ is a maximal isotropic subgroup. For three such we define the Maslov index $\tau(L_1,L_2,L_3)$ as the character of second degree
\[ L_1 \times L_2 \times L_3 \to \C^\times,\quad (x_1,x_2,x_3) \mapsto \zeta(x_1,x_2) \cdot \zeta(x_2,x_3)  \cdot \zeta(x_3,x_1). \]
For every $g \in \Sp(A,\zeta)$ the graph $\Gamma_g \subset A \times A$ of $g : A \to A$ is a Lagrangian subgroup of $A \times A$ with respect to the alternating bicharacter $\zeta^{-1} \times \zeta$. Let $\QW$ be Weissauer's  Witt group of non-degenerate characters of second degree\footnote{Note that Weissauer calls characters of second degree ``quadratic characters'', but we avoid this notion here to not confused it with a sign character.} as in \cite{Weissauer1986}. Recall that a (non-degenerate) character of second degree of a locally compact abelian group $A$ is a continuous map $f:A \to \mb{S}^1 \subset \bC^\times$ such that $(a,b)\mapsto f(ab)f(a)^{-1}f(b)^{-1}$ is a (non-degenerate) bicharacter of $A \times A$. We will actually use a small modification of Weissauer's group: we let $\QW_{2'}$ denote the Witt group of non-degenerate \emph{even} characters (even means $f(x^{-1})=f(x)$) of second degree whose source is a finite abelian group of \emph{odd} order, modulo Witt equivalences. On $\QW_{2'}$ the assignment $f \mapsto f^2$ is an automorphism. We have the map $\QW_{2'} \to \QW$ under which we pull back the Weil index.

Note that while Weissauer only considers non-degenerate characters of second degree, one can attach to each (possible degenerate) even character of second degree  $f: A \to \mb{S}^1$ of a finite abelian group of odd order $A$ an element $[\ol{f} : A/N \ra \mb{S}^1]$ in $\QW_{2'}$ as follows. Let $N \subset A$ be the radical of the bicharacter $b_f: (a,b)\mapsto f(ab)f(a)^{-1}f(b)^{-1}$. Evenness of $f$ implies $f(x)^2=b_f(x,x)$, hence $f(x)=1$ for $x \in N$ since the values of $f$ have odd order. Therefore $f$ descends to an (again even) character of second degree $\ol{f} : A/N \to \C^\times$. Its bicharacter is the factorization $\overline{b_f}$ of the bicharacter $b_f$ to the quotient $A/N$, and we have
\[ \ol{f}(x) = \overline{b_f}(x, \sqrt{x}), \]
 as in \cite[\S3, page 79]{Weissauer1986}.
for any lift $x \in A$ of $\bar x \in A/N$. 

Continuing with a finite group of odd order $A$ we obtain the 2-cocycle $z \in Z^2(\Sp(A,\zeta),\QW_{2'})$\index[terminology]{QW} as 
\[ z(g,h) = \tau(\Gamma_1,\Gamma_g,\Gamma_{gh}), \]
with the Maslov index for three Lagrangian subgroups of $(A \times A,\zeta^{-1} \times \zeta)$. Define\index[notation]{WM@$\WM$} 
\[ \WM(A,\zeta) = \QW_{2'} \boxtimes_z \Sp(A,\zeta). \]
The pushout along the renormalized Weil index $\tilde\gamma:=\gamma \circ \sqrt{(-)} : \QW_{2'} \to \mu_8$\index[notation]{gamma@$\gamma$} will be called $\WM_\gamma(A,\zeta)$, hence we have a short exact sequence
\[ 1 \ra \mu_8 \ra \WM_\gamma(A,\zeta) \ra \Sp(A,\zeta) .\]
Here we use the Weil index as defined in \cite[4.2~Lemma]{Weissauer1986}. 

\begin{lem}\label{lem:FWeilHeisgamma}
Let $(\ol{V},\bar\eta)$ be a symplectic $\F_p$-vector space and let $\psi : \F_p \to \C^\times$ be a non-trivial character. Let $A=\ol{V}$ and $\zeta=\psi\circ\bar\eta$. 
\begin{enumerate}
	\item The identity $\Sp(\ol{V},\bar\eta)=\Sp(A,\zeta)$ and the renormalized Weil index $\tilde\gamma_\psi : \QW(\F_p) \to \C^\times$ splice to an isomorphism  $\WM_\psi(\ol{V},\bar\eta)\to \WM_\gamma(A,\zeta)$.
	\item The transport of $\mathrm{FWeilHeis}_\psi$ under the resulting map 
	\[ \WM_\psi(\ol{V},\bar\eta) \ltimes H(\ol{V})(\bF_p) \to \WM_\gamma(A,\zeta) \ltimes H(\ol{V})(\bF_p)\]	
	only depends on $\zeta=\psi \circ \bar \eta$ and not on the choice of $\bar \eta$ and $\psi$. \end{enumerate}
\end{lem}
\begin{dfn}\label{dfn:FWeilHeis}
	We write $\mathrm{FWeilHeis}$\index[notation]{FWeilHeis} for the representation of  $\WM_\gamma(A,\zeta)\ltimes H(\ol{V})(\bF_p) $ corresponding to $\mathrm{FWeilHeis}_\psi$ as in Lemma \ref{lem:FWeilHeisgamma}, and we write $\mathrm{FWeil}$\index[notation]{FWeil} for its restriction to $\WM_\gamma(A,\zeta)$.
\end{dfn}

\end{appendices}

\printindex[notation]
\printindex[terminology]

\addcontentsline{toc}{section}{References}
\bibliographystyle{amsalpha}
\bibliography{bibliography.bib}

\end{document}